\documentclass[12pt]{dalthesis}
\usepackage[utf8]{inputenc}

\usepackage{graphicx}
\graphicspath{ {./images/} }
\usepackage[hidelinks]{hyperref}
\hypersetup{
	colorlinks=true,
	linkcolor=blue,
	citecolor=black,
	filecolor=magenta,
	urlcolor=blue, }

\usepackage{caption}
\usepackage[toc,page]{appendix}
\usepackage{cite}
\usepackage{url}
\usepackage{amsfonts,float,lscape,color}
\usepackage{textgreek}
\usepackage{upgreek}

\usepackage{amsmath,amsthm,amssymb,mathrsfs}
\usepackage[dvipsnames, table]{xcolor}
\usepackage{color,soul}
\usepackage{graphicx}
\usepackage[most]{tcolorbox}
\usepackage{subcaption}
\usepackage{tabularx}
\usepackage{longtable}
\usepackage{tabularray}

\usepackage{bbm}
\usepackage{mathtools}
\usepackage{cleveref}

\usepackage{faktor}
\usepackage{xfrac}

	\makeatletter
	\newcommand{\linkdest}[1]{\Hy@raisedlink{\hypertarget{#1}{}}}
	\makeatother
\newcommand{\iddef}[1]{\linkdest{#1}{#1}}
\newcommand{\idref}[1]{\hyperlink{#1}{#1}}

\usepackage[percent]{overpic} 
\usepackage{pgfplots}
\pgfplotsset{compat=newest}

\usetikzlibrary{intersections,calc, hobby}
\pgfplotsset{my style/.append style={axis x line=middle, axis y line=
		middle, xlabel={$x$}, ylabel={$y$}, axis equal }}
\usepackage{graphicx,psfrag}
\usepackage{epsf} 
\usepackage[all]{xy}
\usepackage{tikz}
\usepackage{tikz-cd}
\usetikzlibrary{calc,patterns,angles,quotes} 
\usepackage{tkz-fct}  

\usepgfplotslibrary{polar}
\usepgfplotslibrary{fillbetween}

\usepackage{caption}

\tikzstyle{vertex}=[shape=circle, minimum size=2mm, inner sep=0, fill]
\tikzstyle{opendot}=[shape=circle, minimum size=2mm, inner sep=0, fill=white, draw]

\pgfplotsset{my style/.append style={axis x line=middle, axis y line=
		middle, axis equal }}

\usepackage{calc}

\usetikzlibrary{calc}

\usepackage{pgf}
\usepackage{pgffor}

\usetikzlibrary{patterns, patterns.meta}
\usetikzlibrary{shadings}
\usetikzlibrary{shadings}
\usetikzlibrary{arrows,decorations,backgrounds}
\usetikzlibrary{hobby}  
\usetikzlibrary{matrix}

\newcommand{\ud}{\,\mathrm{d}}

\theoremstyle{plain}
\newtheorem{theorem}{Theorem}[section]
\newtheorem*{theorem*}{Theorem}
\newtheorem{prop}[theorem]{Proposition}
\newtheorem{lemma}[theorem]{Lemma}
\newtheorem{cor}[theorem]{Corollary}

\theoremstyle{definition} 
\newtheorem{defn}[theorem]{Definition}
\newtheorem{ex}[theorem]{Example}

\theoremstyle{remark}
\newtheorem{rem}[theorem]{Remark}
\newtheorem{notation}{Notation}[section]

\theoremstyle{plain}
\newtheorem*{namedthm}{\namedthmname}
\newcounter{namedthm}

\makeatletter
\newenvironment{named}[1]
{\def\namedthmname{#1}%
	\refstepcounter{namedthm}%
	\namedthm\def\@currentlabel{#1}}
{\endnamedthm}
\makeatother

\newcommand{\M}{{\fontfamily{pzc}\selectfont \text{m}}}

\usepackage{upgreek}

\makeatletter
\def\namedlabel#1#2{\begingroup
	#2%
	\def\@currentlabel{#2}%
	\phantomsection\label{#1}\endgroup
}
\makeatother

\begin{document}
	
	\title{Subspaces of $L^2(\mathbb{R}^{\lowercase{n}})$ Invariant Under Crystallographic Shifts}
	\author{Tom Potter}

	\phd 
	
	\degree{Doctor of Philosophy}
	\degreeinitial{Ph.D.}
	\faculty{Science}
	\dept{Mathematics and Statistics}
	
	\defencemonth{September}\defenceyear{2025}

	\nolistoftables
	\nolistoffigures
	
	\frontmatter

	\begin{abstract}
	In this thesis we consider crystal groups in dimension $n$ and their natural unitary representation on $L^2(\mathbb{R}^n)$. We show that this representation is unitarily equivalent to a direct integral of factor representations, and use this to characterize the subspaces of $L^2(\mathbb{R}^n)$ invariant under crystal symmetry shifts. Finally, by giving an explicit unitary equivalence of the natural crystal group representation, we find the \textit{central decomposition} guaranteed by direct integral theory. 
	\end{abstract}

\chapter*{List of Abbreviations and Symbols Used}
\addcontentsline{toc}{chapter}{List of Abbreviations and Symbols Used}  
	{	\centering
		\small
		\DefTblrTemplate{firsthead, middlehead,lasthead}{default}{}
		\begin{longtblr}{width=\textwidth,
				colspec={X[.2]X[.7]X[r,.1]},
			row{odd} = {gray!10!white},
			row{1} ={gray!20!white, 6mm},
			rows={m, rowsep=0pt}
			}
				\textbf{Symbol}  & \textbf{Description} & \textbf{Page} \\
$L^2(\mathbb{R}^n)$ & Lebesgue space of complex-valued square-integrable functions & \pageref{Rn}\\
$T_y$ & translation by $y$ & \pageref{translation}\\
$\mathbb{R}^n$ & $n$-dimensional Euclidean space & \pageref{Rn}\\
$\mathbb{T}^n$ & $n$-dimensional torus & \pageref{ntorus}\\
$\mathrm{GL}(n,\mathbb{R})$ & general linear group of degree $n$& \pageref{GLnr}\\
$[x,L]$ & an element of $\mathrm{Aff}(n,\mathbb{R})$ & \pageref{xL}\\
$\mathrm{O}(n)$ & group of orthogonal transformations of degree $n$ & \pageref{On}\\
$\mathrm{I}(n)$ & group of isometries of $\mathbb{R}^n$ & \pageref{In}\\
$\Gamma$ & crystallographic group & \pageref{crystalgroup}\\
$Q$ & quotient map from $\mathrm{Aff}(n,\mathbb{R})$ to $\mathrm{GL}(n,\mathbb{R})$& \pageref{Q}\\
$\Pi$ & the point group of $\Gamma$& \pageref{pi}\\
$\mathrm{T}$ & the translation subgroup of $\Gamma$& \pageref{T}\\
$\tilde{q}$ & canonical quotient map from $\Gamma$ to $\Gamma/\mathrm{T}$ & \pageref{tildeq}\\
$Q^*$ & map identifying $\Gamma/\mathrm{T}$ with $\Pi$ & \pageref{Qstar}\\
$\mathcal{L}$ & Bravais lattice associated with $\Gamma$ & \pageref{L}\\
$p4m$ & symmorphic wallpaper group described in Example \ref{p4m} & \pageref{p4m}\\
$pg$ & nonsymmorphic wallpaper group described in Example \ref{pgexample} & \pageref{pgexample}\\
$\gamma$ & a cross-section for an equivalence relation & \pageref{crosssec}\\
$\mathcal{L}^*$ & dual lattice of $\mathcal{L}$ & \pageref{duallattice}\\
$\widehat{\mathcal{L}}$ & the group of continuous characters on $\mathcal{L}$ & \pageref{contcharonL}\\
$\chi_y$ & a continuous character on $\mathbb{R}^n$ & \pageref{chiy}\\
$\widehat{\mathbb{R}^n}$ &the character group of $\mathbb{R}^n$ & \pageref{Rnchar}\\
$\mathcal{L}^\perp$ & the orthogonal subgroup of $\widehat{\mathbb{R}^n}$ & \pageref{orthsubgroup}\\
$\Gamma^*$ & the semidirect product of $\mathcal{L}^*$ with $\Pi$ & \pageref{Gammastar}\\
$R$ & a fixed fundamental domain for $\Gamma^*$ & \pageref{Rdomain}\\
$\Omega_{\mathrm{T}^*}$ & a fixed fundamental domain for $\mathrm{T}^*$ & \pageref{OmegaTstar}\\
$\chi_\nu^\mathcal{L}$ & the character of $\mathcal{L}$ associated with $\nu \in \Omega_{\mathrm{T}^*}$ & \pageref{XLdef}\\
$X$& the map parametrizing $\widehat{\mathcal{L}}$ & \pageref{Xmappingdef}\\
$\chi^\mathcal{L}$ & a character of $\mathcal{L}$ & \pageref{chiL}\\
$\chi^\mathrm{T}$ & a character of $\mathrm{T}$ & \pageref{chiL}\\
$\cdot_{i}$ & one of the twelve group actions associated with $\Gamma$ & \pageref{TableOfActions}\\
$\vartheta$  & the map identifying $\widehat{\mathrm{T}}$ with $\widehat{\mathcal{L}}$ & \pageref{vartheta}\\
$\chi_\nu^{\mathrm{T}}$  & the continuous character on $\mathrm{T}$ associated with $\nu \in \Omega_{\mathrm{T}^*}$ & \pageref{chiT}\\
$\upsilon$ & the map parametrizing $\widehat{\mathrm{T}}$ & \pageref{upsilon}\\
$\mathscr{V}^*$  & dual space of a vector space $\mathscr{V}$ & \pageref{dualvector}\\
$\sigma(\mathscr{V},\mathscr{F})$  & the weak topology induced by a family of linear functionals & \pageref{weaktopinduced}\\
$\mathscr{B}(\mathscr{H})$ & bounded linear operators on $\mathscr{H}$ & \pageref{BH}\\
$\mathscr{B}(\mathscr{H},\mathscr{K})$ & bounded linear transformations from $\mathscr{H}$ to $\mathscr{K}$ & \pageref{boundedHK}\\
$T^*$ & the Hilbert space adjoint of $T$ & \pageref{adjoint}\\
$\mathrm{U}(\mathscr{H})$ & the group of unitary operators on $\mathscr{H}$ & \pageref{adjoint}\\
$(\pi,\mathscr{H}_\pi)$ & a unitary representation of a group on $\mathscr{H}_\pi$ &\pageref{unitary representation}\\
$\{\mathscr{M},\mathscr{H}\}$  & a von Neumann algebra on $\mathscr{H}$ &\pageref{vnalg}\\
 $\mathscr{U}^{-}$ & von Neumann algebra generated by $\mathscr{U}$ &\pageref{vngen}\\
 $\mathscr{U}^{\prime}$ & the commutant $\mathscr{U}$ &\pageref{commutant} \\
  $\mathscr{U}^{\prime\prime }$ & the double commutant $\mathscr{U}$ &\pageref{doublecommutant}\\
  $\mathscr{Z}(\mathscr{M})$ & the center of $\mathscr{M}$ &\pageref{center}\\
  $\mathscr{H}_1\otimes \cdots\otimes \mathscr{H}_n$ & the tensor product of the Hilbert space $\mathscr{H}_1,\ldots,\mathscr{H}_n$ &\pageref{Hilberttensor}\\
$T_1\otimes\cdots \otimes T_n$ & the tensor product of the operators $T_1,\ldots,T_n$ &\pageref{operatortensor}\\
$\mathscr{M}_1 \overline{\otimes}  \cdots\overline{\otimes} \mathscr{M}_n$& the tensor product of the von Neumann algebras $\mathscr{M}_1,\ldots,\mathscr{M}_n$ &\pageref{vntensor}\\
$\mathrm{Hom}_G(\pi,\sigma)$ & the intertwining space of $\pi$ and $\sigma$ &\pageref{intspace}\\
$\pi\sim \sigma$ & unitary equivalence of $\pi$ and $\sigma$ &\pageref{equivrep}\\
$P_K$ & orthogonal projection of $P$ onto $K$ & \pageref{PK}\\
$\textstyle{\sum}_{i\in I} \oplus \mathscr{H}_i$ & direct sum of Hilbert spaces $\mathscr{H}_i$ & \pageref{directsumrep}\\
$\textstyle{\sum}_{i\in I} \oplus\pi_i$ & direct sum of representations $\pi_i$ & \pageref{directsumrep}\\
$m\rho$& amplification of $\rho$ by m & \pageref{mrho} \\
$\pi_1\otimes \pi_2$ & inner tensor product of representations & \pageref{innertensor} \\
$\int_X^\oplus \mathscr{H}_x \ud\mu(x)$ & direct integral of the Hilbert spaces $\mathscr{H}_x$ & \pageref{directinthilbspace} \\
$\int_X^\oplus \xi(x)\ud\mu(x)$ & a measurable field in $\int_X^\oplus \mathscr{H}_x \ud\mu(x)$  & \pageref{directfield} \\
$\int_X^\oplus T(x)\ud\mu(x)$& direct integral of the operators $T(x)$   & \pageref{directintops} \\
 $\|T\|_\infty$ &  essential supremum of the operator field $\{T(x)\}_{x\in X}$ & \pageref{esssup}\\
$\mathscr{A}$ & the diagonal algebra & \pageref{diag} \\
$\int_X^\oplus  \mathscr{M}_x \ud \mu(x)$ & direct integral of the von Neumann algebras $\mathscr{M}_x$ & \pageref{directVN} \\
$\int_X^\oplus \pi_x(g) \ud\mu(x)$ & direct integral of the representations $\pi_x$ & \pageref{directreps} \\
$\mathscr{P}(G)$ & functions of positive type & \pageref{fnpostype} \\
$\mathscr{P}_1(G)$ & functions of positive type and norm 1 & \pageref{fnpostype} \\
$\widehat{G}$ & the unitary dual space of $G$ & \pageref{dualspace}\\
$C(X)$ & space of continuous functions on $X$ & \pageref{CX} \\
$C_0(X)$ & space of continuous functions vanishing at infinity & \pageref{CX} \\
$\widehat{f}$ & the Fourier transform of $f$ & \pageref{FourierT} \\
$\mathcal{F}$ & the Plancherel transform & \pageref{FT0}\\
 $\mu_{G}$ & Haar measure on LCA group $G$  & \pageref{dualmeas0}\\
$\mu_{\widehat{G}}$ & dual measure of $\mu_G$  & \pageref{dualmeas0}\\
$x\mapsto J(x)$ & measurable range function  & \pageref{range}\\
$L^2(X;\mathscr{H})$ & $\mathscr{H}$-valued Lebesgue $L^2$ space  & \pageref{vectorL2} \\
$(X,\mathcal{M},\mu)$ & measure space with measure $\mu$ & \pageref{measspace}\\
$\overline{\mathcal{M}}$ & completion of $\sigma$-algebra $\mathcal{M}$ & \pageref{measspace}\\
$\overline{\mu}$ & completion of measure $\mu$ & \pageref{mucomplete}\\
$\bigvee_{i\in I} v_i$ & closed linear span of family of vectors $\{v_i\}$ & \pageref{bigvee}\\
$J_W$ & measurable range function associated with $W$ & \pageref{mblerangefun}\\
$\pi$ & used to denote the natural crystallographic representation (from Ch. 4 onwards) & \pageref{pi0}\\
$\hat{\pi}$ & the transformed crystallographic representation & \pageref{FT}\\
$\mathrm{ind}_H^G \sigma$ & representation induced from $\sigma$ & \pageref{induced}\\
$l^2(\Pi)$ & the square integrable functions on $\Pi$ & \pageref{l2pi}\\
$U^{\chi_\nu^\mathrm{T}}$ & representation equivalent to $\mathrm{ind}_H^G \sigma$ & \pageref{uchinotation0}\\
$U^\nu$ & abbreviation for $U^{\chi_\nu^\mathrm{T}}$  & \pageref{uchinotation}\\
$\alpha(L,M)$ & the cocycle of $L$ and $M$ & \pageref{cocycle0}\\
$G_x$ & the stabilizer of $x$ & \pageref{stab}\\
$\Gamma_\nu$ & the stabilizer of $\Gamma$ under the action $\idref{5}$ & \pageref{GammaandPiStab}\\
$\Pi_\nu$ & the stabilizer of $\Pi$ under the action $\idref{6}$ & \pageref{GammaandPiStab}\\
$\widehat{\Gamma_{\nu}}^\nu$ & equivalence classes of irreducible representations of $\Gamma_\nu$ which restrict to a multiple of $\chi_\nu^{\mathrm{T}}$ on $\mathrm{T}$ & \pageref{GammaandPiStab}\\
$\Omega_0$ & the set of free points in $\Omega$ & \pageref{Omega0}\\
$\mathcal{B}_X$ & Borel $\sigma$-algebra on $X$ & \pageref{L2identification}\\
$\Phi$ & the map parametrizing $\widehat{\mathbb{R}^n}$ over $\mathbb{R}$& \pageref{PhiMap}\\
$\mathscr{L}$ & Lebesgue $\sigma$-algebra on $\mathbb{R}^n$ & \pageref{Lebesguesigalg0}\\
$\widehat{\mathscr{L}}$ & Lebesgue $\sigma$-algebra on $\widehat{\mathbb{R}^n}$ & \pageref{Lebesguesigalg0}\\
$\M$ & the measure on $\widehat{\mathbb{R}^n}$ & \pageref{Lebesguesigalg0}\\
$\Upsilon$& unitary transformation from	$L^2(\widehat{\mathbb{R}^n})$ onto  $L^2\big({R};l^2(\Pi)\otimes l^2(\mathcal{L}^*)\big)$ & \pageref{Upsilondef}\\
$l^2(\mathcal{L}^*)$ & the square integrable functions on $\mathcal{L}^*$ & \pageref{l2lstar}\\
$\tau$ & the unitary $U\pi U^{-1}$ on $L^2\big({R};l^2(\Pi)\otimes l^2(\mathcal{L}^*)\big)$ &\pageref{taudef}\\
$U$ & unitary from $L^2(\mathbb{R}^n)$ onto $L^2\big({R};l^2(\Pi)\otimes l^2(\mathcal{L}^*)\big)$ &	\pageref{Udef}\\

$\mathcal{C}$ & homeomorphism from $\mathbb{R}^n/\mathcal{L}^*$ onto $\widehat{\mathrm{T}}$ & \pageref{C}\\
$\widehat{q}$ & quotient map from $\widehat{\mathbb{R}^n}$ to $\widehat{\mathbb{R}^n}/\mathcal{L}^\perp$ &\pageref{qhat}\\
$\upsilon_0$ & the restriction of $\upsilon$ to $\Pi R$ & \pageref{upsilon0}\\
$\widehat{\mathrm{T}}_0$ & image of $\Pi R$ under $\upsilon$ & \pageref{T0} \\
$|\Pi|$ & order of the group $\Pi$ & \pageref{order}\\
$p$ & unitary transformation from $l^2(\Pi)\otimes l^2(\mathcal{L}^*)$ onto $l^2(\Pi\times \mathcal{L}^*)$ &\pageref{p} \\
 $\mathrm{P}$ & unitary transformation from  $L^2\big(R;l^2(\Pi)\otimes l^2(\mathcal{L}^*)\big)$ onto $L^2\big(R;l^2(\Pi\times\mathcal{L}^*)\big)$ & \pageref{P}\\
$\kappa_1$ &unitary transformation from  $L^2\big(R;l^2(\Pi)\otimes l^2(\mathcal{L}^*)\big)$ onto $L^2(\Pi\times R \times \mathcal{L}^*)$ & \pageref{kappa1}\\
$\Lambda$ & bijective mapping from $\Pi \times R$ onto $\Pi R$ & \pageref{Lambda} \\
	$\widetilde{\Lambda}$ & unitary transformation from $L^2(\Pi R\times \mathcal{L}^*)$ onto $L^2(\Pi \times R\times \mathcal{L}^*)$ & \pageref{tildeLambda}\\
$\kappa_2$ & unitary transformation from $L^2(\Pi\times R\times \mathcal{L}^*)$ onto $L^2(\Pi R\times \mathcal{L}^*)$ & \pageref{kappa2}\\
$\widetilde{\upsilon_0\times I}$ & unitary transformation from $L^2(\widehat{\mathrm{T}}_0 \times I)$ onto $L^2(\Pi R \times \mathcal{L}^*)$ & \pageref{upsilon0I}\\
$\tilde\mu_{\widehat{\mathrm{T}}}$ & scaled measure on $\widehat{\mathrm{T}}$ & \pageref{tildemuthat}\\
$\kappa_3$ & unitary transformation from $L^2(\Pi R \times \mathcal{L}^*)$ onto $L^2(\widehat{\mathrm{T}}_0 \times I)$ & \pageref{kappa3}\\
$\kappa_4$ &unitary transformation from $L^2(\widehat{\mathrm{T}}_0\times\mathcal{L}^*)$ onto
$L^2\big(\widehat{\mathrm{T}}_0; l^2(\mathcal{L}^*)\big)$ & \pageref{kappa4}\\
$\kappa_5$ & unitary map identifying $L^2\big(\widehat{\mathrm{T}}_0; l^2(\mathcal{L}^*)\big)$ with	$L^2\big(\widehat{\mathrm{T}}; l^2(\mathcal{L}^*)\big)$  & \pageref{kappa5}\\
$\kappa$ & unitary mapping from $L^2\big(R;l^2(\Pi)\otimes l^2(\mathcal{L}^*)\big)$ onto $L^2\big(\widehat{\mathrm{T}}; l^2(\mathcal{L}^*)\big)$ & \pageref{kappadef}\\
$E_{[y,\mathrm{id}]}$ & evaluation character on $\widehat{\mathrm{T}}$ & \pageref{evalcharonThat}\\
$\mathrm{M}$ & unitary representation of $\mathrm{T}$ on $L^2\big(\widehat{\mathrm{T}}; l^2(\mathcal{L}^*)\big)$ & \pageref{Mrep}\\
$\mathscr{S}_{[z,M]}$ & direct summand of $L^2\big({R};l^2(\Pi)\otimes l^2(\mathcal{L}^*)\big)$&\pageref{Szmdef}\\
$u$ & unitary transformation from $L^2({R}) \otimes l^2(\Pi)\otimes l^2(\mathcal{L}^*)$ onto $L^2\big({R};l^2(\Pi)\otimes l^2(\mathcal{L}^*)\big)$ & \pageref{udef0}\\
$u_{a,g}$ & an element of $L^2\big({R};l^2(\Pi)\otimes l^2(\mathcal{L}^*)\big)$ & \pageref{udef2}\\
$\mathscr{L}_{[z,M]}$ & direct summand of $L^2\big({R};l^2(\Pi)\otimes l^2(\mathcal{L}^*)\big)$ equal to $\mathscr{S}_{[z,M]}$ &\pageref{Lzmdef}\\
$I_{[w,N]}$ & bounded linear transformation from $L^2\big({R};l^2(\Pi)\otimes l^2(\mathcal{L}^*)\big)$ onto $L^2({R})$ & \pageref{Izm}\\
$\mathcal{D}$ & determining set for $L^1(X)$ & \pageref{detset}\\
 $\mathcal{C}(X)$ & the collection of all closed subsets of $X$ & \pageref{CX1}\\
 $\mathcal{K}(G)$ & the set of all closed subgroups of $G$& \pageref{KG}\\
  $\mathcal{S}(G)$ &  the set of all subgroup-representation pairs & \pageref{SG}\\
  $\mathcal{S}_0(G)$ & the set of \textit{cataloguing pairs} for $\widehat{G}$ & \pageref{SG}\\
 $I$ & the inducing map& \pageref{Imap} \\
 $\mathscr{A}(G)$& the set of all pairs $(K, \sigma)$ where $K \in \mathcal{K}(G)$ and $\sigma$ is an irreducible unitary representation of $K$ & \pageref{AG}\\
 $\widetilde{\mathcal{S}}_0(G)$ & equivalence classes in $\mathcal{S}_0(G)$ & \pageref{tildeSOG} \\
$S(\nu)$& irreducible unitary representations of $\Gamma_\nu$ 
 whose restriction to $\mathrm{T}$ is a multiple of $\chi_\nu^\mathrm{T}$ & \pageref{Snu}\\
 $\Psi$ & Injective Borel mapping from $R$ into $\widehat{\Gamma}$ & \pageref{Psi}\\
 $\widehat{\Gamma}_0$ & image in $\widehat{\Gamma}$ of $R$ under $\Psi$&\pageref{Gamma0hat} \\
$\mathcal{E}(X)$& set of all complex-valued functions $f$ defined on a closed subset of $X$&\pageref{EX}\\
$\mathcal{E}_s(G)$& the topological subspace of $\mathcal{E}(G)$ consisting of all $f\in \mathcal{E}(G)$ with domain $D(f)$ a closed subgroup of $G$&\pageref{EsG}\\
$\widetilde{\Psi}$ & unitary transformation from  $L^2\big(R;l^2(\Pi)\otimes l^2(\mathcal{L}^*)\big)$ onto $L^2\big(\widehat{\Gamma};l^2(\Pi)\otimes l^2(\mathcal{L}^*)\big)$ & \pageref{tildePsi}\\
$\rho$ & unitary representation $\rho$ of $\widehat{\Gamma}$ on $L^2\big(\widehat{\Gamma};l^2(\Pi)\otimes l^2(\mathcal{L}^*)\big)$ & \pageref{rho1}\\
$\pi_1|_{\mathrm{T}}$ & map identifying $\mathrm{T}$ with $\mathcal{L}$ & \pageref{pi1}\\
	$\mathcal{M}\vert_E$ & relative $\sigma$-algebra on $E\subset X$ & \pageref{relsigalg}\\
	$\mathcal{P}(Y)$ & power set of $Y$ & \pageref{PY}\\
	$\mu \circ f^{-1}$ & image of measure $\mu$ under $f$ & \pageref{imegameas}\\

				\end{longtblr}
			}

	\begin{acknowledgements}
		First and foremost, the author would like to express his gratitude to his supervisor, Dr. Keith Taylor. Keith met with me regularly, well into his retirement years, providing helpful critiques, stylistic comments, mathematical insights, and guidance, which helped shape this work and me as a mathematician. In terms of intellectual acknowledgements, Keith gave me this problem, directed me to relevant literature, and suggested the tools that would be needed: namely, that the representation would be equivalent to a direct integral of induced representations, thereby allowing us to find the invariant subspaces. Keith also provided critical insights at various times. But more importantly, I am deeply grateful to Keith for his patience, support, and sincere commitment to me as his student.
		
I am also very grateful to the committee---Dr. Yemon Choi (external examiner), Dr. Theo Johnson-Freyd (reader), and Dr. Keith Johnson (reader)---for reading and providing corrections and thoughtful suggestions. These substantially improved the final draft of this work, especially Yemon's detailed list of suggestions and typos.
		
		I would also like to thank the many friends, colleagues, and well-wishers during this period of my life, for their support and encouragement. The encouragement of faculty, administrators, and staff members has helped me to feel supported in the department. I am particularly grateful to Dr. Dorette Pronk for collaborating with me on math outreach projects, for giving feedback on my C.V. and job application documents, and for always believing in me as a student. I would like to thank Dr. Sarah Chisholm for providing feedback on my teaching and creating opportunities in the department to focus on math pedagogy. I am also grateful for the longtime support and encouragement of professors in the department, especially Dr. Bob Par\'e, Dr. Karl Dilcher, and Dr. Chelluri Sastri. I would also like to thank Dr. John Baxter of the English department, for allowing me to sit in on numerous of his upper-level English classes which deepened my appreciation of English literature, challenged my mind in a different direction, and were a lot of fun.
		
		I am deeply grateful to the Bragg family and Eastlink for supporting the NS Math Circles program, which I was closely involved with for 6 years of my graduate student life. This unique opportunity made a huge impression on me. The experience of working with schoolkids to foster a love of math, through interactive and fun activities, is something that few graduate students have the privilege to experience. Eastlink’s funding made it possible, and it was such a pleasure to share the excitement of doing math with young people.
		
		I am also deeply grateful to the family of the late Dr. Patrick Lett, whose generous bursary for late-stage math graduate students helped alleviate the financial stress of doing a Ph.D. and helped make it possible.
		
		Last but not least, I am deeply grateful to my parents, who each supported me tremendously in their own ways during this long journey.
		
	\end{acknowledgements}

	\mainmatter
	
\clearpage

	\chapter{Introduction} 
\section{Background and Literature Review}
This thesis is concerned with  an arbitrary
crystallographic group $\Gamma$ in dimension $n$ and the Hilbert space, $L^2(\mathbb{R}^n)$,
consisting of square-integrable complex-valued functions on $\mathbb{R}^n$. The
natural action of $\Gamma$ on $\mathbb{R}^n$ as isometric affine transformations
determines a unitary representation, $\pi$, of the group $\Gamma$ on $L^2(\mathbb{R}^n)$. We answer the
question of which closed subspaces of $L^2(\mathbb{R}^n)$ are 
$\pi$-invariant. That is, we will characterize those closed
subspaces $V$ of $L^2(\mathbb{R}^n)$ with the property that, if
$f\in V$, then the shifts of $f$ by elements of $\Gamma$
are all in $V$, as well.

One of the motivations for describing the $\pi$-invariant closed
subspaces of $L^2(\mathbb{R}^n)$ comes from the theory of wavelets. This theory was initiated in the 1980's and wavelets quickly
had an impact on signal and image processing. The book of
Daubechies \cite{Daubechies} provides a 
comprehensive introduction to wavelets. Let us briefly sketch enough of 
the theory of wavelets to introduce the concepts of interest to us.
In
dimension one, there are two basic unitary maps on $L^2(\mathbb{R})$
that are useful in classic wavelet theory: The {\em shift operator} $T_1$ is defined by $T_1f(x)=f(x-1)$, for
a.e. $x\in\mathbb{R}$ and all $f\in L^2(\mathbb{R})$. The {\em dilation by $2$ operator}
$D_2$ is given by $D_2f(x)=2^{k/2}f(2^{k}x)$, for
a.e. $x\in\mathbb{R}$ and all $f\in L^2(\mathbb{R})$. A classical wavelet is a
function $\psi\in L^2(\mathbb{R})$ such that 
$\{D_2^jT_1^k\psi : j,k\in\mathbb{Z}\}$ is an orthonormal basis of
$L^2(\mathbb{R})$. Then, for any $f\in L^2(\mathbb{R})$, 
\[
f=\sum_{j=1}^\infty\sum_{k=1}^\infty\langle f,D_2^jT_1^k\psi
\rangle D_2^jT_1^k\psi.
\]
So $f$ can be recovered from the countable set
$W_f=\{\langle f,D_2^jT_1^k\psi\rangle:j,k\in\mathbb{Z}\}$ of doubly
indexed complex numbers. 
Moreover, for $\varepsilon>0$, $W_{f,\varepsilon}=
\{\alpha\in W_f:|\alpha|>\varepsilon\}$ is a finite set from
which $f$ can be approximated.  We will use the term {\em Discrete Wavelet
	Transform} (DWT) for the map $f\to W_f$

The concept of a multiresolution analysis (MRA) is integral to
the development of 
efficient algorithms for the DWT. An MRA is a bilateral sequence 
\[
\cdots V_{-1}\subset V_0\subset V_1\subset V_2 \cdots
\]
of closed subspaces of $L^2(\mathbb{R})$. Their union must be dense
in $L^2(\mathbb{R})$ and their intersection must be trivial. Also, any $V_k$
is just $V_0$ scaled by $2^k$ in the sense that $V_k=D_2^kV_0$.
But the property of an MRA that is of special interest in this thesis is that the central subspace $V_0$ is invariant under shifts by integers.
That is, for any $k\in\mathbb{Z}$, if $f\in V_0$, then $T_kf\in V_0$, where
$T_kf(x)=T_1^kf(x)=f(x-k)$, for a.e. $x\in\mathbb{R}$. Such a closed subspace of $L^2(\mathbb{R})$ is called {\em shift-invariant}.
In dimension one, there
are just two crystallographic groups (up to changes of scale).
When $\Gamma=\mathbb{Z}$, the abelian crystallographic group, the
natural representation $\pi$ is simply given by $\pi(k)=T_k$,
for $k\in\mathbb{Z}$. Thus, shift-invariant is the same as 
$\pi$-invariant for the abelian crystallographic group.

There are many variations of wavelet theory, 
especially in higher dimensions,
with two dimensional versions now widely used in image processing. 
The concept of an MRA is important in higher dimensions with the
central closed subspace now being shift-invariant with the shifts
coming from a full-rank lattice such as $\mathbb{Z}^n$. 

In \cite{KeithJosh}, MacArthur and Taylor introduced a
theory of wavelets with shifts by a crystallographic group acting naturally on $\mathbb{R}^n$. They also defined the concept of a
$\Gamma$-MRA, where $\Gamma$ is a crystallographic group.
In this case, the closed subspace $V_0$, in the $\Gamma$-MRA,
is invariant under the action of $\Gamma$. Hence, the interest in characterizing such closed subspaces.

Motivated by wavelet theory, Bownik \cite{BownikRoss} used the concept of a range function, adapted from Helson \cite{Helson}, to give a characterization of the closed subspaces of $L^2(\mathbb{R}^n)$ that are shift-invariant for shifts by $\mathbb{Z}^n$. With the
appearance of \cite{KeithJosh}, a natural question to
consider was to generalize the characterization of shift-invariant
closed subspaces of $L^2(\mathbb{R}^n)$ to invariance under the action of
shifts by a crystallographic group. Our approach to this question is
to use the theory of unitary representations applied to
the natural representation $\pi$.

A substantial part of this thesis is devoted to deriving an 
explicit decomposition of $\pi$ as a direct integral of distinct factor
representations (Theorem \ref{DecompTheorem}). Each of these factor representations is a
countably infinite multiple of a distinct irreducible unitary 
representation of $\Gamma$. Besides its use to characterize
the $\pi$-invariant closed subspaces of $L^2(\mathbb{R}^n)$\label{Rn}, this direct integral decomposition should be of independent interest.

We use Mackey's theory of induced representations (see \cite{TaylorKaniuth}, \cite{FollandAbstract}, or \cite{Mackey})
to identify a subset $R$ of ``frequency'' space
and a map, $r\to U^r$, of $R$ into the
set of irreducible unitary representations of $\Gamma$. Each
$U^r$ is expressed with an explicit formula \eqref{Uchi1} as a
unitary operator on the Hilbert space $l^2(\Pi)$, where
$\Pi$ is the {\em point group} of $\Gamma$. That is, $\Pi$
is the quotient group of $\Gamma$ modulo the subgroup of
pure translations. Thus, each $U^r$ acts on the same finite-dimensional Hilbert space and they are mutually
inequivalent. Moreover, $\{U^r:r\in R\}$
has enough irreducible unitary representations to provide
the direct integral decomposition of $\pi$ discussed
above. A Hilbert space isomorphism $U$ is found
that ``transforms'' $L^2(\mathbb{R}^n)$ into a direct integral over
$R$ of Hilbert spaces, $\mathscr{H}_r$, where
$\mathscr{H}_r=l^2(\Pi)\otimes\mathscr{H}$, for
each $r\in R$, with
$\mathscr{H}$ a fixed separable infinite-dimensional Hilbert
space. The unitary map we find also transforms $\pi$ into
a unitary representation, $\tau$ that is 
the direct integral of the factor representations $U^r\otimes I$, where $I$ is the identity operator on $\mathscr{H}$. The
measure for these direct integrals is just Lebesgue measure
restricted to $R$.

We use basic results from the theory of von Neumann algebras
applied to this direct integral decomposition to describe
the $\tau$-invariant closed subspaces of the direct 
integral Hilbert space. We can then transfer this 
description back to a characterization of the 
$\pi$-invariant closed subspaces of $L^2(\mathbb{R}^n)$.

In order to present our results in more detail, we will
summarize what was previously known about 
shift-invariant closed subspaces. While doing that,
we introduce some additional terminology. We start with full
translation-invariance.

A closed subspace $V$ of $L^2(\mathbb{R}^n)$ is called \textit{translation-invariant} if $V$ is invariant under translation by elements of $\mathbb{R}^n$; i.e., $f\in V$ implies $T_y f \in V$, for all $y\in\mathbb{R}^n$, where $T_y f(x) = f(x-y)$,\label{translation} for a.e. $x\in\mathbb{R}^n$. Translation-invariant
closed subspaces can be characterized using the Fourier
transform. For $f\in L^1(\mathbb{R}^n)$, the Fourier transform
$\widehat{f}$ is defined by
\[
\widehat{f}(\xi)=\int_{\mathbb{R}^n}f(x)e^{2\pi i\xi\cdot x}dx,
\]
for $\xi\in\mathbb{R}^n$. The Fourier transform is an isometry
as a map from $L^1(\mathbb{R}^n)\cap L^2(\mathbb{R}^n)$ into $L^2(\mathbb{R}^n)$
and extends to a unitary map $\mathcal{F}$
on $L^2(\mathbb{R}^n)$. We will often write $\widehat{f}$ instead
of $\mathcal{F}f$, for $f\in L^2(\mathbb{R}^n)$. We note that 
\[
\widehat{T_yf}=E_y\widehat{f}, \text{ for }f\in L^2(\mathbb{R}^n),
y\in\mathbb{R}^n,
\]
where $E_y$ is the modulation operator given by
$E_yg(\xi)=e^{2\pi i\xi\cdot y}g(\xi)$, for a.e. $\xi\in\mathbb{R}^n$.
This enables the identification of some closed subspaces
of $L^2(\mathbb{R}^n)$ that are translation-invariant. For
any measurable subset $A$ of $\mathbb{R}^n$, let
\[
M_A=\{f\in L^2(\mathbb{R}^n): \widehat{f}(\xi)=0 
\text{ for a.e. }\xi\in\mathbb{R}^n\setminus A\}.
\]
Clearly, the property of vanishing off of $A$ is not 
altered by any modulation operator, so $M_A$ is a
translation-invariant closed subspace of $L^2(\mathbb{R}^n)$. The following theorem, attributed to Wiener, is proved in \cite{Rudin}, pp.~188--190.
\begin{tcolorbox}
\begin{theorem*}\label{Wiener}
	The closed, translation-invariant subspaces of $L^2(\mathbb{R}^n)$ are precisely the sets $M_A$. Furthermore, $M_A = M_B$ if and only if the symmetric difference, 
	$A \triangle B$, is a null set.
\end{theorem*}
\end{tcolorbox}
Much of the pioneering theory of shift-invariant subspaces can be found in Helson's \textit{Lectures on Invariant Subspaces}\cite{Helson}, although he proves everything in terms of modulation-invariance (\cite{Helson}, pp.~7, 38, and 57--60). Helson's work makes heavy use of \textit{range functions}. A range function is a mapping\label{ntorus} 
\[
J: \mathbb{T}^n \to \{\text{closed subspaces of}\ \mathscr{H}\},
\]
where $\mathscr{H}$ is a separable Hilbert space. Much of this work was later adapted by Marcin Bownik \cite{Bownik}, where the theory is expressed explicitly in terms of translation (as opposed to modulation), and where the implications of the ideas in \cite{Helson} are clarified. For this, Bownik takes range functions to be functions from $\mathbb{T}^n$ into the collection of closed subspaces of $l^2(\mathbb{Z}^n)$, and defines an isometric isomorphism $\mathscr{T}: L^2(\mathbb{R}^n) \to L^2(\mathbb{T}^n;l^2(\mathbb{Z}^n))$ by 
\[
\mathscr{T}f(x) = (\widehat{f}(x+k))_{k\in\mathbb{Z}^n};
\]
this map is often called a \textit{fiberization mapping}. The above equation may also be expressed component-wise:
\[
[\mathscr{T}f(x)](k) = \widehat{f}(x+k).  
\]
The reader will find that our work makes use of a similar (albeit more complicated) isometric isomorphism. An example of the usefulness of the fiberization technique is the following theorem, which is Proposition 1.5 in \cite{Bownik}.
\begin{tcolorbox}
\begin{named}{Theorem}\label{Bownik1.5}
	A closed subspace $V\subset L^2(\mathbb{R}^n)$ is shift invariant (that is, invariant under $T_k$ for all $k\in \mathbb{Z}^n)$ if and only if 
	\begin{equation}\label{Helsonshift}
		V = \{f\in L^2(\mathbb{R}^n): \mathscr{T}f(x) \in J(x)\ \text{for a.e.}\ x\in\mathbb{T}^n\},
	\end{equation}
	where $J$ is a (measurable) range function. The correspondence between $V$ and $J$ is one-to-one (identifying range functions that are equal a.e.). There is also a 1-1 correspondence between shift-invariant closed subspaces of $L^2(\mathbb{R}^n)$ and subspaces $W$ of $L^2(\mathbb{T}^n;l^2(\mathbb{Z}^n))$ that are closed under multiplication by exponentials, i.e.:
	\[
	x\mapsto \Phi(x) \in W \implies x\mapsto e^{-2\pi i x \cdot k}\Phi(x) \in W \ \ \ \text{for all}\ k\in \mathbb{Z}^n.
	\]
\end{named}
\end{tcolorbox}
In his paper \cite{Bownik}, Bownik goes on to investigate properties of shift-invariant closed subspaces that are important in
wavelet analysis. He also gives a decomposition of shift-invariant subspaces as an orthogonal sum, which allows him to prove a characterization of operators that commute with shifts. In 2015, Bownik and Kenneth Ross in \cite{BownikRoss} generalized the above theorem to the setting of LCA groups. They consider shifts by closed, co-compact subgroups $H\subset G$, where $G$ is a locally compact abelian group. In the case where $G=\mathbb{R}^n$, the results of \cite{BownikRoss} recover the above theorems by taking $H=\mathbb{R}^n$ and $H=\mathbb{Z}^n$, respectively; but \cite{BownikRoss} also yields the ``in-between'' cases where $H = P(\mathbb{Z}^k \times \mathbb{R}^{n-k})$ for $P\in \text{GL}_n(\mathbb{R})$. We note that some of the techniques in Bownik and Ross's paper inspired the author in this thesis. For example, the technique for constructing the range functions is shown in Theorem 2.4 of \cite{BownikRoss} (this technique traces back to Theorem 8 in Helson \cite{Helson}, p.~59; Helson credits T.P. Srinivasan with the proof). Moreover, in Chapter \ref{OrthConditionChapter}, the author inherited the idea of using \textit{determining sets}, and the technique for showing that certain functions form a determining set (Lemma 3.5 in \cite{BownikRoss}). Bownik and Ross also mention briefly that any translation-invariant subspace of $L^2(\mathbb{R}^n)$ is isomorphic (via the fiberization map) to a space which can be expressed as a \textit{direct integral}. The author---at his supervisor's suggestion---pursued this idea in this thesis. It yields a more powerful and versatile approach to the theory of translation-invariant spaces.

Finally, we mention that Benjamin Manning, in his 2012 Ph.D. thesis (\cite{Manning}) shows that $\Gamma$-invariant subspaces have an orthogonal decomposition with certain defining properties. Using the Plancherel transform defined in Taylor (1989) (\cite{Cstarcrystal}, p.~513), Manning defines an operator-valued Plancherel transform on functions on the crystal group $\Gamma$. He uses this, as well as the dimension function and a \textit{bracket}, which resembles an inner product on the range of a fiberization-type mapping, to show existence of an orthonormal basis for $\Gamma$-invariant subspaces. From this, Manning gets the existence of multiresolution analyses, and gives examples of crystallographic Daubechies-type wavelets. Manning's thesis only deals with symmorphic crystal groups, and does not yield an explicit description of crystallographic-invariant subspaces. But it may be fruitful to explore the connections between these two approaches.

	\section{Examples}\label{Examples}
	
	In this section we see how our main result works in the case a few simple two-dimensional examples. Any undefined terminology will be defined in the body of the thesis.
	
	\begin{ex}\label{Intropgex}
		Consider the \textit{nonsymmorphic} wallpaper group $\Gamma = pg$, which is described in Example \ref{pgexample}; see the figure below for a pattern illustrating its symmetries.
		\begin{figure}[H]
			\centering
			\begin{tikzpicture}[scale=.8]
				\foreach \i [evaluate=\i as \y using \i]  in {0,2}
				\foreach \x in {-4,-2,...,4}
				{
					\shade[thick, lower right = white!90!black, upper left = white!40!black] (\x, \y) -- (\x+2, \y) -- (\x+2,\y+1) -- (\x, \y+1);		
					\draw[line width=1.2pt] (\x, \y) -- (\x+2, \y) -- (\x+2,\y+1) -- (\x, \y+1) -- (\x,\y) -- (\x+2,\y);
					\draw[line width=1.2pt] (\x+.8, \y+.3) -- (\x+.8, \y+.7) -- (\x+1.2,\y+.7); 
				}
				\foreach \i [evaluate=\i as \y using \i]  in {1,3}
				\foreach \x in {-4,-2,...,4}
				{
					\shade[thick, lower left = white!40!black, upper right = white!90!black]  (\x+1, \y) -- (\x+3, \y) -- (\x+3,\y+1) -- (\x+1, \y+1);	
					\draw[line width=1.2pt] (\x+1, \y) -- (\x+3, \y) -- (\x+3,\y+1) -- (\x+1, \y+1) -- (\x+1,\y) -- (\x+3,\y);	
					\draw[line width=1.2pt] (\x+1.8,\y+.7) -- (\x+1.8, \y+.3) -- (\x+2.2, \y+.3); 
				}		
			\end{tikzpicture}
			\caption{A pattern illustrating the symmetries of the wallpaper group $pg$.}
		\end{figure}
		We note that the point group $\Pi$ for $pg$ consists of just the identity element, $\mathrm{id}$, and a vertical reflection, $\sigma$. In Example \ref{pgcross}, we describe a \textit{cross-section} (a map which selects one representative from each equivalence class) for the equivalence relation on $\Gamma$ which identifies elements in the same fiber over the point group. With this cross-section and with $\mathrm{T}$ the group of translations in $\Gamma$, we have
		\[
		pg = \gamma(\text{id}) \mathrm{T} \cup \gamma(\sigma)\mathrm{T}.
		\]
		For simplicity, assume that the wallpaper pattern is arranged and scaled so that the translation lattice $\mathcal{L}$ is just $\mathbb{Z}^2$; then the dual lattice $\mathcal{L}^*$ is also $\mathbb{Z}^2$. Also, the primitive translation vectors are just the standard basis vectors in this case. Thus we can take $\gamma(\mathrm{id}) = [0,\mathrm{id}]$ and $\gamma(\sigma) = [(\tfrac12,0),\sigma]$. Moreover, 
		\[
		\Gamma^*\cong \mathcal{L}^*\rtimes \Pi = \mathbb{Z}^2 \rtimes \Pi,
		\]
		with translation subgroup denoted by $\mathrm{T}^*$. We can choose a fundamental set $\Omega_{\mathrm{T}^*}$ for $\mathrm{T}^*$ to be 
		\[
		\Omega_{\mathrm{T}^*} =\left[-\tfrac12, \tfrac12\right) \times \left[-\tfrac12, \tfrac12\right);
		\]
		and we can choose a fundamental set $\Omega$ for $\Gamma^*$ to be 
		\[
		\Omega = \left[ -\tfrac12, \tfrac12\right) \times \left[-\tfrac12, 0 \right].
		\]
		We illustrate these in the figure below.
		\begin{figure}[H]
			\pgfplotsset{
				axis lines=center, 
				grid,
				grid style={very thin, gray},
				xmin=-4,    xmax=4,   xtick={-4,...,4},
				x tick style={draw=none}, 
				ymin=-4,    ymax=4,   
				y tick style={draw=none},   
				x tick style={draw=none},  
				ytick={-4,...,4},
				restrict y to domain=-4:4, 
				ticklabel style={font=\scriptsize, fill=white, inner sep=2pt},
				domain=-1:11, samples=100,
				no marks, 
				every axis plot post/.append style={ultra thick,},
			}
			\centering
			\begin{tikzpicture}[scale=1.2]
				\begin{axis}[thick,
					xticklabels={-2,,-1,,0,,1,,2},
					yticklabels={-2,,-1,,0,,1,,2},
					xlabel={$x$},
					ylabel={$y$},
						every axis x label/.style={
						at={(ticklabel* cs:1.02)},
						anchor=west,
					},
					every axis y label/.style={
						at={(ticklabel* cs:1.02)},
						anchor=south,
					},
					]
					\fill[line width=1pt, color={rgb: red, 3; blue, 1}, opacity=.35] (-1, -1) -- (1,-1) -- (1,0) -- (-1,0) -- (-1,-1); 
					\draw[line width=1.2pt, color={rgb: red, 1; blue, 1}, opacity=1] (-1,1) -- (-1,-1) -- (1,-1);
					\draw[line width=1.2pt, color={rgb: red, 3; blue, 1}] (-1,0) -- (-1,-1) -- (1,-1);
					\draw[very thin, gray] (1, -1) -- (1,1);
					\draw[very thin, gray] (-1, 1) -- (1,1);
					\draw[line width = 0pt, pattern=north west lines, pattern color={rgb: red, 1; blue, 1}] (-1,-1) rectangle (1,1);
					\node[color=black, anchor= west] at (1,1) {\scriptsize $\Omega_{\mathrm{T}^*}$};
					\node[color=black] at (.5,-.5) {\scriptsize $\Omega$};		
					\draw[line width=1.2pt, color={rgb: red, 3; blue, 1}] (-1, 0) -- (1,0);
					\node[outer sep=0pt,circle, draw={rgb: red, 1; blue, 1}, thin, fill=white, inner sep=1pt] at (-1,1) {};
					\node[outer sep=0pt,circle, draw={rgb: red, 3; blue, 1}, thin, fill=white, inner sep=1pt] at (1,0) {};
					\node[outer sep=0pt,circle, draw={rgb: red, 1; blue, 1}, thin, fill=white, inner sep=1pt] at (1,-1) {};
					\node[outer sep=0pt,circle, draw={rgb: red, 3; blue, 1}, thin, fill=white, inner sep=1pt] at (1,-1) {};
				\end{axis}
			\end{tikzpicture}
			\caption{The fundamental set $\Omega_{\mathrm{T}^*}=\left[-\tfrac12, \tfrac12\right) \times \left[-\tfrac12, \tfrac12\right)$ for $\mathrm{T}^*$; and the fundamental set $\Omega= \left[ -\tfrac12, \tfrac12\right) \times \left[-\tfrac12, 0 \right]$ for $\Gamma^*$.}
		\end{figure}
		Another set of importance to us is the set $\Omega_0$ of \textit{free points} in $\Omega$, which satisfies 
		\begin{equation}\label{Omega0intro}
		\Omega_0 = \{\omega \in\Omega :  L\cdot_{\idref{12}}\omega =\omega \implies L=\mathrm{id} \}.
		\end{equation}
	Thus $\Omega_0$ is the set of points in $\Omega$ whose stabilizer under $\cdot_{\idref{12}}$ is trivial. See \S\ref{PionOmegaTStar} for more information about the action $\cdot_{\idref{12}}$ of $\Pi$ on $\Omega_{\mathrm{T}^*}$. Note that in this example, with the choice of $\Omega_{\mathrm{T}^*}$ and $\Omega$ above, we have
		\[
		\Omega_0 = \left[ -\tfrac12, \tfrac12\right) \times \left(-\tfrac12, 0 \right).
		\]
The fundamental domain $R$ for $\Gamma^*$---as discussed in \S\ref{FundRegionConstruction}---is of special important to us: $\pi$ will be shown to be equivalent to a direct integral of multiples of irreducibles $U^r$ over $R$. The other sets discussed so far are, in fact, chosen after $R$ is constructed. In the present example, the fundamental domain $R$ for $\Gamma^*$ can be chosen to be the open rectangle 
\[
R = \left( -\tfrac12, \tfrac12\right) \times \left(-\tfrac12, 0 \right).
\]
We illustrate $\Omega_0$ and $R$ below:
		\begin{figure}[H]
		\begin{subfigure}{0.5\textwidth}
			\pgfplotsset{
				axis lines=center, 
				grid,
				grid style={very thin, gray},
				xmin=-4,    xmax=4,   xtick={-4,...,4},
				x tick style={draw=none}, 
				ymin=-4,    ymax=4,   
				y tick style={draw=none},   
				x tick style={draw=none},  
				ytick={-4,...,4},
				restrict y to domain=-4:4, 
				ticklabel style={font=\scriptsize, fill=white, inner sep=2pt},
				domain=-1:11, samples=100,
				no marks, 
				every axis plot post/.append style={ultra thick,},
			}
			\centering
		\begin{tikzpicture}[scale=.9]
			\begin{axis}[thick,
				xticklabels={-2,,-1,,0,,1,,2},
				yticklabels={-2,,-1,,0,,1,,2},
				xlabel={$x$},
				ylabel={$y$},
			every axis x label/.style={
				at={(ticklabel* cs:1.02)},
				anchor=west,
			},
			every axis y label/.style={
				at={(ticklabel* cs:1.02)},
				anchor=south,
			},
				]
				\fill[line width=1pt, color={rgb: red, 3; blue, 1}, opacity=.5] (-1, -1) -- (1,-1) -- (1,0) -- (-1,0) -- (-1,-1); 
				\draw[line width=1.2pt, color={rgb: red, 3; blue, 1}] (-1,0) -- (-1,-1);
				\draw[very thin, gray] (1, -1) -- (1,1);
				\draw[very thin, gray] (-1, 1) -- (1,1);
				\draw[very thin, gray] (-1, -1) -- (1,-1);
				\node[color=black] at (.5,-.5) {\scriptsize $\Omega_0$};		
				\node[outer sep=0pt,circle, draw={rgb: red, 3; blue, 1}, thin, fill=white, inner sep=1pt] at (-1,-1) {};
				\node[outer sep=0pt,circle, draw={rgb: red, 3; blue, 1}, thin, fill=white, inner sep=1pt] at (-1,0) {};
			\end{axis}
		\end{tikzpicture}
			\caption{}
		\end{subfigure}%
		\begin{subfigure}{0.5\textwidth}
			\pgfplotsset{
				axis lines=center, 
				grid,
				grid style={very thin, gray},
				xmin=-4,    xmax=4,   xtick={-4,...,4},
				x tick style={draw=none}, 
				ymin=-4,    ymax=4,   
				y tick style={draw=none},   
				x tick style={draw=none},  
				ytick={-4,...,4},
				restrict y to domain=-4:4, 
				ticklabel style={font=\scriptsize, fill=white, inner sep=2pt},
				domain=-1:11, samples=100,
				no marks, 
				every axis plot post/.append style={ultra thick,},
			}
			\centering
			\begin{tikzpicture}[scale=.9]
			\begin{axis}[thick,
				xticklabels={-2,,-1,,0,,1,,2},
				yticklabels={-2,,-1,,0,,1,,2},
				xlabel={$x$},
				ylabel={$y$},
				every axis x label/.style={
					at={(ticklabel* cs:1.02)},
					anchor=west,
				},
				every axis y label/.style={
					at={(ticklabel* cs:1.02)},
					anchor=south,
				},
				]
				\fill[line width=1pt, color=RoyalBlue, opacity=.5] (-1, -1) -- (1,-1) -- (1,0) -- (-1,0) -- (-1,-1); 
				\draw[very thin, gray] (1, -1) -- (1,1);
				\draw[very thin, gray] (-1, -1) -- (-1,0);
				\draw[very thin, gray] (-1, 1) -- (1,1);
				\draw[very thin, gray] (-1, -1) -- (1,-1);
				\node[color=black] at (.5,-.5) {\scriptsize $R$};		
			\end{axis}
		\end{tikzpicture}
			\caption{}
		\end{subfigure}
		\caption{(a) The set $\Omega_0= \left[ -\tfrac12, \tfrac12\right) \times \left(-\tfrac12, 0 \right)$. (b) The fundamental domain $R = \left( -\tfrac12, \tfrac12\right) \times \left(-\tfrac12, 0 \right)$ for $\Gamma^*$.}
	\end{figure}
	
Suppose now that $V \subset L^2(\mathbb{R}^2)$ is invariant under shifts by elements of $\Gamma$. In other words, if $f\in V$ and $[x,L]\in pg$, then $f\big([x,L]^{-1}\cdot_{\idref{1}}y\big) \in V$. By Theorem \ref{BigInvTheorem2}, there exists a measurable range function $r \mapsto K_r$ from $R$ into the closed subspaces of $l^2(\mathbb{Z}^2)$, uniquely determined almost everywhere by $V$, such that $V$ consists of all $f \in L^2(\mathbb{R}^2)$ satisfying
		\[
		\widehat{f}(\chi_{[z,M]\cdot_{\idref{1}}r}) = \chi_z(x_M)\langle F(r),\delta_M \otimes \delta_{z}\rangle
		\]
		for a.e. $r \in R$ and all $[z,M] \in \Gamma^*$, for some $F\in \int_R^\oplus l^2(\Pi)\otimes K_r \ud r$. In this equation, the direct integral is regarded as a subspace of $L^2(R;l^2(\Pi)\otimes l^2(\mathcal{L}^*))$; $\widehat{f}$ denotes taking the Plancherel transform of $f$; $\chi_{[z,M]\cdot_{\idref{1}} r}$ denotes a continuous character\footnote{i.e.: $\chi_y(u) := e^{2\pi i y\cdot u}$ for all $u\in \mathbb{R}^n$.} on $\mathbb{R}^2$; and $x_M$ comes from the cross-section, $\gamma(M) = [x_M,M]$. Thus, writing $z \in \mathbb{Z}^2$ as $z=(z_1,z_2)$, the above equation becomes
		\begin{align}
			\widehat{f}(\chi_{[z,M]\cdot_{\idref{1}}r})	&= e^{2\pi i z\cdot x_M}\langle F(r),\delta_M \otimes \delta_{z}\rangle \nonumber \\
			&=  \left\{ \begin{array}{ll}
				e^{2\pi i \tfrac12 z_1} \langle F(r),\delta_\sigma \otimes \delta_{z}\rangle &\text{if } M=\sigma \nonumber \\
				\langle F(r),\delta_{\mathrm{id}} \otimes \delta_{z}\rangle &\text{if } M=\mathrm{id}
			\end{array} \right. \nonumber  \\
&	= \left\{ \begin{array}{ll}
				e^{\pi i z_1} \langle F(r),\delta_\sigma \otimes \delta_{z}\rangle &\text{if } M=\sigma \\
				\langle F(r),\delta_{\mathrm{id}}  \otimes \delta_{z}\rangle &\text{if } M=\mathrm{id}
			\end{array} \right. \nonumber \\
		&=	\left\{ \begin{array}{ll}
				\langle F(r),\delta_\sigma \otimes \delta_{z}\rangle &\text{if } M=\sigma\ \text{and}\ z_1\ \text{even}\\
				-\langle F(r),\delta_\sigma \otimes \delta_{z}\rangle &\text{if } M=\sigma\ \text{and}\ z_1\ \text{odd}\\
				\langle F(r),\delta_{\mathrm{id}}  \otimes \delta_{z}\rangle &\text{if } M=\mathrm{id}
			\end{array} \right. 	\label{pgexintro}
		\end{align}
		for a.e. $r\in R$ and all $z\in \mathbb{Z}^2$. Since 
		\[
		\{\chi_{[z,M]\cdot_{\idref{1}}r}: r\in R, [z,M] \in \Gamma^*\}
		\]
		is almost all of $\widehat{\mathbb{R}^2}$, equation \eqref{pgexintro} defines $\widehat{f}$ almost everywhere, and hence as an element of $L^2(\widehat{\mathbb{R}^2})$. Therefore, $V$ is invariant under shifts by $pg$ if and only if there is some measurable range function $r \mapsto K_r$ from $R$ into the closed subspaces of $l^2(\mathbb{Z}^2)$, uniquely determined almost everywhere by $V$, such that $V$ consists of all $f \in L^2(\mathbb{R}^2)$ satisfying equation \eqref{pgexintro} for a.e. $r\in R$ and all $[z,M] \in \Gamma^*$, for some $F\in \int_R^\oplus l^2(\Pi)\otimes K_r \ud r$. We shall illustrate this example in more detail in section \ref{Gammapgex}.
	\end{ex}
	
	\begin{ex}\label{Introexp1}
		Consider the simplest of all the wallpaper groups $\Gamma = p1 \cong \{\mathrm{id}\}\rtimes \mathcal{L}$, which consists of only translations. This wallpaper group is \textit{symmorphic}. The figure below shows a pattern illustrating its symmetry.
		
		\begin{figure}[H]
			\centering
			\begin{tikzpicture}[scale=0.7]
				\foreach \i [evaluate=\i as \y using \i]  in {0,2,4,6}
				\foreach \x in {-4,-2,...,4}
				{
					\shade[thick, upper left = white!30!black] (\x, \y) -- (\x+2, \y) -- (\x+2,\y+2) -- (\x, \y+2);		
					\draw[line width=1.2pt] (\x, \y) -- (\x+2, \y) -- (\x+2,\y+2) -- (\x, \y+2) -- (\x,\y) -- (\x+2,\y);
					\draw[line width=1.2pt] (\x+.8, \y+.8) -- (\x+.8, \y+1.2) -- (\x+1.2,\y+1.2); 
				}
			\end{tikzpicture}
			\caption{A pattern illustrating the symmetry of the wallpaper group $p1$; this group has only translation-symmetry.}
		\end{figure}
		We note that the point group $\Pi$ for $p1$ consists of just the identity element. Thus the cross-section is trivial: $\gamma(\mathrm{id}) = [0,\mathrm{id}]$; i.e., $x_{\mathrm{id}} = 0$. For simplicity, assume that the wallpaper pattern is arranged and scaled so that the translation lattice $\mathcal{L}$ is just $\mathbb{Z}^2$; then the dual lattice $\mathcal{L}^*$ is also $\mathbb{Z}^2$. Also, the primitive translation vectors are just the standard basis vectors in this case. We have $\Gamma^*\cong \mathcal{L}^*\rtimes \{\mathrm{id}\}$ and the translation subgroup $\mathrm{T}^*$ coincides with $\Gamma^*$ in this case. We can choose a fundamental set $\Omega_{\mathrm{T}^*}$ for $\mathrm{T}^*$ to be 
		\[
		\Omega_{\mathrm{T}^*} =\left[-\tfrac12, \tfrac12\right) \times \left[-\tfrac12, \tfrac12\right). 
		\]
		Note that $\Omega = \Omega_{\mathrm{T}^*}$ is of course a fundamental set for $\Gamma^*$ as well. We note that for this wallpaper group, the set $\Omega_0$ of points in $\Omega$ with trivial stabilizer coincides with $\Omega$, since the point group is trivial. Moreover, for this example, a fundamental domain $R$ for $\Gamma^*$ is given by the open square
		\[
		R = \left(-\tfrac12, \tfrac12\right) \times \left(-\tfrac12, \tfrac12\right)
		\]
		We illustrate $\Omega_0$ and $R$ in the figure below:
			\begin{figure}[H]
			\begin{subfigure}{0.5\textwidth}
				\pgfplotsset{
					axis lines=center, 
					grid,
					grid style={very thin, gray},
					xmin=-4,    xmax=4,   xtick={-4,...,4},
					x tick style={draw=none}, 
					ymin=-4,    ymax=4,   
					y tick style={draw=none},   
					x tick style={draw=none},  
					ytick={-4,...,4},
					restrict y to domain=-4:4, 
					ticklabel style={font=\scriptsize, fill=white, inner sep=2pt},
					domain=-1:11, samples=100,
					no marks, 
					every axis plot post/.append style={ultra thick,},
				}
				\centering
			\begin{tikzpicture}[scale=.9]
				\begin{axis}[thick,
					xticklabels={-2,,-1,,0,,1,,2},
					yticklabels={-2,,-1,,0,,1,,2},
					xlabel={$x$},
					ylabel={$y$},
				every axis x label/.style={
					at={(ticklabel* cs:1.02)},
					anchor=west,
				},
				every axis y label/.style={
					at={(ticklabel* cs:1.02)},
					anchor=south,
				},
					]
					\draw[line width=1.2pt, color={rgb: red, 1; blue, 1}] (-1,1) -- (-1,-1) -- (1,-1);
					\draw[line width=1.2pt, color={rgb: red, 3; blue, 1}] (-1,1) -- (-1,-1) -- (1,-1);
					\draw[very thin, gray] (1, -1) -- (1,1);
					\draw[very thin, gray] (-1, 1) -- (1,1);
					\draw[line width = 0pt, pattern=north west lines, pattern color={rgb: red, 1; blue, 1}] (-1,-1) rectangle (1,1);
					\fill[line width=1pt, color={rgb: red, 3; blue, 1}, opacity=.35] (-1, -1) -- (1,-1) -- (1,1) -- (-1,1) -- (-1,-1);
					\node[color=black, anchor= west] at (1,1) {\scriptsize $\Omega_0=\Omega =\Omega_{\mathrm{T}^*}$};		
					\node[outer sep=0pt,circle, draw={rgb: red, 1; blue, 1}, thin, fill=white, inner sep=1pt] at (1,-1) {};
					\node[outer sep=0pt,circle, draw={rgb: red, 1; blue, 1}, thin, fill=white, inner sep=1pt] at (-1,1) {};
					\node[outer sep=0pt,circle, draw={rgb: red, 3; blue, 1}, thin, fill=white, inner sep=1pt] at (1,-1) {};
					\node[outer sep=0pt,circle, draw={rgb: red, 3; blue, 1}, thin, fill=white, inner sep=1pt] at (-1,1) {};
				\end{axis}
			\end{tikzpicture}
				\caption{}
			\end{subfigure}%
			\begin{subfigure}{0.5\textwidth}
				\pgfplotsset{
					axis lines=center, 
					grid,
					grid style={very thin, gray},
					xmin=-4,    xmax=4,   xtick={-4,...,4},
					x tick style={draw=none}, 
					ymin=-4,    ymax=4,   
					y tick style={draw=none},   
					x tick style={draw=none},  
					ytick={-4,...,4},
					restrict y to domain=-4:4, 
					ticklabel style={font=\scriptsize, fill=white, inner sep=2pt},
					domain=-1:11, samples=100,
					no marks, 
					every axis plot post/.append style={ultra thick,},
				}
				\centering
				\begin{tikzpicture}[scale=.9]
					\begin{axis}[thick,
						xticklabels={-2,,-1,,0,,1,,2},
						yticklabels={-2,,-1,,0,,1,,2},
						xlabel={$x$},
						ylabel={$y$},
						every axis x label/.style={
							at={(ticklabel* cs:1.02)},
							anchor=west,
						},
						every axis y label/.style={
							at={(ticklabel* cs:1.02)},
							anchor=south,
						},
						]
						\draw[very thin, gray] (1, -1) -- (1,1);
						\draw[very thin, gray] (-1, 1) -- (1,1);
						\draw[very thin, gray] (1, -1) -- (-1,-1) -- (-1,1);
						\fill[line width=1pt, color=RoyalBlue, opacity=.5] (-1, -1) -- (1,-1) -- (1,1) -- (-1,1) -- (-1,-1);
						\node[color=black] at (.5,.5) {\scriptsize $R$};		
					\end{axis}
				\end{tikzpicture}
				\caption{}
			\end{subfigure}
			\caption{(a) The fundamental set $\Omega=\Omega_{\mathrm{T}^*}=\left[-\tfrac12, \tfrac12\right) \times \left[-\tfrac12, \tfrac12\right)$ for $\mathrm{T}^*=\Gamma^*$. We note that $\Omega_0 = \Omega$. (b) The fundamental domain $R = \left( -\tfrac12, \tfrac12\right) \times \left(-\tfrac12, \tfrac12 \right)$ for $\Gamma^*$.}
		\end{figure}
		Suppose now that $V \subset L^2(\mathbb{R}^2)$ is invariant under shifts by elements of $\Gamma$. In other words, if $f\in V$ and $[x,\mathrm{id}]\in p1$, then $f\big([x,\mathrm{id}]^{-1}\cdot_{\idref{1}}y\big) = f(y-x)\in V$. \textit{Note that in the case $\Gamma=p1$, the notion of $\Gamma$-invariance coincides with the notion of shift-invariance described above.} We have from Theorem \ref{BigInvTheorem2} that there exists a measurable range function $r \mapsto K_r$ from $R$ into the closed subspaces of $l^2(\mathbb{Z}^2)$, uniquely determined almost everywhere by $V$, such that $V$ consists of all $f \in L^2(\mathbb{R}^2)$ satisfying
		\[
		\widehat{f}(\chi_{[z,\mathrm{id}]\cdot_{\idref{1}}r}) = \chi_z(x_\mathrm{id})\langle F(r),\delta_\mathrm{id} \otimes \delta_{z}\rangle,
		\]
		for a.e. $r\in R$ and all $z\in \mathbb{Z}^2$, for some $F\in \int_{R}^\oplus l^2(\{\mathrm{id}\})\otimes K_r \ud r$. Since $\chi_z(0)= 1$, the above equations says
		\begin{equation}\label{p1ex}
			\widehat{f}(\chi_{r+z})=\langle F(r),\delta_\mathrm{id} \otimes \delta_{z}\rangle,
		\end{equation}
		for a.e. $r\in R$ and all $z\in \mathbb{Z}^2$. Since 
		\[
		\{\chi_{r+z}: r\in R, z \in \mathbb{Z}^2\}
		\]
		is almost all of $\widehat{\mathbb{R}^2}$, equation \eqref{p1ex} defines $\widehat{f}$ almost everywhere, and hence as an element of $L^2(\widehat{\mathbb{R}^2})$. Therefore, $V$ is invariant under shifts by $p1$ if and only if there is some measurable range function $r \mapsto K_r$ from $R$ into the closed subspaces of $l^2(\mathbb{Z}^2)$, uniquely determined almost everywhere by $V$, such that $V$ consists of all $f \in L^2(\mathbb{R}^2)$ satisfying equation \eqref{p1ex} for a.e. $r\in R$ and all $z\in \mathbb{Z}^2$, for some $F\in \int_R^\oplus l^2(\{\mathrm{id}\})\otimes K_r \ud r$. By identifying $L^2\big(R;l^2(\{\mathrm{id}\})\otimes l^2(\mathbb{Z}^2)\big)$ with $L^2\big(\mathbb{T}^n;l^2(\mathbb{Z}^2)\big)$, we recover the well-known characterization described in the Theorem on page \pageref{Bownik1.5}, with $\mathscr{T}f=F \in \int_{R}^\oplus K_r \ud r$.
	\end{ex}

	\begin{ex}\label{Introexp4m}
		We now consider a symmorphic wallpaper group with more symmetries, namely $\Gamma = p4m$. This group is discussed in greater detail in Example \ref{p4m}. The figure below shows a pattern illustrating its symmetries.
		\begin{figure}[H]
			\centering
			\begin{tikzpicture}[scale=.8]
				\foreach \i [evaluate=\i as \y using 2*\i]  in {0,1,...,4}
				\foreach \x in {-4,-2,...,4}
				{
					\fill[line width=1.2pt, color=gray!10] (\x,\y) --(\x+2,\y) -- (\x+2,\y+ 2) -- (\x,\y+2) -- (\x,\y) -- (\x+2,\y); 
					\shade[thick, lower left = black!50!white, upper right = black!50!white] (\x+1, \y+1) -- (\x+2, \y+1) -- (\x+2,\y+2) -- (\x+1, \y+1)-- (\x+2, \y+1);		
					\shade[thick, lower left = black!50!white, upper right = black!50!white] (\x+1, \y+1) -- (\x+1, \y+2) -- (\x+2,\y+2) -- (\x+1, \y+1)-- (\x+1, \y+2);
					\shade[thick, lower right = black!50!white, upper left = black!50!white] (\x+1, \y+1) -- (\x+1, \y+2) -- (\x,\y+2) -- (\x+1, \y+1)-- (\x+1, \y+2);
					\shade[thick, lower right = black!50!white, upper left = black!50!white] (\x+1, \y+1) -- (\x, \y+2) -- (\x,\y+1) -- (\x+1, \y+1)-- (\x, \y+2);
					\shade[thick, upper right = black!50!white, lower left = black!50!white] (\x+1, \y+1) -- (\x, \y+1) -- (\x,\y) -- (\x+1, \y+1)-- (\x, \y+1);
					\shade[thick, upper right = black!50!white, lower left = black!50!white] (\x+1, \y+1) -- (\x, \y) -- (\x+1,\y) -- (\x+1, \y+1)-- (\x, \y);
					\shade[thick, upper left = black!50!white, lower right= black!50!white] (\x+1, \y+1) -- (\x+1, \y) -- (\x+2,\y) -- (\x+1, \y+1)-- (\x+1, \y);
					\shade[thick, upper left = black!50!white, lower right= black!50!white] (\x+1, \y+1) -- (\x+2, \y) -- (\x+2,\y+1) -- (\x+1, \y+1)-- (\x+2, \y);
					\draw[line width=1.2pt, color = black][-] (\x,\y) -- (\x+2,\y+2);
					\draw[line width=1.2pt, color = black][-] (\x,\y+2) -- (\x+2,\y); 
					\draw[line width=1.2pt, color = black][-] (\x+1,\y) -- (\x+1,\y+2);
					\draw[line width=1.2pt, color = black][-] (\x,\y+1) -- (\x+2,\y+1);
					\draw[line width=1.2pt, color = black][-] (\x,\y) -- (\x+2,\y); 
					\draw[line width=1.2pt, color = black][-] (\x,\y) -- (\x,\y+2); 
					\draw[line width=1.2pt, color = black][-] (-4,0) -- (-4,10) -- (6,10) -- (6,0) -- (-4,0);
					\draw[line width=1.2pt, color = black][-] (\x+1.5,\y+1.2) -- (\x+1.8,\y+1.2) -- (\x+1.8,\y+1.5);
					\draw[line width=1.2pt, color = black][-] (\x+.2,\y+1.5) -- (\x+.2,\y+1.2) -- (\x+.5,\y+1.2);
					\draw[line width=1.2pt, color = black][-] (\x+1.5,\y+1.8) -- (\x+1.2,\y+1.8) -- (\x+1.2,\y+1.5);
					\draw[line width=1.2pt, color = black][-] (\x+.5,\y+1.8) -- (\x+.8,\y+1.8) -- (\x+.8,\y+1.5);
					\draw[line width=1.2pt, color = black][-] (\x+.5,\y+.2) -- (\x+.8,\y+.2) -- (\x+.8,\y+.5);
					\draw[line width=1.2pt, color = black][-] (\x+.2,\y+.5) -- (\x+.2,\y+.8) -- (\x+.5,\y+.8);
					\draw[line width=1.2pt, color = black][-] (\x+1.5,\y+.8) -- (\x+1.8,\y+.8) -- (\x+1.8,\y+.5);
					\draw[line width=1.2pt, color = black][-] (\x+1.2,\y+.5) -- (\x+1.2,\y+.2) -- (\x+1.5,\y+.2);
					\draw[line width=1.2pt, color = black][-] (-4,0.1) -- (-4,0) -- (-3.8,0);
				}
			\end{tikzpicture}
			\caption{A pattern illustrating the symmetries of the wallpaper group $p4m$.}
		\end{figure}
		We note that the point group $\Pi$ for $p4m$ is $D_8$, the dihedral group of order 8. That is, 
		\[
		\Pi  = \langle r,s | r^4 = s^2 = \mathrm{id},\ rs = sr^{-1}\rangle,
		\]
		where 
		\[
		\mathrm{id} = \left[ \begin{array}{cc}
			1 & 0 \\
			0 & 1
		\end{array} \right], \ \ \
		r =  \left[ \begin{array}{cc}
			0 & -1 \\
			1 & 0
		\end{array} \right], \ \ \ \text{and}\ \ \ 
		s = \left[ \begin{array}{cc}
			0 & 1 \\
			1 & 0
		\end{array} \right].
		\]
		Since the group is symmorphic, the cross-section $\gamma$ can be chosen such that 
		\[
		\gamma(M) = [0,M]\ \ \ \text{for all}\ \ \ M\in \Pi;
		\] 
		that is, $x_M = 0$ for all $M$, where $[x_M,M] = \gamma(M)$. For simplicity, assume that the wallpaper pattern is arranged and scaled so that the translation lattice $\mathcal{L}$ is just $\mathbb{Z}^2$; then the dual lattice $\mathcal{L}^*$ is also $\mathbb{Z}^2$. Also, the primitive translation vectors are just the standard basis vectors in this case. Moreover, 
		\[
		\Gamma^*\cong \mathcal{L}^*\rtimes \Pi = \mathbb{Z}^2 \rtimes \Pi,
		\]
with translation subgroup denoted by $\mathrm{T}^*$. We can choose a fundamental set $\Omega_{\mathrm{T}^*}$ for $\mathrm{T}^*$ to be 
		\[
		\Omega_{\mathrm{T}^*} =\left[-\tfrac12, \tfrac12\right) \times \left[-\tfrac12, \tfrac12\right);
		\]
		and we can choose a fundamental set $\Omega$ for $\Gamma^*$ to be 
		\[
		\Omega = \{(x,y):-1 \leq x \leq 0,  x \leq y \leq 0\}.
		\]
	In this case, the set $\Omega_0$ and $R$ are just the interior of $\Omega$: i.e.: 
	\[
	\Omega_0 = R = \{(x,y): -1 < x < 0, x < y < 0\}.
	\]
We illustrate these sets in the figure below.
	\begin{figure}[H]
	\begin{subfigure}{0.5\textwidth}
			\pgfplotsset{
				axis lines=center, 
				grid,
				grid style={very thin, gray},
				xmin=-4,    xmax=4,   xtick={-4,...,4},
				x tick style={draw=none}, 
				ymin=-4,    ymax=4,   
				y tick style={draw=none},   
				x tick style={draw=none},  
				ytick={-4,...,4},
				restrict y to domain=-4:4, 
				ticklabel style={font=\scriptsize, fill=white, inner sep=2pt},
				domain=-1:11, samples=100,
				no marks, 
				every axis plot post/.append style={ultra thick,},
			}
				\centering
			\begin{tikzpicture}[scale=0.9]
				\begin{axis}[thick,
					xticklabels={-2,,-1,,0,,1,,2},
					yticklabels={-2,,-1,,0,,1,,2},
					xlabel={$x$},
					ylabel={$y$},
				every axis x label/.style={
					at={(ticklabel* cs:1.02)},
					anchor=west,
				},
				every axis y label/.style={
					at={(ticklabel* cs:1.02)},
					anchor=south,
				},
					]
					\draw[line width=1.2pt, color={rgb: red, 1; blue, 1}, opacity=1] (-1,1) -- (-1,-1) -- (1,-1);
						\draw[line width=1.2pt, color={rgb: red, 3; blue, 1}] (-1,-1) -- (0,0) -- (-1,0) -- (-1,-1);
					\draw[very thin, gray] (1, -1) -- (1,1);
					\draw[very thin, gray] (-1, 1) -- (1,1);
					\draw[line width = 0pt, pattern=north west lines, pattern color={rgb: red, 1; blue, 1}] (-1,-1) rectangle (1,1);
					\fill[line width=1pt, color={rgb: red, 3; blue, 1}, opacity=.35] (-1, 0) -- (-1,-1) -- (0,0) -- (-1,0);
					\node[color=black, anchor= west] at (1,1) {\scriptsize $\Omega_{\mathrm{T}^*}$};
					\node[color=black, anchor= west] at (-1.05,-.32) {\scriptsize $\Omega$};				
					\node[outer sep=0pt,circle, draw={rgb: red, 1; blue, 1}, thin, fill=white, inner sep=1pt] at (1,-1) {};
					\node[outer sep=0pt,circle, draw={rgb: red, 1; blue, 1}, thin, fill=white, inner sep=1pt] at (-1,1) {};
				\end{axis}
			\end{tikzpicture}
			\caption{}
		\end{subfigure}%
			\begin{subfigure}{0.5\textwidth}
				\pgfplotsset{
					axis lines=center, 
					grid,
					grid style={very thin, gray},
					xmin=-4,    xmax=4,   xtick={-4,...,4},
					x tick style={draw=none}, 
					ymin=-4,    ymax=4,   
					y tick style={draw=none},   
					x tick style={draw=none},  
					ytick={-4,...,4},
					restrict y to domain=-4:4, 
					ticklabel style={font=\scriptsize, fill=white, inner sep=2pt},
					domain=-1:11, samples=100,
					no marks, 
					every axis plot post/.append style={ultra thick,},
				}
					\centering
				\begin{tikzpicture}[scale=0.9]
					\begin{axis}[thick,
						xticklabels={-2,,-1,,0,,1,,2},
						yticklabels={-2,,-1,,0,,1,,2},
						xlabel={$x$},
						ylabel={$y$},
					every axis x label/.style={
						at={(ticklabel* cs:1.02)},
						anchor=west,
					},
					every axis y label/.style={
						at={(ticklabel* cs:1.02)},
						anchor=south,
					},
						]
						\draw[very thin, gray] (1, -1) -- (1,1);
						\draw[very thin, gray] (-1, 1) -- (1,1);
						\draw[very thin, gray] (-1, -1) -- (0,-1);
						\draw[very thin, gray] (-1, -1) -- (-1,0);
						\fill[line width=1pt, color=RoyalBlue, opacity=.5] (-1, 0) -- (-1,-1) -- (0,0) -- (-1,0);
						\node[color=black] at (.1,-.7) {\scriptsize $\Omega_0=R$};				
					\end{axis}
				\end{tikzpicture}
				\caption{}
			\end{subfigure}
				\caption{(a) The fundamental set $\Omega_{\mathrm{T}^*}=\left[-\tfrac12, \tfrac12\right) \times \left[-\tfrac12, \tfrac12\right)$ for $\mathrm{T}^*$; and the fundamental set $\Omega = \{(x,y): x \leq y \leq 0, -1 \leq x \leq 0\}$ for $\Gamma^*$. (b) The set $\Omega_0$ of points in $\Omega$ with trivial stabilizer, which coincides with the fundamental domain $R$ for $\Gamma^*$.}
		\end{figure}

		Suppose now that $V \subset L^2(\mathbb{R}^2)$ is invariant under shifts by elements of $\Gamma$. In other words, if $f\in V$ and $[x,L]\in p4m$, then $f\big([x,L]^{-1}\cdot_{\idref{1}}y\big) \in V$. By Theorem \ref{BigInvTheorem2}, there exists a measurable range function $r \mapsto K_r$ from $R$ into the closed subspaces of $l^2(\mathbb{Z}^2)$, uniquely determined almost everywhere by $V$, such that $V$ consists of all $f \in L^2(\mathbb{R}^2)$ satisfying
		\[
		\widehat{f}(\chi_{[z,M]\cdot_{\idref{1}}r}) = \chi_z(x_M)\langle F(r),\delta_M \otimes \delta_{z}\rangle 
		\]
		for a.e. $r\in R$ and all $[z,M] \in \Gamma^*$, for some $F\in \int_{R}^\oplus l^2(\{\mathrm{id}\})\otimes K_r \ud r$. Since 
		\[
		\chi_z(x_M) = \chi_z(0)= 1,
		\]
		the above equation becomes
		\begin{equation}\label{p4meq}
			\widehat{f}(\chi_{[z,M]\cdot_{\idref{1}}r})=\langle F(r),\delta_M \otimes \delta_{z}\rangle.
		\end{equation}
		Since 
		\[
		\{\chi_{[z,M]\cdot_{\idref{1}}r}: r\in R, [z,M] \in \Gamma^*\}
		\]
		is almost all of $\widehat{\mathbb{R}^2}$, equation \eqref{p4meq} defines $\widehat{f}$ almost everywhere, and hence as an element of $L^2(\widehat{\mathbb{R}^2})$. Therefore $V$ is invariant under shifts by $p4m$ if and only if there is some measurable range function $r \mapsto K_r$ from $R$ into the closed subspaces of $l^2(\mathbb{Z}^2)$, uniquely determined almost everywhere by $V$, such that $V$ consists of all $f \in L^2(\mathbb{R}^2)$ satisfying equation \eqref{p4meq} for a.e. $r \in R$ and all $[z,M] \in \Gamma^*$, for some $F\in \int_R^\oplus l^2(\Pi)\otimes K_r \ud r$.
	\end{ex}
	
	\section{Overview}
	
	In Chapter \ref{Preliminaries}, we define crystallographic groups: these are symmetry groups of repeating patterns in space. In dimension 2, these are the familiar wallpaper groups, and in dimension 3, these are the symmetries of crystal lattices. We introduce the point group $\Pi$ and the translation subgroup $\mathrm{T}$, and discuss symmorphic and non-symmorphic crystal groups. We introduce \textit{fundamental domains} and \textit{fundamental sets}---concepts essential to our study---and we construct a fundamental domain $R$ for $\Gamma^*$. Furthermore, we introduce the dual lattice $\mathcal{L}^*$ associated with a crystallographic group, and list the different group actions relevant to our study. Finally, we construct fundamental sets $\Omega_{\mathrm{T}^*}$ and $\Omega$, for $\mathrm{T}^*$ and $\Gamma^*$, respectively. The set $\Omega$ is shown to satisfy $m(\Omega-R) =0$.
	
	In Chapter \ref{Prerequisite}, we lay the theoretical foundations for the work of later chapters. This includes defining von Neumann algebras and unitary representations, and providing the reader with some basic facts about these objects. We also introduce the theory of direct integrals, and include a brief introduction to the essential tools of abstract Fourier analysis. We then describe the \textit{central decomposition} of a unitary representation. In the last section, we introduce range functions and show how any measurable range function $J: X \to \{\text{closed subspaces of}\ \mathscr{H}\}$ corresponds to a measurable field of Hilbert spaces $\{J(x): x\in X\}$. We show how to embed the direct integral of these $J(x)$ in $L^2(X;\mathscr{H})$. We then introduce the \textit{orthogonality condition} for any closed subspace of $L^2(X;\mathscr{H})$. We prove a key theorem which states that any closed subspace satisfying the orthogonality condition can be expressed as a direct integral of closed subspaces of $\mathscr{H}$, and moreover that the orthogonal projection onto this subspace decomposes as a direct integral of projections.
	
	In Chapter \ref{CrystalShiftandInduced}, we introduce the fundamental object of our study, the natural crystallographic shift representation $\pi$. We then introduce the concept of \textit{induced representations}, and use a theorem from \cite{TaylorKaniuth} to completely describe (up to unitary equivalence) all the irreducible representations of a crystallographic group, in terms of induced representations. We show that the irreducible representations are obtained by inducing from the stabilizer subgroup $\Gamma_{\omega}$ of the character $\chi_\omega^\mathrm{T}$, where $\omega\in\Omega$. We introduce a set $\Omega_0$ on which these stabilizer subgroups are just the translation subgroup $\mathrm{T}$, and such that $R\subset \Omega_0 \subset \Omega$. It follows that for $\omega \in \Omega_0$, the induced representations are representations on $l^2(\Pi)$, where $\Pi$ is the point group of $\Gamma$. Moreover, for $\omega \in \Omega_0$, these induced representations have an explicit formula.
	
	In Chapter \ref{DecompChapter}, we find a unitary operator $U$ transforming $\pi$ into a representation $\tau$ on $L^2(R; l^2(\Pi)\otimes l^2(\mathcal{L}^*))$. Moreover, $\tau$ is a direct integral of factor representations over $R$.
	
	In Chapter \ref{OrthConditionChapter}, we show that any $\tau$-invariant subspace of $L^2(R; l^2(\Pi)\otimes l^2(\mathcal{L}^*))$ satisfies the orthogonality condition. By the key theorem in Chapter 3, this ensures that any $\tau$-invariant subspace is a direct integral of closed subspaces of $l^2(\Pi)\otimes l^2(\mathcal{L}^*)$, so that our characterization is complete.
	
	In Chapter \ref{InvSubspaceChapter}, we use our characterization of $\tau$-invariant subspaces, together with the map found in Chapter \ref{DecompChapter}, to characterize the $\pi$-invariant subspaces of $L^2(\mathbb{R}^n)$, thereby achieving the primary goal of this thesis. In addition, we show that the representation $\tau$ is unitarily equivalent to a representation $\rho$ which is a direct integral of factor representations over the unitary dual space $\widehat{\Gamma}$ of $\Gamma$, with the property that the center of the von Neumann algebra generated by $\rho$ is the diagonal algebra. Thus we give the \textit{central decomposition of $\pi$}.
	
Another goal of this thesis was to recast certain techniques from the area of shift-invariant subspaces in more abstract terms. The techniques used in this thesis can be used to prove known results such as the theorems mentioned in the first part of this introduction.
	
	\chapter[Crystallographic Groups, Fundamental Domains,\\ and Group  Actions]{Crystallographic Groups, Fundamental Domains, and Group Actions}\label{Preliminaries}
	\section{Crystallographic groups}
	
	Let $\mathrm{GL}(n,\mathbb{R})$ denote the group of invertible linear transformations of $\mathbb{R}^n$, with identity element $\text{id}$. We shall follow the notation used in \cite{KeithJosh}. For $x\in\mathbb{R}^n$ and $L \in \mathrm{GL}(n,\mathbb{R})$\label{GLnr}, we let $[x,L]$ denote the affine transformation on $\mathbb{R}^n$ given by $[x,L]z = L(z+x)$\label{xL}. These transformations form a group, called the \textit{affine group}:
	\[
	\mathrm{Aff}(n,\mathbb{R}) = \{[x,L]: x\in\mathbb{R}^n, L \in \mathrm{GL}(n,\mathbb{R})\}
	\]
	with group multiplication
	\begin{equation}\label{AffnMult}
		[x,L][y,M] = [M^{-1}x + y, LM];
	\end{equation}
	inverse
	\[
	[x,L]^{-1} = [-Lx,L^{-1}];
	\]
	and identity $[0,\text{id}]$, the identity transformation on $\mathbb{R}^n$. 
	\begin{rem}
		A more standard definition is
		\[
		\mathrm{Aff}(n,\mathbb{R}) = \mathbb{R}^n\rtimes \mathrm{GL}(n,\mathbb{R}),
		\]
		(where elements of this group are regarded as affine transformations of $\mathbb{R}^n$). \footnote{The (naturally isomorphic) group we have defined here is in fact
		\[
		\big(\mathbb{R}^n\rtimes \mathrm{GL}(n,\mathbb{R})^{\mathrm{op}} \big)^{\mathrm{op}}, 
		\]
		where the ``op'' superscript denotes the \textit{opposite group}. Given a group, $G = (G,*)$, its \textit{opposite group} is the group $G^{\mathrm{op}} = (G,*^\prime)$, with underlying set $G$, and binary operation $*^\prime$ given by $g *^\prime h = h * g$. We will occasionally identify groups with their opposites, writing $=$ in place of the more accurate $\cong$.}
	\end{rem}
	The affine group is the group of all invertible affine transformations on $\mathbb{R}^n$. The group of translations in $\mathrm{Aff}(n,\mathbb{R})$, explicitly
	\[
	\mathrm{T}\big(\mathrm{Aff}(n,\mathbb{R})\big) :=\{[x,\text{id}]: x\in\mathbb{R}^n\}, 
	\]
	is a normal subgroup of $\mathrm{Aff}(n,\mathbb{R})$. Let $Q:  \mathrm{Aff}(n,\mathbb{R}) \to \mathrm{GL}(n,\mathbb{R})$ be the map $Q([x,L]) = L$; then $Q$\label{Q} is a homomorphism onto $\mathrm{GL}(n,\mathbb{R})$ with kernel $\mathrm{T}\big(\mathrm{Aff}(n,\mathbb{R})\big)$, so that
	\[
	\faktor{\mathrm{Aff}(n,\mathbb{R})}{\mathrm{T}\big(\mathrm{Aff}(n,\mathbb{R})\big)} \cong \mathrm{GL}(n,\mathbb{R}), 
	\]
	by the first isomorphism theorem. 
	
	Let $\mathrm{O}(n)$\label{On} denote the group of orthogonal transformations of $\mathbb{R}^n$ and let 
	\[
	\mathrm{I}(\mathbb{R}^n) = \{[x,L]\in\mathrm{Aff}(n,\mathbb{R}): L\in\mathrm{O}(n)\}
	\]
	denote the group of isometries of $\mathbb{R}^n$\label{In}. Then $\mathrm{I}(\mathbb{R}^n)$ inherits a topology from $C(\mathbb{R}^n)$, the space of all continuous maps from $\mathbb{R}^n$ to itself.\footnote{$C(\mathbb{R}^n)$ has the topology of uniform convergence on compact sets.} Let $\Gamma$ be a discrete subgroup of $\mathrm{I}(\mathbb{R}^n)$. Elements of $\Gamma$ are affine transformations, and hence act on $\mathbb{R}^n$ by 
	\begin{equation}\label{sec:firstaction}
		[x,L]\cdot_{\idref{1}} u= [x,L]u = L(u+x)
	\end{equation}
	(see \S \ref{groupactions}). We define the set
	\[
	\mathbb{R}^n/\Gamma = \{[u] :u\in\mathbb{R}^n\},
	\]
	where $[u] = \Gamma u= \big\{[x,L]\cdot_{\idref{1}} u:[x,L]\in \Gamma\big\}$ is the $\Gamma$-orbit of $u$ in $\mathbb{R}^n$. We give $\mathbb{R}^n/\Gamma$ the quotient topology, i.e., the topology that makes $p: \mathbb{R}^n \to \mathbb{R}^n/\Gamma$, defined by $p(u) = [u]$, a \textit{quotient map}; it consists precisely of the sets whose preimages under $p$ are open in $\mathbb{R}^n$. With this topology, we call $\mathbb{R}^n/\Gamma$ the \textit{orbit space} of the action of $\Gamma$ on $\mathbb{R}^n$ (\cite{Ratcliffe}, p.~231). The orbit space is an example of a \textit{quotient space}. In general, given a group $G$ and a set $X$, we say that an action of $G$ on $X$ is \textit{cocompact} when the orbit space $X/G$ is compact.
	\begin{defn}
		A subgroup $\Gamma$\label{crystalgroup} of $\mathrm{I}(\mathbb{R}^n)$ is called a \textit{crystallographic group} (\textit{crystal group} for short) if it is discrete and its action \eqref{sec:firstaction} is cocompact.
	\end{defn}
We note that the cocompactness requirement in the above definition is met precisely when the space of cosets $\mathrm{I}(\mathbb{R}^n)/\Gamma$ is compact (Proposition 1.9 in \cite{Szczepanski}).
	\section{The Point Group and Translation Subgroup}
	
	Let $\Gamma$ be a crystal group, and let \label{pi}
	\[
	\Pi = Q(\Gamma) = \{L\in \mathrm{O}(n) : [x,L] \in \Gamma\ \text{for some}\ x\in\mathbb{R}^n\}.
	\]
	It can be shown that $\Pi$ is a finite subgroup of $\mathrm{O}(n)$, called the \textit{point group} (see ~\cite{Ratcliffe},~p.~311). The group of translations in $\Gamma$, denoted $\mathrm{T}(\Gamma)$ or $\mathrm{T}$ for short, is given by\label{T}
	\[
	\mathrm{T}= \Gamma\cap \mathrm{T}\big(\mathrm{Aff}(n,\mathbb{R})\big) = \{[x,L]\in\Gamma:L=\text{id}\}.
	\]
	This is an abelian normal subgroup of $\Gamma$, and is the kernel of $Q\vert_{\Gamma}$; thus $\Gamma/\mathrm{T} \cong \Pi$, by the first isomorphism theorem. For potential later reference, we explicitly define maps
	\begin{equation}\label{qmaps}
		\Gamma \overset{\tilde{q}}{\longrightarrow} \Gamma/\mathrm{T} \overset{Q^*}{\longrightarrow} \Pi,
	\end{equation}
	where $\tilde{q}$\label{tildeq} is the canonical quotient map, and $Q^*([x,L]\mathrm{T}) = Q([x,L])$, so that $Q\vert_{\Gamma} = Q^*\circ \tilde{q}$. To see that $Q^*$\label{Qstar} is well-defined, note that if $[y,M]\mathrm{T}=[x,L]\mathrm{T}$, then $[x,L]^{-1}[y,M]\in\mathrm{T}$, from which it follows that $L=M$. We note here that a discrete group $\Gamma \subset \mathrm{I}(\mathbb{R}^n)$ is a crystal group if and only if the subgroup $\mathrm{T}$ of translations has finite index in $\Gamma$ and is generated by $n$ linearly independent vectors in $\mathbb{R}^n$ (\cite{Ratcliffe},~p.~309). 
	
	We define a \textit{lattice} of $\mathbb{R}^n$ to be a subgroup generated by $n$ linearly independent vectors of $\mathbb{R}^n$. Equivalently, a lattice is a subset of $\mathbb{R}^n$ of the form $\mathcal{L} = B\mathbb{Z}^n$ for some invertible $n \times n$ matrix $B$. We note that every lattice of $\mathbb{R}^n$ is a discrete subgroup of $\mathbb{R}^n$ (Corollary 2, p.~158, \cite{Ratcliffe}). Some authors use the term lattice more generally, where $B$ may be noninvertible, and say that the lattice is \textit{full rank} when $B$ is invertible. If $\Gamma \subset \mathrm{I}(\mathbb{R}^n)$ is a crystal group, its translation subgroup $\mathrm{T}$ determines a full rank lattice $\mathcal{L}$\label{L} of $\mathbb{R}^n$ known as the \textit{Bravais lattice}, after the French physicist Auguste Bravais. Namely, $\mathcal{L} = \{x\in\mathbb{R}^n: [x,\mathrm{id}] \in \mathrm{T}\}$. This is sometimes called the \textit{translation lattice} or \textit{vector lattice} of $\Gamma$. The translation subgroup $\mathrm{T}$ consists of translation operations by vectors of the form
	\begin{equation*}
		k_1 \mathbf{b}_1 + k_2 \mathbf{b}_2+ \cdots + k_n \mathbf{b}_n , \qquad k_i \in \mathbb{Z},\ i=1,2,\cdots, n,
	\end{equation*}
	where the $\mathbf{b}_i$ are linearly independent by the condition of $\mathrm{T}$ being of full rank. Thus the Bravais lattice $\mathcal{L}$ associated with $\Gamma$ is the subgroup
	\begin{equation}\label{translationvectors}
		\mathcal{L}=\{k_1 \mathbf{b}_1 + k_2 \mathbf{b}_2+ \cdots + k_n \mathbf{b}_n : k_i \in \mathbb{Z},\ i=1,2,\cdots, n\}.
	\end{equation}
	The $\mathbf{b}_i$ are called \textit{primitive translation vectors} for the lattice, and give a \textit{lattice basis} for the Bravais lattice. Letting $B = [\begin{array}{cccc}
		\mathbf{b}_1 &  \mathbf{b}_2 & \dots &  \mathbf{b}_n
	\end{array} ]$, we see that $\mathcal{L}= B\mathbb{Z}^n$.
	\begin{rem} The set of primitive translations vectors that generate a given lattice is not unique.
	\end{rem}
	
	We note that $\mathcal{L}$ is the image of $\mathrm{T}$ under the map $\pi_1: \mathrm{Aff}(n,\mathbb{R})\to \mathbb{R}^n$ which maps $[x,L]$ to $x$. Indeed, it is clear from the definition of $\mathcal{L}$ that $\pi_1\vert_{\mathrm{T}}$ is a bijective group homomorphism from $\mathrm{T}$ to $\mathcal{L}$. 
\begin{defn}
The \textit{group of symmetries} (or \textit{symmetry group}) of a lattice $\mathcal{L}\subset \mathbb{R}^n$ is the group $\mathrm{Sym}(\mathcal{L})$ of all isometries of $\mathbb{R}^n$ that leave $\mathcal{L}$ invariant, with multiplication given by composition.
\end{defn} 
\begin{rem}
Each full rank lattice $\mathcal{L}$ arises as the Bravais lattice of some crystal group, e.g., $\mathcal{L}\rtimes \{\mathrm{id}\}$.
\end{rem}
\begin{prop}
$\mathrm{Sym}(\mathcal{L})$ is the unique symmorphic crystal group $\bar\Gamma = \mathcal{L}\rtimes \bar\Pi$ having $\mathcal{L}$ as its Bravais lattice and with the property that $\bar\Pi$ has maximal order among all symmorphic crystal groups $\mathcal{L}\rtimes \Pi$ with Bravais lattice $\mathcal{L}$.
\end{prop}
\begin{proof}
Let $G=\mathrm{Sym}(\mathcal{L})$. Let 
\[
\mathrm{T}(G)=	\mathrm{T}\big(\mathrm{Aff}(n,\mathbb{R})\big) \cap G.
\]
Then $\mathrm{T}(G)$ is a normal subgroup of $G$, being the kernel of $Q\vert_G$. Moreover, $\mathrm{T}(G) = \{ [y,\mathrm{id}]: y\in\mathcal{L}\}$; for if $[x,\mathrm{id}] \in \mathrm{T}(G)$ for some $x\in \mathbb{R}^n$, then $[x,\mathrm{id}]y \in \mathcal{L}$ for each $y\in \mathcal{L}$, whence $x\in\mathcal{L}$. 

We first show that $G$ splits as $\mathcal{L}\rtimes Q(G)$, where $Q$ is the homomorphism described on page \pageref{Q}. Otherwise, $[0,L]\notin G$ for some $L \in Q(G)$. Then since $L \in Q(G)$, there exists $u \in \mathbb{R}^n$ such that $[u,L] \in G$. Then by the normality of $\mathrm{T}(G)$ in $G$, we have
\[
[Ly,\mathrm{id}] = [u,L][y,\mathrm{id}][u,L]^{-1} \in \mathrm{T}(G)\quad \text{for each}\ y\in \mathcal{L}.
\]
Thus $[0,L]y = Ly \in\mathcal{L}$ for all $y \in \mathcal{L}$; i.e. $[0,L]$ is an isometry of $\mathbb{R}^n$ which leaves $\mathcal{L}$ invariant. But this implies that $[0,L]\in G$, a contradiction. This proves that $G = \mathcal{L}\rtimes Q(G)$. This argument also shows that $Q(G)$ is a subgroup of $\mathrm{O}(n)$ that leaves $\mathcal{L}$ invariant.

Finally, since $\bar{\Pi}$ is the point group of $\bar\Gamma$, and $\bar\Gamma$ has Bravais lattice $\mathcal{L}$, it follows that $\bar\Pi$ leaves $\mathcal{L}$ invariant. Thus $\bar\Gamma = \mathcal{L}\rtimes \bar\Pi$ leaves $\mathcal{L}$ invariant, and hence $\bar\Gamma \subset G$. Therefore, $\bar{\Pi} = Q(\bar\Gamma) \subset Q(G)$. If this containment is proper, then $\mathcal{L}\rtimes Q(G)$ is a symmorphic crystal group whose point group is larger than $\bar{\Pi}$, a contradiction. Thus $Q(G) = \bar\Pi$, whence $G = \mathcal{L}\rtimes Q(G) = \mathcal{L}\rtimes \bar{\Pi}$. 
\end{proof} 

Bravais lattices are considered equivalent if their symmetry groups are isomorphic\footnote{Two $n$-dimensional crystallographic groups are isomorphic if and only if they are conjugate in $\mathrm{Aff}(n,\mathbb{R})$ (Corollary 3, p.~318 of \cite{Ratcliffe}).}; in dimension 2 there are five symmetrically distinct Bravais lattices, and in dimension 3 there are fourteen. 
	
	In crystallography, a \textit{unit cell} for a lattice is a subset of $\mathbb{R}^n$ whose translates by the lattice elements fill the space without gaps or overlaps. A unit cell whose volume is minimal is called a \textit{primitive unit cell}. A primitive unit cell is not unique; but one possible choice is given by the parallelotope 
	\[
	x_1 \mathbf{b}_1 + x_2 \mathbf{b}_2+ \cdots + x_n \mathbf{b}_n , \qquad x_i \in [0,1),\ i=1,2,\cdots, n,
	\]
	where the $\mathbf{b}_i$ are primitive translation vectors. This primitive unit cell is an example of a \textit{fundamental domain} for the translation subgroup of $\Gamma$; we will discuss fundamental domains later.
	
	\section{Symmorphic and Nonsymmorphic Crystal Groups}
	\begin{defn}
		$\Gamma$ is called \textit{symmorphic} if $[0,L] \in \Gamma$ for all $L\in\Pi$. $\Gamma$ is \textit{nonsymmorphic} if there exists $L\in\Pi$ such that $[0,L] \notin \Gamma$. Some authors use the term \textit{split} for symmorphic, which refers to the splitting of a short exact sequence
		\[
		[0,\mathrm{id}] \to \mathrm{T}  \overset{\iota}{\longrightarrow} \Gamma  \overset{Q}{\longrightarrow} \Pi  \to \mathrm{id}.
		\]
		Thus $\Gamma$ is symmorphic if and only if $\Gamma \cong \mathrm{T} \rtimes \Pi$.
	\end{defn}
	Crystal groups are symmetry groups of crystal patterns; a \textit{crystal pattern} is a set $P$ of points in $\mathbb{R}^n$ for which the translations leaving $P$ invariant form a lattice (\cite{CrystalSummer}, p.~3). In dimension 2, these groups are often referred to as \textit{wallpaper groups}. It has been shown that two crystal groups are isomorphic precisely when they are conjugate in $\mathrm{Aff}(n,\mathbb{R})$. Up to isomorphism, there are seventeen wallpaper groups, four of which are nonsymmorphic. In dimension 3 there are 219 isomorphism classes of crystal groups, 73 of which are nonsymmorphic. Crystallographers often classify crystal groups by a slightly more rigid notion of equivalence: namely, equivalence up to conjugation by \textit{orientation-preserving} affine transformation. See \cite{CrystalSummer},~p.~35 for more on this distinction.
	
	\begin{ex}\label{p4m}
		Consider the symmorphic wallpaper group $\Gamma=p4m$. Its point group is $D_8$, the dihedral group of order eight. It contains two rotation centres of order four, as well as reflections along horizontal, vertical, and diagonal axes. It also contains glide reflections along axes parallel to the diagonals, and rotations of order two at the intersection of the glide reflection axes. See the figure below for a pattern whose symmetries are described by the group. See also Example \ref{p4mcodedpicture}; Wikipedia \cite{enwiki:1221937434} also has helpful pictures illustrating the symmetries in the unit cell for all 17 wallpaper groups.
	\end{ex}
	\begin{figure}[H]
		\centering
		\includegraphics[width=12cm]{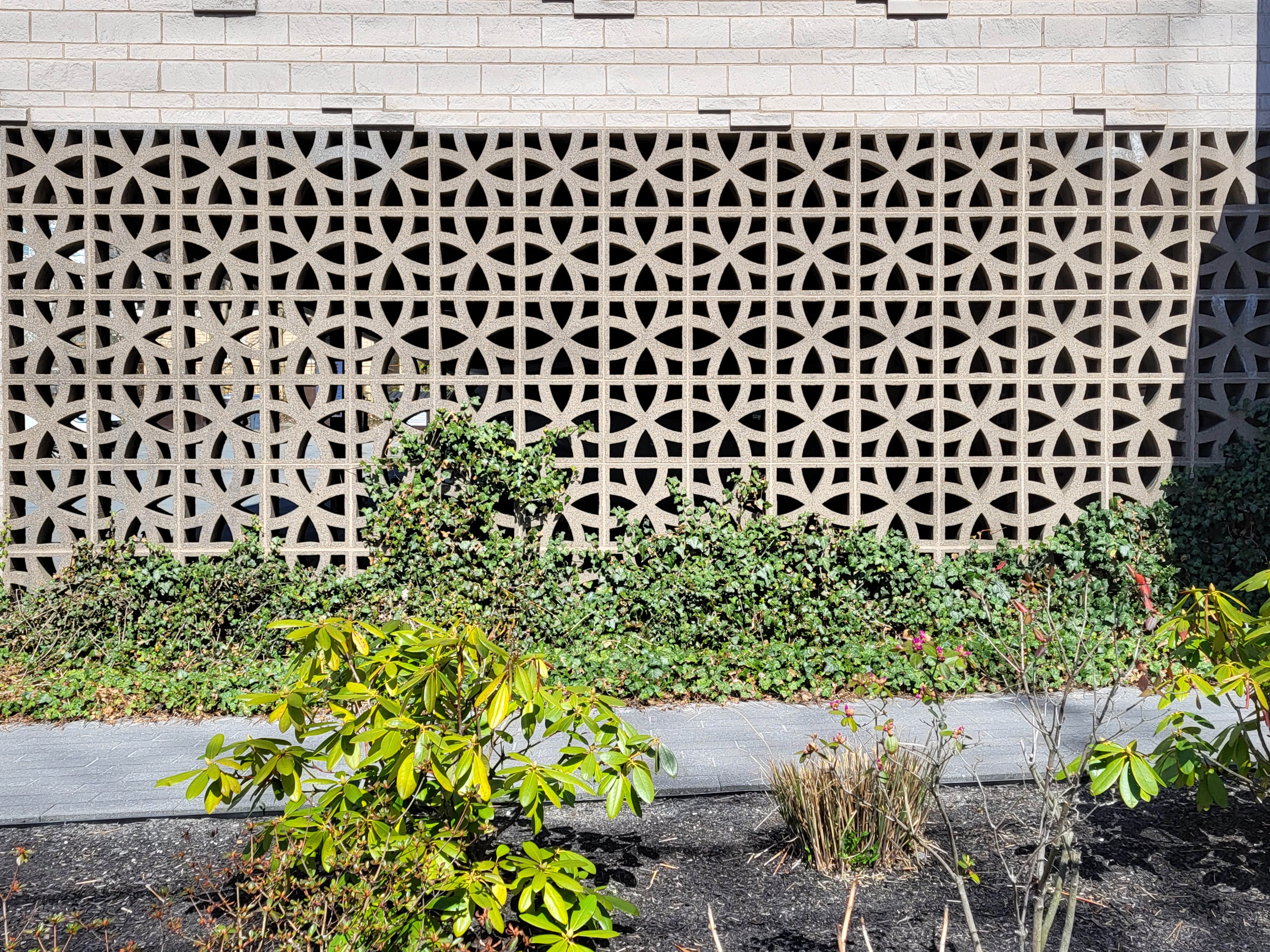}
		\caption{A pattern illustrating the wallpaper group $p4m$, photographed on Oxford Street in Halifax, Nova Scotia, near Dalhousie University.}
	\end{figure}
	
	\begin{ex}\label{pgexample}
		Consider the wallpaper group $pg$. It contains no rotations or reflections, but has \textit{glide reflections} along parallel axes. A glide reflection is a reflection together with a translation in the direction parallel to the axis of reflection. See the figure below for a pattern whose symmetries are described by the group. The purple vectors illustrate a choice of primitive translation vectors, and the hatched square illustrates a choice of primitive unit cell. With the choice of primitive translation vectors in the pattern illustrated below, the glide reflection of $pg$ is expressed as $[(k+\tfrac12)\mathbf{b}_1,\sigma]$, where $k\in\mathbb{Z}$, and $\sigma$ is the element of $\mathrm{O}(n)$ fixing $\mathbf{b}_1$ and sending $\mathbf{b}_2$ to $-\mathbf{b}_2$. $pg$ is an example of a nonsymmorphic crystal group, because $[0,\sigma] \notin \Gamma$.
		
		\begin{figure}[H]
			\centering
			\begin{tikzpicture}
				\foreach \i [evaluate=\i as \y using \i]  in {0,2}
				\foreach \x in {-4,-2,...,4}
				{
					\shade[thick, lower right = white!90!black, upper left = white!40!black] (\x, \y) -- (\x+2, \y) -- (\x+2,\y+1) -- (\x, \y+1);		
					\draw[line width=1.2pt] (\x, \y) -- (\x+2, \y) -- (\x+2,\y+1) -- (\x, \y+1) -- (\x,\y) -- (\x+2,\y);
					\draw[line width=1.2pt] (\x+.8, \y+.3) -- (\x+.8, \y+.7) -- (\x+1.2,\y+.7); 
				}
				\foreach \i [evaluate=\i as \y using \i]  in {1,3}
				\foreach \x in {-4,-2,...,4}
				{
					\shade[thick, lower left = white!40!black, upper right = white!90!black]  (\x+1, \y) -- (\x+3, \y) -- (\x+3,\y+1) -- (\x+1, \y+1);	
					\draw[line width=1.2pt] (\x+1, \y) -- (\x+3, \y) -- (\x+3,\y+1) -- (\x+1, \y+1) -- (\x+1,\y) -- (\x+3,\y);	
					\draw[line width=1.2pt] (\x+1.8,\y+.7) -- (\x+1.8, \y+.3) -- (\x+2.2, \y+.3); 
				}
				\draw[pattern=north west lines, pattern color={rgb: red, 2; blue, 4}] (1,1) rectangle (3,3);		
				\draw[color={rgb:  red, 2; blue, 4}, line width=2pt][->] (1,2) -- (1,1) -- (3,1);
				\node[color={rgb: red, 2; blue, 4}, anchor=north west] at (2.9,1.1) {$\mathbf{b}_1$};
				\node[color={rgb: red, 2; blue, 4}, anchor=south east] at (.95,2.55) {$\mathbf{b}_2$};
				\draw[color={rgb: red, 2; blue, 4}, line width=2pt][->] (1,1) -- (1,3);		
			\end{tikzpicture}
			\caption{A pattern illustrating the symmetries of the wallpaper group $pg$, with primitive translation vectors $\mathbf{b}_1$ and $\mathbf{b}_2$, and a choice of primitive unit cell (hatched).}
		\end{figure}
	\end{ex}
	
	\section{Cross-Sections and Transversals}
	\begin{defn}
		Let $\sim$ be an equivalence relation on a set $X$, and let $X/\sim$ denote the set of equivalence classes. A \textit{transversal} for $\sim$ is a set $\mathcal{X}\subseteq X$ that meets every equivalence class in exactly one point. A \textit{cross-section} for $\sim$ is a function $\gamma: X/\sim \to X$ such that, for each equivalence class $\mathcal{E} \in X/\sim$, $\gamma(\mathcal{E}) \in \mathcal{E}$. When $H$ is a subgroup of $G$, we may speak of a cross-section (transversal, resp.) of $G/H$, which we define to be a cross-section (transversal, resp.) for the equivalence relation whose classes are the left cosets of $H$.
	\end{defn}
	\begin{rem}
		If we have an action of a group $G$ on a space $X$, then there is an equivalence relation on $X$ whose classes are the orbits $\mathcal{O} \in X/G$. A transversal for this equivalence relation will be referred to as a \textit{transversal associated with the group action}.
	\end{rem}
	\begin{rem}
		Note that if $\gamma$\label{crosssec} is a cross-section for $\sim$, then $\mathrm{image}(\gamma)$ is a transversal. Indeed, $\mathrm{image}(\gamma)$ intersects the equivalence classes $\mathcal{E}$ in the singleton $\gamma(\mathcal{E})$. Conversely, if $\mathcal{X}$ is a transversal for $\sim$, then $\gamma(\mathcal{E}) = \mathcal{E}\cap \mathcal{X}$ gives a cross-section.
	\end{rem}
	\begin{rem}
		Given an equivalence relation $\sim$ on $X$, we may consider the quotient map $q:X\to X/\sim$ given by $q(x)= [x]$, where $[x]$ denotes the equivalence class containing $x$. Thus $x\sim y \iff q(x) = q(y)$. Then $\gamma$ is a cross-section for $\sim$ if and only if $q\circ \gamma = \mathrm{id}_{X/\sim}$. Indeed, for $\mathcal{E} \in X/\sim$,
		\[
		\mathcal{E} = [\gamma(\mathcal{E})] \iff \gamma(\mathcal{E})\in \mathcal{E}.
		\]
	\end{rem}
	
	We know from algebra that $\Gamma$ is partitioned by the cosets of $\mathrm{T}$. This gives an equivalence relation $\sim$ on $\Gamma$, where the equivalence classes are the left cosets of $\mathrm{T}$. Thus, in this case, $\Gamma/\sim$ is just the quotient group $\Gamma/\mathrm{T}$. Since $\mathrm{T}$ is the kernel of the homomorphism $Q$, the left cosets of $\mathrm{T}$ are the fibers of $Q$ over elements of $\Pi$. Define a map $\gamma: \Pi \to \Gamma$ as follows: for each $L\in \Pi$, we let $\gamma(L)= [x_L,L]$ be a coset representative for the coset that is the fiber of $Q$ over $L$. With the previously defined maps
		\begin{equation*}
		\Gamma \overset{\tilde{q}}{\longrightarrow} \Gamma/\mathrm{T} \overset{Q^*}{\longrightarrow} \Pi,
	\end{equation*}
we define a cross-section $\tilde{\gamma}$ for $\sim$ by $\tilde\gamma = \gamma\circ Q^*$. Thus, by the above remark, $\tilde q \circ \tilde\gamma  = \mathrm{id}_{\Gamma/\mathrm{T}}$. It follows that $\gamma = \tilde{\gamma}\circ (Q^*)^{-1}$ and $Q\circ \gamma$ is the identity automorphism of $\Pi$. By a slight abuse of language we refer to $\gamma$ as a cross-section for $\Gamma$. The following diagram illustrates the maps involved:
	\[
	\begin{tikzcd} 
		\Gamma \arrow[r, shift left=.37ex, "\tilde{q}"]  \arrow[to=1-3, bend left=30, "Q"] 
		& \arrow[l,shift left=.37ex, "\tilde\gamma"]  \Gamma/\mathrm{T} \arrow[r, "Q^*"] 
		& \Pi 	\arrow[to=1-1, bend left=30, "\gamma"].
	\end{tikzcd}
	\]
	\begin{tcolorbox}
		\begin{notation}\label{crosssection}
		In the sequel, we assume that for each crystallographic group $\Gamma$ we have fixed a cross-section $\gamma: \Pi \to \Gamma$ with $\gamma(\mathrm{id}) = [0,\mathrm{id}]$. We reserve the Greek $\gamma$ for this cross-section. In the sequel we write $x_L$ for the fixed element of $\mathbb{R}^n$ such that $\gamma(L)= [x_L,L]$.
	\end{notation}
	\end{tcolorbox}
	\begin{rem}
	If $\Gamma$ is symmorphic, then we can simply take $\gamma(L) = [0,L]$.
		\end{rem}
	\begin{ex}\label{pgcross}
		Continuing the example of the group $pg$, we can take $\gamma(\mathrm{id}) = [0, \mathrm{id}]$ and $\gamma(\sigma) = [\tfrac12 \mathbf{b}_1 ,\sigma]$. We have 
		\begin{align*}
			pg &= \gamma(\text{id}) \mathrm{T} \cup \gamma(\sigma)\mathrm{T}\\
			&=\{[k_1\mathbf{b}_1 +k_2\mathbf{b}_2,\text{id}] : k_1, k_2 \in \mathbb{Z} \} \cup \{[(k_1+\tfrac12)\mathbf{b}_1 + k_2 \mathbf{b}_2,\sigma] : k_1, k_2\in\mathbb{Z}\}.
		\end{align*}
	\end{ex}
	
	\section{Fundamental Regions, Domains, and Sets}\label{FundRegionConstruction}
	
	We follow the treatment of Ratcliffe \cite{Ratcliffe} in this section, narrowing the definitions and proofs to fit our purpose of showing the existence of a fundamental domain for a crystallographic group. The proofs in this section are essentially taken from there.
	
	\begin{defn} A \textit{fundamental domain} for a group $\Gamma \subset \mathrm{I}(\mathbb{R}^n)$ is a subset $\mathcal{R} \subset \mathbb{R}^n$ such that  
		\newcounter{ctr}
		\begin{list}{(\roman{ctr})}{\usecounter{ctr}}
			\setcounter{ctr}{0}
			\item $\mathcal{R}$ is open in $\mathbb{R}^n$;
			\item the members of $\{[x,L]\mathcal{R}: [x,L]\in\Gamma\}$ are mutually disjoint;
			\item $\mathbb{R}^n = \bigcup \{[x,L]\overline{\mathcal{R}}: [x,L]\in\Gamma\}$;
			\item $\mathcal{R}$ is connected.
		\end{list}
		A subset $\mathcal{R}$ satisfying (i) to (iii) is called a \textit{fundamental region} for $\Gamma$.
	\end{defn}
	\begin{ex}\label{p4mcodedpicture}
		Consider the two-dimensional symmorphic wallpaper group $\Gamma=p4m$ described in Example \ref{p4m}. The figure below shows a pattern whose symmetries are described by the group. The orange triangle (boundary excluded) represents a fundamental domain for $\Gamma$.
		\begin{figure}[H]
			\centering
			\begin{tikzpicture}
				\foreach \i [evaluate=\i as \y using 2*\i]  in {0,1,...,4}
				\foreach \x in {-4,-2,...,4}
				{
					\draw[line width=.5pt, color=orange, opacity=.2, dotted][-] (1,5) -- (2,5) -- (2,6) -- (1,5) -- (2,5);
					\fill[line width=1.2pt, color=gray!10] (\x,\y) --(\x+2,\y) -- (\x+2,\y+ 2) -- (\x,\y+2) -- (\x,\y) -- (\x+2,\y); 
					\fill[line width=2pt, color=orange, opacity=.03] (1, 5) -- (2, 5) -- (2,6) -- (1, 5)-- (2, 5);	
					\shade[thick, lower left = black!50!white, upper right = black!50!white] (\x+1, \y+1) -- (\x+2, \y+1) -- (\x+2,\y+2) -- (\x+1, \y+1)-- (\x+2, \y+1);		
					\shade[thick, lower left = black!50!white, upper right = black!50!white] (\x+1, \y+1) -- (\x+1, \y+2) -- (\x+2,\y+2) -- (\x+1, \y+1)-- (\x+1, \y+2);
					\shade[thick, lower right = black!50!white, upper left = black!50!white] (\x+1, \y+1) -- (\x+1, \y+2) -- (\x,\y+2) -- (\x+1, \y+1)-- (\x+1, \y+2);
					\shade[thick, lower right = black!50!white, upper left = black!50!white] (\x+1, \y+1) -- (\x, \y+2) -- (\x,\y+1) -- (\x+1, \y+1)-- (\x, \y+2);
					\shade[thick, upper right = black!50!white, lower left = black!50!white] (\x+1, \y+1) -- (\x, \y+1) -- (\x,\y) -- (\x+1, \y+1)-- (\x, \y+1);
					\shade[thick, upper right = black!50!white, lower left = black!50!white] (\x+1, \y+1) -- (\x, \y) -- (\x+1,\y) -- (\x+1, \y+1)-- (\x, \y);
					\shade[thick, upper left = black!50!white, lower right= black!50!white] (\x+1, \y+1) -- (\x+1, \y) -- (\x+2,\y) -- (\x+1, \y+1)-- (\x+1, \y);
					\shade[thick, upper left = black!50!white, lower right= black!50!white] (\x+1, \y+1) -- (\x+2, \y) -- (\x+2,\y+1) -- (\x+1, \y+1)-- (\x+2, \y);
					\draw[line width=1.2pt, color = black][-] (\x,\y) -- (\x+2,\y+2);
					\draw[line width=1.2pt, color = black][-] (\x,\y+2) -- (\x+2,\y); 
					\draw[line width=1.2pt, color = black][-] (\x+1,\y) -- (\x+1,\y+2);
					\draw[line width=1.2pt, color = black][-] (\x,\y+1) -- (\x+2,\y+1);
					\draw[line width=1.2pt, color = black][-] (\x,\y) -- (\x+2,\y); 
					\draw[line width=1.2pt, color = black][-] (\x,\y) -- (\x,\y+2); 
					\draw[line width=1.2pt, color = black][-] (-4,0) -- (-4,10) -- (6,10) -- (6,0) -- (-4,0);
					\draw[line width=1.2pt, color = black][-] (\x+1.5,\y+1.2) -- (\x+1.8,\y+1.2) -- (\x+1.8,\y+1.5);
					\draw[line width=1.2pt, color = black][-] (\x+.2,\y+1.5) -- (\x+.2,\y+1.2) -- (\x+.5,\y+1.2);
					\draw[line width=1.2pt, color = black][-] (\x+1.5,\y+1.8) -- (\x+1.2,\y+1.8) -- (\x+1.2,\y+1.5);
					\draw[line width=1.2pt, color = black][-] (\x+.5,\y+1.8) -- (\x+.8,\y+1.8) -- (\x+.8,\y+1.5);
					\draw[line width=1.2pt, color = black][-] (\x+.5,\y+.2) -- (\x+.8,\y+.2) -- (\x+.8,\y+.5);
					\draw[line width=1.2pt, color = black][-] (\x+.2,\y+.5) -- (\x+.2,\y+.8) -- (\x+.5,\y+.8);
					\draw[line width=1.2pt, color = black][-] (\x+1.5,\y+.8) -- (\x+1.8,\y+.8) -- (\x+1.8,\y+.5);
					\draw[line width=1.2pt, color = black][-] (\x+1.2,\y+.5) -- (\x+1.2,\y+.2) -- (\x+1.5,\y+.2);
					\draw[line width=1.2pt, color = black][-] (-4,0.1) -- (-4,0) -- (-3.8,0);
					\node at (1.6,5.35) {${\scriptscriptstyle \mathcal{R}}$};
				}
			\end{tikzpicture}
			\caption{A pattern illustrating the symmetries of the wallpaper group $p4m$, with a choice of fundamental domain $\mathcal{R}$ shaded orange ($\mathcal{R}$ does not include the boundary of this triangle since it is open).}
		\end{figure}
	\end{ex}
	
	\begin{defn}
		A subset $\mathfrak{S} \subset \mathbb{R}^n$ is called a \textit{fundamental set} for $\Gamma \subset \mathrm{I}(\mathbb{R}^n)$ if and only if $\mathfrak{S}$ meets each $\Gamma$-orbit in a singleton. Thus, a fundamental set for $\Gamma$ is a transversal for the action of $\Gamma$ on $\mathbb{R}^n$.
	\end{defn}
	
	\begin{lemma}\label{rigidity}
		If $U\subset\mathbb{R}^n$ is open and nonempty, and if there is $[x,L]\in\mathrm{Aff}(n,\mathbb{R})$ such that $[x,L]u=u$ for all $u\in U$ (that is, $[x,L]$ fixes $U$), then $[x,L] = [0,\mathrm{id}]$.
	\end{lemma}
	\begin{proof}
		Let $U\subset\mathbb{R}^n$ be open and nonempty, with $[x,L]$ fixing $U$. Then there is $y\in \mathbb{R}^n$ such that $0\in [y,\mathrm{id}]U$, i.e., $[y,\mathrm{id}]U$ is a neighbourhood of $0$. Then there is a basis $\{v_1,\dots,v_n\}$ of $\mathbb{R}^n$ contained in $[y,\mathrm{id}]U$. Since $U$ is fixed by $[x,L]$, it follows that $[y,\mathrm{id}]U$ is fixed by $[y,\mathrm{id}][x,L][-y,\mathrm{id}]$. Hence
		\[
		[y,\mathrm{id}][x,L][-y,\mathrm{id}]v_i = v_i\ \ \ \text{for all}\ i=1,\dots, n,
		\]
		whence
		\[
		[y,\mathrm{id}][x,L][-y,\mathrm{id}] = [0,\mathrm{id}].
		\]
		It follows that $[x,L] = [0,\mathrm{id}]$.
	\end{proof}
	\begin{lemma}\label{Rboundarydisjoint}
		Let $\Gamma$ be a subgroup of $\mathrm{I}(\mathbb{R}^n)$, and let $\mathcal{R}$ be a fundamental region for $\Gamma$. For each $r\in \mathcal{R}$, the orbit $\Gamma r$ is disjoint from the boundary $\partial \mathcal{R}$.
	\end{lemma}
	\begin{proof}
		Let $r\in \mathcal{\mathcal{R}}$, and let $[x,L]r \in \Gamma r$. We consider two cases: if $[x,L] =  [0,\mathrm{id}]$, then $r \in \mathcal{R}$, which is disjoint from $\partial \mathcal{R}$, since $\mathcal{R}$ is open. If, on the other hand, $[x,L]\neq [0,\mathrm{id}]$, then $[x,L]r \in [x,L]\mathcal{R}$, which is open and disjoint from $\mathcal{R}$. Thus there exists an open ball $B([x,L]r,\delta) \subset \mathcal{\mathcal{R}}^c$. If $z\in \partial \mathcal{R} \cap \mathcal{R}^c$, then $z$ is an accumulation point of $\mathcal{R}$, by Theorem 3.1.6 in \cite{OSearcoid}. Thus there is no ball $B(z,\delta) \subset \mathcal{\mathcal{R}}^c$ when $z\in \partial \mathcal{\mathcal{R}} \cap \mathcal{R}^c$. Hence $[x,L]r \notin \partial \mathcal{\mathcal{R}}$. Therefore no element of $\Gamma r$ is in $\partial \mathcal{R}$, for $r\in \mathcal{R}$.
	\end{proof}
	\begin{lemma}\label{FundSetExist}
		Let $\Gamma$ be a subgroup of $\mathrm{I}(\mathbb{R}^n)$. An open subset  $\mathcal{R}\subset\mathbb{R}^n$ is a fundamental region for $\Gamma$ if and only if there is a fundamental set $\mathfrak{S}$ for $\Gamma$ with $\mathcal{R} \subset \mathfrak{S} \subset\overline{\mathcal{R}}$.
	\end{lemma}
	\begin{proof}
		Suppose first that $\mathcal{R}$ is a fundamental region for $\Gamma$. Since the members of $\{[x,L]\mathcal{R} : [x,L]\in \Gamma\}$ are mutually disjoint, it follows that $\mathcal{R}$ contains at most one element from each $\Gamma$-orbit. Indeed, suppose $r_0\in \mathcal{O} \in \mathbb{R}^n/\Gamma$ with $r_0\in \mathcal{R}$. Any other element $z$ of $\mathcal{O}$ is of the form $z=[x,L]r_0 \in [x,L]\mathcal{R}$ with $[x,L]\neq [0,\mathrm{id}]$, and hence $z\notin \mathcal{R}$. On the other hand, $\overline{\mathcal{R}}$ must contain at least one element from each $\Gamma$-orbit, because $\mathbb{R}^n = \bigcup \{[x,L]\overline{\mathcal{R}}: [x,L]\in\Gamma\}$; and hence any $w_0\in\mathbb{R}^n$ is in $[x,L]\overline{\mathcal{R}}$, some $[x,L] \in \Gamma$. Thus $[x,L]^{-1}w_0\in\overline{\mathcal{R}}$, and by definition, $[x,L]^{-1}w_0$ is in the $\Gamma$-orbit of $w_0$. By the axiom of choice, we can choose an element $w_{\mathcal{O}} \in \mathcal{O}\cap\overline{\mathcal{R}}$ for each $\mathcal{O} \in \mathbb{R}^n/\Gamma$. Letting $\mathfrak{S}$ be the set of these chosen elements, we see that $\mathcal{R}\subset \mathfrak{S}$. To see this, note that since $\mathcal{R}$ is open, it follows that $\overline{\mathcal{R}}=\mathcal{R}\, \dot\cup\, \partial \mathcal{R}$, and hence $\mathcal{O}\cap\overline{\mathcal{R}} = (\mathcal{O}\cap \mathcal{R}) \, \dot\cup\, (\mathcal{O}\cap\partial \mathcal{R})$ for each $\mathcal{O} \in \mathbb{R}^n/\Gamma$. If $r_1\in \mathcal{R}$, then since $\mathcal{R}$ contains at most one element of $\Gamma r_1$, and $\partial \mathcal{R}$ is disjoint from $\Gamma r_1$ by Lemma \ref{Rboundarydisjoint}, any other element  $w_1\in \Gamma r_1$ with $w_1\neq r_1$ is outside $\overline{\mathcal{R}}=\mathcal{R} \cup \partial \mathcal{R}$. Thus $r_1$ is the unique element of $\overline{\mathcal{R}}\cap \Gamma r_1$, and therefore $r_1\in \mathfrak{S}$.
		
		Suppose conversely that $\mathfrak{S}$ is a fundamental set for $\Gamma$ in $\mathbb{R}^n$ with $\mathcal{R}\subset \mathfrak{S} \subset \overline{\mathcal{R}}$. We first show that the members of $\{[x,L]\mathcal{R}: [x,L]\in\Gamma\}$ are mutually disjoint. Suppose there are elements $[x,L], [y,M]\in \Gamma$ such that $[x,L]\mathcal{R}\cap[y,M]\mathcal{R} \neq \emptyset$; i.e., $\mathcal{R}\cap [x,L]^{-1}[y,M]\mathcal{R} \neq \emptyset$. Then there are $r_0,r_1 \in \mathcal{R}$ such that $r_0 = [x,L]^{-1}[y,M]r_1$, so that $r_0$ is the in $\Gamma$-orbit of $r_1$. Since $\mathfrak{S}$ contains only one element from each orbit, it follows that $r_0=r_1$ and $[y,M]^{-1}[x,L]r_0 = r_0$. But $r_0$ was an arbitrary element of $\mathcal{R}\cap [x,L]^{-1}[y,M]\mathcal{R}$, so that $\mathcal{R}\cap [x,L]^{-1}[y,M]\mathcal{R}$ is a nonempty open set fixed by $[y,M]^{-1}[x,L]$. Hence $[y,M]^{-1}[x,L] = [0,\mathrm{id}]$, by Lemma \ref{rigidity}; i.e., $[x,L] = [y,M]$. 
		
		To show that the union of the $[x,L]\overline{\mathcal{R}}$ is $\mathbb{R}^n$, observe that since $\mathfrak{S}$ is a fundamental set, and $\mathfrak{S}\subset \overline{\mathcal{R}}$, we have
		\[
		\mathbb{R}^n = \bigcup \{[x,L]\mathfrak{S}: [x,L] \in \Gamma\} \subset \bigcup\{[x,L]\overline{\mathcal{R}}: [x,L] \in \Gamma\}.
		\]
	\end{proof}
	For a crystal group $\Gamma$, it follows from Lemma \ref{rigidity} that there exists a point $a\in\mathbb{R}^n$ whose stabilizer $\Gamma_a$ is trivial (Theorem 6.6.12 in \cite{Ratcliffe}). Given such an $a$ and $[x,L]\neq [0,\mathrm{id}]$ in $\Gamma$, we define the open half-space
	\[
	H_{[x,L]}(a) = \{y\in \mathbb{R}^n: d(y,a)<d(y,[x,L]a)\}.
	\]
	The open half-space $H_{[x,L]}(a)$ contains $a$ and its boundary is the hyperplane orthogonal to the line segment joining $a$ and $[x,L]a$ and equidistant from these two points. In other words, its boundary is the perpendicular bisector of this line segment. 
	\begin{defn}
The \textit{Dirichlet domain} $D(a)$ with center $a$ for a crystal group $\Gamma$ is defined to be 
\[
\boxed{D(a) = \bigcap\{H_{[x,L]}(a): [x,L]\in \Gamma\ \text{with}\ [x,L]\neq [0,\mathrm{id}]\}}\,.
\]
	\end{defn}
	
	\begin{theorem}\label{FunDomainExist}
		Let $\Gamma$ be a crystal group and let $a \in \mathbb{R}^n$ have trivial stabilizer. Then $D(a)$ as defined above is a fundamental domain for $\Gamma$.
	\end{theorem}
	\begin{proof}
		\textbf{Step 1: We show that $\boldsymbol{D(a)}$ is open.} Let $\overline{B}(a,\delta)$ be a closed ball containing $a$. Since $\Gamma$ is discrete, it follows that the orbits of $\Gamma$ are discrete subsets of $\mathbb{R}^n$ (p.~163, \cite{Ratcliffe}). It then follows that $\overline{B}(a,\delta)$ contains only finitely many points of any orbit $\mathcal{O}\in\mathbb{R}^n/\Gamma$, because any infinite subset of a compact set must have an accumulation point in the set. Let 
		\[
		K_{[x,L]} = \mathbb{R}^n - H_{[x,L]}(a)
		\]
		for each $[x,L]\neq [0,\mathrm{id}]$ in $\Gamma$. Then $K_{[x,L]}$ is closed, and 
		\begin{equation}\label{distance1}
		d(a,K_{[x,L]}) = \tfrac12 d(a,[x,L]a),
		\end{equation}
		because $\partial K_{[x,L]}$ is the bisector of $a$ and $[x,L]a$.  Let
		\[
		\mathcal{S} = \{K_{[x,L]}: [x,L]\neq [0,\mathrm{id}]\}.
		\]
		Let $w_0\in \mathbb{R}^n$, and choose $\delta$ be such that $w_0\in B(a,\delta)$. We show that $B(a,\delta)$ meets only finitely many $K_{[x,L]}\in\mathcal{S}$. Suppose that $B(a,\delta)$ meets $K_{[x,L]}$. Then by \eqref{distance1}, 
\[		
			d(a,[x,L]a) = 2	d(a,K_{[x,L]}) < 2\delta.
	\]
	Hence, if $B(a,\delta)\cap K_{[x,L]}\neq \emptyset$, then $B(a,2\delta)$ contains $[x,L]a$. Consequently, as $B(a,2\delta)$ contains only finitely many points of the orbit $\Gamma a$, it follows that the ball $B(a,\delta)$ meets only finitely many of the sets $K_{[x,L]}$. Therefore $B(a,\delta)$ is a neighbourhood of $w_0$ that meets only finite many sets in $\mathcal{S}$. It follows that $D(a)$ is open; for if $w_1\in D(a)$, we have shown there exists a neighbourhood $U$ of $w_1$ such that $U$ meets only finitely many sets in $\mathcal{S}$. Labeling these sets $K_1, \dots, K_n$, we have that $V = U - \cup_{i=1}^n K_i$ is an open neighbourhood of $w_1$, with
		\begin{align*}
			V &\subset \mathbb{R}^n - \cup_{i=1}^n K_i \\
			&\subset \mathbb{R}^n - \cup\{K_{[x,L]}: [x,L]\neq [0,\mathrm{id}]\} \\
			&= D(a).
		\end{align*}
		
		\textbf{Step 2: We show existence of a fundamental set $\boldsymbol{\mathfrak{S}}$ for $\boldsymbol{\Gamma}$ such that $\boldsymbol{D(a) \subset \mathfrak{S} \subset \overline{D}(a)}$, and $\boldsymbol{D(a)}$ is connected.} The set $\mathfrak{S}$ is obtained as follows: for each orbit $\Gamma w_0$, select a point $w_1 \in \{u\in \Gamma w_0: d(a,\Gamma w_0) = d(a,u)\}$; thus $w_1$ is among the points of $\Gamma w_0$ which are nearest $a$. Let $\mathfrak{S}$ be the union of these chosen points. Then $\mathfrak{S}$ is a fundamental set for $\Gamma$. If $w_2\in D(a)$ and $[x,L]\neq [0,\mathrm{id}]$ in $\Gamma$, then since $w_2\in H_{[x,L]^{-1}}(a)$, we have
		\[
		d(w_2,a) < d(w_2,[x,L]^{-1}a) = d([x,L]w_2,a);
		\]
		therefore $w_2$ is the unique point of $\Gamma w_2$ nearest $a$. Thus $D(a)\subset \mathfrak{S}$.
		
		Now let $z \in \mathfrak{S} - \{a\}$, let $\ell(a,z)$ be the (closed) line segment in $\mathbb{R}^n$ joining $a$ to $z$, and let $w_3$ be a point in the interior of $\ell(a,z)$; thus 
		\begin{equation}\label{firstseg}
			d(w_3,a) = d(z,a) - d(z,w_3).
		\end{equation}	
		For $[x,L]\neq [0,\mathrm{id}]$ in $\Gamma$ we have
		\[
		d(z,a) \leq d([x,L]^{-1}z, a) = d(z,[x,L]a),
		\]
		whence 
		\begin{equation}\label{segs0}
			d(w_3,a) \leq d(z,[x,L]a) - d(z,w_3) \leq d(w_3,[x,L]a).
		\end{equation}
		Comparing with \eqref{firstseg} we see that this holds with equality precisely when 
		\begin{equation}\label{segs}
			d(z,a) =  d(z,[x,L]a) = d(z,w_3) + d(w_3,[x,L]a).
		\end{equation}
		If \eqref{segs} holds, then $z$, $w_3$, and $[x,L]a$ must be collinear. In this case, the line segments $\ell(z,[x,L]a)$ and $\ell(z,a)$ both extend $\ell(w_3,z)$ and have the same length; hence they coincide. Thus $a=[x,L]a$, and since the stabilizer of $a$ is trivial, we get $[x,L] = [0,\mathrm{id}]$, a contradiction. Therefore there must be inequality in \eqref{segs0}; i.e., $d(w_3,a) < d(w_3,[x,L]a)$. Hence $w_3\in H_{[x,L]}(a)$ for all $[x,L]\neq [0,\mathrm{id}]$, so that $w_3\in D(a)$. Since $w_3$ was an arbitrary point in the interior of the line segment $\ell(a,z)$, we have $z\in \overline{D}(a)$. This proves that $\mathfrak{S}\subset \overline{D}(a)$ and hence $D(a)$ is a fundamental region. Moreover, if $z \in D(a)$, then the whole line segment $\ell(z,a) \subset D(a)$; this shows that $D(a)$ is connected, and hence a fundamental domain.
	\end{proof}
	\begin{defn}
		A collection $\mathcal{S}$ of subsets of $\mathbb{R}^n$ is called \textit{locally finite} if, for each point $w\in\mathbb{R}^n$, there is an open neighbourhood $U$ of $w$ such that $U$ meets only finitely many members of $\mathcal{S}$. A fundamental region $\mathcal{R}$ for a group $\Gamma \subset \mathrm{I}(\mathbb{R}^n)$ is said to be locally finite if $\{[x,L]\overline{\mathcal{R}}: [x,L] \in \Gamma\}$ is a locally finite collection of subsets of $\mathbb{R}^n$. 
	\end{defn}
	
	\begin{rem}\label{fundDomainFact1}
		 The fundamental domain $D(a)$ in the above theorem is a convex set, being the intersection of open half-spaces of $\mathbb{R}^n$.
		 \end{rem}
	\begin{rem}\label{fundDomainFact2} The fundamental domain $D(a)$ is locally finite. The present argument is from \cite{Ratcliffe}, p.~242,: let $\delta>0$; if $B(a,\delta)$ meets $[x,L]\overline{D}(a)$ for some $[x,L]\in\Gamma$, then it is easy to show that there exists $w\in D(a)$ such that $[x,L]w\in B(a,\delta)$. We have
			\begin{align*}
				d(a,[x,L]a) &\leq d(a,[x,L]w) + d([x,L]w,[x,L]a)\\
				&<\delta+d(w,a)\\
				&\leq \delta+d(w,[x,L]^{-1}a)\\
				&=\delta+d([x,L]w,a)\\
				&<2\delta.
			\end{align*}
			But this is possible for only finitely many $[x,L]\in \Gamma$, for each $\Gamma$-orbit in $\mathbb{R}^n$ is discrete (from 5.3.5  in \cite{Ratcliffe} and Lemma 5, p.~160 of the same). Thus $D(a)$ is locally finite.
	\end{rem}
	
	\section{The Dual Lattice}
	
	Let $\Gamma \subset \mathrm{I}(\mathbb{R}^n)$ be a crystallographic group, and let $\mathcal{L}$ denote its Bravais lattice. 
	\begin{defn}
	The \textit{dual lattice} of $\mathcal{L}$\label{duallattice}, denoted $\mathcal{L}^*$, is defined by
	\[
	\mathcal{L}^* = \{z \in \mathbb{R}^n: x\cdot z \in \mathbb{Z} \ \text{for all}\ x \in \mathcal{L}\},
	\]
	where $x\cdot z = z^T x$ for $x,z\in\mathbb{R}^n$ (here the superscript $T$ denotes matrix transpose). 
	\end{defn}
	\begin{rem}
		We note that the dual of $\mathcal{L}$ is defined more generally to be the set of $f\in (\mathbb{R}^n)^*$---the vector space dual of $\mathbb{R}^n$---such that $f(x) \in\mathbb{Z}$ for all $x\in \mathcal{L}$; this agrees with our definition above if we identify $(\mathbb{R}^n)^*$ with $\mathbb{R}^n$ via the dot product.\footnote{To be precise, each $f\in(\mathbb{R}^n)^*$ is a real linear functional on $\mathbb{R}^n$, and hence can be expressed as a dot product: $f(x) = x_f \cdot x$ for some $x_f\in \mathbb{R}^n$. Then $f \leftrightarrow x_f$ identifies $(\mathbb{R}^n)^*$ with $\mathbb{R}^n$.}
	\end{rem}
	
	If $\mathcal{L} = B\mathbb{Z}^n$, then it can be shown that $\mathcal{L}^*=(B^T)^{-1}\mathbb{Z}^n$. Indeed, 
	\begin{equation}\label{Lstarzn}
		\begin{aligned}
			z\in\mathcal{L}^* &\iff x\cdot z \in \mathbb{Z}\ \ \ \text{for all}\ x\in\mathcal{L}\\
			&\iff (Bk)^Tz= Bk\cdot z \in\mathbb{Z}\ \ \ \text{for all}\ k\in\mathbb{Z}^n\\
			&\iff k\cdot (B^T z) = k^T B^T  z \in \mathbb{Z}\ \ \ \text{for all}\ k\in\mathbb{Z}^n\\
			&\iff B^Tz\in(\mathbb{Z}^n)^* = \mathbb{Z}^n\\
			&\iff z\in (B^T)^{-1}\mathbb{Z}^n.
		\end{aligned}
	\end{equation}
	
	We also require a different type of dual in our work. Let $\widehat{\mathcal{L}}$\label{contcharonL} denote the set of \textit{continuous characters} on $\mathcal{L}$---that is, the set of continuous group homomorphisms of $\mathcal{L}$ into the circle group $\mathbb{T}$. Then $\widehat{\mathcal{L}}$ is an abelian group under pointwise multiplication, and is called the \textit{character group} of $\mathcal{L}$. Likewise, we may consider $\mathbb{R}^n$ as an LCA (locally compact, abelian) group, and it is well known that its character group is $\{x\mapsto \chi_y(u) :=e^{2\pi i y\cdot u} :y\in\mathbb{R}^n\}$\label{chiy}, and hence isomorphic with $\mathbb{R}^n$. We record this notation here, as we will use it frequently in the sequel:
	\begin{tcolorbox}
		\begin{defn}\label{Rnchar}
			The character group of $\mathbb{R}^n$ is $\widehat{\mathbb{R}^n} = \{\chi_y:y \in \mathbb{R}^n\}$, where
			\[
			\chi_y(u) := e^{2\pi i y\cdot u}.
			\]
		\end{defn}
	\end{tcolorbox}
	The isomorphism between $\mathbb{R}^n$ and its dual leads some authors to simply identify them, writing $y \in \widehat{\mathbb{R}^n}$ in place of $\chi_y$. For us, it will be more clear and convenient to keep the $\chi_y$ notation, and even to make the mapping explicit: $\Phi: \mathbb{R}^n \to \widehat{\mathbb{R}^n}$ given by $\Phi(y) = \chi_y$. 
	
	We now describe $\widehat{\mathcal{L}}$, considering $\mathcal{L}$ as a subgroup of $\mathbb{R}^n$. To this end, we define the \textit{orthogonal subgroup} $\mathcal{L}^\perp$ of $\widehat{\mathbb{R}^n}$ to be\label{orthsubgroup}
	\[
	\boxed{\mathcal{L}^\perp= \{\chi_y \in \widehat{\mathbb{R}^n}: \chi_y(x) = 1\ \ \text{for all}\ x\in\mathcal{L}\}}\,.
	\]
	Observing that $\chi_y(x) =1$ for all $x\in\mathcal{L}$ if and only if $x\cdot y \in\mathbb{Z}$ for all $x\in\mathcal{L}$, which by definition occurs when $y \in \mathcal{L}^*$, we deduce that $\mathcal{L}^\perp \cong \mathcal{L}^*$. By Proposition 1.86 in \cite{TaylorKaniuth}, $\widehat{\mathcal{L}} \cong \widehat{\mathbb{R}^n}/\mathcal{L}^\perp$, and hence, identifying $\widehat{\mathbb{R}^n}$ with $\mathbb{R}^n$, we have 
	\[
	\boxed{\widehat{\mathcal{L}} \cong \mathbb{R}^n/\mathcal{L}^*}\,.
	\]
	Note that $\Pi$ acts on $\mathcal{L}^*$ by matrix multiplication, which follows from the lattice-invariance of $\mathcal{L}$ under multiplication by elements of $\Pi$ (see Action \idref{7} below). Indeed, let $z \in \mathcal{L}^*$ and $L\in \Pi$; then for all $x\in\mathcal{L}$ we have
	\[
	x\cdot Lz = L^T x \cdot z = L^{-1}x\cdot z \in\mathbb{Z};
	\]
	thus $Lz \in \mathcal{L}^*$. Therefore we can define a group $\Gamma^*$ as follows: 
 \begin{tcolorbox}
\begin{defn}\label{Gammastar}
 We define the group
\[
\Gamma^*:= \{[z,L]: z\in\mathcal{L}^*, L\in \Pi\},
	\] 
	which has group multiplication
	\[
	[z,L][w,M] = [M^{-1}z + w, LM];
	\]
	inverse
	\[
	[z,L]^{-1} = [-Lx,L^{-1}];
	\]
	and identity $[0,\mathrm{id}]$.
\end{defn}
	\end{tcolorbox}
 Thus\footnote{Technically, $\Gamma^* = (\mathcal{L}^*\rtimes \Pi^{\mathrm{op}})^{\mathrm{op}}$.}
	\[
	\Gamma^* = \mathcal{L}^*\rtimes \Pi \subset\mathrm{I}(\mathbb{R}^n),
	\]
	with the action of $\Pi^{\mathrm{op}}$ on $\mathcal{L}^*$ given by matrix multiplication, and where the elements of this group are regarded as affine transformations of $\mathbb{R}^n$ via the usual action: $[z,L]u = L(u+z)$. In fact, $\Gamma^*$ is a crystal group, because its translation subgroup is of full rank and finite index $|\Pi|$. 
\begin{rem}
We note that, while $(\mathcal{L}^*)^* = \mathcal{L}$, in general $(\Gamma^*)^* \neq \Gamma$, unless $\Gamma$ is symmorphic.
\end{rem}
\begin{tcolorbox}
	\begin{notation}
			Let the translation subgroup of $\Gamma^*$ be denoted by $\mathrm{T}^*$, i.e., $\mathrm{T}^*=\{[z,\mathrm{id}]:z\in\mathcal{L}^*\}$.
	\end{notation}
\end{tcolorbox}
 We know from Theorem \ref{FunDomainExist} that there exists a fundamental domain $R$ for $\Gamma^*$, and $R$ is convex and locally finite, by Remarks \ref{fundDomainFact1} and \ref{fundDomainFact2}. 
 \begin{tcolorbox}
 	\begin{defn}\label{Rdomain}
 		In the sequel, we let $R$ be a fixed fundamental domain for $\Gamma^*$.
 	\end{defn}
 \end{tcolorbox}
	\begin{lemma}\label{translatesopen}
		Let $U \subset\mathbb{R}^n$ be open. Then $LU = \{Lu: u \in U\}$ is open in $\mathbb{R}^n$ for each $L\in \mathrm{O}(n)$. In particular, if $R$ is a fundamental domain for $\Gamma^*$, then $[0,L]R$ is open for each $L\in \Pi$.
	\end{lemma}
	\begin{proof}
		Let $L\in \mathrm{O}(n)$ and let $z \in LU$, so that $L^{-1}z \in U$. Since $U$ is open in $\mathbb{R}^n$, there exists $\delta>0$ such that $L^{-1}z\in B(L^{-1}z,\delta) \subset U$. It is easily verified that $LB(L^{-1}z,\delta) = B(z,\delta)$. Hence, $z \in B(z,\delta) \subset LU$, and thus $LU$ is open.
	\end{proof}
	\begin{lemma}
		Let $\Pi R=\cup\{[0,L]R: L\in \Pi\}$. Then $\overline{\Pi R} = \Pi \overline{R}$, where
		\[
		\Pi \overline{R}	= \bigcup_{L\in \Pi}[0,L]\overline{R}.
		\]
	\end{lemma}
	\begin{proof}
		It is easy to see that $\Pi \overline{R} \subset \overline{\Pi R}$, for if $z\in \Pi \overline{R}$, then $z\in [0,L]\overline{R}$, for some $L\in \Pi$. Thus $z= Lw$ for some $w\in \overline{R}$. Then $w=\lim_{n\to\infty} r_n$ for some sequence $\{r_n\}_{n=1}^\infty \subset R$. Since multiplication by $L$ is continuous, it follows that $z=Lw = \lim_{n\to\infty} L r_n \in \overline{\Pi R}$. We now show that $\Pi \overline{R}$ is closed. To see this, first note that 
		\begin{equation}\label{complement}
			\mathbb{R}^n - [0,L]\overline{R} = [0,L](\mathbb{R}^n -\overline{R}).
		\end{equation}
		Indeed, if $u\notin \overline{R}$, then $[0,L]u\notin [0,L]\overline{R}$ (since $[0,L]$ is invertible). Hence 
		\[
		[0,L](\mathbb{R}^n -\overline{R}) \subset	\mathbb{R}^n - [0,L]\overline{R}.
		\]
		On the other hand, if $v \in \mathbb{R}^n-[0,L]\overline{R}$, then writing $v = [0,L](L^{-1}v)$, we have $L^{-1}v \notin \overline{R}$, and hence $v \in [0,L](\mathbb{R}^n -\overline{R})$. Thus we have established \eqref{complement}. It follows that
		\[
		\mathbb{R}^n - \Pi \overline{R} =\bigcap_{L\in \Pi} \mathbb{R}^n - [0,L]\overline{R} = \bigcap_{L\in \Pi}  [0,L](\mathbb{R}^n -\overline{R}).
		\]
		Moreover, $[0,L](\mathbb{R}^n -\overline{R})$ is open for each $L\in \Pi$, by Lemma \ref{translatesopen}, and therefore $\mathbb{R}^n - \Pi \overline{R}$ is open, being a finite intersection of open sets. Thus we have shown that $\Pi \overline{R}$ is closed. Now since $\Pi R \subset \Pi\overline{R}$, we have
		$\overline{\Pi R} \subset  \Pi\overline{R}$, and hence $\overline{\Pi R} = \Pi\overline{R}$.
		
	\end{proof}
	\begin{lemma}\label{Pirfundreg}
		The open set $\Pi R=\cup\{[0,L]R: L\in \Pi\}$ is a fundamental region for $\mathrm{T}^*$.
	\end{lemma}
	\begin{proof}
		Since $\Pi R$ is a finite union of open sets, it is open. We show that the members of $\{[y,\mathrm{id}]\Pi R: [y,\mathrm{id}]\in \mathrm{T}^*\}$ are disjoint. Indeed, let $u\in [y,\mathrm{id}]\Pi R \cap [z,\mathrm{id}]\Pi R$. Then there exists $L, M \in \Pi$ and $r_1, r_2 \in R$ such that 
		\[
		u= [y,\mathrm{id}][0,L]r_1 = [z,\mathrm{id}][0,M]r_2.
		\]
		Thus $u\in [y,L]R \cap [z,M]R$. Since $R$ is a fundamental domain for $\Gamma^*$, this forces $[y,L] = [z,M]$, so that $y=z$. Thus $[y,\mathrm{id}]\Pi R \cap [z,\mathrm{id}]\Pi R = \emptyset$ whenever $y\neq z$. Lastly, we show that $\bigcup_{[y,\mathrm{id}]\in\mathrm{T}^*}[y,\mathrm{id}]\Pi\overline{R} = \mathbb{R}^n$. For this, first note that
		\[
		[y,\mathrm{id}]\Pi\overline{R} = \bigcup_{L\in \Pi} [y,\mathrm{id}] [0,L]\overline{R} = \bigcup_{L\in \Pi} [y,L]\overline{R}.
		\]
		Thus
		\begin{align*}
			\bigcup_{[y,\mathrm{id}]\in\mathrm{T}^*}[y,\mathrm{id}]\Pi\overline{R} &= 	\bigcup_{y\in\mathcal{L}^*}\Big[\bigcup_{L\in \Pi} [y,L]\overline{R}\Big]\\
			&=	\bigcup_{[y,L]\in\Gamma^*} [y,L]\overline{R}\\
			&=\mathbb{R}^n,
		\end{align*}
		because $R$ is a fundamental region for $\Gamma^*$.
	\end{proof}
	\begin{cor}\label{paramdomaindef}
		It follows from the above Lemma and \ref{FundSetExist} that there exists a fundamental set $\Omega_{\mathrm{T}^*}$ for $\mathrm{T}^*$ such that $\Pi R \subset \Omega_{\mathrm{T}^*}\subset\Pi \overline{R}$.
	\end{cor}
	\begin{tcolorbox}
		\begin{defn}\label{OmegaTstar}
			In the sequel, $\Omega_{\mathrm{T}^*}$ will be a fixed fundamental set for $\mathrm{T}^*$ satisfying $\Pi R \subset \Omega_{\mathrm{T}^*}\subset\Pi \overline{R}$.
		\end{defn}
	\end{tcolorbox}
	\begin{lemma}\label{QuotientTrans}
		$\Omega_{\mathrm{T}^*}$ is a transversal for $\mathbb{R}^n/\mathcal{L}^*$---that is, a transversal for the equivalence relation on $\mathbb{R}^n$ whose equivalence classes are the left cosets of $\mathcal{L}^*$. In other words, $\Omega_{\mathrm{T}^*}$ meets each coset of $\mathcal{L}^*$ exactly once.
	\end{lemma}
	\begin{proof}
		We have shown that $\Omega_{\mathrm{T}^*}$ meets each $\mathrm{T}^*$-orbit precisely once, so it is sufficient to show that these orbits coincide with the cosets of $\mathcal{L}^*$ in $\mathbb{R}^n$. Let $\mathcal{O}$ be a $\mathrm{T}^*$-orbit; then 
		\begin{align*}
			\mathcal{O} =	\mathrm{T}^* u &= \{[z,\mathrm{id}]u : [z,\mathrm{id}]\in\mathrm{T}^*\}\ \ \text{for some}\ u\in \mathbb{R}^n\\
			&= \{[z,\mathrm{id}]u : z\in\mathcal{L}^*\}\\
			&= \{z+u:z\in\mathcal{L}^*\}\\
			&= u+\mathcal{L}^*.
		\end{align*}
		Thus $\Omega_{\mathrm{T}^*}$ is a transversal for $\mathbb{R}^n/\mathcal{L}^*$. 
	\end{proof}
	\begin{tcolorbox}
		\begin{defn}\label{XLdef}
			For each $\nu\in \Omega_{\mathrm{T}^*}$ and each $x\in\mathcal{L}$, we define 
			\[
			\chi_\nu^{\mathcal{L}}= \chi_\nu\vert_{\mathcal{L}}.
			\]
		\end{defn}
	\end{tcolorbox}
	Thus,
	\begin{equation}\label{characterofL}
		\chi_\nu^{\mathcal{L}}(x) = e^{2\pi i \nu\cdot x}\ \ \ \text{for each}\ x\in \mathcal{L}.
	\end{equation}
	It is immediate that each $\chi_\nu^{\mathcal{L}}$ is a continuous character, that is, an element of $\widehat{\mathcal{L}}$.
	\begin{defn}\label{Xmappingdef}
		Define $X: \Omega_{\mathrm{T}^*}\to\widehat{\mathcal{L}}$ by 
		\[
	\boxed{	X(\nu) = \chi_\nu^{\mathcal{L}}}\,.
		\]
	\end{defn}
	\begin{lemma}\label{Xmapping}
		$X$ is injective.
	\end{lemma}
	\begin{proof}
		If $\chi_{\nu_1}^{\mathcal{L}}=\chi_{\nu_2}^{\mathcal{L}}$ then 
		\[
		e^{2\pi i \nu_1 \cdot x} = e^{2\pi i \nu_2\cdot x} \ \ \ \text{for all}\ x\in\mathcal{L},
		\]
		which implies that $\nu_2-\nu_1\in \mathcal{L}^*$, i.e., $\nu_2\in \nu_1+\mathcal{L}^*$. But $\Omega_{\mathrm{T}^*}$ intersects each coset of $\mathcal{L}^*$ in a singleton by Lemma \ref{QuotientTrans}, so that $\nu_1=\nu_2$. Thus $X$ is injective. 
	\end{proof}
	We now show that $X$ is in fact a bijection, and hence gives a parametrization of $\widehat{\mathcal{L}}$: 
	\begin{equation}\label{ParamOfLhat}
		\boxed{\text{Every element of }\widehat{\mathcal{L}}\text{ is of the form }\chi_\nu^{\mathcal{L}}\text{ for a unique }\nu\in \Omega_{\mathrm{T}^*}}
	\end{equation}
	Indeed, $X$ is the composition of several identifications, most of which were discussed previously. We now make these explicit: $X$ is composed of the mappings $X_1, X_2$, and $X_3$ shown in the following commutative diagram:
	\[
	\begin{tikzcd} 
		\Omega_{\mathrm{T}^*} \arrow[to=1-4, bend right=30, "X"{name=D}] \arrow[r, "X_1"{name=U}]
		& \mathbb{R}^n/\mathcal{L}^* \arrow[r, "X_2"{name=U}]
		& \widehat{\mathbb{R}^n}/\mathcal{L}^\perp  \arrow[r, "X_3"{name=U}]
		&  \widehat{\mathcal{L}} ,
	\end{tikzcd}
	\]
	where we define 
	\[
	X_1(\nu) = \nu+\mathcal{L}^*,
	\]
	\[
	X_2(\nu+\mathcal{L}^*) = \chi_\nu\mathcal{L}^\perp,
	\]
	and 
	\[
	X_3(\chi_\nu\mathcal{L}^\perp) = \chi_\nu\vert_{\mathcal{L}}.
	\]
	\begin{rem}\label{qrestrictedbij}
		The map $X_1$ is just $q\vert_{\Omega_{\mathrm{T}^*}}$, where $q:\mathbb{R}^n \to  \mathbb{R}^n/\mathcal{L}^*$ is the quotient map. Note that if $\nu_1, \nu_2 \in \Omega_{\mathrm{T}^*}$, and $q(\nu_1) = q(\nu_2)$, then $\nu_1 \in \nu_2 +\mathcal{L}^*$. But $\Omega_{\mathrm{T}^*}$ is a transversal for $\mathbb{R}^n/\mathcal{L}^*$, so this implies $\nu_1 = \nu_2$, and hence $q\vert_{\Omega_{\mathrm{T}^*}}$ is injective. On the other hand, given $\nu_1 + \mathcal{L}^* \in \mathbb{R}^n/\mathcal{L}^*$, if follows from $\Omega_{\mathrm{T}^*}$ being a transversal for $\mathbb{R}^n/\mathcal{L}^*$ that there exists $\nu_2 \in\Omega_{\mathrm{T}^*}$ such that $\nu_2 \in \nu_1 + \mathcal{L}^*$. Hence $q(\nu_2) =  \nu_1 + \mathcal{L}^*$, and thus $q\vert_{\Omega_{\mathrm{T}^*}}$ is surjective. Thus we have shown that $q\vert_{\Omega_{\mathrm{T}^*}}$ is a bijection.
	\end{rem}
	The map $X_3$ is obtained by applying the first isomorphism theorem to the map $\chi_\nu \mapsto \chi_\nu\vert_{\mathcal{L}}$ from $\widehat{\mathbb{R}^n}$ to $\widehat{\mathcal{L}}$. Thus, it is a well-defined bijection. Indeed, $X_3$ is an isomorphism of topological groups, by Proposition 1.86 in \cite{TaylorKaniuth}. It is easily verified that $X_2$ is well-defined and bijective. It is immediate that these three maps compose to give $X$.
	
	\section{Table of Relevant Group Actions}\label{groupactions}
	Here we provide a table of the relevant group actions involved in the crystal group setting. In Appendix \ref{GroupActions}, we verify that these satisfy the group action criteria; we include this verification for reference purposes. While this table may seem imposing, several of the actions simplify. In particular, we note that
\begin{itemize}
	\item Action 7 is just the restriction of Action 2 to $\mathcal{L}$.
\item  Actions 8 through 11 are merely transferring Actions 3 through 6 from $\mathrm{T}$ and $\widehat{\mathrm{T}}$ to $\mathcal{L}$ and $\widehat{\mathcal{L}}$, respectively, via a natural bijection (see Remark \ref{actionfromiso}). 
\item Action 12 merely applies a matrix $L \in \Pi$ to a vector $\nu \in \Omega_{\mathrm{T}^*}$, then shifts the resulting vector $L\nu$ by the unique lattice vector that pushes $L\nu$ back into $\Omega_{\mathrm{T}^*}$.
\item It will be shown that $\widehat{\mathrm{T}}$ and $\widehat{\mathcal{L}}$ are parametrized by $\Omega_{\mathrm{T}^*}$, and that Actions 5, 6, 10 and 11 have simplified forms. These are listed in Table \ref{SimpTableOfActions}.
\end{itemize}
\begin{tcolorbox}
	\begin{notation}\label{chiL}
		We conventionally denote characters of abelian groups by $\chi$. We often write $\chi^{\mathcal{L}}$ or $\chi^{\mathrm{T}}$ to indicate or emphasize the group.
	\end{notation}
\end{tcolorbox}
	\begin{table}[H]
		\centering
		{
			\small
		\begin{tblr}{width=\textwidth,
					colspec={X[c,.6]X[c,.4]X[c]X[c,2.2]},
					row{1} ={cyan!20, 6mm},
					row{2-Z}={cyan!5, 6mm},
				}
				\textbf{Name}  & \textbf{Symbol} & \textbf{Group and Set} & \textbf{Definition} \\
				\textbf{Action \idref{1}  }& $\cdot_{\idref{1}}$ & $\Gamma$ acting on $\mathbb{R}^n$ &  $[x,L]\cdot_{\iddef{1}}u=L(u+x)$  \\ 
				\textbf{Action \idref{2}  }& $\cdot_{\idref{2}}$ & $\Pi$ acting on $\mathbb{R}^n$ & $L\cdot_{\iddef{2}}u=Lu$ \\
				\textbf{Action \idref{3} }& $\cdot_{\idref{3}}$ & $\Gamma$ acting on $\mathrm{T}$ & $[x,L]\cdot_{\iddef{3}}[y,\mathrm{id}]=[x,L][y,\mathrm{id}][x,L]^{-1}$ \\
				\textbf{Action \idref{4} }& $\cdot_{\idref{4}}$ & $\Pi$ acting on $\mathrm{T}$ & $L\cdot_{\iddef{4}}[y,\mathrm{id}]=\gamma(L)\cdot_{\idref{3}}[y,\mathrm{id}]$ \\
				\textbf{Action \idref{5} }& $\cdot_{\idref{5}}$ & $\Gamma$ acting on $\widehat{\mathrm{T}}$ &	$([x,L]\cdot_{\iddef{5}} \chi^{\mathrm{T}})([y,\text{id}]) = \chi^{\mathrm{T}}([x,L]^{-1}\cdot_{\idref{3}} [y,\text{id}])$ \\
				\textbf{Action \idref{6} }& $\cdot_{\idref{6}}$ & $\Pi$ acting on $\widehat{\mathrm{T}}$ &	$	L\cdot_{\iddef{6}} \chi^{\mathrm{T}} = \gamma(L)\cdot_{\idref{5}} \chi^{\mathrm{T}}$ \\
				\textbf{Action \idref{7} }& $\cdot_{\idref{7}}$ & $\Pi$ acting on $\mathcal{L}$ &	$L\cdot_{\iddef{7}} y = L\cdot_{\idref{2}} y = Ly$ \\
				\textbf{Action \idref{8} }& $\cdot_{\idref{8}}$ & $\Gamma$ acting on $\mathcal{L}$ &	$[x,L]\cdot_{\iddef{8}} y = \pi_1|_{\mathrm{T}}\big([x,L]\cdot_{\idref{3}}(\pi_1|_{\mathrm{T}}^{-1}(y)\big)$ \\
				\textbf{Action \idref{9} }& $\cdot_{\idref{9}}$ & $\Pi$ acting on $\mathcal{L}$ &	$L\cdot_{\iddef{9}} y = \gamma(L)\cdot_{\idref{8}}  y$ \\
				\textbf{Action \idref{10} }& $\cdot_{\idref{10}}$ & $\Gamma$ acting on $\widehat{\mathcal{L}}$ &	$	([x,L]\cdot_{\iddef{10}} \chi^{\mathcal{L}})(y) = \chi^{\mathcal{L}}([x,L]^{-1}\cdot_{\idref{8}} y)$ \\
				\textbf{Action \idref{11} }& $\cdot_{\idref{11}}$ & $\Pi$ acting on $\widehat{\mathcal{L}}$ &	$L\cdot_{\iddef{11}} \chi^{\mathcal{L}} = \gamma(L)\cdot_{\idref{10}} \chi^{\mathcal{L}}$ \\
				\textbf{Action \idref{12} }& $\cdot_{\idref{12}}$ & $\Pi$ acting on $\Omega_{\mathrm{T}^*}$ &	$L\cdot_{\iddef{12}} \nu = (q\vert_{\Omega_{\mathrm{T}^*}})^{-1}\big(q(L\nu)\big)$ 
			\end{tblr}
			\captionof{table}{The relevant group actions for a crystal group $\Gamma$.}
			\label{TableOfActions}
		}
	\end{table}

	\section[The Transversal \texorpdfstring{$\Omega$}{Omega} for the Action of \texorpdfstring{$\Pi$}{Pi} on \texorpdfstring{$\Omega_{\mathrm{T}^*}$}{Parameter Domain}]{The Transversal \texorpdfstring{$\boldsymbol{\Omega}$}{Omega} for the Action of \texorpdfstring{$\boldsymbol{\Pi}$}{Pi} on \texorpdfstring{$\boldsymbol{\Omega_{\mathrm{T}^*}}$}{Parameter Domain}}\label{fundsetOmega}
	
	We now construct a transversal $\Omega$ for the action $\cdot_{\idref{12}}$ of $\Pi$ on $\Omega_{\mathrm{T}^*}$. This is a transversal for the equivalence relation whose equivalence classes are the $\Pi\cdot_{\idref{12}}$-orbits in $\Omega_{\mathrm{T}^*}$. So we must construct a set $\Omega$ which meets each of these orbits exactly once. Let $R$ be as defined in \ref{Rdomain}.
	\begin{lemma}\label{action12onR}
		We have $\Pi\cdot_{\idref{12}} r = \Pi r$ for each $r\in R$, and hence $\Pi\cdot_{\idref{12}} R = \Pi R$, where $\Pi\cdot_{\idref{12}} R = \bigcup_{L\in \Pi} \{ L\cdot_{\idref{12}}r: r\in R\}$. In particular, $\Pi\cdot_{\idref{12}} R$ is a fundamental region for $\mathrm{T}^*$ (by Lemma \ref{Pirfundreg}). Moreover, $\Omega_{\mathrm{T}^*} = \Pi\cdot_{\idref{12}} \overline{R}$, where $\Pi\cdot_{\idref{12}} \overline{R} = \bigcup_{L\in \Pi} \{ L\cdot_{\idref{12}}u: u\in \overline{R}\}$.
	\end{lemma}
	\begin{proof}
		Note that for $r\in R$, $Lr \in LR \subset \Omega_{\mathrm{T}^*}$. Thus $Lr$ is the unique element in $(Lr+\mathcal{L}^*)\cap \Omega_{\mathrm{T}^*}$, and hence
		\begin{equation}\label{12actionOnR}
			L \cdot_{\idref{12}}r= q\vert_{\Omega_{\mathrm{T}^*}}^{-1}\big(q(Lr)\big)= q\vert_{\Omega_{\mathrm{T}^*}}^{-1}(Lr +\mathcal{L}^*) = Lr,
		\end{equation}
	where $q:\mathbb{R}^n \to  \mathbb{R}^n/\mathcal{L}^*$ is the quotient map, and $q\vert_{\Omega_{\mathrm{T}^*}}$ is its restriction to $\Omega_{\mathrm{T}^*}$. It follows that  $\Pi\cdot_{\idref{12}} r = \Pi r$. Now we note that, since $\Omega_{\mathrm{T}^*}\subset \Pi \overline{R}$, and since $q\vert_{\Omega_{\mathrm{T}^*}}$ is a bijection---as shown in Remark \ref{qrestrictedbij}---we have
		\[
		\Omega_{\mathrm{T}^*} = q\vert_{\Omega_{\mathrm{T}^*}}^{-1}\big(q(\Omega_{\mathrm{T}^*})\big) \subset q\vert_{\Omega_{\mathrm{T}^*}}^{-1}\big(q( \Pi \overline{R})\big) \subset\Omega_{\mathrm{T}^*}.
		\]
		Thus 
		\begin{align*}
			\Omega_{\mathrm{T}^*} &= q\vert_{\Omega_{\mathrm{T}^*}}^{-1}\big(q( \Pi \overline{R})\big) \\
			&= \{q\vert_{\Omega_{\mathrm{T}^*}}^{-1}\big(q( Lu)\big) :L\in \Pi, u\in \overline{R}\}\\
			&= \{ L \cdot_{\idref{12}}u: L\in \Pi, u\in \overline{R}\}\\
			&= \Pi \cdot_{\idref{12}}\overline{R}.
		\end{align*}
	\end{proof}
	\begin{lemma}
		For each $r\in R$, $\Pi\cdot_{\idref{12}} r$ is disjoint from $\partial R$.
	\end{lemma}
	\begin{proof}
		By Lemma \ref{action12onR}, we have $\Pi\cdot_{\idref{12}} r = \Pi u$ for each $r\in R$. Thus
		\[
		\Pi\cdot_{\idref{12}} r = \{[0,L]r: L\in \Pi\} \subset \Gamma^* r
		\]
		for each $r \in R$. But $\Gamma^* r$ is disjoint from $\partial R$, by Lemma \ref{Rboundarydisjoint}, and hence $\Pi\cdot_{\idref{12}} r$ is disjoint from $\partial R$.
	\end{proof}
	\begin{prop}\label{ConstructOmega}
		There exists a transversal $\Omega$ for the action $\cdot_{\idref{12}}$ of $\Pi$ on $\Omega_{\mathrm{T}^*}$, with $R\subset \Omega \subset\overline{R}$.
	\end{prop}
	\begin{proof}
		We first note that $R$ contains at most one element from each $\Pi\cdot_{\idref{12}}$-orbit. Indeed, if $r\in R$ and $w \in \Pi\cdot_{\idref{12}} r$, then there exists $L\in \Pi$ such that 
		\[
		w = L \cdot_{\idref{12}} r \in L\cdot_{\idref{12}} R = [0,L]R,
		\]
		which is disjoint from $R$ unless $L=\mathrm{id}$, which implies $w=r$. We next note that $\overline{R}$ meets each $\Pi\cdot_{\idref{12}}$-orbit in $\Omega_{\mathrm{T}^*}$ at least once. Indeed, if $\nu \in \Omega_{\mathrm{T}^*}$, then it follows from Lemma \ref{action12onR} that $\nu = L \cdot_{\idref{12}}u$ for some $L\in \Pi$ and $u\in \overline{R}$. But then $L^{-1}\cdot_{\idref{12}}\nu = u \in \overline{R}$, so that $u\in \Pi\cdot_{\idref{12}} \nu \cap \overline{R}$. It follows now from the axiom of choice that, for each orbit $\mathcal{O} \in \Omega_{\mathrm{T}^*}/\Pi$, we can choose an element $u_{\mathcal{O}} \in \mathcal{O}\cap\overline{R}$. Letting $\Omega$ be the set of these chosen elements, we see that $R\subset \Omega$. To see this, note that since $R$ is open, it follows that $\overline{R}=R\, \dot\cup\, \partial R$, and hence $\mathcal{O}\cap\overline{R} = (\mathcal{O}\cap R) \, \dot\cup\, (\mathcal{O}\cap\partial R)$ for each $\mathcal{O} \in \Omega_{\mathrm{T}^*}/\Pi$. If $r_1\in R$, then since $R$ contains at most one element of $\Pi\cdot_{\idref{12}} r_1$, and $\partial R$ is disjoint from $\Pi\cdot_{\idref{12}} r_1$, any other element $w_1\in \Pi\cdot_{\idref{12}} r_1$ with $w_1\neq r_1$ is outside $\overline{R}=R \cup \partial R$. Thus $r_1$ is the unique element of $\overline{R}\cap \Pi\cdot_{\idref{12}} r_1$, and so $r_1\in \Omega$.
	\end{proof}
	\begin{tcolorbox}
		\begin{defn}\label{Omega}
		In the sequel, $\Omega$ will be a fixed transversal $\Omega$ for the action $\cdot_{\idref{12}}$ of $\Pi$ on $\Omega_{\mathrm{T}^*}$ satisfying $R\subset \Omega \subset\overline{R}$.
	\end{defn}
	\end{tcolorbox}

 We show that the transversal $\Omega$ constructed above is a fundamental set for the action of $\Gamma^*$ on $\mathbb{R}^n$. 
	\begin{prop}\label{Omegafundset}
		The set $\Omega$ constructed above is a fundamental set for $\Gamma^*$.
	\end{prop}
	\begin{proof}
		Let $\mathcal{O}\in \mathbb{R}^n/\Gamma^*$. We first show that $\Gamma\cap \mathcal{O}\neq\emptyset$. Let $u\in\mathcal{O}$. Then since $\Omega_{\mathrm{T}^*}$ is a fundamental set for $\mathrm{T}^*$, there exists $[y,\mathrm{id}]\in\mathrm{T}^*$ such that $[y,\mathrm{id}]u \in \Omega_{\mathrm{T}^*}$. And then since $\Omega$ is a transversal for the action of $\Pi$ on $\Omega_{\mathrm{T}^*}$, there exists $L\in \Pi$ such that $L\cdot_{\idref{12}} [y,\mathrm{id}]u \in \Omega$. Since 
		\[
		L\cdot_{\idref{12}} [y,\mathrm{id}]u = q\vert_{\Omega_{\mathrm{T}^*}}^{-1}\big(q(L[y,\mathrm{id}]u)\big) = q\vert_{\Omega_{\mathrm{T}^*}}^{-1}\big(L[y,\mathrm{id}]u+\mathcal{L}^* \big) \in L[y,\mathrm{id}]u+\mathcal{L}^*,
		\]
		there exists some $z\in\mathcal{L}^*$ such that $L\cdot_{\idref{12}} [y,\mathrm{id}]u = L[y,\mathrm{id}]u+z$. In other words, $L\cdot_{\idref{12}} [y,\mathrm{id}]u = [L,y+L^{-1}z]u$, so that $L\cdot_{\idref{12}} [y,\mathrm{id}]u \in \Gamma^*u=\mathcal{O}$. Thus $\mathcal{O}\cap \Omega\neq \emptyset$, so that $\Omega$ meets each $\Gamma^*$-orbit at least once. We now show that $\Omega$ meets each $\Gamma^*$-orbit at most once. Let $\omega_1, \omega_2 \in \Omega$. Suppose that $\Gamma^*\omega_1 = \Gamma^*\omega_2$, so that $\omega_2 = [y,L]\omega_1$ for some $[y,L]\in\Gamma^*$; i.e., $\omega_2=  L\omega_1 + Ly$. In particular, since $\omega_2 \in \Omega \subset \Omega_{\mathrm{T}^*}$, we have $\omega_2 \in (L\omega_1 + \mathcal{L}^*) \cap \Omega_{\mathrm{T}^*}$. Note that
		\[
		L\cdot_{\idref{12}} \omega_1 =  q\vert_{\Omega_{\mathrm{T}^*}}^{-1}\big(q(L\omega_1 )\big) \in (L\omega_1 + \mathcal{L}^*) \cap \Omega_{\mathrm{T}^*}.
		\]
		But since $\Omega_{\mathrm{T}^*}$ is a transversal for $\mathbb{R}^n/\mathcal{L}^*$, it follows that $L\cdot_{\idref{12}} \omega_1 = \omega_2$. Thus $\omega_2 \in \Pi\cdot_{\idref{12}} \omega_1$, which contradicts $\Omega$ being a transversal for the action of $\Pi$ on $\Omega_{\mathrm{T}^*}$, unless $\omega_1 = \omega_2$. Therefore $\Omega$ meets each $\Gamma^*$-orbit exactly once, and hence is a fundamental set for $\Gamma^*$.
	\end{proof}
	
		\section{Important Relationships Between Actions}\label{PionOmegaTStar}
	We show how the action of $\Pi$ on $\widehat{\mathcal{L}}$, described by \idref{11} in the above table, is reflected in an action on the parameter domain $\Omega_{\mathrm{T}^*}$. The parametrization \ref{ParamOfLhat} allows us to give a concrete description of the action of $\Pi$ on $\widehat{\mathcal{L}}$ in terms of an action of $\Pi$ on the parameter space $\Omega_{\mathrm{T}^*}$. We describe this action presently. Recall that $\Pi$ acts on $\widehat{\mathcal{L}}$ via the fixed cross-section $\gamma$: 
	\[
	L\cdot_{\idref{11}}  \chi^{\mathcal{L}}= \gamma(L)\cdot_{\idref{10}} \chi^{\mathcal{L}}.
	\]
	Let $\gamma(L) =[x_L,L]$ for all $L\in \Pi$, as in \ref{crosssection}. Let $\nu\in \Omega_{\mathrm{T}^*}$ and $y\in\mathcal{L}$. Let $\vartheta: \widehat{\mathrm{T}} \overset{\cong}{\to}\widehat{\mathcal{L}}$\label{vartheta} be defined by $\vartheta(\chi^{\mathrm{T}}) = \chi^{\mathrm{T}} \circ f^{-1}$, where $f = \pi_1|_{\mathrm{T}}$. We have
	\begin{align*}
		(L\cdot_{\idref{11}}  \chi_\nu^{\mathcal{L}})(y)&= \big(\gamma(L)\cdot_{\idref{10}} \chi_\nu^{\mathcal{L}}\big)(y)\\
		&=\vartheta^{-1}(\chi_\nu^{\mathcal{L}})([x_L,L]\cdot_{\idref{3}}[y,\mathrm{id}])\ \text{as in equation}\ \eqref{converttentothree}  \\
		&=\vartheta^{-1}(\chi_\nu^{\mathcal{L}})([L^{-1}y,\mathrm{id}])\ \text{as in equation}\ \eqref{conjbygamma} \\
		&= \chi_\nu^{\mathcal{L}}(L^{-1}y)\\
		&=e^{2\pi i \nu\cdot L^{-1}y}\\
		&=e^{2\pi i L\nu\cdot y}.
	\end{align*}
	Note that since $L\cdot_{\idref{12}}\nu$ belongs to $L\nu + \mathcal{L}^*$, we have 
	\[
	e^{2\pi i L\nu\cdot y} = e^{2\pi i (L\cdot_{\idref{12}}\nu)\cdot y}.
	\]
	It therefore follows from the above equalities that 
	\[
	(L\cdot_{\idref{11}}  \chi_\nu^{\mathcal{L}})(y) = e^{2\pi i (L\cdot_{\idref{12}}\nu)\cdot y} = \chi_{L\cdot_{\idref{12}}\nu}^{\mathcal{L}}(y)
	\]
	for all $\nu\in \Omega_{\mathrm{T}^*}$ and $y\in\mathcal{L}$.
	Thus we see that 
	\begin{equation}\label{piaction}
		\boxed{L\cdot_{\idref{11}} \chi_\nu^{\mathcal{L}} = \chi_{L\cdot_{\idref{12}}\nu}^{\mathcal{L}}}\, .
	\end{equation}
	Equation \eqref{piaction} shows how the action of $\Pi$ on $\widehat{\mathcal{L}}$ is reflected in the action on the parameter domain.

		Having in \S\ref{fundsetOmega} constructed a transversal $\Omega$ for the action $\cdot_{\idref{12}}$ of $\Pi$ on $\Omega_{\mathrm{T}^*}$, it follows from equation \eqref{piaction} that $X(\Omega) =\{\chi_\omega^{\mathcal{L}} : \omega\in \Omega\}$ is a transversal for the action of $\Pi$ on $\widehat{\mathcal{L}}$. We now examine how the parametrization \eqref{ParamOfLhat} of $\widehat{\mathcal{L}}$ can be extended, via $\vartheta$, to give a parametrization of $\widehat{\mathrm{T}}$. For $\nu\in\Omega_{\mathrm{T}^*}$ and $[y,\mathrm{id}]\in\mathrm{T}$, we obtain a character $\chi_\nu^{\mathrm{T}}$ of $\mathrm{T}$:
	\begin{tcolorbox}
		\begin{defn} \label{chiT}
			For  $\nu\in\Omega_{\mathrm{T}^*}$ and $[y,\mathrm{id}]\in\mathrm{T}$, we define
			\begin{equation}\label{ChiomegaTdef1}
				\chi_\nu^{\mathrm{T}}([y,\mathrm{id}])=\chi_\nu(y) .
			\end{equation}
		\end{defn}
	\end{tcolorbox}
	Thus,
	\begin{equation}\label{ChiomegaTcor1}
		\chi_\nu^{\mathrm{T}}([y,\mathrm{id}]) = e^{2\pi i \nu\cdot y}
	\end{equation}
	for all $\nu\in\Omega_{\mathrm{T}^*}$ and $[y,\mathrm{id}]\in\mathrm{T}$.
	Letting $f=\pi_1\vert_{\mathrm{T}}: \mathrm{T} \to \mathcal{L}$ and $\vartheta: \widehat{\mathrm{T}} \to\widehat{\mathcal{L}}$, we see that
	\[
	\chi_\nu^{\mathrm{T}}([y,\mathrm{id}]) = \chi_\nu^{\mathcal{L}}\big(f([y,\mathrm{id}])\big) = \vartheta^{-1}(\chi_\nu^{\mathcal{L}})([y,\mathrm{id}]);
	\] 
	that is, $\chi_\nu^{\mathrm{T}} = \vartheta^{-1}(\chi_\nu^{\mathcal{L}})$.
	It is clear from \eqref{ChiomegaTcor1} that the $\chi_\nu^{\mathrm{T}}$ defined above are indeed characters of $\mathrm{T}$, and moreover,
	\begin{equation}\label{ParamOfThat}
		\boxed{\text{Every element of }\widehat{\mathrm{T}}\text{ is of the form }\chi_\nu^{\mathrm{T}}\text{ for a unique }\nu\in \Omega_{\mathrm{T}^*}}\,.
	\end{equation}
	Furthermore, we have
	\begin{align*}
		L\cdot_{\idref{6}} \chi_\nu^{\mathrm{T}}([y,\mathrm{id}]) &= \big(L\cdot_{\idref{6}} \vartheta^{-1}(\chi_\nu^{\mathcal{L}})\big)([y,\mathrm{id}]) \\
		&= (L\cdot_{\idref{11}}\chi_\nu^{\mathcal{L}})(y)\ \ \text{by equation}\ \eqref{forpionlhat}\\
		&=(\chi_{L\cdot_{\idref{12}}\nu}^\mathcal{L})(y)\ \ \text{by equation}\ \eqref{piaction}\\
		&= \chi_{L\cdot_{\idref{12}}\nu}^\mathcal{L}\big(f([y,\mathrm{id}])\big)\\
		&= \chi_{L\cdot_{\idref{12}}\nu}^{\mathrm{T}}([y,\mathrm{id}]).
	\end{align*}
	Thus, 
	\begin{equation}\label{actionsixtwelverel}
		\boxed{	L\cdot_{\idref{6}} \chi_\nu^{\mathrm{T}}=\chi_{L\cdot_{\idref{12}}\nu}^{\mathrm{T}} \ \ \text{for all}\ \nu \in \Omega_{\mathrm{T}^*}}\,.
	\end{equation}
	Hence by \eqref{12actionOnR}, we have
	\begin{equation}\label{pionthatreflectparamdomain}
		L\cdot_{\idref{6}} \chi_r^{\mathrm{T}}=\chi_{Lr}^{\mathrm{T}} \ \ \ \text{whenever}\ r \in R.
	\end{equation}
	We saw in Proposition \ref{ConstructOmega} that $\Omega$ is a transversal for the action $\cdot_{\idref{12}}$ of $\Pi$ on $\Omega_{\mathrm{T}^*}$. Consequently, 
	\begin{equation*}
		(\vartheta^{-1}\circ X)(\Omega) = \{\chi_\omega^{\mathrm{T}}: \omega\in \Omega\}
	\end{equation*}
	is a transversal for the action of $\Pi$ on $\widehat{\mathrm{T}}$.
	\begin{tcolorbox}
		\begin{defn}\label{upsilon}
			Define $\upsilon = \vartheta^{-1}\circ X$. 
		\end{defn}
	\end{tcolorbox}
	That is, $\upsilon$ is the composition
	\[
	\begin{tikzcd} 
		\Omega_{\mathrm{T}^*}\arrow[to=1-3, bend right=30, "\upsilon"{name=D}] \arrow[r, "X"{name=U}]
		& \widehat{\mathcal{L}} \arrow[r, "\vartheta^{-1}"{name=vartheta}]
		& \widehat{\mathrm{T}} \, .
	\end{tikzcd}
	\]
	Thus $\upsilon(\nu) = \chi_\nu^\mathrm{T}$ for all $\nu \in \Omega_{\mathrm{T}^*}$, and 
	\begin{tcolorbox}
	\begin{equation}\label{transveralpiThat}
	\upsilon(\Omega) = \{\chi_\omega^{\mathrm{T}}: \omega\in \Omega\}\ \text{is a transversal for the action of}\ \Pi\ \text{on}\ \widehat{\mathrm{T}}.
	\end{equation}
	\end{tcolorbox}		
		Finally, we note that it was shown in Remark \ref{moreaboutactiongammaonThat} that
	\begin{equation}
		[x,L]\cdot_{\idref{5}} \chi^{\mathrm{T}} = 	\gamma(L)\cdot_{\idref{5}} \chi^{\mathrm{T}}
	\end{equation}
	for all $[x,L] \in \Gamma$ and any choice of cross-section  $\gamma: \Pi \to \Gamma$. Putting this together with the above formulas, we get the following:
	\begin{tcolorbox}
	\begin{prop}\label{actionsrel}
	Let $\gamma$ be the cross-section fixed in \ref{crosssection}. For any $[x,L] \in \Gamma$, we have
		\begin{equation}
	[x,L]\cdot_{\idref{5}} \chi_\nu^{\mathrm{T}} = 	\gamma(L)\cdot_{\idref{5}}\chi_\nu^{\mathrm{T}}  =	L\cdot_{\idref{6}} \chi_\nu^{\mathrm{T}}=\chi_{L\cdot_{\idref{12}}\nu}^{\mathrm{T}}.
	\end{equation}
		\end{prop}
	\end{tcolorbox}
This proposition, together with \eqref{convertten} and \eqref{piaction}, yields the following table of group actions with simplified forms.
			\begin{table}[H]
		\centering
		{
			\small
			\begin{tblr}{width=\textwidth,
					colspec={X[c,.6]X[c,.4]X[c]X[c,2.2]},
					row{1} ={cyan!20, 6mm},
					row{2-Z}={cyan!5, 6mm},
				}
				\textbf{Name}  & \textbf{Symbol} & \textbf{Group and Set} & \textbf{Simplified Form} \\
				\textbf{Action \idref{5} }& $\cdot_{\idref{5}}$ & $\Gamma$ acting on $\widehat{\mathrm{T}}$ &	$[x,L]\cdot_{\idref{5}} \chi_\nu^{\mathrm{T}} = \chi_{L\cdot_{\idref{12}}\nu}^{\mathrm{T}}$ \\
				\textbf{Action \idref{6} }& $\cdot_{\idref{6}}$ & $\Pi$ acting on $\widehat{\mathrm{T}}$ &	$	L\cdot_{\idref{6}} \chi_\nu^{\mathrm{T}} = \chi_{L\cdot_{\idref{12}}\nu}^{\mathrm{T}}$ \\
				\textbf{Action \idref{10} }& $\cdot_{\idref{10}}$ & $\Gamma$ acting on $\widehat{\mathcal{L}}$ &	$	[x,L]\cdot_{\idref{10}} \chi_\nu^{\mathcal{L}} = \chi_{L\cdot_{\idref{12}}\nu}^{\mathcal{L}}$ \\
				\textbf{Action \idref{11} }& $\cdot_{\idref{11}}$ & $\Pi$ acting on $\widehat{\mathcal{L}}$ &	$L\cdot_{\idref{11}} \chi_\nu^{\mathcal{L}} = \chi_{L\cdot_{\idref{12}}\nu}^{\mathcal{L}}$ 
			\end{tblr}
			\captionof{table}{Group actions that simplify.}
			\label{SimpTableOfActions}
		}
	\end{table}
	We end this section with the following important observation. 
	\begin{prop}\label{differencecrosssec}
Let $\gamma, \gamma^\prime: \Pi \to \Gamma$ be two cross-sections; let $\gamma(L)  = [x_L,L]$ and $\gamma^\prime(L)=[x_L^\prime, L]$ for all $L\in \Pi$. Then $x_L - x_L^\prime \in\mathcal{L}$.
	\end{prop}
	\begin{proof}
By definition, $\mathrm{T} = \{[y,\mathrm{id}]: y\in \mathcal{L}\}$. Since $[x_L,L]$ and $[x_L^\prime, L]$ are in the same left coset of $\mathrm{T}$ in $\Gamma$ (namely, the fiber of $Q$ over $L$), we have $[x_L,L]\in  [x_L^\prime, L] \mathrm{T}$. Thus 
\[
[x_L^\prime,L]^{-1}[x_L,L]\in \mathrm{T};
\]
that is, $[x_L-x_L^\prime,\mathrm{id}] \in \mathrm{T}$, whence $x_L - x_L^\prime\in\mathcal{L}$.
	\end{proof}
	
	\newpage
	
	\chapter[Prerequisite Theory: Unitary Representations and\\  von Neumann  Algebras]{Prerequisite Theory: Unitary Representations and von Neumann Algebras}\label{Prerequisite}
	\section[Weak-operator topology, Strong-operator topology,\\ and Unitary Representations]{Weak-operator topology, Strong-operator topology, and Unitary Representations}
	
	The following proposition will help us to clearly define two important topologies. It is a slight modification of Theorem 5.14 in \cite{FollandReal}. Any relevant terminology in this section can be found in \cite{FollandReal} and \cite{KadRing}.
	
	\begin{prop}\label{TopfromSemi}
		Given a family $\{p_\alpha:\alpha\in A\}$ of seminorms on a vector space $\mathscr{V}$, we can give $\mathscr{V}$ the topology $\mathscr{T}$ generated by the sets
		\[
		B_\alpha(x,\epsilon) = \{y\in\mathscr{V}: p_\alpha(y-x)<\epsilon \};
		\]
		i.e.,
		\[
		\{B_\alpha(x,\epsilon): \alpha \in A, x\in \mathscr{V}, \epsilon>0\}
		\]
		is a subbasis for $\mathscr{T}$.
		\begin{list}{(\arabic{ctr})}{\usecounter{ctr}}\label{weaktobasis}
			\setcounter{ctr}{0}
			\item Let $\mathscr{F}$ denote the collection of all finite subsets of $A$. For $F \in \mathscr{F}$, define
			\[
			I_F(x,\epsilon)= \bigcap_{\alpha\in F} B_\alpha(x,\epsilon).
			\]
			Then $\{I_F(x,\epsilon): F\in \mathscr{F}, \epsilon>0\}$ is a neighbourhood basis at $x$, and hence
			\[
			\{I_F(x,\epsilon): F\in \mathscr{F}, x\in \mathscr{V}, \epsilon>0\}
			\]
			is a basis for the topology $\mathscr{T}$.
			\item \label{convergesemi} If $\{x_i\}_{i\in I}$ is a net in $\mathscr{V}$, then $x_i \to x$ if and only if $p_\alpha(x_i-x) \to 0$ for all $\alpha \in A$.
			\item $(\mathscr{V}, \mathscr{T})$ is a locally convex topological vector space. 
		\end{list}
	\end{prop}
	The first part of this proposition is proved as follows: Let $x\in \mathscr{V}$ and let $\mathcal{U}$ be an open set containing $x$. Since the $B_\alpha(x,\epsilon)$ sets are a subbasis for $\mathscr{T}$, $\mathcal{U}$ consists of a union of finite intersections of these (Proposition 4.4 in \cite{FollandReal}). Thus there exist $x_i, \alpha_i, \epsilon_i$ such that $x \in \cap_{i=1}^n B_{\alpha_i}(x_i,\epsilon_i)\subset\mathcal{U}$. For each $i$, let $\delta_i= \epsilon_i-p_{\alpha_i}(x-x_i)$. Letting $\delta = \min \delta_i$ and $F = \{\alpha_1,\dots,\alpha_n\}$ we have $I_F(x,\delta) \subset \bigcap_{i=1}^n B_{\alpha_i} (x_i, \epsilon_i)$. Indeed, if $y \in I_F(x,\delta)$, then for each $i=1,\dots,n$ we have
	\begin{align*}
		p_{\alpha_i}(y-x_i) & \leq p_{\alpha_i}(y-x)+p_{\alpha_i}(x-x_i)\\
		&< \delta_i+ p_{\alpha_i}(x-x_i)\\
		&= \epsilon_i,
	\end{align*}
	so that $I_F(x,\delta) \subset B_{\alpha_i} (x_i, \epsilon_i)$ for each $i=1,\dots,n$. Thus $x\in I_F(x,\delta)\subset\mathcal{U}$, which proves the first assertion. The remaining assertions are proved in \cite{FollandReal}.
	\begin{defn} Let $\mathscr{V}$ be a complex vector space.
		\begin{list}{(\arabic{ctr})}{\usecounter{ctr}}
			\item A family  $\mathscr{F}$ of linear functionals on $\mathscr{V}$ is said to \textit{separate points of $\mathscr{V}$} if $x\neq 0$ implies there exists $\rho\in\mathscr{F}$ such that $\rho(x)\neq 0$. 
			\item A family  $\{p_\alpha\}_{\alpha\in A}$ of seminorms on $\mathscr{V}$ is said to \textit{separate points of $\mathscr{V}$} if $x\neq 0$ implies there exists $\alpha \in A$ such that $p_\alpha(x)\neq 0$. 
			\end{list}
	 \end{defn}
	 \begin{lemma}\label{SepptsHaus}
	 		Let $\mathscr{V}$ be a vector space with a topology defined by a family $\{p_\alpha\}_{\alpha\in A}$ of seminorms, as in Proposition \ref{TopfromSemi}. Then $(\mathscr{V},\mathscr{T})$ is Hausdorff if and only if $\{p_\alpha\}_{\alpha\in A}$ separates points of $\mathscr{V}$.
	 \end{lemma}
	\begin{proof}
		This is Proposition 5.16 in \cite{FollandReal}.
	\end{proof}
	\begin{lemma}\label{Ex46Folland}
	Suppose that $\mathscr{V}$ is a vector space, $\mathscr{W}$ is a normed vector space, and $\{T_\alpha\}_{\alpha \in A}$ is a collection of linear maps from $\mathscr{V}$ to $\mathscr{W}$. Each $T_\alpha$ defines a seminorm on $\mathscr{V}$ by $x \mapsto \|Tx\|$. By Proposition \ref{TopfromSemi}, these seminorms generate a topology $\mathscr{T}^\prime$ on $\mathscr{V}$. If $\mathscr{T}$ denotes the weak topology on $\mathscr{V}$ generated by the collection $\{T_\alpha\}_{\alpha \in A}$---that is, the weakest (coarsest) topology which makes each $T_\alpha$ continuous---then $\mathscr{T} = \mathscr{T}^\prime$.
	\end{lemma}
	\begin{proof}
		This is Exercise 46 on p.~170 of \cite{FollandReal}.
				\end{proof}
	\begin{rem}
 Suppose that $\mathscr{F}$ is a family of linear functionals on $\mathscr{V}$; that is, $\mathscr{F}$ is a family of linear maps from $\mathscr{V}$ into $\mathbb{C}$. Then the weak topology on $\mathscr{V}$ generated by $\mathscr{F}$ is also called the \textit{weak topology induced on $\mathscr{V}$ by $\mathscr{F}$}; it is Hausdorff whenever $\mathscr{F}$ separates points of $\mathscr{V}$, by Lemmas \ref{SepptsHaus} and \ref{Ex46Folland}.
  \end{rem}
  \begin{tcolorbox}
  \begin{notation}
   Suppose that $\mathscr{F}$ is a family of linear functionals on $\mathscr{V}$. The weak topology induced on $\mathscr{V}$ by $\mathscr{F}$ is denoted by $\sigma(\mathscr{V},\mathscr{F})$.\label{weaktopinduced}
  \end{notation}
\end{tcolorbox}
\begin{defn}\label{dualvector}
	Let $\mathscr{V}$ be a normed vector space. The \textit{dual space} of $\mathscr{V}$---denoted $\mathscr{V}^*$---is the space of all continuous linear functionals on $\mathscr{V}$.
	\end{defn}

Let $\mathscr{V}$ be a normed vector space, and define $\mathrm{Ev}_x: \mathscr{V}^* \to \mathbb{C}$ by $\mathrm{Ev}_x(\rho) = \rho(x)$. Then $\{\mathrm{Ev}_x\}_{x\in \mathscr{V}}$ is a family of continuous linear functionals on $\mathscr{V}^*$ which separates points.
	\begin{defn}
We define the \textit{weak-* topology on $\mathscr{V}^*$} to be the weak topology generated by $\{\mathrm{Ev}_x\}_{x\in \mathscr{V}}$. That is, the weak-* topology on $\mathscr{V}^*$ is $\sigma(\mathscr{V}^*,\{\mathrm{Ev}_x\}_{x\in \mathscr{V}})$.
	\end{defn}
	The weak-* topology on $\mathscr{V}$ is Hausdorff, by Lemma \ref{SepptsHaus}.
	\begin{notation}\label{BH}
Let $\mathscr{H}$ is a complex Hilbert space and let $\mathscr{B}(\mathscr{H})$ denote the space of bounded linear operators on $\mathscr{H}$.
\end{notation}
  	\begin{defn}
 We define the \textit{strong-operator topology} on $\mathscr{B}(\mathscr{H})$ to be the weak topology generated by the linear maps $\{T\mapsto Tx: x\in \mathscr{H}\}$.
	 	\end{defn}
		\begin{rem}
	Since $p_x(T) = \|Tx\|$ defines seminorm on $\mathscr{B}(\mathscr{H})$, for each $x\in \mathscr{H}$, and the collection $\{p_x:x\in\mathscr{H}\}$ separates points of $\mathscr{B}(\mathscr{H})$, it follows from Lemmas \ref{SepptsHaus} and \ref{Ex46Folland} that the strong-operator topology on $\mathscr{B}(\mathscr{H})$ is generated by the seminorms $p_x$, and is Hausdorff.
		\end{rem}	
	\begin{rem}	
	By (\ref{convergesemi}) of Proposition \ref{TopfromSemi}, a net $\{T_i\}_{i\in I}$ in $\mathscr{B}(\mathscr{H})$ converges to $T$ in the strong-operator topology if and only if 
	\[
	|p_x(T_i-T)| = \|T_ix - Tx\| \to 0\ \  \text{for all}\ x\in \mathscr{H}. 
	\]
	In other words, $\{T_i\}$ converges to $T$ in the strong-operator topology if and only if $\{T_i\}$ converges to $T$ pointwise in the norm topology on $\mathscr{H}$. 
	\end{rem}
		Let $\mathscr{H}$ is a complex Hilbert space and consider the family $\mathscr{F}_{\mathrm{w}}$ of linear functionals $\rho_{(x,y)}: \mathscr{B}(\mathscr{H}) \to \mathbb{C}$ defined by 
		\[
		\rho_{(x,y)}(T) = \langle Tx,y\rangle\ \ \text{for}\ x,y\in \mathscr{H}, T\in \mathscr{B}(\mathscr{H}).
		\]
	\begin{defn}
			 	We define the \textit{weak-operator topology} on $\mathscr{B}(\mathscr{H})$ to be the weak topology induced by $\mathscr{F}_{\mathrm{w}}$. 
	\end{defn}
		It is easy to show that $\mathscr{F}_{\mathrm{w}}$ separates points of $\mathscr{B}(\mathscr{H})$, and hence the weak-operator topology is Hausdorff.
	\begin{lemma}[Lemma 13.3 in Munkres]
		Let $\mathscr{T}_1$ and $\mathscr{T}_2$ be a topologies on $X$ with bases $\mathcal{B}_1$ and $\mathcal{B}_2$, respectively. The following are equivalent:
		\begin{list}{(\alph{ctr})}{\usecounter{ctr}}\label{weaktop}
			\setcounter{ctr}{0}
			\item  $\mathscr{T}_1 \subset \mathscr{T}_2$ (i.e., $\mathscr{T}_2$ is finer than $\mathscr{T}_1$).
			\item For each $\mathcal{U}\in\mathcal{B}_1$ and each $x\in \mathcal{U}$, there exists $\mathcal{V} \in \mathcal{B}_2$ such that $x\in \mathcal{V}\subset \mathcal{U}$.
		\end{list}
	\end{lemma}
	\begin{prop}
		The weak-operator topology is in fact induced by the family of linear functionals $\{\rho_{(x,x)}:x\in\mathscr{H}\}$. That is, if we let
		\begin{itemize}
			\item $\mathscr{F}_1$ denote the family of linear functionals $\{\rho_{(x,y)}:x,y\in\mathscr{H}\}$; and
			\item $\mathscr{F}_2$ denote the family of linear functionals $\{\rho_{(x,x)}:x\in\mathscr{H}\}$; 
		\end{itemize}
		then $\sigma(\mathscr{V},\mathscr{F}_1) = \sigma(\mathscr{V},\mathscr{F}_2)$. 
			\end{prop}
		\begin{proof}To see this, first note that $\mathscr{F}_2\subset \mathscr{F}_1$, from which it follows that $\sigma(\mathscr{V},\mathscr{F}_2) \subset \sigma(\mathscr{V},\mathscr{F}_1)$. To show the reverse inclusion, let $F_1$ denote the family of seminorms $\{x\mapsto |\rho_{(x_j,y_j)}|: j=1,\dots,n\}$, and let $A\in I_{F_1}(A_0,\epsilon) \in\sigma(\mathscr{V},\mathscr{F}_1)$. Applying the \textit{polarization identity}, we have
			\begin{align}\label{pol}
				\langle Tx_i,y_i\rangle  &= \tfrac14 \big( \langle T(x_j+y_j), x_j+y_j\rangle - \langle T(x_j-y_j), x_j-y_j\rangle \nonumber\\
				&\qquad \quad + i \langle T(x_j+i y_j),x_j+i y_j\rangle - i\langle T(x_j-i y_j), x_j-i y_j\rangle \big)
			\end{align} 
			for all $j=1,\dots,n$ and all $T\in\mathscr{B}(\mathscr{H})$. Letting $F_2$ denote the family of seminorms
			\[
			\big\{x\mapsto |\rho_{(x_j+c y_j,x_j+c y_j)}|: j\in\{1,\dots,n\}, c\in \{\pm 1,\pm i\}\big\},
			\]
			it follows from \eqref{pol} with $T=A-A_0$ that $A\in I_{F_2}(A,\epsilon) \subset I_{F_1}(A_0,\epsilon)$. Thus $\sigma(\mathscr{V},\mathscr{F}_1) \subset \sigma(\mathscr{V},\mathscr{F}_2)$ by the preceding lemma. 
\end{proof}
It follows from the above proposition that a net $\{T_i\}_{i\in I}$ in $\mathscr{B}(\mathscr{H})$ converges to $T$ in the weak-operator topology if and only if 
	\[
	p_{x,x}(T_i-T)=|\rho_{x,x}(T_i-T)| = |\langle (T_i - T)x, x\rangle| \to 0\ \  \text{for all}\ x\in \mathscr{H}. 
	\]
	In general, the weak-operator topology is strictly coarser than the strong-operator topology which, in turn, is coarser than the norm topology on $\mathscr{B}(\mathscr{H})$. But note:
	\begin{rem}[Theorem 5.1.2 in \cite{KadRing}]\label{convexclosure}
		The weak- and strong-operator closures of a convex subset of $\mathscr{B}(\mathscr{H})$ coincide. 
	\end{rem}
	\begin{prop}[Theorem 2.4.2 in \cite{KadRing}]
		Suppose that $\mathscr{H}$ and $\mathscr{K}$ are Hilbert spaces. Let $\mathscr{B}(\mathscr{H},\mathscr{K})$\label{boundedHK} denote the space of bounded linear operators from $\mathscr{H}$ into $\mathscr{K}$. If $T\in \mathscr{B}(\mathscr{H},\mathscr{K})$, there is a unique element $T^*$ of $\mathscr{B}(\mathscr{K},\mathscr{H})$ such that
		\[
		\langle T^*x,y\rangle = \langle x,Ty\rangle \ \ \ \text{for all}\ x\in \mathscr{K}, y\in \mathscr{H}.
		\]
	\end{prop}
	\begin{defn}\label{adjoint}
	The operator $T^*$ in the above proposition is called the \textit{(Hilbert) adjoint of $T$}.
	\end{defn}
	\begin{prop}[Proposition 2.4.5 in \cite{KadRing}]
	If $T\in \mathscr{B}(\mathscr{H},\mathscr{K})$, then $T$ is an isometric isomorphism from $\mathscr{H}$ onto $\mathscr{K}$ if and only if it is invertible, with inverse $T^{-1} = T^*$. In this case, $T$ is called a \textit{unitary map}.
	\end{prop}
	In the case where $\mathscr{H}=\mathscr{K}$, the above proposition says that $U\in\mathscr{B}(\mathscr{H})$ is an isomorphism from $\mathscr{H}$ onto itself if and only if $UU^* = U^*U = I$, where $I$ denotes the identity operator. 
	\begin{defn}\label{UH}
		Any bounded linear operator $U\in\mathscr{B}(\mathscr{H})$ which satisfies $UU^* = U^*U = I$ is called a \textit{unitary operator}. The collection of unitary operators on $\mathscr{H}$ forms a group under composition of operators; we denote this group by $\mathrm{U}(\mathscr{H})$.
	\end{defn}

	\begin{defn}\label{contunitaryrepdef}
		Let $G$ be a locally compact group. We define a \textit{continuous unitary representation} of $G$ to be a pair $(\pi,\mathscr{H}_\pi)$\label{unitary representation}, where $\mathscr{H}_\pi$ is a Hilbert space, and $\pi$ is a homomorphism from $G$ into $\mathrm{U}(\mathscr{H}_\pi)$ such that $\pi$ on continuous with respect to the strong-operator topology on $\mathrm{U}(\mathscr{H}_\pi)$ inherited from $\mathscr{B}(\mathscr{H}_\pi)$.
	\end{defn}
	 We will often refer to a continuous unitary representation of $G$ as simply a \textit{representation}.
	 \begin{defn}
	 	The \textit{dimension} of a unitary representation $(\pi,\mathscr{H}_\pi)$ of $G$ is the dimension of the Hilbert space $\mathscr{H}_\pi$.
	 \end{defn}
	\begin{rem}
		In the literature, one often finds the required continuity in the above definition to be with respect to the weak-operator topology. This is not a different definition, because the two topologies coincide on $\mathrm{U}(\mathscr{H}_\pi)$: see \cite{FollandAbstract}, pp.~73--74.
	\end{rem}
	
	\section{Essential facts about von Neumann Algebras}
	
	Having established the topological prerequisites, we are now ready to describe a class of operator algebras that will play a crucial role for us. In the definitions that follow, $\mathscr{H}$ will be a complex Hilbert space.
	\begin{defn} 
A \textit{normed $*$-algebra} is a normed algebra $\mathscr{A}$ over $\mathbb{C}$ with norm $\|\cdot\|$ that satisfies $\|ab\|\leq \|a\|\|b\|$ for all $a,b \in \mathscr{A}$, and an involution $a\mapsto a^*$ which satisfies $\|a^*\|=\|a\|$. If $\mathscr{A}$ is complete with respect to this norm then we call $\mathscr{A}$ a \textit{Banach $*$-algebra}.
	\end{defn}
	\begin{defn}
	An \textit{abstract $\mathrm{C}^*$-algebra}, which is a Banach $*$-algebra $A$ whose norm satisfies the \textit{$\mathrm{C}^*$-condition}: $\|a^*a\|=\|a\|^2$. See \cite{TaylorKaniuth} p.~35 for more on abstract $\mathrm{C}^*$-algebras.
		\end{defn}
		\begin{defn} 
		A subset $\mathscr{U}$ of $\mathscr{B}(\mathscr{H})$ is called \textit{self-adjoint} if $T^*\in  \mathscr{U}$ whenever $T\in\mathscr{U}$, where $T^*$ denotes the Hilbert adjoint of the linear operator $T$.
	\end{defn}	
		\begin{defn} 
		A subalgebra $\mathfrak{U}$ of $\mathscr{B}(\mathscr{H})$ is called a \textit{$*$-subalgebra} if it is self-adjoint. $\mathscr{B}(\mathscr{H})$ and its closed $*$-subalgebras are examples of Banach $*$-algebras.
			\end{defn}	
			\begin{defn} 
			 A \textit{concrete $\mathrm{C}^*$-algebra} is a norm-closed $*$-subalgebra of $\mathscr{B}(\mathscr{H})$.
			 	\end{defn}
			 		\begin{defn}\label{vnalg}
		A \textit{von Neumann algebra} is a $*$-subalgebra $\mathscr{M}$ of $\mathscr{B}(\mathscr{H})$ which contains the identity $I$ and is closed in the weak-operator (and hence strong-operator, norm) topology. We sometimes write $\{\mathscr{M},\mathscr{H}\}$ for this von Neumann algebra, to emphasize that $\mathscr{M}$ is a subalgebra of $\mathscr{B}(\mathscr{H})$. In cases where there is no possible confusion, we often just write $\mathscr{M}$.
			\end{defn}
				\begin{defn}\label{vngen}
		If $\mathscr{U}$ is an arbitrary subset of $\mathscr{B}(\mathscr{H})$, then there is a smallest von Neumann algebra $\mathscr{U}^{-}$ on $\mathscr{H}$ containing $\mathscr{U}$; it is called the \textit{von Neumann algebra generated by $\mathscr{U}$} (A.11 in \cite{DixmierC}).
			\end{defn}

	\begin{defn}[The Commutant]\label{commutant} 
		Let $\mathscr{U} \subset \mathscr{B}(\mathscr{H})$. We define the \textit{commutant} of $\mathscr{U}$---denoted $\mathscr{U}^\prime$---to be the collection of bounded operators on $\mathscr{H}$ commuting with all the operators in $\mathscr{U}$. That is
		\[
		\mathscr{U}^\prime=\{T\in \mathscr{B}(\mathscr{H}): ST=TS \ \ \text{for all}\ S\in \mathscr{U}\}.
		\]
	\end{defn}
	\begin{rem}	\label{WOTClosed} Let $\mathscr{U}\subset \mathscr{B}(\mathscr{H})$. Then
		\begin{list}{(\arabic{ctr})}{\usecounter{ctr}}
			\setcounter{ctr}{0}
			\item if $\mathscr{V} \subset \mathscr{B}(\mathscr{H})$ with $\mathscr{V} \subset \mathscr{U}$, then $\mathscr{U}^\prime\subset \mathscr{V}^\prime$.
			\item We denote the double commutant $(\mathscr{U}^\prime)^\prime$ by $\mathscr{U}^{\prime\prime}$. It follows from the definition that $\mathscr{U}\subset \mathscr{U}^{\prime\prime}$.\label{doublecommutant}
			\item  The commutant $\mathscr{U}^\prime$ is weak-operator closed. Indeed, suppose that $\{T_n\}\subset\mathscr{U}$ and $T_n \to T$ in the weak-operator sense. Then $ST_n \to ST$ in the weak-operator sense, for all $S\in\mathscr{U}$; indeed, for every $x,y\in\mathscr{H}$,
			\[
			\langle ST_n x,y\rangle = \langle T_nx,S^*y\rangle \to \langle Tx,S^*y\rangle = \langle STx,y\rangle.
			\]
			Therefore,
			\[
			\langle TSx,y\rangle = \lim_{n\to\infty} \langle T_nSx,y\rangle 	=\lim_{n\to\infty} \langle ST_n x,y\rangle = \langle STx,y\rangle.
			\]
			Since $x,y \in\mathscr{H}$ are arbitrary, it follows that $TS=ST$; and since $S\in \mathscr{U}$ was arbitrary, it follows that $T\in\mathscr{U}^\prime$.
			\item If $\mathscr{U}$ is self-adjoint, then $\mathscr{U}^\prime$ is a von Neumann algebra (A.5 in \cite{DixmierC}).
		\end{list}
	\end{rem}
	
	\begin{theorem}[Double Commutant---5.3.1 in \cite{KadRing}] If $\mathfrak{U}$ is a $*$-subalgebra of $\mathscr{B}(\mathscr{H})$ containing the identity operator, then the strong-operator closure of $\mathfrak{U}$ coincides with $\mathfrak{U}^{\prime\prime}$.
	\end{theorem}
	\begin{rem}Let $\mathfrak{U}$ be a $*$-subalgebra of $\mathscr{B}(\mathscr{H})$ containing the identity.
		\begin{list}{(\arabic{ctr})}{\usecounter{ctr}}\label{doubleComCor}
			\setcounter{ctr}{0}
			\item We showed in \ref{WOTClosed} that $\mathfrak{U}^{\prime\prime}$ is weakly-operator closed. Hence we have 
			\[
			\mathfrak{U}\subset \overline{\mathfrak{U}}^{\|\cdot\|} \subset \overline{\mathfrak{U}}^{\mathrm{SOT}} \subset \overline{\mathfrak{U}}^{\mathrm{WOT}} \subset \mathfrak{U}^{\prime\prime}
			\] 
			Thus the above theorem shows that the weak-operator and strong-operator closures of $\mathfrak{U}$ coincide (a fact we already knew by Remark \ref{convexclosure}).
			\item It also follows from the above theorem that $\mathfrak{U}$ is a von Neumann algebra if and only if $\mathfrak{U}=\mathfrak{U}^{\prime\prime}$. Some authors take this property as the definition of a von Neumann algebra.
		\end{list}
	\end{rem}
	\begin{cor} Another corollary of the above theorem is that if $\mathscr{U} \subset \mathscr{B}(\mathscr{H})$ is self-adjoint, then $\mathscr{U}^{\prime\prime}$ is the von Neumann algebra generated by $\mathscr{U}$. That is, 
		\[
		\mathscr{U}^{\prime\prime} = \mathscr{U}^{-}
		\]
	\end{cor}
	
	\begin{defn}
		An element $E\in \mathscr{B}(\mathscr{H})$ which satisfies $E^2 = E =E^*$ is called a \textit{projection}. That is, a projection is a self-adjoint idempotent in $\mathscr{B}(\mathscr{H})$.
	\end{defn}
	It is well-known that a projection $E\in \mathscr{B}(\mathscr{H})$ is the \textit{orthogonal projection} onto the closed subspace $\{\xi \in\mathscr{H}: E\xi = \xi\}$---see \S 2.5 in \cite{KadRing}.
	\begin{prop}\label{ProjDense}
		Let $\{\mathscr{M},\mathscr{H}\}$ be a von Neuamnn algebra, where $\mathscr{H}$ is a complex Hilbert space. Then the projections in $\mathscr{M}$ form a total subset of $\mathscr{M}$ with respect to the norm topology. Thus any von Neumann algebra over a complex Hilbert space is generated by its projections.
	\end{prop}
	\begin{proof} Since $\mathscr{H}$ is a complex Hilbert space, we can express any operator $T$ in $\mathscr{M}$ as a linear combination of self-adjoint operators in the usual way (see \cite{KadRing}, p.~105). Thus, we may assume for simplicity that $A$ is self-adjoint. Furthermore, letting $\mathscr{R}\subset\mathscr{M}$ denote the (abelian) von Neumann algebra generated by $T$ and $I$, we have a \textit{spectral resolution} of $T$; that is, a family $\{E_\lambda\}$ of projections, indexed by $\mathbb{R}$, in $\mathscr{R}$, satisfying certain properties (see 5.2.2 in \cite{KadRing}). Among these is the following: for each partition $\{\lambda_0,\dots,\lambda_n\}$ such that $\lambda_0<-\|T\|$ and $\lambda_n = \|T\|$, we form the Riemann sum $\sum_{j=1}^n \lambda_j^\prime E(j)$, where $E(j) = E_{\lambda_{j}} -E_{\lambda_{j-1}}$, and $\lambda_j^\prime$ is taken to be any point in $[\lambda_{j-1}, \lambda_j]\cap \mathrm{sp}(T)$ when this intersection is nonempty, and $\lambda_{j-1}$ otherwise; then the $E(j)$ are mutually orthogonal projections, and the Riemann sum approximates $T$ in norm; and this approximation converges to $T$ in norm as we take finer partitions.
	\end{proof}
\begin{defn}\label{center}
		 The \textit{center} of a von Neumann algebra $\mathscr{M}$---which we denote $\mathscr{Z}(\mathscr{M})$---is the set of all operators in $\mathscr{M}$ that commute with all other operators in $\mathscr{M}$: that is, $\mathscr{Z}(\mathscr{M}) = \mathscr{M}\cap \mathscr{M}^\prime$.
		 	\end{defn}
\begin{rem}		
 $\mathscr{Z}(\mathscr{M})$ is itself a von Neumann algebra, and  $\mathscr{Z}(\mathscr{M}) =  \mathscr{Z}(\mathscr{M}^\prime)$.
\end{rem}

	\begin{defn} A \textit{factor} is a von Neumann algebra $\mathscr{M}$ with trivial center $\mathscr{Z}(\mathscr{M})= \mathbb{C}I$.
	\end{defn}
	\begin{defn}
		We sometimes write $\mathbb{C}_\mathscr{H}$ for the von Neumann algebra $\mathbb{C}I_\mathscr{H}$ consisting of scalar multiples of the identity operator $I_\mathscr{H}$ on $\mathscr{H}$.
	\end{defn}
	\begin{lemma}\label{closureofconjugate}
		Let $\mathfrak{U}$ be a $*$-subalgebra of $\mathscr{B}(\mathscr{H})$, and let $T \in \mathscr{B}(\mathscr{H},\mathscr{K})$ be a unitary map. Then
		\[
		T\overline{\mathfrak{U}} T^{-1}  = \overline{T\mathfrak{U}T^{-1}},
		\]
		where the bar over the algebras indicates closure in the weak-operator topology.
	\end{lemma}
	\begin{proof}
		Suppose that $A \in \overline{\mathfrak{U}}$ and let $\{A_n\} \subset \mathfrak{U}$ with $A_n \to A$ in the weak-operator sense. Then for any $x,y\in\mathscr{K}$, we have
		\begin{align*}
			\langle (TA_nT^{-1} - TAT^{-1})x,y\rangle &= \langle (T(A_n - A)T^{-1})x,y\rangle\\
			&= \langle (A_n-A)(T^{-1}x),T^{-1}y\rangle\\
			&\to 0.
		\end{align*}
		Thus $T\overline{\mathfrak{U}} T^{-1} \subset \overline{T\mathfrak{U}T^{-1}}$. Conversely, if $TA_nT^{-1} \to TAT^{-1}$ in the weak-operator topology on $\mathscr{B}(\mathscr{K})$, then for $x,y\in\mathscr{H}$, we have $x= T^{-1}x_0$ and $y=T^{-1}y_0$ for some $x_0,y_0 \in \mathscr{K}$, and 
		\[
		\langle (A_n-A)x,y\rangle =	\langle (A_n-A)(T^{-1}x_0),T^{-1}y_0\rangle = \langle (TA_nT^{-1} - TAT^{-1})x_0,y_0\rangle  \to 0,
		\]
		whence $\overline{T\mathfrak{U}T^{-1}} \subset T\overline{\mathfrak{U}} T^{-1}$.
	\end{proof}
	\begin{defn}
		Two von Neumann algebras, $\{\mathscr{M}_1,\mathscr{H}_1\}, \{\mathscr{M}_2,\mathscr{H}_2\}$ are said to be \textit{unitarily equivalent} if there exists a unitary map $U: \mathscr{H}_1\to\mathscr{H}_2$ such that
		\[
		\mathscr{M}_2 = \{UTU^{-1}: T\in\mathscr{M}_1\};
		\]
		i.e.: $\mathscr{M}_2 = U\mathscr{M}_1U^{-1}$. 
	\end{defn}
	
	\begin{lemma}\label{centerlemma}
		Suppose that $\{\mathscr{M}_1,\mathscr{H}_1\}$ and $\{\mathscr{M}_2,\mathscr{H}_2\}$ are unitarily equivalent von Neumann algebras, with intertwining operator $U$: i.e., $\mathscr{M}_2 = U\mathscr{M}_1U^{-1}$. Then 
		\[
		\mathscr{Z}(\mathscr{M}_2) = U\mathscr{Z}(\mathscr{M}_1)U^{-1}.
		\]
	\end{lemma}
	\begin{proof}
		We have
		\begin{align*}
			\mathscr{Z}(\mathscr{M}_2) &=\mathscr{Z}(U\mathscr{M}_1U^{-1})\\
			&= \{A \in\mathscr{M}_2: AB=BA\ \text{for all}\ B\in U\mathscr{M}_1U^{-1}\}\\
			&= \{A \in\mathscr{M}_2: AUTU^{-1} = UTU^{-1}A\ \text{for all}\ T\in \mathscr{M}_1\}\\
			&=\{A \in\mathscr{M}_2: U^{-1}AUT= TU^{-1}AU\ \text{for all}\ T\in \mathscr{M}_1\}\\
			&=\{A \in\mathscr{M}_2: U^{-1}AU\in \mathscr{Z}(\mathscr{M}_1)\}\\
			&= U\mathscr{Z}(\mathscr{M}_1)U^{-1}.
		\end{align*}
	\end{proof}
	\begin{defn}[Tensor product of von Neumann algebras]\label{vnTensorDef} Let 
		\[
		\{\mathscr{M}_1,\mathscr{H}_1\}, \dots,\{\mathscr{M}_n,\mathscr{H}_n\}
		\]
		be von Neumann algebras, and let $\mathscr{H}=\mathscr{H}_1\otimes \cdots\otimes \mathscr{H}_n$\label{Hilberttensor}. Let $\mathscr{M}_0$ denote the $*$-subalgebra of $\mathscr{B}(\mathscr{H})$ consisting of all finite sums of linear operators of the form $T_1\otimes\cdots \otimes T_n$\label{operatortensor}, where $T_i\in \mathscr{M}_i$ for $i=1,\dots,n$. We define the \textit{(von Neumann algebra) tensor product}, $\mathscr{M}_1 \overline{\otimes}  \cdots\overline{\otimes} \mathscr{M}_n$\label{vntensor}, to be the von Neumann algebra $\mathscr{M}_0^{-}$ generated by $\mathscr{M}_0$.
	\end{defn}
	
	\begin{prop}[Elementary Properties]\label{VNTensorElem}\hspace{2cm}
		\begin{list}{(\arabic{ctr})}{\usecounter{ctr}}
			\setcounter{ctr}{0}
			\item If $\{\mathscr{M}_i,\mathscr{H}_i\}_{i=1}^n$ are von Neumann algebras and $\mathscr{U}_i \subset \mathscr{B}(\mathscr{H}_i)$ with $\mathscr{M}_i = \mathscr{U}_i^{-}$ for $i=1,\dots, n$, then $\mathscr{M}_1 \overline{\otimes}  \cdots\overline{\otimes} \mathscr{M}_n$ is the von Neumann algebra generated by the set 
			\[
			\mathscr{U} = \{T_1\otimes\cdots \otimes T_n: T_1 \in\mathscr{U}_1,\dots, T_n \in\mathscr{U}_n\}.
			\]
			\item  The tensor product of von Neumann algebras is associative: i.e.:
			\[
			(\mathscr{M}_1\overline{\otimes}\cdots \overline{\otimes} \mathscr{M}_m)\overline{\otimes} (\mathscr{M}_{m+1}\overline{\otimes}\cdots\overline{\otimes} \mathscr{M}_n) = \mathscr{M}_1\overline{\otimes} \cdots \overline{\otimes}\mathscr{M}_n.
			\]
			\item Suppose that $\{\mathscr{M},\mathscr{H}\}$ is a von Neumann algebra, and $\mathbb{C}_\mathscr{K}$ is the von Neumann algebra consisting of scalar multiples of the identity operator $I_\mathscr{K}$ on $\mathscr{K}$. Then 
			\[
			\mathscr{M}\overline{\otimes}\mathbb{C}_\mathscr{K} = \{T\otimes I_\mathscr{K}: T\in \mathscr{M}\}.
			\]
		\end{list}
	\end{prop}
	
	\begin{theorem}[The Commutation Theorem]\label{Commutation} Let $\mathscr{M}_1$ and $\mathscr{M}_2$ be von Neumann algebras. Then 
		\[
		(\mathscr{M}_1 \overline{\otimes}\mathscr{M}_2)^\prime = \mathscr{M}_1^\prime \overline{\otimes}\mathscr{M}_2^\prime.
		\]
		I.e., the commutant of the tensor product is the tensor product of the commutants. By associativity, the theorem holds for the tensor product of $n$ von Neumann algebras.
	\end{theorem}
	
	\begin{rem} This is a hard theorem. It was first proved by M. Tomita in 1967, using the theory of unbounded operators. An elementary proof that does not rely on unbounded operator theory was given by Rieffel and van Daele in 1975 (\cite{RieffelvanDaele}). See also \cite{KadRing}, pp.~821--828.
	\end{rem}

	\section[Representations of Normed \texorpdfstring{$*$-algebras}{Star Algebras} and Locally Compact \\ Groups]{Representations of Normed \texorpdfstring{$*$-algebras}{Star Algebras} and Locally Compact Groups}
	
	Some of the references we cite use representations of $\mathrm{C}^*$-algebras, so we discuss them briefly here. More generally, we define representations of $*$-algebras as follows:
	\begin{defn}
		 Given a normed $*$-algebra $\mathscr{A}$, we define a \textit{$*$-representation of $\mathscr{A}$} (or simply a \textit{representation of $\mathscr{A}$}) to be a pair $(\pi,\mathscr{H}_{\pi})$, where $\mathscr{H}_{\pi}$ is a Hilbert space and $\pi$ is a homomorphism of $\mathscr{A}$ into $\mathscr{B}(\mathscr{H}_{\pi})$ such that $\pi(a^*)=\pi(a)^*$ for all $a\in \mathscr{A}$.
	\end{defn}
		As with unitary representations of locally compact groups, \textit{dimension} of a representation $(\pi,\mathscr{H}_\pi)$ of $\mathscr{A}$ is just the dimension of $\mathscr{H}_\pi$. We often just write $\pi$ for $(\pi,\mathscr{H}_{\pi})$ when the Hilbert space is understood. It follows from spectral theory that any $*$-representation $\pi$ of $\mathscr{A}$ is continuous (see Theorem 4.1.8 in \cite{KadRing}).		
\begin{defn}
		We say that a representation $\pi$ of $\mathscr{A}$ is \textit{nondegenerate} if $\{\pi(a)\xi: a\in \mathscr{A}, \xi\in \mathscr{H}_{\pi}\}$ is total in $\mathscr{H}_{\pi}$; that is, if the closure of its linear span is all of $\mathscr{H}_{\pi}$. 
	\end{defn}
\begin{rem}
	The condition of $\{\pi(a)\xi: a\in \mathscr{A}, \xi\in \mathscr{H}_{\pi}\}$ being total in $\mathscr{H}_{\pi}$ is equivalent to requiring that for any nonzero $\xi\in\mathscr{H}_{\pi}$, there exists an $a\in \mathscr{A}$ such that $\pi(a)\xi \neq 0$ (see Proposition I.9.2 in \cite{Tak}). All the $*$-representations of normed $*$-algebras that we encounter in this thesis will be nondegenerate.
	\end{rem}
	\begin{rem}
		Some of the references we cite will mention the nondegeneracy condition, at least in their proofs, which is why we include it here. In particular, if $I \in \pi(\mathscr{A})$, then $\pi$ is nondegenerate, and moreover, $\pi(\mathscr{A})^{\prime\prime}$ is the von Neumann algebra generated by $\pi(\mathscr{A})$, by the Double Commutant Theorem.
	\end{rem}
	\begin{defn}
		Let $(\pi,\mathscr{H}_\pi)$ be a representation of a normed $*$-algebra $\mathscr{A}$. Then $\pi$ is called \textit{cyclic} with \textit{cyclic vector} $\xi$ if the set $\pi(\mathscr{A})\xi=\{\pi(a)\xi:a\in\mathscr{A}\}$ is dense in $\mathscr{H}_\pi$.
	\end{defn}
		\begin{defn}
		Let $(\pi,\mathscr{H}_\pi)$ be a unitary representation of a locally compact group $G$. Then $\pi$ is called \textit{cyclic} with \textit{cyclic vector} $\xi$ if the set $\{\pi(x)\xi:x\in G\}$ is total in $\mathscr{H}_\pi$; that is, if the closed linear span of this set is $\mathscr{H}_\pi$.
	\end{defn}
\begin{defn}
	Given two $*$-representations $\pi$ and $\sigma$ of a normed $*$-algebra $\mathscr{A}$, we define the \textit{intertwining space} of $\pi$ and $\sigma$ to be 
			\[
			\mathrm{Hom}_\mathscr{A}(\pi,\sigma) = \{T\in\mathscr{B}(\mathscr{H}_\pi, \mathscr{H}_\sigma): T\pi(a)=\sigma(a)T, \ \ \text{for all}\ a\in \mathscr{A}\}.
			\]
			\end{defn}
			Note that $\mathrm{Hom}_\mathscr{A}(\pi,\pi) = \pi(\mathscr{A})^\prime$. We make an analogous definition for representations of a locally compact group:
			\begin{defn}
			 Likewise, given two continuous unitary representations $\pi$ and $\sigma$ of a locally compact group $G$, we define the \textit{intertwining space} $\mathrm{Hom}_G(\pi,\sigma)$\label{intspace} in a directly analogous way, and note that $\mathrm{Hom}_G(\pi,\pi) = \pi(G)^\prime$, just as in the normed $*$-algebra case.
	\end{defn}
	We note that all definitions that follow can be expressed for $*$-representations of normed $*$-algebras, but we state them for unitary representations of locally compact groups---these being the objects most frequently discussed in this thesis. Recall that unitary representations of locally compact groups were defined in \ref{contunitaryrepdef}.
	
	\begin{defn}[Unitary Equivalence] Two representations $(\pi,\mathscr{H}_\pi)$ and $(\sigma,\mathscr{H}_\sigma)$ of a locally compact group $G$ are said to be \textit{unitarily equivalent} (or simply \textit{equivalent}) if there is an unitary operator $U \in \mathrm{Hom}_G(\pi,\sigma)$. We sometimes write $\pi\sim \sigma$\label{equivrep} to express unitary equivalence.
	\end{defn}
	
	\begin{defn}[Invariant Subspace; Irreducible Representation]\label{invarDef} \hspace{2cm}
		\begin{list}{(\arabic{ctr})}{\usecounter{ctr}}
			\setcounter{ctr}{0}
			\item  Let $(\pi,\mathscr{H}_{\pi})$ be a unitary representation of a locally compact group $G$. A closed subspace $K\subset \mathscr{H}$ is called \textit{$\pi$-invariant} if $\pi(x)\xi\in K$, for all $\xi\in\mathscr{H}_{\pi}$ and $x\in G$. 
			\item We say that $\pi$ is \textit{irreducible} if $\{0\}$ and $\mathscr{H}_{\pi}$ are the only $\pi$-invariant closed subspaces of $\mathscr{H}_{\pi}$.
		\end{list}
	\end{defn}
	\begin{rem}\label{EquivInvar}
		Let $(\pi,\mathscr{H}_\pi)$ and  $(\sigma,\mathscr{H}_\sigma)$ be unitary representations of $G$; let $K\subset \mathscr{H}_\pi$ be a closed subspace; and let and let $TK = \{T\xi:\xi\in K\}$. We will often use the fact that if $T\in \mathrm{Hom}_G(\pi,\sigma)$, then $K$ is $\pi$-invariant if and only if $TK$ is $\sigma$-invariant.
	\end{rem}
	
	\begin{prop}[Cf. II.6.1.4 in \cite{Blackadar}]\label{invorth} Let $(\pi,\mathscr{H})$ be a unitary representation, and let $K$ be a closed subspace of $\mathscr{H}$. We write $P_K$\label{PK} for the orthogonal projection onto $K$. The following are equivalent:
		\newcounter{ctr1}
		\begin{list}{(\roman{ctr1})}{\usecounter{ctr1}}
			\setcounter{ctr1}{0}
			\item $K$ is $\pi$-invariant;
			\item $K^\perp$, the orthogonal complement of $K$, is $\pi$-invariant;
			\item $P_K \in \pi(G)^\prime$.
		\end{list}
	\end{prop}
	\begin{proof}
		The equivalence of $(\mathrm{i})$ and $(\mathrm{ii})$ follows easily from the fact that $\pi(G)$ is a self-adjoint subset of $\mathscr{B}(\mathscr{H})$ (because $\pi(x)^* = \pi(x)^{-1} =\pi(x^{-1}))$. To prove that $(\mathrm{i})\implies (\mathrm{iii})$, suppose that $K$ is $\pi$-invariant. Let $x\in G$. We need to show that
		\[
		\pi(x)P_K = P_K\pi(x).
		\]
		For every $\xi \in \mathscr{H}$, we have
		\[
		P_K\pi(x)\xi = P_K[\underbrace{\pi(x)P_K\xi}_{\in K} + \underbrace{\pi(x)(I-P_K)\xi}_{\in K^\perp}] = \pi(x)P_K\xi,
		\]
		as desired. To see that $(\mathrm{iii})\implies (\mathrm{i})$, note that if $\xi\in K$, then 
		\[
		\pi(x)\xi = \pi(x)P_K\xi = P_K\pi(x)\xi \in K.
		\]
	\end{proof}
	\begin{cor}\label{IrrRepCrit}
		A representation $\pi$ of $G$ on $\mathscr{H}$ is irreducible if and only if $\pi(G)^\prime$ consists of scalar multiples of the identity, i.e.: $\pi(G)^\prime = \mathbb{C}I = \mathbb{C}_\mathscr{H}$. It follows that if $G$ is abelian, then every irreducible representation of $G$ is one-dimensional. See \cite{FollandAbstract}, p.~77.
	\end{cor}
	
	\begin{defn}\label{directsumrep}
 Let $\mathscr{H} =  \textstyle{\sum}_{i\in I} \oplus \mathscr{H}_i$\footnote{Direct sums of Hilbert spaces and bounded linear transformations are discussed in \cite{KadRing}, \S~2.6., which also treats tensor products of Hilbert spaces and bounded linear transformations.}, and let $\{(\pi_i,\mathscr{H}_i):i=1,\dots,n\}$ be unitary representations of $G$. If $\pi(x) = \textstyle{\sum}_{i\in I}  \pi_i(x)$ for all $x\in G$, then we write $\pi = \textstyle{\sum}_{i\in I}  \pi_i$ and say that $\pi$ is the \textit{direct sum} of the $\pi_i$.
	\end{defn} 
		\begin{defn}
		If $\mathscr{H}_i = \mathscr{H}_0$ and $\pi_i = \rho$ for all $i\in I$, for some fixed representation $\rho$, then we say that $\pi$ is a \textit{multiple} or \textit{amplification} of $\rho$ by $\mathrm{card}(I)$. When  $\mathrm{card}(I) = m$, we sometimes write $m\rho$\label{mrho} for the amplification of $\rho$ by $m$.
	
	\end{defn}
	
	\begin{defn}
	Let $(\pi,\mathscr{H})$ be a unitary representation of $G$. If $K$ is a $\pi$-invariant subspace of $\mathscr{H}$, then $\pi\vert_K$ gives a representation of $G$ on $K$. Here $\pi\vert_K$ denotes the restriction of $\pi$ to $K$---that is, $\pi\vert_K (x) = \pi(x)\vert_K$. Any representation obtained in this way is called a \textit{subrepresentation} of $G$.
	\end{defn}
	If $K$ is a $\pi$-invariant subspace of $\mathscr{H}$, then it follows from Proposition \ref{invorth} that $\mathscr{H} = K\oplus K^\perp$, and that $\pi = \pi\vert_K \oplus \pi\vert_{K^\perp}$.
	\begin{defn}[Inner tensor product]\label{innertensor} 
		Given two representations $(\pi_1,\mathscr{H}_1)$ and $(\pi_2,\mathscr{H}_2)$ of a locally compact group $G$, we define their \textit{inner tensor product} $\pi_1\otimes \pi_2$ to be the unitary representation on $\mathscr{H}_1\otimes \mathscr{H}_2$ given by
		\[
		(\pi_1\otimes \pi_2)(x) = \pi_1(x) \otimes \pi_2(x),\ \ \ \text{for all}\ x\in G.
		\]
	\end{defn}
	In the special case where one of the representations in the above definition is the trivial representation $I$, we have that $\pi\otimes I$ is isomorphic to a direct sum of $d$ copies of $\pi$, where $d = \dim\mathscr{H}$. Thus $\pi\otimes I$ is the amplification of $\pi$ by $d$.
	The following theorem will be used at a later point:
	\begin{theorem}\label{commutantfact}
		Suppose that $\rho$ is an irreducible representation of $G$ on $\mathscr{H}_1$ and $I$ is the trivial representation on $\mathscr{H}_2$. Suppose $P$ is a projection in $[(\rho\otimes I)(G)]^\prime$. Then $P = I \otimes P_V$ for some closed subspace $V$ of $\mathscr{H}_2$.
	\end{theorem}
	
	\begin{proof}
		Note that
		\begin{align*}
			[(\rho\otimes I)(G)]^\prime &= \{T \in \mathscr{B}(\mathscr{H}_1\otimes \mathscr{H}_2): T(\rho\otimes I)(x) = (\rho\otimes I)(x)T \ \ \text{for all}\ x\in G\}\\
			&= \{T \in \mathscr{B}(\mathscr{H}_1\otimes \mathscr{H}_2): T(\rho(x)\otimes I) = (\rho(x)\otimes I)T \ \ \text{for all}\ x\in G\}\\
			&= \mathcal{T}^\prime,
		\end{align*}
		where $\mathcal{T} = \{\rho(x)\otimes I: x\in G\}$. Now note that
		\[
		\mathcal{T}^\prime = (\mathcal{T}^{\prime\prime})^\prime =  (\mathcal{T}^{-})^\prime.
		\]
		Letting $\mathcal{S}= \{\rho(x): x\in G\}$, it follows from the elementary property (1) in (\ref{VNTensorElem}) that
		\begin{align*}
			\mathcal{T}^{-} &= \{\rho(x)\otimes I: x\in G\}^- \\
			&= \{S\otimes I: S\in \mathcal{S} \}^- \\
			&= \mathcal{S}^- \,\overline{\otimes}\, \mathbb{C}_{\mathscr{H}_2}.
		\end{align*}
		Then it follows from the Commutation Theorem (\ref{Commutation}) that 
		\[
		(\mathcal{T}^{-})^\prime = (\mathcal{S}^-)^\prime \,\overline{\otimes}\, \mathbb{C}_{\mathscr{H}_2}^\prime = \mathcal{S}^\prime \,\overline{\otimes} \,\mathscr{B}(\mathscr{H}_2) = [\rho(G)]^\prime \,\overline{\otimes}\, \mathscr{B}(\mathscr{H}_2).
		\]
		Thus 
		\[
		[(\rho\otimes I)(G)]^\prime =  [\rho(G)]^\prime \,\overline{\otimes}\, \mathscr{B}(\mathscr{H}_2).
		\]
		Now since $\rho$ is irreducible, we have $[\rho(G)]^\prime = \mathbb{C}_{\mathscr{H}_1}$, by Corollary \ref{IrrRepCrit}, so that
		\begin{equation*}
		[(\rho\otimes I)(G)]^\prime = \mathbb{C}_{\mathscr{H}_1}\,\overline{\otimes}\, \mathscr{B}(\mathscr{H}_2)
		\end{equation*}
		In particular, if $P$ is a projection in $[(\rho\otimes I)(G)]^\prime$, then $P$ is the orthogonal projection onto its range space $W=\{Px: x\in  \mathscr{H}_1\otimes \mathscr{H}_2\}$. And $P \in  \mathbb{C}_{\mathscr{H}_1}\,\overline{\otimes}\, \mathscr{B}(\mathscr{H}_2)$ implies---by property (3) in \ref{VNTensorElem}---that $P = I \otimes T$ for some $T \in \mathscr{B}(\mathscr{H}_2)$. And $P$ being a projection then implies that $T = P_V$ for some closed subspace $V \subset \mathscr{H}_2$. Thus $P = I \otimes P_V$, as claimed.
	\end{proof}
	\begin{cor}\label{CorforTensorSubspaces}
		Suppose that $\rho$ is an irreducible representation of $G$ on $\mathscr{H}_1$ and $I$ is the trivial representation on $\mathscr{H}_2$. Suppose $P$ is a projection in $[(\rho\otimes I)(G)]^\prime$. Then the invariant subspaces for $\rho\otimes I$ are precisely $\{\mathscr{H}_1 \otimes V: V\subset \mathscr{H}_2\ \text{a closed subspace}\}$, so that the subrepresentations of $\rho\otimes I$ are precisely $\{\rho \otimes (I\vert_V): V\subset \mathscr{H}_2\ \text{a closed subspace} \}$.
	\end{cor}
	The powerful Commutation theorem is not necessary to prove Theorem \ref{commutantfact}. Indeed, the elementary lemma below can be used:
	\begin{lemma}[Theorem 7.11 in \cite{FollandAbstract}]\label{FollandCommuteLemma}
		Suppose $\mathscr{H}_1$ and $\mathscr{H}_2$ are Hilbert spaces, and $\mathcal{S}$ is a subset of $\mathscr{B}(\mathscr{H}_1)$ such that the only bounded operators on $\mathscr{H}_1$ that commute with every $S\in\mathcal{S}$ are scalar multiples of $I$. Then the bounded operators on $\mathscr{H}_1 \otimes \mathscr{H}_2$ that commute with $S\otimes I$ for every $S \in \mathcal{S}$ are precisely those of the form $I\otimes T$ where $T\in \mathscr{B}(\mathscr{H}_2)$.
	\end{lemma}
	
	\begin{proof}[Elementary Proof of \ref{commutantfact}]
		As before, we have
		\begin{align*}
			[(\rho\otimes I)(G)]^\prime &= \{T \in \mathscr{B}(\mathscr{H}_1\otimes \mathscr{H}_2): T(\rho(x)\otimes I) = (\rho(x)\otimes I)T \ \ \text{for all}\ x\in G\}\\
			&= \{T \in \mathscr{B}(\mathscr{H}_1\otimes \mathscr{H}_2): T(S\otimes I) = (S\otimes I)T \ \ \text{for all}\ S\in \mathcal{S}\},
		\end{align*}
		where $\mathcal{S}= \{\rho(x): x\in G\}$. Since $\rho$ is irreducible, the only bounded operators on $\mathscr{H}_1$ that commute with every $S\in\mathcal{S}$ are scalar multiples of the identity. It follows by Lemma \ref{FollandCommuteLemma} that 
		\begin{equation}\label{tensorvnalg}
		[(\rho\otimes I)(G)]^\prime=  \{I\otimes T: T\in \mathscr{B}(\mathscr{H}_2)\}= \mathbb{C}_{\mathscr{H}_1} \,\overline{\otimes}\, \mathscr{B}(\mathscr{H}_2),
	\end{equation}
		and the rest follows as before.
	\end{proof}
	\begin{defn}[Factor Representations]
		A unitary representation $\pi$ of $G$ on $\mathscr{H}$ is called a \textit{factor representation} (or sometimes a \textit{primary representation}) if $\pi(G)^-$ is a factor. I.e., $\pi$ is a factor representation if $\mathscr{Z}(\pi(G)^-) = \mathbb{C}_{\mathscr{H}}$.
	\end{defn}
		\begin{cor}\label{tensorrepfactor}
		Let $\rho$ be an irreducible representation of $G$ on $\mathscr{H}_1$, and let $I$ be the trivial representation on $\mathscr{H}_2$. Then $\rho\otimes I$ is a factor representation of $G$ on $\mathscr{H}_1\otimes \mathscr{H}_2$.
	\end{cor}
	\begin{proof}
		It follows from \eqref{tensorvnalg} and Examples 11.1.4 and 11.2.1 in \cite{KadRing} that
		\[
		[(\rho\otimes I)(G)]^- = 	[(\rho\otimes I)(G)]^{\prime\prime} = [\mathbb{C}_{\mathscr{H}_1} \,\overline{\otimes}\, \mathscr{B}(\mathscr{H}_2)]^\prime= \mathscr{B}(\mathscr{H}_1)\,\overline{\otimes}\, \mathbb{C}_{\mathscr{H}_2}.
		\]
	And
		\[
			[[(\rho\otimes I)(G)]^-]^\prime = [(\rho\otimes I)(G)]^\prime = \mathbb{C}_{\mathscr{H}_1} \,\overline{\otimes}\, \mathscr{B}(\mathscr{H}_2).
		\]
Therefore,
\begin{align*}
\mathscr{Z}([(\rho\otimes I)(G)]^-) &=[(\rho\otimes I)(G)]^-\cap 	[[(\rho\otimes I)(G)]^-]^\prime\\
&= \big(\mathscr{B}(\mathscr{H}_1)\,\overline{\otimes}\, \mathbb{C}_{\mathscr{H}_2}\big)\cap \big(\mathbb{C}_{\mathscr{H}_1} \,\overline{\otimes}\, \mathscr{B}(\mathscr{H}_2)\big)\\
&= \{S\otimes I_{\mathscr{H}_2}: S\in \mathscr{B}(\mathscr{H}_1) \} \cap \{I_{\mathscr{H}_1} \otimes T: T\in \mathscr{B}(\mathscr{H}_2)\}\\
&= \mathbb{C}_{\mathscr{H}_1\otimes \mathscr{H}_2}.
\end{align*}
	\end{proof}
	\section{Direct Integral Theory}
	We now introduce the relevant theory of direct integrals. The direct integral construct is a generalization of the concept of direct sum of Hilbert spaces. It involves indexing over a measure space $(X,\mathcal{M},\mu)$ and allows us to decompose representations into factor representations in a canonical way. To do this, we must find the right notion of measurability for a function landing in different Hilbert spaces. Different authors take slightly different approaches to this; we shall mostly follow Takesaki \cite{Tak}, Folland \cite{FollandAbstract}, and Dixmier \cite{DixmierC}. In this section we will refer the reader to these references for the proofs.
	
	\begin{defn}
		 Let $(X,\mathcal{M})$ is a measurable space equipped with a $\sigma$-finite, (positive) measure $\mu$. We call a family $\{\mathscr{H}_x\}_{x\in X}$ of Hilbert spaces a \textit{field of Hilbert spaces} over $X$. An element $\xi \in \Pi_{x\in X} \mathscr{H}_x$ is called a \textit{vector field} on $X$; some authors call this a \textit{section} of the field $\{\mathscr{H}_x\}_{x\in X}$.
		\end{defn}
		\begin{defn}
			A \textit{measurable field of Hilbert spaces} on $(X,\mathcal{M}, \mu)$ is defined to be a field of Hilbert spaces $\{\mathscr{H}_x\}_{x\in X}$ together with a subspace $\mathcal{E}$ of the product vector space $\Pi_{x\in X} \mathscr{H}_x$ with the following properties:
		\newcounter{ctr2}
		\begin{list}{(\roman{ctr2})}{\usecounter{ctr2}}
			\setcounter{ctr2}{0}
			\item For any $\xi \in \mathcal{E}$, the function $x\in X \mapsto \|\xi(x)\|$ is $\mu$-measurable.
			\item For any $\eta \in \Pi \mathscr{H}_x$, if the function $x\in X \mapsto \langle \xi(x),\eta(x)\rangle \in \mathbb{C}$ is $\mu$-measurable for every $\xi\in\mathcal{E}$, then $\eta\in\mathcal{E}$. Vector fields in $\mathcal{E}$ are called \textit{measurable}.
			\item There exists a countable sequence $\{\xi_n\}_{n=1}^\infty \subset \mathcal{E}$ such that for each $x\in X$, $\mathscr{H}_x =\overline{\text{span}}\{\xi_n(x)\}$. The countable collection $\{\xi_n\}_{n=1}^\infty$ is called a \textit{fundamental sequence of measurable vector fields.}
			\end{list}
		\end{defn}
	\begin{rem}	\label{directintpolarmeas}\hspace{2cm}
		\begin{list}{(\alph{ctr})}{\usecounter{ctr}}
			\setcounter{ctr}{0}
			\item Condition (iii) above implies that each $\mathscr{H}_x$ is separable.
			\item It follows from the polarization identity and (i) above that the function $x\mapsto \langle\xi(x), \eta(x)\rangle_x$ is $\mu$-measurable for all $\xi, \eta\in \mathcal{E}$.
		\end{list}
	\end{rem}
	\begin{rem}\label{FollandMeasField} Folland \cite{FollandAbstract} starts with a simpler definition for a measurable field of Hilbert spaces: he requires only the existence of a sequence $\{\xi_n\}_{n=1}^\infty$ of vector fields such that
		\begin{list}{(\arabic{ctr})}{\usecounter{ctr}} 
			\setcounter{ctr}{0}
			\item $x\mapsto \langle \xi_m(x),\xi_n(x)\rangle_x$ is measurable for all $m,n\in\mathbb{N}$.
			\item For each $x\in X$, $\mathscr{H}_x=\overline{\text{span}}\{\xi_n(x)\}$.
		\end{list}
		The existence of such a sequence implies the existence of a subspace $\mathcal{E}\subset \Pi_{x\in X} \mathscr{H}_x$ satisfying (i)--(iii) in the above definition; namely, we can take
		\begin{equation}\label{Follandmeasfieldscriteria}
			\mathcal{E} = \{\xi\in \Pi_{x\in X} \mathscr{H}_x: x\mapsto \langle \xi(x),\xi_n(x) \rangle_x \ \text{is}\ \mu\text{-measurable for every}\ n \}.
		\end{equation}
		This is Lemma 8.10 in \cite{Tak}. Thus
		
		\begin{tcolorbox}
			\begin{center}
				A sequence of vector fields $\{\xi_n\}_{n\in \mathbb{N}}$ satisfying (1) and (2) ensures that $\{\mathscr{H}_x\}_{x\in X}$ is a measurable field of Hilbert spaces with fundamental sequence $\{\xi_n\}$ and measurable vector fields given by \eqref{Follandmeasfieldscriteria}.
			\end{center}
		\end{tcolorbox}
	\end{rem}
	
	For the remainder of this section, we assume that $(X,\mathcal{M})$ is a measure space equipped with a $\sigma$-finite, (positive) Borel measure $\mu$. We also assume that $\{\mathscr{H}_x\}$ is a measurable field of Hilbert spaces on $(X,\mathcal{M},\mu)$, without specifying the family $\mathcal{E}$ of measurable vector fields unless necessary to avoid confusion.
	
	\begin{defn}[Direct integral of Hilbert Spaces]
		Let $\{\mathscr{H}_x\}$ be a measurable field of Hilbert spaces on $(X,\mathcal{M},\mu)$ and let $\mathscr{H}_0$ be the collection of measurable vector fields $\xi$ such that 
		\begin{equation}\label{directnorm}
			\|\xi\| = \Big(\int_X \|\xi(x)\|_x^2\ud\mu(x) \Big)^{1/2} <\infty.
		\end{equation}
		Then $\mathscr{H}_0$ is a vector space with addition and scalar multiplication defined pointwise. The sesquilinear form 
		\begin{equation}\label{sesqui}
			\langle \xi, \eta\rangle = \int_{X} \langle \xi(x),\eta(x)\rangle_x \ud \mu(x),\ \ \ \text{for}\ \xi, \eta\in \mathscr{H}_0
		\end{equation}
		is a pre-inner product (sometimes called a semi-inner product) on $\mathscr{H}_0$. Identifying vector fields which are equal $\mu$-a.e. makes \eqref{sesqui} into an inner product on $\mathscr{H}_0$. Thus, defining an equivalence relation $\sim$ on $\mathscr{H}_0$ by $\xi\sim \eta \iff \xi=\eta$ a.e., we get an inner product space $\mathscr{H} = \mathscr{H}_0/\sim$. An argument similar to the usual proof that $L^2(\mu)$ is complete shows that $\mathscr{H}$ is complete, and hence, a Hilbert space. We call $\mathscr{H}$ the \textit{direct integral} of the measurable field of Hilbert spaces $\mathscr{H}_x$ with respect to $\mu$, and write\label{directinthilbspace}
		\[
		\mathscr{H}= \int_X^\oplus \mathscr{H}_x \ud\mu(x).
		\]
		\begin{rem}It is common to abuse notation and write $\xi\in\mathscr{H}$ to mean a representative of the equivalence class of square-integrable, measurable vector fields equal to $\xi$ almost everywhere. When we wish to be more precise, we sometimes use brackets $[\xi]$ to denote this equivalence class.
		\end{rem}
		\begin{notation}
	We write each $\xi\in\mathscr{H}$ as $\xi=\int_X^\oplus \xi(x)\ud\mu(x)$\label{directfield}. When the measure $\mu$ is understood, we sometimes omit it from the integral to simplify notation.
		\end{notation} 
	\end{defn}
	
	\begin{ex}\label{constantField}
		Given a field of identical, separable Hilbert spaces---say $\mathscr{H}_x = \mathscr{H}$ for all $x$---the orthonormal basis $\{e_n\}_{n=1}^\infty$ of $\mathscr{H}$ gives a fundamental sequence of measurable vector fields, by taking $e_n(x) = e_n$ for all $x$. This makes $\{\mathscr{H}_x\}_{x\in X}$ a measurable field of Hilbert spaces. The collection of measurable vector fields is given by 
		\[
		\mathcal{E}= \{\xi\in \Pi_{x\in X} \mathscr{H}_x: x\mapsto \langle \xi(x),e_n(x) \rangle_x \ \text{is}\ \mu\text{-measurable for every}\ n \}.
		\] 
		In this case, $\{\mathscr{H}_x:x\in X\}$ is said to be a \textit{constant field} of Hilbert spaces. Moreover, we have
		\[
		\int_X^\oplus \mathscr{H}_x \ud\mu(x) = L^2(X;\mathscr{H}),
		\]
		where the right-hand side is the measurable $\mathscr{H}$-valued square integrable functions, to be defined in Section \ref{MeasRangeSection}. We will make frequent use of this construction in Chapter \ref{DecompChapter}.
	\end{ex}
See Lemma 8.12 in \cite{Tak} for the following.
	\begin{lemma}\label{basissequence}
Let $\{\mathscr{H}_x\}$ be a measurable field of Hilbert spaces on $(X,\mathcal{M},\mu)$. For each $x \in X$, let $d(x) = \dim \mathscr{H}_x \in [0,\infty]$. Then $\{x\in X: d(x) = m\}$ is a measurable set for $m=0,1,2,\ldots ,\infty$. Moreover, there exists a fundamental sequence $\{\xi_n\}_{n=1}^\infty$ of measurable vector fields such that
\begin{list}{(\roman{ctr})}{\usecounter{ctr}}
	\item For each $x\in X$, $\{\xi_1(x),\xi_2(x),\ldots,\xi_{d(x)}(x)\}$ is an orthonormal basis for $\mathscr{H}_x$.
	\item $\xi_n(x) = 0$ for $n>d(x)$, if $d(x)<\infty$.
	\end{list}
	\end{lemma}
	\begin{defn}
	 Let $\{\mathscr{H}_x\}$ and $\{\mathscr{K}_x\}$ be two measurable fields of Hilbert spaces on $(X,\mathcal{M},\mu)$. An operator field $T \in \Pi_{x\in X} \mathscr{B}(\mathscr{H}_x,\mathscr{K}_x)$ is called \textit{measurable} if for any measurable vector field $\xi \in \{\mathscr{H}_x\}$, the vector field $x\mapsto T(x)\xi(x) \in \mathscr{K}_x$ is measurable. 
	\end{defn}
 Given a measurable operator field $T \in \Pi_{x\in X} \mathscr{B}(\mathscr{H}_x,\mathscr{K}_x)$, we define the norm\label{esssup}
	\[
\|T\|_\infty = \underset{x\in X}{\mathrm{ess}\,\sup} \|T(x)\|.
	\] 
If $\|T\|_\infty < \infty$, then $\|T(x)\xi(x)\|_x \leq \|T\|_\infty \|\xi(x)\|_x$ a.e., so that $T$ gives rise to a bounded linear operator---also denoted by $T$---from $\int_X^\oplus \mathscr{H}_x \ud x \to \int_X^\oplus \mathscr{K}_x \ud x$ in the obvious way: for $\xi \in \int_X^\oplus \mathscr{H}_x \ud x$, we define
	\begin{equation}\label{directintegraloperators}
	T\xi = \int_X^\oplus T(x)\xi(x)\ud x \in \int_X^\oplus \mathscr{K}_x \ud x.
	\end{equation}
Thus, $T\xi$ is the (equivalence class of) the measurable vector field $x\mapsto T(x)\xi(x)$. 
\begin{defn}\label{directintops}
The operator \eqref{directintegraloperators} is called the \textit{direct integral} of the bounded measurable field of operators $\{T(x)\}_{x\in X}$, and we write
			\[
			T = \int_X^\oplus T(x)\ud\mu(x).
			\]
	\end{defn}
	\begin{defn}
	The operators $T: \int_X^\oplus \mathscr{H}_x \ud x \to \int_X^\oplus \mathscr{K}_x \ud x$ which can be written as $\int_X^\oplus T(x)\ud\mu(x)$ for some measurable field $\{T(x)\}_{x\in X}$ such that $\|T\|_\infty < \infty$, are called \textit{decomposable}.
	\end{defn}
	\begin{defn}
	 If $\mathscr{H} = \int_X^\oplus \mathscr{H}_x \ud x$ and $T\in\mathscr{B}(\mathscr{H})$ is decomposable with $T(x)$ a scalar for all $x$, then we say that $T$ is a \textit{diagonal operator}; the algebra of all diagonal operators on $\mathscr{H}$ is called the \textit{diagonal algebra}.
	\end{defn}
		\begin{rem}\label{diagonalalgrem}
In the sequel, we usually denote the diagonal algebra by $\mathscr{A}$.\label{diag} The diagonal algebra is a commutative von Neumann algebra (\cite{DixmierC}, p.~388). Given $T\in \mathscr{A}$, we have $T= \int_X^\oplus T(x)\ud\mu(x)$, where $T(x)$ is a scalar for every $x \in X$. Thus we have $T(x) = g(x)I_x$ for some function $g: X\to \mathbb{C}$, where $I_x$ denotes the identity operator on $\mathscr{H}_x$. Letting $S=\{x\in X: \mathscr{H}_x = \{0\}\}$ and letting $f =g\mathbbm{1}_{X-S}$, we have $T(x) = f(x)I_x$ for all $x\in X$. Since $T$ is a measurable field, it follows that $f$ is measurable. To see this, let $\xi_1$ be the first measurable field in the sequence described in Lemma \ref{basissequence}. We have $\|\xi_1(x)\| = 1$ on $X-S$. Since $x\mapsto T(x)\xi_1(x) = f(x)\xi_1(x)$ is a measurable field, we have that
\[
x\mapsto \langle f(x)\xi_1(x), \xi_1(x)\rangle = f(x)
\]
is measurable. Thus we have $T(x) = M_f(x)$ for all $x$, where $M_f(x) = f(x)I_x$ for some measurable $f: X\to \mathbb{C}$, and setting $M_f = \int_X^\oplus M_f(x)\ud\mu(x)$, we have $T = M_f$. Thus
\begin{equation}\label{MF}
(M_f\xi)(x) = f(x)\xi(x) \ \text{for all}\ \xi \in \int_X^\oplus \mathscr{H}_x\ud\mu(x).
\end{equation}
Moreover, one readily sees that $\|M_f(x)\| = |f(x)|$, so that $\|f\|_\infty = \|M_f\|_\infty = \|T\|_\infty < \infty$. Hence,
\[
\boxed{\mathscr{A} =  \{M_f: f \in L^\infty(X)\}}\,,
\]
where $M_f$ is given by \eqref{MF}.
	\end{rem}
	\begin{lemma}\label{directintegralprops}
		Direct integrals of measurable fields of operators respect all the usual algebraic operations: Given two measurable, essentially bounded fields of operators, $\{S(x)\}_{x\in X}$ and  $\{T(x)\}_{x\in X}$ on $\{\mathscr{H}_x\}_{x\in X}$, we have
		\begin{list}{(\arabic{ctr})}{\usecounter{ctr}}
			\setcounter{ctr}{0}
			\item $\int_X^\oplus S(x) \ud x  + \int_X^\oplus T(x) \ud x = \int_X^\oplus S(x) + T(x) \ud x$.
			\item $\Big(\int_X^\oplus S(x)\ud x\Big) \Big(\int_X^\oplus T(x)\ud x\Big) = \int_X^\oplus S(x)T(x)\ud x$.
			\item $\Big(\int_X^\oplus T(x)\ud x\Big)^* = \int_X^\oplus T(x)^* \ud x$. 
			\item Let $I(x)$ denote the identity operator on $\mathscr{H}_x$ for each $x$. Then 
			\[
			I = \int_X^\oplus I(x)\ud\mu(x)
			\] 
			is the identity operator on $\int_X^\oplus \mathscr{H}_x\ud x$.
		\end{list}
	\end{lemma}
	\begin{proof}
		Straightforward from the definitions.
	\end{proof}
	\begin{lemma}[Proposition 7.25 and Corollary 7.27 of \cite{FollandAbstract}]\label{aeequalops} \hspace{2cm}
		\begin{list}{(\arabic{ctr})}{\usecounter{ctr}}
			\setcounter{ctr}{0}
			\item If  $\{T(x)\}_{x\in X}$ is an essentially bounded measurable field of operators, then $\|\int_X^\oplus T(x)\ud x\| = \|T\|_\infty$.
			\item If $\int_X^\oplus S(x)\ud x = \int_X^\oplus T(x)\ud x$, then $S(x) = T(x)$ a.e..
		\end{list}
	\end{lemma}
	
	\begin{prop} Let $(X,\mathcal{M},\mu)$ be a $\sigma$-finite Borel measure space and let $\{\mathscr{H}_x\}_{x\in x}$ and $\{\mathscr{K}_x\}_{x\in x}$ be two measurable fields of separable Hilbert spaces on $(X,\mathcal{M},\mu)$. If $\mathscr{H}_x = \mathscr{K}_x$ for a.e. $x\in X$, then $\int_X^\oplus \mathscr{H}_x\ud\mu(x) \cong \int_X^\oplus \mathscr{K}_x\ud\mu(x)$.
	\end{prop}
	
	\begin{defn}[Measurable fields of von Neumann Algebras]\hspace{.2cm}
		
		Let $\{\mathscr{H}_x\}$ be a measurable field of Hilbert spaces on $(X,\mathcal{M},\mu)$. A field
		\[
		\{\mathscr{M}_x, \mathscr{H}_x\}_{x\in X}
		\]
		of von Neumann algebras is said to be \textit{measurable} if there exists a countable family $\{T_n\}_{n\in\mathbb{N}}$ of measurable fields of operators over $X$ such that $\mathscr{M}_x$ is generated by  $\{T_n(x)\}_{n\in\mathbb{N}}$ for a.e. $x$.
	\end{defn}
	\begin{lemma}
		Let  $\{\mathscr{M}_x, \mathscr{H}_x\}_{x\in X}$ be a measurable field of von Neumann algebras over $X$. Then the field $\{\mathscr{M}_x^\prime, \mathscr{H}_x\}_{x\in X}$ is a measurable, where $\mathscr{M}_x^\prime$ denotes the commutant of $\mathscr{M}_x$, as usual.
	\end{lemma}
	\begin{theorem}[Theorem 8.18 in \cite{Tak}]\label{818Tak} Let  $\{\mathscr{M}_x, \mathscr{H}_x\}_{x\in X}$ be a measurable field of von Neumann algebras on $(X,\mathcal{M},\mu)$. Let $\mathscr{H}=\int_X^\oplus \mathscr{H}_x\ud x$.
		\begin{list}{(\arabic{ctr})}{\usecounter{ctr}}
			\setcounter{ctr}{0}
			\item
			Let $\{\mathscr{M},\mathscr{H}\}$ denote the decomposable fields of operators in $\Pi_x  \{\mathscr{M}_x,\mathscr{H}_x\}$; i.e., the operators of the form $T=\int_X^\oplus T(x)\ud x$ on $\mathscr{H}$, with $T(x) \in \mathscr{M}_x$ for every $x\in X$ and $\|T\|_\infty < \infty$. Then $\mathscr{M}$ is a von Neumann algebra on $\mathscr{H}$; we write
			\begin{equation*}
				\{\mathscr{M},\mathscr{H}\} = \int_X^\oplus \{\mathscr{M}_x,\mathscr{H}_x\} \ud x.
			\end{equation*}
			We often abbreviate this as
			\begin{equation}\label{directVN}
				\mathscr{M} = \int_X^\oplus  \mathscr{M}_x \ud x.
			\end{equation}
			\item The commutant of $\mathscr{M}$ is 
			\[
			\mathscr{M}^\prime =  \int_X^\oplus \mathscr{M}_x^\prime \ud x.
			\]
			\item The diagonal algebra $\mathscr{A}$ is contained in the center $\mathscr{Z}=\mathscr{M}\cap \mathscr{M}^\prime$ of $\mathscr{M}$.
		\end{list}
	\end{theorem}
	\begin{defn}[Direct integrals of von Neumann algebras]
		Given  $\{\mathscr{M}_x, \mathscr{H}_x\}_{x\in X}$ as in the above theorem, the von Neumann algebra $\mathscr{M}$ in \eqref{directVN} is called the \textit{direct integral of} $\{\mathscr{M}_x\}_{x\in X}$.
	\end{defn}
	\begin{prop}[Corollary 8.20 in \cite{Tak}]\label{820Tak}
		Let $\mathscr{M} = \int_X^\oplus  \mathscr{M}_x \ud x$ be a direct integral of von Neumann algebras. Then the center $\mathscr{Z}$ of $\mathscr{M}$ is also expressed as a direct integral
		\[
		\mathscr{Z} = \int_X^\oplus \mathscr{Z}_x \ud x,
		\]
		with $\mathscr{Z}_x =  \mathscr{M}_x \cap \mathscr{M}_x^\prime$. In particular, $\mathscr{Z}$ coincides with the diagonal algebra $\mathscr{A}$ if an only if $\mathscr{M}_x$ is a factor for almost every $x\in X$.
	\end{prop}
	
	For what follows, we assume that $G$ be a locally compact group.
	\begin{defn}
	Let $\{\mathscr{H}_x\}_{x\in X}$ be a measurable field of Hilbert spaces over $X$, and for each $x\in X$, let $\pi_x$ be a unitary representation of $G$ on $\mathscr{H}_x$. We say that $\{\pi_x\}_{x\in X}$ is a \textit{measurable field of representations} of $G$ if, for each $g\in G$, the operator field $x\mapsto \pi_x(g) \in\mathscr{B}(\mathscr{H}_x)$ is measurable. 
	\end{defn}
Given a measurable field of representations $\{\pi_x\}_{x\in X}$ of $G$, we have $\|\pi_x(g)\| = 1$ for all $x\in X$ and all $g\in G$; thus, $\underset{x\in X}{\mathrm{ess}\,\sup} \|\pi_x(x)\|=1$ for each $g\in G$, and hence we can form the direct integrals\label{directreps}
	\[
	\pi(g) = \int_X^\oplus \pi_x(g) \ud\mu(x),\ \ \ g\in G.
	\]
	It follows from the properties \ref{directintegralprops}, and the dominated convergence theorem that $g\mapsto \pi(g)$ defines a unitary representation of $G$ on $\int_X^\oplus \mathscr{H}_x\ud \mu(x)$. 
	\begin{defn}
 Let $\{\pi_x\}_{x\in X}$ be a measurable field of representations of $G$, and let
			\[
			\pi(g) = \int_X^\oplus \pi_x(g) \ud\mu(x),\ \ \text{for } g\in G.
			\]
			We call $\pi$ the \textit{direct integral of the representations $\pi_x$}, and write $\pi = \int_X^\oplus \pi_x \ud \mu(x)$.
			\end{defn}
 If $\mathscr{A}$ is a $\mathrm{C}^*$-algebra, then we define a measurable field of representations $\{\pi_x\}_{x\in X}$ of $\mathscr{A}$, and the direct integral of these representations, in precisely the same way.

	\section{Functions of Positive Type}\label{postype}
	If $G$ is a locally compact group, then $L^1(G)$ is a Banach $*$-algebra: its multiplication is convolution, and its involution is given by\footnote{$\Delta$ in this formula denotes the modular function---see \cite{TaylorKaniuth}, p.~3 for the definition of the modular function. Discrete groups are \textit{unimodular}, meaning that $\Delta \equiv 1$. In particular, the crystallographic groups are unimodular.} 
	\[
	f^*(x) = \Delta(x^{-1})\overline{f(x^{-1})}\ \ \text{for}\ f \in L^1(G) \ \ \text{and}\ x\in G.
	\]
A linear functional $\Lambda$ on $L^1(G)$ is called \textit{positive} if $\Lambda(f^* * f)\geq 0$ for all $f\in L^1(G)$. A function $\phi \in L^\infty(G)$ is said to be of \textit{positive type} if the associated bounded linear functional it defines on $L^1(G)$ is positive, i.e.:
\[
\int_G (f^* * f)(x)\phi(x)\ud x \geq 0 \ \ \text{for all}\ f\in L^1(G).
\]
\begin{rem}
Given a unitary representation $(\pi,\mathscr{H}_\pi)$ of $G$, and $\xi \in \mathscr{H}_\pi$, it is easily verified (see p.~30 in \cite{TaylorKaniuth}) that the function $\phi_{\xi,\xi}$ given by 
\begin{equation}\label{fnpostypeassoc}
\phi_{\xi,\xi}(x) = \langle \pi(x)\xi,\xi\rangle \ \ \text{for all}\ x\in G,
\end{equation}
is of positive type.
\end{rem}
		\begin{defn}
Let $(\pi,\mathscr{H}_\pi)$ be a unitary representation of $G$. A function $\phi$ of positive type is said to be \textit{associated with $\pi$} if $\phi = \phi_{\xi,\xi}$ for some $\xi \in \mathscr{H}_\pi$, where $\phi_{\xi,\xi}$ is given by \eqref{fnpostypeassoc}. Given a set $S$ of representations of $G$, we say that a function $\phi$ of positive type is \textit{associated with S} if $\phi$ is associated with $\sigma$ for some $\sigma \in S$. 
	\end{defn}
	\begin{rem} 
		Equivalent representations have the same functions of positive type associated with them.
		\end{rem}
\begin{notation}
Let $G$ be a locally compact group. We let $\mathscr{P}(G)$\label{fnpostype} denote the set of all continuous functions of positive type, and $\mathscr{P}_1(G) = \{\phi \in \mathscr{P}(G): \|\phi\|_\infty = 1\}$.
\end{notation}
	
	\section{Essential Theory of Fourier Analysis}\label{EssentialFourier}
	In this section we introduce some of the essential theory of Fourier analysis. The theorems in this section are given without proof; we refer the reader to \cite{FollandAbstract}, \cite{TaylorKaniuth}, and \cite{HewittRoss} for proofs.
	
		For any topological space $X$, we let $\mathcal{B}_X$ denote the Borel $\sigma$-algebra on $X$. A measure $\mu$ on $(X,\mathcal{B}_X)$ is called a \textit{Borel measure}; i.e., $\mu$ is a measure whose domain is the Borel sets. A Borel measure $\mu$ on $(X,\mathcal{B}_X)$ is said to be \textit{Radon measure} if it satisfies the following:
	\begin{list}{(\roman{ctr})}{\usecounter{ctr}}
		\setcounter{ctr}{0}
		\item $\mu$ is finite on all compact sets;
		\item $\mu$ is \textit{outer regular} on all Borel sets, i.e.:
		\[
		\mu(E) = \inf\{\mu(U):U\supset E,\ U\ \text{open}\}
		\]
		for all $E\in\mathcal{B}_X$;
		\item $\mu$ is \textit{inner regular} on all open sets, i.e.:
		\[
		\mu(U) = \sup\{\mu(K):K\subset U,\ K\ \text{compact}\}
		\]
		for all open $U\subset X$.
	\end{list}
	\begin{defn}
		Let $G$ be a locally compact group. A Borel measure $\mu$ on $G$ is called \textit{left-invariant} if $\mu(xE)=\mu(E)$ for any $x\in G$ and Borel subset $E$ of $G$. There exists a nonzero left-invariant Radon measure $\mu_G$ on $G$ satisfying:
		\begin{list}{(\alph{ctr1})}{\usecounter{ctr1}}
			\setcounter{ctr1}{0}
			\item $\mu_G(U)>0$ for any nonempty open subset $U$ of $G$.
			\item If $\nu$ is any nonzero left-invariant Radon measure on $G$, then there is a constant $c>0$ such that $\nu = c\mu_G$.
		\end{list}
		Any measure from the family $\{c\mu_G: c>0\}$ is a \textit{left Haar measure} on $G$. There is analogous notion of \textit{right Haar measure}, and when $G$ is abelian, discrete, or compact, the left Haar measure is also right-invariant. In each of these circumstances, we simply use the term \textit{Haar measure}.
	\end{defn}

	\begin{defn}\label{dualspace}
	 For a locally compact group $G$, the \textit{(unitary) dual space} of $G$, denoted $\widehat{G}$, is the set of equivalence classes of irreducible representations of $G$.  
	 	\end{defn}
	 	\begin{rem}\label{DualTop}
	 There is a topology on $\widehat{G}$ called the \textit{Fell topology} or \textit{dual space topology}---see \cite{TaylorKaniuth}, \S~1.6.
	 	\end{rem}
			 When $G$ is abelian, all its irreducible representations are $1$-dimensional, and hence are equivalent to representations on $\mathbb{C}$. Any unitary operator on $\mathbb{C}$ can be described by multiplication by an element of the circle group, $\mathbb{T}$. Identifying the unitary operators on $\mathbb{C}$ with the elements of $\mathbb{T}$ they correspond to, we see that $\widehat{G}$ may be identified with the set of \textit{continuous characters} of $G$---that is, the set of continuous group homomorphisms into the circle group. The set of continuous characters of $G$ forms an abelian group under pointwise multiplication. Thus when $G$ is abelian, $\widehat{G}$ is a group.
		\begin{defn}\label{chargroup}
The dual space $\widehat{G}$ of an abelian group $G$ is called the \textit{character group} or \textit{Pontryagin dual} of $G$.
	\end{defn}
	\begin{rem}
		 When $G$ is abelian, the Fell topology on $\widehat{G}$ coincides with the topology of uniform convergence on compact subsets of $G$---see \cite{FollandAbstract}, p.~88---and with this topology, $\widehat{G}$ is locally compact.
	\end{rem}
		\begin{rem}\label{convergenceinGhat}
			When $G$ is an abelian group, the characters of $G$ are contained in the subset $\mathscr{P}_1(G)$ of $L^\infty(G)$---the set of continuous functions of positive type with norm 1.\footnote{It follows from Theorem 3.25 in \cite{FollandAbstract} that the continuous characters of $G$ are precisely extreme points of $\mathscr{P}_1(G)$.} As the dual space of $L^1(G)$\footnote{That $L^1(G)^* \cong L^\infty(G)$ holds whenever the Haar measure on $G$ is $\sigma$-finite---see \cite{FollandReal}, Exercise 25 on p.~192, and \cite{FollandAbstract} pp.~50--51.}, we give $L^\infty(G)$ the weak-$*$ topology, and give $\mathscr{P}_1(G)$ relative weak-$*$ topology inherited from $L^\infty(G)$. It can be shown that the topology of uniform convergence on compact subsets of $G$ agrees with the restriction of the weak-$*$ topology to $\mathscr{P}_1(G)$ (this is Theorem 3.31 in \cite{FollandAbstract}). Since the weak-$*$ topology is Hausdorff, $\widehat{G}$ is Hausdorff.
		\end{rem}
		\begin{rem}
	Any 1-dimensional representation is cyclic. In particular, if $G$ is abelian, then the continuous characters of $G$ are cyclic with cyclic vector $1\in\mathbb{C}$. Moreover, if $\chi$ is a continuous character of $G$, then $\chi$ is a function of positive type associated with itself, since $\chi=\chi_{1,1}$, where $\chi_{1,1}(x)=\langle \pi(x)1,1\rangle$ for all $x\in G$.
		\end{rem}
		\begin{lemma}\label{dualiso}
			Let $A$ and $B$ be locally compact abelian groups. If $A$ is topologically isomorphic with $B$, then $\widehat{A}$ is topologically isomorphic with $\widehat{B}$.
		\end{lemma}
		\begin{proof}
	This follows from Exercises 7.6 and 7.9 in \cite{Deitmar}.		
			\end{proof}
	\begin{defn}\label{CX} Let $X$ be a locally compact Hausdorff space, and let
		\[
		C(X)= \{f:X\to \mathbb{C}: f\ \text{is continuous}\}.
		\]
		A function $f\in C(X)$ is said to \textit{vanish at infinity} if $\{x\in G: |F(x)|>\epsilon\}$ is compact for all $\epsilon>0$. We let
		\[
		C_0(X) = \{f\in C(X): f\ \text{vanishes at infinity}\}.
		\]
	\end{defn}
	\begin{defn}\label{FourierT} Let $G$ be a locally compact abelian group with Haar measure $\mu$. For $f\in L^1(G)$, the \textit{Fourier transform} of $f$ is the function $\widehat{f}: \widehat{G} \to \mathbb{C}$ defined by 
			\begin{equation}
			\boxed{	\widehat{f}(\chi) = \int_G f(x)\chi(x)\ud\mu(x)\ \ \ \text{for all}\ \chi \in \widehat{G}}\, .
			\end{equation}
		\end{defn}
		\begin{prop}[Proposition 1.79 in \cite{TaylorKaniuth}]
We may regard $C_0(\widehat{G})$ as a commutative $\mathrm{C}^*$-algebra algebra: operations are defined pointwise, and norm is $\|\cdot\|_\infty$. Then $f\mapsto \widehat{f}$ is an injective $*$-homomorphism of $L^1(G)$ onto a dense subset of $C_0(\widehat{G})$.
\end{prop}
\begin{tcolorbox}
\begin{defn}
The linear operator $\mathcal{F}: L^1(G) \to C_0(\widehat{G})$ given by $\mathcal{F}(f)=\widehat{f}$ is called the \textit{Fourier transform}.
\end{defn}
\end{tcolorbox}
\begin{rem}	 The fact that $\mathcal{F}(L^1(G)) \subset C_0(\widehat{G})$ is the abstract form of the Riemann-Lebesgue lemma.
	\end{rem}
			\begin{rem}
			Some authors instead define $\widehat{f}(\chi) = \int_G f(x)\overline{\chi(x)}\ud x$. This is merely a matter of convention.
			\end{rem}
	
	\begin{prop}[Theorem 4.22 in \cite{FollandAbstract}]	\hspace{2cm}
		\begin{list}{(\arabic{ctr})}{\usecounter{ctr}}
			\setcounter{ctr}{0}
			\item Let $G$ be a compact group, and let $M(G)$ denote the space of complex Radon measures on $G$. This a Banach algebra, called the \textit{measure algebra} of $G$, with convolution as its product (\cite{FollandAbstract}, pp.~54--55.)
			\item Given $\mu \in M(\widehat{G})$, we define $\phi_\mu$ on $G$ by 
			\[
			\phi_\mu(x) = \int_{\widehat{G}} \chi(x) \ud\mu(\chi) \ \ \ \text{for all } x\in G.
			\]
			Then $\phi_\mu$ is a bounded continuous function on $G$.
			\item Let $\mathcal{B}^1  = \{\phi_\mu: \mu \in M(\widehat{G})\} \cap L^1(G)$. It is known that $\mathcal{B}^1$ is dense in $L^1(G)$. 
			\item $\mathcal{F}(\mathcal{B}^1) \subset L^1(\widehat{G})$. Moreover, the Haar measure $\mu_{\widehat{G}}$ on $\widehat{G}$ can be normalized so that 
			\begin{equation}\label{inversion}
				f(x) = \int_{\widehat{G}} \widehat{f}(\chi) \overline{\chi(x)} \ud\mu_{\widehat{G}}(\chi)\ \ \ \text{for all}\ x\in G \ \text{and all}\ f\in\mathcal{B}^1.
			\end{equation}
		\end{list}
	\end{prop}

	\begin{defn}\label{dualmeas0} If $G$ is a locally compact abelian group with left Haar measure $\mu_{G}$, then the Haar measure $\mu_{\widehat{G}}$ on $\widehat{G}$ such that \eqref{inversion} holds is called the \textit{dual measure} of $\mu_{G}$.
	\end{defn}	
	
	\begin{rem}\label{scalingduality}
		If the dual of $\mu_{G}$ is $\mu_{\widehat{G}}$, then the dual of $c\mu_{G}$ is $c^{-1} \mu_{\widehat{G}}$.	
	\end{rem}	
	\begin{prop}[Proposition 1.88 in \cite{TaylorKaniuth}]
$G$ is compact if and only if $\widehat{G}$ is discrete; $G$ is discrete if and only if $\widehat{G}$ is compact.
	\end{prop}
			\begin{prop}[Proposition 4.25 in \cite{FollandAbstract}]
If $G$ is compact with Haar measure $\mu_{G}$ such that $\mu_{G}(G) =1$, then the dual measure $\mu_{\widehat{G}}$ on $\widehat{G}$ is the counting measure. If $G$ is discrete with Haar measure chosen to be counting measure $\nu_G$, then the dual measure $\nu_{\widehat{G}}$ on $\widehat{G}$ satisfies $\nu_{\widehat{G}}(\widehat{G}) =1$.
		\end{prop}

	\begin{theorem}[Plancherel Theorem---4.26 in \cite{FollandAbstract}]\label{Plancherel} Let $G$ be a locally compact abelian group. The Fourier transform on $L^1(G)\cap L^2(G)$ extends uniquely to a unitary map from $L^2(G)$ to $L^2(\widehat{G})$.
	\end{theorem}
	\begin{tcolorbox}
	\begin{defn}
		The unitary map from $L^2(G)$ to $L^2(\widehat{G})$ in the above theorem is called the \textit{Plancherel transform}.
		\end{defn}
	\end{tcolorbox}
		\begin{rem}\label{FT0}
	The Plancherel transform is sometimes referred to loosely as the Fourier transform. We use the calligraphic $\mathcal{F}$ for both.
		\end{rem}
	\begin{cor}[4.27 in \cite{FollandAbstract}]\label{dualbasis} If $G$ is compact and abelian with Haar measure $\mu_{G}$ such that and $\mu_{G}(G)=1$, then $\widehat{G}$ is an orthonormal basis for $L^2(G,\mu_{G})$.
	\end{cor}

	\begin{theorem}[Pontryagin Duality Theorem---4.32 in \cite{FollandAbstract}]\label{PontDual} Let $G$ be a locally compact abelian group. Each $x\in G$ defines a character $E_x$ on $\widehat{G}$ by $E_x(\chi) = \chi(x)$. The map $\mathbf{Ev}: x\mapsto E_x$ is a group homomorphism, and $\mathbf{Ev}: G \to \widehat{\widehat{G}}$ is an isomorphism of topological groups.
	\end{theorem}
	
	\begin{cor}\label{FourierInv2} If $f\in L^1(G)$ and $\widehat{f}\in L^1(\widehat{G})$ then 
		\[
		f(x) = \int_{\widehat{G}} \widehat{f}(\chi) \overline{\chi(x) }\ud \mu_{\widehat{G}}(\chi) \ \ \ \text{for a.e.}\ x\in G.
		\]
		If $f$ is continuous then this equation holds for every $x$.
	\end{cor}
	\section{The Central Decomposition}
	\begin{defn}
		 A unitary representation $\pi$ of a locally compact group $G$ is said to be of \textit{Type I} if the von Neumann algebra generated by $\pi(G)$ is of Type I (see A.35 in \cite{DixmierC} for a definition of Type I von Neumann algebra). See 5.4.2 and 13.4.1 in \cite{DixmierC}.
		\end{defn}
		\begin{defn}
		A locally compact group $G$ is said to be of \textit{Type I} if every unitary representation of $G$ is of Type I. 
	\end{defn}
	\begin{rem}\label{factorrepform}
A group is Type I if and only if, for each factor representation $\rho$ of $G$, the factor generated by $\rho(G)$ is of Type I (see \cite{DixmierC}, 13.9.4). This, in turn, is equivalent to requiring that every factor representation $\rho$ of $G$ is equivalent to a representation of the form $\tilde{\rho}\otimes I$, where $\tilde{\rho}$ is irreducible. (See \cite{FollandAbstract}, pp.~229 and 236.) Together with Corollary \ref{tensorrepfactor}, we have
\begin{tcolorbox}
	The factor representations of a Type I group $G$ are precisely those which are unitarily equivalent to $\tilde{\rho}\otimes I$ for some irreducible representation $\tilde{\rho}$ of $G$.
\end{tcolorbox}
	\end{rem}
 For any second countable group $G$, the dual space $\widehat{G}$ is equipped with a useful $\sigma$-algebra $\mathcal{M}$, called the \textit{Mackey Borel Structure}, introduced by George Mackey. 
\begin{defn}
	A measurable space is called \textit{standard} if it is measurably isomorphic to a Borel subset of a complete and separable metric space. A measure $\mu$ on $(\widehat{G},\mathcal{M})$ is called \textit{standard} if there is a subset $A\subset \widehat{G}$ such that $\mu(\widehat{G}-A)=0$ and such that the relative Mackey Borel structure on $A$ inherited from $\widehat{G}$ is standard.
	\end{defn}
\begin{rem}
	 When $G$ is Type I, the Mackey Borel Structure coincides with the Borel $\sigma$-algebra for the Fell topology, and the Mackey Borel $\sigma$-algebra on $\widehat{G}$ is standard (and hence any measure on it is also). See \cite{FollandAbstract}, pp.~228--230, 249 for more details.
\end{rem}
	\begin{theorem}[Central Decomposition (simplified to Type I setting)]\label{centraldecomptheorem}
		Let $G$ be a Type I, second countable, locally compact group. Let $(\pi, \mathscr{H})$ be a unitary representation of $G$, with $\mathscr{H}$ separable. Let $\widehat{G}$ be given the Mackey Borel $\sigma$-algebra. There is a (standard) measure $\mu$ on $\widehat{G}$, a measurable field of Hilbert spaces $\{\mathscr{H}_\zeta\}_{\zeta\in\widehat{G}}$, and a measurable field of (Type I) factor representations $\{\rho_\zeta\}_{\zeta\in \widehat{G}}$ on $\{\mathscr{H}_\zeta\}_{\zeta\in\widehat{G}}$ with the following properties:
		\begin{list}{(\arabic{ctr})}{\usecounter{ctr}}
			\setcounter{ctr}{0}
			\item we have $\rho_{\zeta}\sim \tilde{\zeta}\otimes I$ for some $\tilde{\zeta} \in \zeta$, for $\mu$-a.e. $\zeta \in \widehat{G}$.
			\item there is a unitary equivalence $U$ between $\pi$ and $\rho =\int_{\widehat{G}}^\oplus \rho_\zeta\ud\mu(\zeta)$. 
			\item $U$ transforms the center of $\pi(G)^-$ into the algebra $\mathscr{A}$ of diagonal operators on $\int_{\widehat{G}}^\oplus \mathscr{H}_\zeta\ud\mu(\zeta)$. Equivalently, by Lemmas \ref{closureofconjugate} and \ref{centerlemma}, $\mathscr{Z}(\rho(G)^-)= \mathscr{A}$.
			\item If $\mu^\prime$ and $\{\rho_\zeta^\prime\}_{\zeta\in \widehat{G}}$ satisfy the same properties, then $\mu$ is equivalent to $\mu^\prime$, and $\rho_\zeta$ is equivalent to $\rho_\zeta^\prime$ for $\mu$-a.e. $\zeta \in \widehat{G}$.
		\end{list}
	\end{theorem} 
	\begin{rem}
		The above theorem adapts Theorem 7.34 in \cite{FollandAbstract} to the Type I setting. The full theorem gives a unique disintegration of an arbitrary unitary representation of a second countable group on a separable Hilbert space as a direct integral of (uniquely determined a.e.) factor representations. 
		
		This theorem is adapted from Ernest \cite{Ernest}, which builds on a theorem of von Neumann (see 7.29 and 7.30(a) in \cite{FollandAbstract}). The main advantage of this decomposition is the uniqueness, which is gained by considering factor representations instead of irreducibles. For a proof of this theorem, see \cite{DixmierC}---Theorem 8.4.2., together with the Remarks 18.6.1--3, and 18.7.1--6.
	\end{rem}
	
	\section{Decompositions with Measurable Range Functions}\label{MeasRangeSection}
	In this section we lay out the relationship between range functions and direct integrals. We briefly discuss vector-valued functions and operators, and introduce the vector-valued Lebesgue space $L^2(X;\mathscr{H})$. We show that this space can be viewed as a direct integral over $X$ of the constant field $\{\mathscr{H}_x \equiv \mathscr{H}\}$. We then define what it means for a range function to be measurable, and in Theorem \ref{bigequiv} we describe the relationship between the measurability of a range function \label{range}
		\[
	J: X \to \{\text{closed subspaces of } \mathscr{H}\}
	\]
	and the existence of a measurable field of subspaces of $\mathscr{H}$. In Lemma \ref{embedDirectInt} we show how the direct integral $\int_X^\oplus J(x)\ud x$ of these subspaces may be regarded as a closed subspace of $L^2(X;\mathscr{H})$. In Corollary \ref{criteriaForDirectIntSubspace} we show that $\int_X^\oplus J(x)\ud x$ consists precisely of $\xi\in L^2(X;\mathscr{H})$ such that $\xi(x) \in J(x)$ a.e.. In the last part of this section we prove some necessary lemmas and then show how a closed subspace $W \subset L^2(X;\mathscr{H})$ allows us to construct range functions. We prove that if $W$ satisfies an \textit{orthogonality condition}, then there is a unique (up to a.e.-equality) measurable range function $J_W$ associated with $W$. We show that $W$ is the direct integral of the subspaces in the range of $J_W$.
	
	When dealing with vector-valued functions, we need a suitable notion of measurability. 
	\begin{defn} Let $\mathscr{H}$ and $\mathscr{K}$ be Hilbert spaces. \begin{list}{(\arabic{ctr})}{\usecounter{ctr}}
			\setcounter{ctr}{0}
			\item Given a measure space $(X,\mathcal{M},\mu)$, a function $\xi: X \rightarrow \mathscr{H}$ is called {\it weakly measurable} if for each $b \in \mathscr{H}$, the map $x\mapsto \langle \xi(x), b\rangle$ is a measurable complex-valued function. 
			\item A function $K: X \rightarrow \mathscr{B}(\mathscr{H},\mathscr{K})$ is said to be \textit{weakly measurable} if $K(x)a$ is weakly measurable for each $a\in \mathscr{H}$; that is, $K$ is weakly measurable if $x\mapsto \langle K(x)a,b\rangle$ is measurable for each $a \in \mathscr{H}$ and each $b \in \mathscr{K}$.
		\end{list}
	\end{defn}
	\begin{rem}\label{weakmeasITObases}
		If $\mathscr{H}$ and $\mathscr{K}$ in the above definition are separable Hilbert spaces with orthonormal bases $\{e_n\}_{n\in\mathbb{N}}$ and $\{f_n\}_{n\in\mathbb{N}}$, respectively, then $\xi: X\to \mathscr{H}$ is weakly measurable if and only if the functions $x\mapsto \langle \xi(x),e_n\rangle$ are measurable for all $n$. This follows from the linearity and continuity of the inner product in the second variable, as well as the fact that pointwise limits of measurable functions are measurable. Likewise, $K: X \rightarrow \mathscr{B}(\mathscr{H},\mathscr{K})$ is weakly measurable if and only if the functions $x\mapsto \langle K(x)e_m,f_n\rangle$ are measurable for all $m,n\in\mathbb{N}$. This follows from the continuity of $K(x)$ for each $x$, the continuity of the inner product in each variable, and the fact that measurability is preserved by pointwise limits.
	\end{rem}
	\begin{defn} Given a measure space $(X,\mathcal{M},\mu)$, we define $L^2(X;\mathscr{H})$ to be the set of equivalence classes of weakly measurable $\mathscr{H}$-valued functions $\xi$---where equivalence $\sim$ is given by $\mu$-a.e. equality---satisfying $\int_X \|\xi(x)\|_{\mathscr{H}}^2 \ud \mu(x) <\infty$. That is,
		\[
		L^2(X;\mathscr{H}) = \Big\{\xi:X\to \mathscr{H}\ \text{weakly measurable and}\  \int_X \|\xi(x)\|_{\mathscr{H}}^2 \ud \mu(x) <\infty \Big\}\Big/\sim
		\]
		where $\sim$ denotes the equivalence relation $\xi\sim \eta \iff \xi=\eta\ \mu$-a.e.. We write $[\xi]_{L^2}$, or just $[\xi]$, for the equivalence class of $\xi$. \textbf{N.B.:} We often abuse notation and write $\xi\in L^2(X;\mathscr{H})$ to mean that $\xi$ is a representative of the equivalence class of square-integrable, weakly measurable functions equal to $f$ almost everywhere.
		$L^2(X;\mathscr{H})$\label{vectorL2} is a Hilbert space, with norm given by $\|[\xi]\|^2 =  \int_X \|\xi(x)\|_{\mathscr{H}}^2 \ud \mu(x)$, and inner product given by 
		\[
		\langle [\xi],[\eta]\rangle = \int_X \langle \xi(x),\eta(x)\rangle_{\mathscr{H}} \ud \mu(x).
		\]
	\end{defn}
	\begin{rem}\label{MeasVectorField}
		Recall from Example \ref{constantField} that if $\{\mathscr{H}_x\}_{x\in X}$ is the constant field of separable Hilbert spaces, i.e., $\mathscr{H}_x = \mathscr{H}$ for all $x\in X$, then the orthonormal basis $\{e_n\}_{n\in\mathbb{N}}$ of $\mathscr{H}$ provides a the fundamental sequence of measurable fields, by taking $e_n(x) = e_n$ for all $x$. This makes $\{\mathscr{H}_x\}_{x\in X}$ a measurable field of Hilbert spaces with collection of measurable vector fields is given by 
		\begin{equation}\label{measvectorfields}
			\mathcal{E}= \{\xi\in \Pi_{x\in X} \mathscr{H}_x: x\mapsto \langle \xi(x),e_n(x) \rangle_x \ \text{is}\ \mu\text{-measurable for every}\ n \}.
		\end{equation}
		Since $\{e_n\}_{n\in\mathbb{N}}$ is an orthonormal basis for $\mathscr{H}$, and the inner product is continuous in the second variable, it follows that: 
		\begin{tcolorbox}
			\begin{center}
				A vector field $\xi$ associated with the constant field $\mathscr{H}_x\equiv \mathscr{H}$ is measurable if and only if it is a weakly measurable $\mathscr{H}$-valued function.
			\end{center}
		\end{tcolorbox}
		Thus
		\[
		\int_X^\oplus \mathscr{H}_x \ud\mu(x) = L^2(X;\mathscr{H}).
		\] 
	\end{rem}
	
	\begin{defn} Let $(X,\mathcal{M},\mu)$ be a $\sigma$-finite measure space, and let $\mathscr{H}$ be a separable Hilbert space.
		\begin{list}{(\arabic{ctr})}{\usecounter{ctr}}
			\setcounter{ctr}{0}
			\item 	A \textit{range function} is a mapping $J:X \to \{\text{closed subspaces of}\ \mathscr{H}\}$.
			\item A range function $J$ is said to be measurable if the operator-valued function $P_J:X \to \mathscr{B}(\mathscr{H})$---which takes $x$ the orthogonal projection $P_J(x)$ of $\mathscr{H}$ onto $J(x)$---is weakly measurable.
		\end{list}
	\end{defn}
	\begin{rem}
	A range function gives a field $\{J(x)\}_{x\in X}$ of closed subspaces of $\mathscr{H}$.
	\end{rem}
	Let $(X,\mathcal{M},\mu)$ be a measure space.\label{measspace} Recall that the \textit{completion of $\mathcal{M}$} with respect to $\mu$ is defined as follows: 
	\[
	\overline{\mathcal{M}} = \{E\cup N_0: E\in\mathcal{M},\ \text{and}\ N_0 \subset N\ \text{for some}\ N\in \mathcal{M}\ \text{with}\ \mu(N) = 0\}.
	\]
	Moreover, one can extend $\mu$ to a measure $\overline{\mu}$ on $\overline{\mathcal{M}}$, called the \textit{completion of} $\mu$\label{mucomplete}, by setting
	\[
	\overline{\mu}(E\cup N_0) = \mu(E)
	\]
	for all $E\cup N_0\in \overline{\mathcal{M}}$, where $E\in \mathcal{M}$ and $N_0 \subset N$ for some $N\in \mathcal{M}$ with $\mu(N)=0$. A measure $\mu$ on $(X,\mathcal{M})$ is said to be \textit{complete} if its domain includes all subsets of sets of measure zero. The following standard result will be useful:
	\begin{prop}[Proposition 2.11 in \cite{FollandReal}]\label{MeasAEresult} Let $(X,\mathcal{M},\mu)$ be a measure space. The following are equivalent: 
		\begin{list}{(\arabic{ctr})}{\usecounter{ctr}}
			\setcounter{ctr}{0}
			\item $\mu$ is complete.
			\item If $f = g$ $\mu$-a.e., and $f$ is measurable, then so is $g$.
			\item If $\{f_n\}_{n\in\mathbb{N}}$ is a sequence of measurable functions and $f_n \to f$ $\mu$-a.e., then $f$ is measurable.
		\end{list}
	\end{prop}
	\begin{notation}
	Given a vector space $V$, and a family $\{v_i\}_{i\in I}$ of vectors in $V$, we write
	\[
	\overline{\text{span}} \{v_i: i \in I\}
	\] 
	for the \textit{closed linear span} of $\{v_i\}_{i\in I}$---i.e.: the smallest closed subspace of $V$ containing the $v_i$. We also use the notation\label{bigvee}
	\[
	\bigvee_{i\in I} v_i = \bigvee\{v_i: i\in I\} =	\overline{\text{span}} \{v_i\}.
	\]
	We use these notations interchangeably.
	\end{notation}
	\begin{theorem}\label{bigequiv}
		Let $(X,\mathcal{M},\mu)$ be a complete, $\sigma$-finite measure space and let $\mathscr{H}$ be a separable Hilbert space. Let $\{e_n\}$ be an orthonormal basis for $\mathscr{H}$, so that we may regard $L^2(X;\mathscr{H})$ as the direct integral of the constant field of Hilbert spaces $\{\mathscr{H}_x \equiv \mathscr{H}\}$ with fundamental sequence $\{e_n(x)\equiv e_n\}$, and with measurable vector fields given by \eqref{measvectorfields}.
		For each $x\in X$, let $J(x)$ be a closed subspace of $\mathscr{H}$. Let $P_J:X \to \mathscr{B}(\mathscr{H})$ denote the operator-valued function taking $x$ to the orthogonal projection $P_J(x)$ of $\mathscr{H}$ onto $J(x)$. The following are equivalent:
		\begin{list}{(\arabic{ctr})}{\usecounter{ctr}}
			\setcounter{ctr}{0}
			\item $J: x\mapsto J(x)$ is a measurable range function.
			\item $x \mapsto \langle P_J(x)a,b\rangle$ is measurable for all $a,b \in \mathscr{H}$ (i.e., $P_J$ is a weakly measurable operator-valued function).
			\item $x \mapsto \langle P_J(x)e_n,e_m\rangle$ is measurable for all $m,n\in\mathbb{N}$.
			\item  $\{J(x)\}_{x\in X}$ is a measurable field of Hilbert spaces with fundamental sequence $\{f_n(x) := P_J(x)e_n\}$ and measurable vector fields given by 
			\begin{equation}\label{EJ}
				\mathcal{E}_J = \{\xi\in \Pi_{x\in X} J(x): x\mapsto \langle \xi(x),f_n(x) \rangle_x \ \text{is}\ \mu\text{-measurable for every}\ n \}.
			\end{equation}
			\item $x \mapsto \langle P_J(x) \xi(x),b\rangle$ is measurable for all $b\in \mathscr{H}$ and all weakly measurable $\mathscr{H}$-valued functions $\xi$ on $X$.
			\item $x \mapsto \langle P_J(x) \xi(x),b\rangle$ is measurable for all $b\in \mathscr{H}$ and all $\xi\in L^2(X;\mathscr{H})$.
			\item $P_J \in \Pi_{x\in X} \mathscr{B}(\mathscr{H}_x)$ is a measurable operator field.
		\end{list}
	\end{theorem}
	\begin{proof}
		(1) $\iff$ (2): This is by definition. 
		
		(2) $\iff$ (3): by Remark \ref{weakmeasITObases}.
		
		(3) $\implies$ (4): Assume that (3) holds. Since the measurable vector fields $\mathcal{E}$ are given by \eqref{measvectorfields}, (3) implies that $f_n$ is a measurable vector field, for each $n$, where $f_n(x):=P_J(x) e_n$. Moreover, $x\mapsto \langle f_m(x),f_n(x)\rangle$ is measurable for all $m, n\in\mathbb{N}$; indeed,
		\begin{align*}
			\langle f_m(x),f_n(x) \rangle &= \langle P_J(x)e_m, P_J(x)e_n\rangle\\
			&= \langle P_J^*(x) P_J(x)e_m, e_n\rangle\\
			&= \langle P_J^2(x)e_m, e_n\rangle\\
			&= \langle P_J(x)e_m, e_n\rangle,
		\end{align*}
		and $x \mapsto \langle P_J(x)e_m, e_n\rangle$ is measurable, by assumption. Furthermore, it is easy to see that $J(x) = \bigvee\{P_J(x) e_n\}$ for each $x$; indeed, given $a\in\mathscr{H}$, we can write $a$ as a linear combination of the basis elements: $a = \sum_n a_n e_n$, whence 
		\[
		P_J(x) a = \sum_n a_n P_J(x) e_n  \in  \overline{\text{span}}\{f_n(x)\}
		\]
		for each $x\in X$. Therefore $J(x) = P_J(x)\mathscr{H}  = \overline{\text{span}}\{f_n(x)\}$ for each $x\in X$, and hence, by Remark \ref{FollandMeasField}, $\{J(x)\}_{x\in X}$ is a measurable field of Hilbert spaces with fundamental sequence $\{f_n\}_{n\in \mathbb{N}}$, and with measurable vector fields given by
		\begin{equation*}
			\mathcal{E}_J = \{\xi\in \Pi_{x\in X} J(x): x\mapsto \langle \xi(x),f_n(x) \rangle_x \ \text{is}\ \mu\text{-measurable for every}\ n \}.
		\end{equation*}
		
		(4) $\implies$ (3): Assume (4), so that $f_n$ is the fundamental sequence of measurable vector fields. As above, we have  $\langle P_J(x)e_m, e_n\rangle =\langle f_m(x),f_n(x) \rangle$, and $x\mapsto \langle f_m(x),f_n(x) \rangle$ is measurable, by Remark \ref{directintpolarmeas}.
		
		(2) $\iff$ (5): This follows from Proposition \ref{KfWeaklyMeasurable} and the fact that constant functions are weakly measurable (Lemma \ref{LinearCombinationsWeaklyMeasurable}).
		
		(5) $\iff$ (7): By definition, $P_J$ being a measurable operator field means that $x\mapsto P_J(x)\xi(x)$ is a measurable vector field whenever $\xi$ is a measurable vector field. And by Remark \ref{MeasVectorField}, this is equivalent to $x\mapsto P_J(x)\xi(x)$ being weakly measurable whenever $\xi$ is weakly measurable, which is precisely what (5) says. 
		
		(5) $\implies$ (6): This is immediate since $L^2(X;\mathscr{H}) \subset \{\text{Weakly measurable}\ \mathscr{H}-\text{valued functions}\}$.
		
		(6) $\implies$ (5): Assume (6), so that $x \mapsto \langle P_J(x) \xi(x),b\rangle$ is measurable for all $b\in \mathscr{H}$ and all $\xi\in L^2(X;\mathscr{H})$. Let $\xi$ be a weakly measurable $\mathscr{H}$-valued function. It follows from Lemma \ref{StrongMeasConv} and Theorem \ref{WeakStrong} that $\xi$ is the a.e.-limit of a sequence $\{\phi_j\}$ of simple functions. Since $\phi_j \in L^2(X;\mathscr{H})$ for all $j$, we have $x\mapsto \langle P_J(x)\phi_j(x) , b\rangle$ is measurable for each $b\in \mathscr{H}$, by assumption. Moreover,
		\[
		\langle P_J(x)\xi(x) , b\rangle = \lim_{j\to\infty} \langle P_J(x)\phi_j(x) , b\rangle \ \ \text{a.e.},
		\]
		so it follows by Proposition \ref{MeasAEresult} that $x\mapsto \langle P_J(x)\xi(x) , b\rangle$ is measurable for all $b\in\mathscr{H}$.
	\end{proof}
	
	\begin{rem}\label{measopconstantfield}
		With $(X,\mathcal{M},\mu)$, $\mathscr{H}$, and $\{e_n\}$ as in the above theorem, and $T\in \Pi_{x\in X} \mathscr{B}(\mathscr{H})$ an operator field, it is easy to see that if we replace $P_J$ with $T$, then the equivalence of (2), (3), (5), (6), and (7) in the above theorem still holds. In particular, the equivalence of (2) and (7) does not require $X$ to be complete or $\sigma$-finite, and gives:
		\begin{tcolorbox}
			\begin{center}
				An operator field $T$ associated with the constant field $\mathscr{H}_x\equiv \mathscr{H}$ is measurable precisely when it is a weakly measurable $\mathscr{B}(\mathscr{H})$-valued function.
			\end{center}
		\end{tcolorbox}
	\end{rem}
		The preceding theorem shows that a measurable range function 
		\[
		J: X \to \{\text{closed subspaces of } \mathscr{H}\}
		\]
	 yields a measurable field of subspaces $\{J(x)\}_{x\in X}$ with measurable vector fields given by \eqref{EJ}. Thus we may form the direct integral $\int_X^\oplus J(x) \ud x$. We now show how this direct integral can be regarded as a subspace of $L^2(X;\mathscr{H})$.
\begin{lemma} Let $(X,\mathcal{M},\mu)$ be a $\sigma$-finite measure space, and let $\mathscr{H}$ be a separable Hilbert space with orthonormal basis $\{e_n\}_{n\in\mathbb{N}}$. As usual, we regard $L^2(X;\mathscr{H})$ as the direct integral of the constant field of Hilbert spaces $\{\mathscr{H}_x \equiv \mathscr{H}\}$, for which a collection $\mathcal{E}$ of measurable vector fields is given by \eqref{measvectorfields}. Let $J:X \to \{\text{closed subspaces of}\ \mathscr{H}\}$ be a measurable range function, so that $\{J(x):x\in X\}$ is a measurable field of Hilbert spaces with a collection $\mathcal{E}_J$ of measurable vector fields given by \eqref{EJ}. We have $\mathcal{E}_J \subset \mathcal{E}$.
	\end{lemma}
			\begin{proof}
				Suppose that $\xi \in \mathcal{E}_J$. Since $\xi(x) \in J(x)$ for each $x$, we have
				\begin{align*}
					\langle \xi(x),e_n(x) \rangle &= \langle P_J(x)\xi(x),e_n(x)\rangle \\
					&= \langle \xi(x), P_J(x) e_n(x) \rangle\\
					&= \langle \xi(x), f_n(x) \rangle,
				\end{align*}
				and $x\mapsto \langle \xi(x), f_n(x) \rangle$ is measurable since $\xi \in \mathcal{E}_J$. Therefore, $x\mapsto\langle \xi(x),e_n(x) \rangle$ is measurable, whence $\xi\in \mathcal{E}$. 
			\end{proof}
			\begin{notation}
			 Recall that, given a measurable field of Hilbert spaces $\{J(x)\}$ with measurable vector fields $\mathcal{E}_J$, we define an equivalence relation on the square-integrable fields in $\mathcal{E}_J$ as follows: two measurable, square-integrable vector fields are said equivalent precisely when they are equal almost everywhere. We denote this equivalence relation by $\sim_J$, and for $\xi \in \mathcal{E}_J$, we write $[\xi]_J$ for its equivalence class. Recall that the set of all these equivalence classes is the direct integral $\int_X^\oplus J(x)\ud x$. Also recall that two elements of $L^2(X;\mathscr{H})$ were defined to be equivalent if there were equal a.e.; we denote this by $\sim_{L^2}$ and write $[\xi]_{L^2}$ for the equivalence class of an element $\xi\in L^2(X;\mathscr{H})$.
			 \end{notation}
	\begin{lemma}\label{embedDirectInt}
Let $J$, $\mathscr{H}$, and $\mathcal{E}_J$ be as in the preceding lemma. The map $\iota([\xi]_J) = [\xi]_{L^2}$---which sends the element $[\xi]_J\in \int_X^\oplus J(x)\ud x$ to the equivalence class of all square-integrable, weakly measurable vector fields equal a.e. to $\xi$---embeds $\int_X^\oplus J(x)\ud x$ in $L^2(X;\mathscr{H})$. Thus we may regard $\int_X^\oplus J(x)\ud x$ as a subspace of $L^2(X;\mathscr{H})$ by identifying it with its image under $\iota$.
	\end{lemma}
	 \begin{proof} $\iota$ is well-defined, for if $\xi = \eta$ a.e., then $[\xi]_{L^2} = [\eta]_{L^2}$. And since the inner products agree, it follows that $\iota$ is an isometry:
			\[
			\|[\xi]_{L^2}\|^2 = \int_X \|\xi(x)\|^2\ud x=	\|[\xi]_J\|^2.
			\]
			Hence $\iota$ is an isometric isomorphism onto its image in $L^2(X;\mathscr{H})$.
\end{proof}
	\begin{prop}[Proposition 2.1 in \cite{BownikRoss}]\label{Mjsubspace}
		Let $(X,\mathcal{M},\mu)$ be a $\sigma$-finite measure space, $J$ be a range function (not necessarily measurable), and $\mathscr{H}$ be a separable Hilbert space. Then 
		\[
		M_J :=\{\xi\in L^2(X;\mathscr{H}): \xi(x)\in J(x)\ \text{for a.e.}\ x\in X\}
		\]
		is a closed linear subspace of $L^2(X;\mathscr{H})$.
	\end{prop}
	\begin{proof}
		Compare Remark (ii) on p.~284 of \cite{Bownik}.
	\end{proof}
	
	\begin{cor}\label{criteriaForDirectIntSubspace}
		Let $(X,\mathcal{M},\mu)$ be a $\sigma$-finite measure space, and let $\mathscr{H}$ be a separable Hilbert space. Let $J:X \to \{\text{closed subspaces of}\ \mathscr{H}\}$ be a measurable range function, so that $\{J(x):x\in X\}$ is a measurable field of Hilbert spaces with a collection $\mathcal{E}_J$ of measurable vector fields given by \eqref{EJ}. Regarding $\int_X^\oplus J(x)\ud x$ as a closed subspace of $L^2(X;\mathscr{H})$ as in the above remark, we have
		\begin{equation}
			M_{J}:=\{\xi\in L^2(X;\mathscr{H}): \xi(x) \in J(x)\ \text{for a.e.}\ x\in X\} = \int_X^\oplus J(x) \ud x.
		\end{equation}
	\end{cor}
	\begin{proof}
		Suppose that $\xi \in \int_X^\oplus J(x)\ud x$, where we regard this as a subspace of $L^2(X;\mathscr{H})$. Then $\xi \in L^2(X;\mathscr{H})$ by definition, and $[\xi]_{L^2}= \iota([\xi_0]_J)=[\xi_0]_{L^2}$ for some $\xi_0 \in  \int_X^\oplus J(x)\ud x$. It follows that $\xi(x) \in J(x)$ almost everywhere, and hence $\xi\in M_J$.
		
		Suppose conversely that $\xi \in M_J$. There exists a null set $N$ such that $\xi(x) \in J(x)$ on $X-N$. Define $\xi_0$ by
		\[
		\xi_0(x) = \left\{ \begin{array}{ll}
			\xi(x) &\text{if } x\in X-N\\
			0 &\text{if } x\in N.
		\end{array} \right.
		\]
		In other words, $\xi_0 = \mathbbm{1}_{X-N}\xi$, where $\mathbbm{1}_{X-N}$ denotes the indicator function of $X-N$---that is, the function taking value $1$ on $X-N$ and $0$ on $N$. Then $\xi_0 \in \Pi_{x\in X} J(x)$, and $\int_X \|\xi_0(x)\|^2\ud x  = \int_{X-N} \|\xi(x)\|^2\ud x < \infty$. Therefore it only remains to show that $\xi_0 \in \mathcal{E}_J$. We have
		\begin{equation}\label{projtrick}
			\begin{aligned}
				\langle \xi_0(x),f_n(x)\rangle &= \langle \xi_0(x),P_J(x) e_n\rangle\\
				&= \langle P_J(x)\xi_0(x),e_n\rangle \\
				&= \langle \xi_0(x),e_n\rangle  \\
				&=\mathbbm{1}_{X-N}(x)\langle \xi(x),e_n\rangle 
			\end{aligned}
		\end{equation}
		Now $x\mapsto \langle \xi(x),e_n\rangle$ is measurable since $\xi \in \mathcal{E}$ (since $\xi$ is weakly measurable). It follows that $x\mapsto \langle \xi_0(x),f_n(x)\rangle $ is measurable, being the product of two measurable functions. Thus $\xi_0 \in \mathcal{E}_J$. We have $[\xi]_{L^2} = \iota([\xi_0]_J)$ where $\xi_0 \in \int_X^\oplus J(x) \ud x$; thus $\xi \in \int_X^\oplus J(x) \ud x$ when we regard this direct integral as a subspace of $L^2(X;\mathscr{H})$.
	\end{proof}
	Suppose $(X,\mathcal{M},\mu)$, $\mathscr{H}$, and $J$ are as in the preceding corollary. Let $P_J(x)$ denote the orthogonal projection of $\mathscr{H}$ onto $J(x)$, and let $P_J$ be the operator-valued map $x\mapsto P_J(x)$. Then in particular, condition (7) of Theorem \ref{bigequiv} is satisfied, i.e.: $P_J \in \Pi_{x\in X} \mathscr{B}(\mathscr{H}_x, J(x))$ is a measurable operator field. Since $\underset{x\in X}{\mathrm{ess}\,\sup}  \|P_J(x)\| \leq 1<\infty$, we may form the direct integral $P_J = \int_X^\oplus P_J(x)\ud x: \int_{x\in X}^\oplus \mathscr{H}_x \ud x \to \int_{x\in X}^\oplus J(x)\ud x$, as in Definition \ref{directintops}.
	\begin{lemma}\label{projRem}
	Regarding $\int_{x\in X}^\oplus J(x)\ud x$ as a closed subspace of $L^2(X;\mathscr{H})$, we have that $\int_X^\oplus P_J(x)\ud x$ is the orthogonal projection from $L^2(X;\mathscr{H})$ onto $\int_{x\in X}^\oplus J(x)\ud x$. 
	\end{lemma}
	\begin{proof}
		This follows immediately from Lemma \ref{directintegralprops}.
		\end{proof}
	The following proposition will be useful later.
	\begin{prop}\label{measrangefntensorprod}
		Let $\mathscr{H}$ and $\mathscr{K}$ be separable Hilbert spaces, with $\mathscr{H}\neq \{0\}$. Then 
		\[
		K: X \to \{\text{closed subspaces of }\mathscr{K}\}
		\]
		is a measurable range function if and only if $x\mapsto \mathscr{H} \otimes K(x)$ is a measurable range function from $X$ to the closed subspaces of $\mathscr{H}\otimes \mathscr{K}$.
	\end{prop} 
	\begin{proof}
		Let $\mathscr{H}$ and $\mathscr{K}$ have orthonormal bases $\{e_m\}_{m\in \mathbb{N}}$ and $\{f_n\}_{n\in \mathbb{N}}$, respectively. Suppose first that $K$ is a measurable range function. Then by Theorem \ref{bigequiv}, $\{K(x)\}_{x\in X}$ is a measurable field of Hilbert spaces with fundamental sequence $\{P_K(x)f_n\}$ and measurable vector fields given by 
		\begin{equation*}
			\{\xi\in \Pi_{x\in X} K(x): x\mapsto \langle \xi(x),P_K(x)f_n \rangle_x \ \text{is}\ \mu\text{-measurable for every}\ n \}.
		\end{equation*}
		We claim that $\{\mathscr{H} \otimes K(x)\}$ is a measurable field of Hilbert spaces with fundamental sequence $\{e_m \otimes P_K(x) f_n\}_{(m,n)\in\mathbb{N}^2}$, where $P_K(x)$ is the orthogonal projection of $\mathscr{K}$ onto $K(x)$. This is equivalent to showing that the following two conditions are satisfied:
		\begin{list}{(\arabic{ctr})}{\usecounter{ctr}} 
			\setcounter{ctr}{0}
			\item $x\mapsto \langle e_{m_1} \otimes P_K(x) f_{n_1},e_{m_2} \otimes P_K(x) f_{n_2}\rangle_x$ is measurable for all $(m_i,n_i) \in\mathbb{N}^2, i=1,2$.
			\item For each $x\in X$, 
			\[
			\mathscr{H}\otimes K(x) =\bigvee\{e_m \otimes P_K(x) f_n: (m,n)\in \mathbb{N}^2\}.
			\]
		\end{list}
		For the first of these conditions, we have
		\begin{equation}\label{meastens}
			\begin{aligned}
				x&\mapsto \langle e_{m_1} \otimes P_K(x) f_{n_1},e_{m_2} \otimes P_K(x) f_{n_2}\rangle\\
				&= \langle e_{m_1}, e_{m_2} \rangle \langle  P_K(x) f_{n_1} ,  P_K(x) f_{n_2}  \rangle \\
				&= \delta_{m_1}(m_2) \langle  P_K(x) f_{n_1},  P_K(x) f_{n_2} \rangle,
			\end{aligned}
		\end{equation}
		which is measurable on account of $\{P_K(x)f_n\}_{n\in \mathbb{N}}$ being a fundamental sequence for the measurable field $\{K(x)\}$. For the second condition, let $a \in \mathscr{H}$ and $b\in \mathscr{K}$. Then $a = \sum_{m\in \mathbb{N}} a_m e_m$ and $b=\sum_{n\in \mathbb{N}} b_n f_n$, and 
		\[
		P_K(x) b = \sum_{n\in \mathbb{N}} b_n P_K(x) f_n  \in  \overline{\text{span}}\{P_K(x) f_{n}\}
		\]
		for each $x\in X$. It follows that
		\[
		a\otimes P_K(x) b  = \sum_{(m,n)\in\mathbb{N}^2} a_mb_n (e_m \otimes P_K(x) f_n) \in \bigvee_{(m,n)\in\mathbb{N}^2} e_m \otimes P_K(x) f_n,
		\]
		Since each element of $K(x)$ is of the the form $P_K(x) b$ for some $b\in \mathscr{K}$, it follows from the above equation that all simple tensors in $\mathscr{H}\otimes K(x)$ are contained in 
		\[
		\bigvee \{e_m \otimes P_K(x) f_n: (m,n)\in\mathbb{N}^2\};
		\]
		therefore, by linearity, 
		\begin{equation}\label{meastensor}
			\mathscr{H}\otimes K(x) = \bigvee \{e_m \otimes P_K(x) f_n: (m,n)\in\mathbb{N}^2\}.
		\end{equation}
		Thus $x\mapsto \mathscr{H} \otimes K(x)$ is a measurable field of Hilbert spaces with fundamental sequence $\{e_m \otimes P_K(x) f_n\}_{(m,n)\in\mathbb{N}^2}$, which by Theorem \ref{bigequiv}, proves that $x\mapsto \mathscr{H} \otimes K(x)$ is a measurable range function.
		
		For the converse, suppose that $x\mapsto \mathscr{H} \otimes K(x)$ is a measurable range function. Then $\{\mathscr{H} \otimes K(x)\}$ is a measurable field of Hilbert spaces with fundamental sequence $\{P_{\mathscr{H}\otimes K(x)}(e_m \otimes f_n)\}$. Equivalently, $\{e_m \otimes P_K(x) f_n\}_{(m,n)\in\mathbb{N}^2}$ is a fundamental sequence (since $P_{\mathscr{H}\otimes K(x)} = P_{\mathscr{H}}\otimes P_{K(x)}$); in particular, \eqref{meastensor} holds for each $x$. Note that since $\mathscr{H}\neq \{0\}$, the orthonormal basis $\{e_m\}_{m\in\mathbb{N}}$ for $\mathscr{H}$ is nonempty: it contains $e_1$. There exists an isometric isomorphism $\epsilon(x): \{e_1\} \otimes K(x) \to K(x)$ given by  $e_1\otimes f \mapsto f$. To see this, first let $\mathscr{H}_0(x) = \{e_1\otimes f: f\in K(x)\}$. Then $\mathscr{H}_0(x)$ is dense in $\mathscr{H}\otimes K(x)$, and we define $\epsilon(x)$ on $\mathscr{H}_0(x)$ by $e_1\otimes f \mapsto f$. Clearly, $\epsilon(x)$ is an isometric isomorphism from $\mathscr{H}_0(x)$ onto $K(x)$. Given $g\in \{e_1\} \otimes K(x)$, we have $g=\lim_{n\to \infty} e_1\otimes f_n$ for some sequence $\{f_n\}_{n=1}^\infty \subset K(x)$. Then $\{e_1\otimes f_n\}_{n=1}^\infty$ is Cauchy, whence $\{f_n\}_{n=1}^\infty$ is Cauchy in $K(x)$, whence $f_n$ converges to some $f\in K(x)$. Thus $g = \lim_{n\to\infty} e_1\otimes f_n = e_1\otimes f$, and $\mathscr{H}_0(x) = \mathscr{H}\otimes K(x)$.
		
		We  now show that $K(x) = \bigvee_{n\in\mathbb{N}} P_K(x) f_n$. Let $f\in K(x)$. By \eqref{meastensor}, there exist scalars $a_{m,n} \in \mathbb{C}$ such that
		\[
		e_1\otimes f = \sum_{(m,n)\in\mathbb{N}^2} a_{m,n} e_m \otimes P_K(x)f_n
		\]
Let $\langle e_1\rangle$ denote the 1-dimensional space $\{ke_1: k\in \mathbb{C}\}$. Then 
		\begin{align*}
e_1\otimes f &=	(P_{\langle e_1\rangle} \otimes P_K(x))(e_1\otimes f) \\
&= \sum_{(m,n)\in\mathbb{N}^2} a_{m,n} (P_{\langle e_1\rangle} \otimes P_K(x))(e_m\otimes P_K(x)f_n)\\
	&= \sum_{(m,n)\in\mathbb{N}^2} a_{m,n} P_{\langle e_1\rangle}(e_m) \otimes [P_K(x)]^2f_n \\
	&= \sum_{n\in\mathbb{N}} e_1 \otimes a_{1,n} P_K(x)f_n.
		\end{align*} 
It follows that
\begin{align*}
f = \epsilon(x)(e_1\otimes f) &= \epsilon(x)\Big( \sum_{n\in\mathbb{N}} e_1 \otimes a_{1,n} P_K(x)f_n\Big)\\
&= \sum_{n\in\mathbb{N}} \epsilon(x)\big(e_1 \otimes a_{1,n} P_K(x)f_n\big)\\
&= \sum_{n\in\mathbb{N}} a_{1,n} P_K(x)f_n.
\end{align*}
Thus we have $K(x) = \bigvee_{n\in\mathbb{N}} P_K(x) f_n$, as claimed. Now, since
		\[
		x\mapsto \langle e_{m_1} \otimes P_K(x) f_{n_1},e_{m_2} \otimes P_K(x) f_{n_2}\rangle
		\]
		is measurable for all $(m_i,n_i)\in\mathbb{N}^2$, $i=1,2$, we have, in particular, that
		\[ 
		x\mapsto \langle P_K(x) f_{m}, P_K(x) f_{n} \rangle =  \langle e_1 \otimes P_K(x) f_{m},e_1 \otimes P_K(x) f_{n}\rangle
		\]
		is measurable for all $m, n \in \mathbb{N}$. Therefore by Remark \ref{FollandMeasField}, $\{K(x)\}$ is a measurable field of Hilbert spaces with fundamental sequence $\{P_K(x) f_n\}_{n\in\mathbb{N}}$.
	\end{proof}
	\begin{lemma}\label{vectorsep}
	Let $(X,\mathcal{M},\mu)$ be a complete, $\sigma$-finite measure space such that $L^2(X)$ is separable (e.g., $X$ a measurable subset of $\mathbb{R}^n$). Let $\mathscr{H}$ be a separable Hilbert space. Then $L^2(X;\mathscr{H})$ is separable.
	\end{lemma}
	\begin{proof}
		This follows from the fact that $L^2(X;\mathscr{H}) \cong L^2(X)\otimes \mathscr{H}$---see \cite{ReedSimon}, p.~52.
		\end{proof}
	The theorems that follow may be compared with Theorem 2.4 in \cite{BownikRoss}, and parts of the proofs follow the arguments presented there.
	\begin{theorem}\label{mainpointwisedecompthm}
		Let $(X,\mathcal{M},\mu)$ and $\mathscr{H}$ be as in the preceding lemma, and regard $L^2(X;\mathscr{H})$ as the direct integral of the constant field $\{\mathscr{H}_x \equiv \mathscr{H}\}$. Let $\mathcal{A}$ be a countable dense subset of $L^2(X;\mathscr{H})$. 
		\begin{list}{(\arabic{ctr})}{\usecounter{ctr}}
			\setcounter{ctr}{0}
			\item Let $W\subset L^2(X;\mathscr{H})$ be a closed subspace. Then $\{P_W\Phi: \Phi \in \mathcal{A}\}$ is a dense subset of $W$, so that $W=\bigvee \{P_W\Phi: \Phi \in \mathcal{A}\}$.
			\item\label{constructingJWA} If, for each $\Phi \in\mathcal{A}$, we fix a function $f(\Phi)$ in the equivalence class $[P_W\Phi]_{L^2}$, then
			\begin{equation}\label{rangefromsubspace}
				J_W^\mathcal{A}[f](x) = \bigvee\{f(\Phi)(x): \Phi \in\mathcal{A}\}
			\end{equation}
			defines a range function. Moreover, $J_W^\mathcal{A}[f]$ is---up to almost everywhere equivalence---independent of the choice of representatives $f(\Phi)\in[P_W\Phi]_{L^2}$.
		\end{list}
	\end{theorem}
		\begin{proof}
		(1) Clearly $P_W\Phi \in W$  for all $\Phi\in\mathcal{A}$. Let $\xi\in W$. Since $\mathcal{A}$ is dense in $L^2(X;\mathscr{H})$, we can find a sequence $\{\Phi_n\}\subset\mathcal{A}$ converging to $\xi$. Since
		\[
		\|P_W\Phi_n - \xi\| = \|P_W\Phi_n - P_W \xi\| \leq \|\Phi_n - \xi\|,
		\]
		we can find $P_W\Phi_n$ arbitrarily close to $\xi$. Thus  $\{P_W\Phi: \Phi \in \mathcal{A}\}$ is dense in $W$, and hence $W=\bigvee \{P_W\Phi: \Phi \in \mathcal{A}\}$.
		
		(2) For each $\Phi \in\mathcal{A}$, let $f(\Phi)$ and $g(\Phi)$ be everywhere-defined functions in the equivalence class $[P_W\Phi]_{L^2}$, and define $J_W^\mathcal{A}[f]$ and $J_W^\mathcal{A}[g]$ as in \eqref{rangefromsubspace}. We show that these range functions agree almost everywhere. To see this, note that for each $\Phi \in \mathcal{A}$, there exists a null set $N_\Phi(f,g)$ such that $f(\Phi) = g(\Phi)$ on $X-N_\Phi(f,g)$. Letting $N = \bigcup_{\Phi\in\mathcal{A}} N_\Phi(f,g)$, we see that any finite linear combination $\sum_{n=1}^N a_\Phi f(\Phi)(x)$ is in $J_W^\mathcal{A}[g](x)$, for $x\in X-N$. Since $J_W^\mathcal{A}[g](x)$ is closed, it follows that $J_W^\mathcal{A}[f](x)\subset J_W^\mathcal{A}[g](x)$, for $x\in X-N$. And by a symmetric argument, $J_W^\mathcal{A}[f](x)\supset J_W^\mathcal{A}[g](x)$ for $x\in X-N$. Therefore $J_W^\mathcal{A}[f] =J_W^\mathcal{A}[g]$ almost everywhere.
		\end{proof}
		\begin{notation}
			 We write $J_W^\mathcal{A}$\label{JWA} for any of the equivalent range functions constructed as in \eqref{constructingJWA}.
		\end{notation}
	\begin{theorem}\label{mainpointwisedecompthm2}
Let $(X,\mathcal{M},\mu)$, $\mathscr{H}$, and $\mathcal{A}$ be as in the preceding theorem. Let $W\subset L^2(X;\mathscr{H})$ be a closed subspace which satisfies the following orthogonality condition:
			\begin{tcolorbox}
				\begin{center}
					For any $\xi \in L^2(X;\mathscr{H})$
					\begin{equation}\label{orthCond}
						\xi\perp W \implies \xi(x)\perp J_W^\mathcal{A}(x)\ \text{for a.e.}\ x.
					\end{equation}
				\end{center}
			\end{tcolorbox}
			Then $J_W^\mathcal{A}$, as defined in \ref{mainpointwisedecompthm2}, is a measurable range function.
\end{theorem}
		\begin{proof}
			 For each $x\in X$, let $P_J(x)$ denote the orthogonal projection of $\mathscr{H}$ onto $J_W^\mathcal{A}(x)$. Let $P_W$ denote the orthogonal projection of $L^2(X;\mathscr{H})$ onto $W$. Let $\xi \in L^2(X;\mathscr{H})$. Since $\mathcal{A}$ is dense in $L^2(X;\mathscr{H})$, there exists a sequence $\Phi_{n}$ converging to $\xi$ in $L^2(X;\mathscr{H})$; that is,
			\[
			\int_X \|\Phi_n(x) - \xi(x)\|_{\mathscr{H}}^2 \ud x \to 0.
			\]
			Define $g_n(x) := \|\Phi_n(x) - \xi(x)\|_{\mathscr{H}}^2$; then $g_n \to 0$ in $L^1(X;\mathbb{R})$. By a standard result (Corollary 2.32 in \cite{FollandReal}), there exists a subsequence $\{g_{n_j}\}_{j=1}^\infty$ tending to $0$ a.e., so that $\Phi_{n_j}(x) \to \xi(x)$ on the complement of some null set $\mathcal{N}_0$. By continuity of $P_J(x)$ and the inner product, it follows that
			\[
			\langle P_J(x)\Phi_{n_j}(x),b\rangle \to \langle P_J(x)\xi(x),b\rangle\ \ \ \text{for}\ x\in X-\mathcal{N}_0,
			\]
			for any element $b \in \mathscr{H}$. 
			
			Let $I$ denote the identity operator on $L^2(X;\mathscr{H})$---that is, $I\zeta = \zeta$ for all $\zeta \in L^2(X;\mathscr{H})$. Then $(I-P_W)\Phi_{n_j} \perp W$ for all $j=1,2,\dots $. Thus, by the orthogonality condition, there exist null sets $N_j$, $j=1,2,\dots$, such that 
			\[
			[(I-P_W)\Phi_{n_j}](x) \perp J_W^\mathcal{A}(x)\ \ \ \text{for }x\in X-N_j.
			\]
			It follows that
			\[
			P_J(x)[(I-P_W)\Phi_{n_j}](x) = 0\ \ \ \text{for }x\in X-N_j.
			\]
			Thus
			\begin{equation}\label{projectioncalcae}
				P_J(x)\Phi_{n_j}(x) = P_J(x)(P_W\Phi_{n_j})(x)\ \ \ \text{for }x\in X-N_j.
			\end{equation}
			Recall that there is a fixed function $f(\Phi_{n_j})$ that agrees with $P_W\Phi_{n_j}$ a.e., and that $f(\Phi_{n_j})(x) \in J_W^\mathcal{A}(x)$ for all $x$, by construction. Let $N_j^\prime$ be the null set on whose complement $f(\Phi_{n_j})=P_W\Phi_{n_j}$, and let $\mathcal{N}_j = N_j \cup N_j^\prime$. Then $(P_W \Phi_{n_j})(x) \in J_W^\mathcal{A}(x)$ on $X-\mathcal{N}_j$, and it follows from \eqref{projectioncalcae} that
			\begin{equation*}
				P_J(x)\Phi_{n_j}(x) = (P_W\Phi_{n_j})(x)\ \ \ \text{for }x\in X-\mathcal{N}_j.
			\end{equation*}
			But $P_W\Phi_{n_j}\in L^2(X;\mathscr{H})$ so that, in particular, $P_W\Phi_{n_j}$ is weakly measurable. This means that, for any $b\in\mathscr{H}$, $x\mapsto \langle (P_W\Phi_{n_j})(x),b\rangle$ is measurable. Let $N= \bigcup_{j=0}^\infty \mathcal{N}_j$; then $N$ is a null set. Let $\mathbbm{1}_{X-N}$ denote the indicator function of $X-N$---that is, the function taking value $1$ on $X-N$ and $0$ on $N$. Then 
			\[
			x\mapsto  \mathbbm{1}_{X-N}(x)\langle P_J(x)\Phi_{n_j}(x),b\rangle =\mathbbm{1}_{X-N}(x) \langle (P_W\Phi_{n_j})(x),b\rangle
			\]
			is measurable, being the product of two measurable functions. By the continuity of $P_J(x)$ and the inner product,
			\[
			\mathbbm{1}_{X-N}(x)\langle P_J(x)\xi(x),b\rangle = \lim_{j\to\infty} \mathbbm{1}_{X-N}(x)  \langle P_J(x)\Phi_{n_j}(x),b\rangle.
			\]
			Therefore, $x\mapsto \mathbbm{1}_{X-N}(x)\langle P_J(x)\xi(x),b\rangle$ is measurable, being the pointwise limit of a sequence of measurable functions. It follow from Proposition \ref{MeasAEresult} that $x\mapsto \langle P_J(x)\xi(x),b\rangle$ is measurable, as it agrees a.e. with a measurable function. Therefore, by the equivalence of (1) and (6) in Theorem \ref{bigequiv}, $J_W^\mathcal{A}$ is a measurable range function.
			\end{proof}
			With the notation of the preceding theorem, if $W\subset L^2(X;\mathscr{H})$ satisfies the orthogonality condition \eqref{orthCond}, then it follows by Theorem \ref{bigequiv} that $\{J_W^\mathcal{A}(x)\}_{x\in X}$ is a measurable field of Hilbert spaces with fundamental sequence $\{f_n(x) := P_J(x)e_n\}$, where $P_J(x)$ denotes the orthogonal projection from $\mathscr{H}$ onto $J_W^\mathcal{A}(x)$. The measurable vector fields are given by 
			\begin{equation}\label{EJ1}
				\mathcal{E}_{J^\mathcal{A}} = \{\xi\in \Pi_{x\in X} J_W^\mathcal{A}(x): x \mapsto \langle \xi(x),f_n(x) \rangle_x \ \text{is}\ \mu\text{-measurable for every}\ n \}.
			\end{equation}
		\begin{theorem}\label{pointwisedecompcor1}
		Let $(X,\mathcal{M},\mu)$, $\mathscr{H}$, and $\mathcal{A}$ be as in Theorem \ref{mainpointwisedecompthm}. Suppose that $W\subset L^2(X;\mathscr{H})$ is a closed subspace satisfying the orthogonality condition \eqref{orthCond}. Regarding $\int_{X}^\oplus J_W^\mathcal{A}(x)\ud\mu(x)$ as a subspace of $L^2(X;\mathscr{H})$, we have
				\[
				W = \int_{X}^\oplus J_W^\mathcal{A}(x)\ud\mu(x).
				\]
			\end{theorem}
		\begin{proof}
We first show that $W \subset \int_X^\oplus J_W^\mathcal{A}(x)\ud x$, where, as usual, we regard $\int_X^\oplus J_W^\mathcal{A}(x)\ud x$ as a subspace of $L^2(X;\mathscr{H})$ via the embedding $\iota$. Let $\xi \in W$. We seek a function $\xi_0\in \int_X^\oplus J_W^\mathcal{A}(x)\ud x$ such that $\xi = \xi_0$ a.e., so that $\iota([\xi_0]_J) =[\xi]_{L^2}$. Since $\{P_W\Phi:\Phi\in\mathcal{A}\}$ is dense in $W$, there is a sequence $\{\Phi_n\}$ such that $P_W\Phi_n\to \xi$. There exists a subsequence  $\{\Phi_{n_j}\}_{j=1}^\infty$ and a null set $N_0$ such that $P_W\Phi_{n_j}(x) \to \xi(x)$ on $X-N_0$.  Moreover, for each $j=1,2,\dots$, we have fixed a function $f(\Phi_{n_j}) \in [P_W \Phi_{n_j}]_{L^2}$, so that $f(\Phi_{n_j})(x) = P_W\Phi_{n_j}$ on the complement of some null set $N_j$. Note that $f(\Phi_{n_j})(x)\in J_W^\mathcal{A}(x)$ by definition, so that $P_W \Phi_{n_j}(x) \in   J_W^\mathcal{A}(x)$ on $X-N_j$. Let $N=\bigcup_{j=0}^\infty N_j$, and let $\xi_j = P_W \Phi_{n_j}$. Then $N$ is a null set, $\xi_j \to \xi$ pointwise on $X-N$, and $\xi_j(x)\in J_W^\mathcal{A}(x)$ on $X-N$ for all $j=1,2,\dots$. Since $J_W^\mathcal{A}(x)$ is a closed subspace for each $x$, it follows that $\xi(x) \in J_W^\mathcal{A}(x)$ for all $x\in X-N$. We define 
		\[
		\xi_0(x) = \left\{ \begin{array}{ll}
			\xi(x) &\text{if } x\in X-N\\
			0 &\text{if } x\in N.
		\end{array} \right.
		\]
		In other words, $\xi_0 = \mathbbm{1}_{X-N}\xi$, where $\mathbbm{1}_{X-N}$ denotes the indicator function of $X-N$. Then $\xi_0 \in \Pi_{x\in X} J_W^\mathcal{A}(x)$, and $\int_X \|\xi_0(x)\|^2\ud x  = \int_{X-N} \|\xi(x)\|^2\ud x < \infty$. It only remains to show $\xi_0 \in 	\mathcal{E}_{J^\mathcal{A}}$, i.e., that $x\mapsto \langle \xi_0(x),f_n(x)\rangle$ is measurable, where $f_n(x) = P_J(x)e_n$. Just as in \eqref{projtrick}, we have
		\[
		\langle \xi_0(x),f_n(x)\rangle = \mathbbm{1}_{X-N}(x)\langle \xi(x),e_n\rangle. 
		\]
		Now $x\mapsto \langle \xi(x),e_n\rangle$ is measurable since $\xi \in \mathcal{E}$ (since it is weakly measurable). It follows that $x\mapsto \langle \xi_0(x),f_n(x)\rangle $ is measurable, being the product of two measurable function. Thus $\xi_0 \in	\mathcal{E}_{J^\mathcal{A}}$. We have $[\xi]_{L^2} = \iota([\xi_0]_J)$ where $\xi_0 \in \int_X^\oplus J_W^\mathcal{A}(x) \ud x$, as desired.
		
		For the reverse inclusion, recall that any element $\eta \in L^2(X;\mathscr{H})$ can be written as $\eta = \eta_0 + \eta_1$, where $\eta_0 \in W$ and $\eta_1\in W^\perp$. Thus to show that $\int_X^\oplus J_W^\mathcal{A}(x) \ud x - W = \emptyset$, it suffices to show that $\eta \in \int_X^\oplus J_W^\mathcal{A}(x) \ud x \cap W^\perp$ implies $\eta = 0$. Let $\eta \in\int_X^\oplus J_W^\mathcal{A}(x) \ud x \cap W^\perp$. Then $\eta(x) \in J(x)$ a.e., since $\eta \in \int_X^\oplus J_W^\mathcal{A}(x) \ud x$. On the other hand, the orthogonality condition implies that $\eta(x) \perp J_W^\mathcal{A}(x)$ for almost every $x$. Thus $\eta(x)=0$ a.e.; i.e., $\eta =0$ as an element of $L^2(X;\mathscr{H})$.
		\end{proof}
		\begin{cor}\label{ProjCor}
		Let $(X,\mathcal{M},\mu)$, $\mathscr{H}$, and $\mathcal{A}$ be as in Theorem \ref{mainpointwisedecompthm}. Suppose that $W\subset L^2(X;\mathscr{H})$ is a closed subspace satisfying the orthogonality condition \eqref{orthCond}. Then
			\[
			P_W =\int_{X}^\oplus P_{J}(x)\ud x,
			\] 
			where $P_W$ is the orthogonal projection of $L^2(X;\mathscr{H})$ onto $W$ and $P_J(x)$ is the orthogonal projection from $\mathscr{H}$ onto $J_W^\mathcal{A}(x)$.
		\end{cor}
		\begin{proof}
		Since $J_W^\mathcal{A}$ is measurable, $x \mapsto P_{J}(x)$ is a measurable operator field by Theorem \ref{bigequiv}. The result follows immediately from Lemma \ref{projRem}, which says that $\int_{X}^\oplus P_{J}(x)\ud x$ is the orthogonal projection from $L^2(X;\mathscr{H})$ onto $\int_{X}^\oplus J_W^\mathcal{A}(x)\ud x = W$.
			\end{proof}
	Had we chosen a different countable dense collection $\mathcal{B} \subset L^2(X;\mathscr{H})$ in Theorem \ref{mainpointwisedecompthm}, we would have gotten a different range function $J_W^\mathcal{B}$.	
		\begin{prop}
	  If $J_W^\mathcal{A}$ and $J_W^\mathcal{B}$ satisfy the orthogonality condition (\ref{orthCond}), then $J_W^\mathcal{A} = J_W^\mathcal{B}$ almost everywhere.
	\end{prop}	
\begin{proof}
Let $P_J^{\mathcal{A}}(x): \mathscr{H} \to J_W^\mathcal{A}(x)$ and $P_J^{\mathcal{B}}(x): \mathscr{H} \to J_W^\mathcal{B}(x)$ be the orthogonal projections from $\mathscr{H}$ onto $J_W^\mathcal{A}(x)$ and $J_W^\mathcal{B}(x)$, respectively. By Corollary \ref{ProjCor}, we have
	\[
	\int_X^\oplus P_J^{\mathcal{A}}(x)\ud x= P_W=\int_X^\oplus P_J^{\mathcal{B}}(x)\ud x.
	\]
	Hence, by Lemma \ref{aeequalops},  $P_J^{\mathcal{A}}(x) =  P_J^{\mathcal{B}}(x)$ a.e., whence $J_W^\mathcal{A}(x) = J_W^\mathcal{B}(x)$ a.e..
\end{proof}
		It follows from the above Proposition that we may, by a slight abuse of language, make the following definition.
	\begin{defn}\label{mblerangefun}
		Let $W\subset L^2(X;\mathscr{H})$ be a closed subspace satisfying the orthogonality condition \eqref{orthCond}. We define \textit{the measurable range function associated with $W$} to be any of the (a.e. equivalent) functions $J_W^\mathcal{A}$ defined in \eqref{constructingJWA} of Theorem \ref{mainpointwisedecompthm2}. We denote this by $J_W$.
	\end{defn}

	\chapter[Particular Unitary Representations of \texorpdfstring{$\boldsymbol{\Gamma}$}{Gamma}]{Particular Unitary Representations of \texorpdfstring{$\boldsymbol{\Gamma}$}{Gamma}}\label{CrystalShiftandInduced}
	
	\section[The Natural Representation of \texorpdfstring{$\Gamma$}{Gamma} on  \texorpdfstring{$L^2(\mathbb{R}^n)$}{L2}]{The Natural Representation of \texorpdfstring{$\boldsymbol{\Gamma}$}{Gamma} on \texorpdfstring{$\boldsymbol{L^2(\mathbb{R}^n)}$}{L2}}\label{pihatpi}
	
	In this work, we are primarily interested in the unitary representation $\pi$ of $\Gamma$ on $L^2(\mathbb{R}^n)$ defined presently, where $\Gamma$ is a crystal group. Let $[x,L]\in\Gamma$ and $f$ a function in $L^2(\mathbb{R}^n)$ (with the usual understanding that functions equal almost everywhere are identified in $L^2(\mathbb{R}^n)$). Define $\pi[x,L]: \mathbb{R}^n \to \mathbb{C}$ by \label{pi0}
	\[
	[\pi[x,L](f)](u) = f\big([x,L]^{-1}\cdot_{\idref{1}}u\big).
	\]
	Thus,
	\[
	[\pi[x,L](f)](u) = f(L^{-1}u -x).
	\]
Since $L \in \mathrm{O}(n)$, a direct calculation shows that $\|\pi[x,L]\|_2 = \|f\|_2$. Thus $\pi[x,L]$ is an isometry on $L^2(\mathbb{R}^n)$. It is easily verified that $\pi$ is a unitary representation of $\Gamma$ on $L^2(\mathbb{R}^n)$. Recall that the Fourier transform $\mathcal{F}: L^1(\mathbb{R}^n) \to C_0(\widehat{\mathbb{R}^n})$ is given by 
	\[
	\mathcal{F}f(\chi_y) =  \int_{\mathbb{R}^n} f(u) \chi_y(u) \ud u=\int_{\mathbb{R}^n} f(u) e^{2\pi i y \cdot u} \ud u;
	\]
	and moreover, by continuous extension from a dense subset of $L^2(\mathbb{R}^n)$, $\mathcal{F}$ becomes a unitary operator from $L^2(\mathbb{R}^n)$ onto $L^2(\widehat{\mathbb{R}^n})$, called the Plancherel transform (this is Theorem \ref{Plancherel}). Define 
	\[
	\widehat{\pi} [x,L] = \mathcal{F} \pi[x,L] \mathcal{F}^{-1} \ \ \text{for all}\ [x,L]\in\Gamma.
	\]
	For $\phi \in L^2(\widehat{\mathbb{R}^n})$, let $f=  \mathcal{F}^{-1} \phi$. For $y\in \mathbb{R}^n$, we have
	\begin{align*}
		[\widehat{\pi}[x,L](\phi)](\chi_y) &= [\mathcal{F} (\pi[x,L] \mathcal{F}^{-1})(\phi)](\chi_y)\\
		&= \int_{\mathbb{R}^n} \pi[x,L](f)(u) e^{2\pi i y\cdot u}\ud u \\
		&= \int_{\mathbb{R}^n} f(L^{-1} u-x) e^{2\pi i y\cdot u}\ud u,
	\end{align*}	
	with the usual understanding that we are dealing with a.e. equivalence classes. Letting $\psi(u) = f(L^{-1}u-x) e^{2\pi i y\cdot u}$ be the integrand, we may replace $\psi$ with $\psi\circ L$ and use Theorem 2.44 in \cite{FollandReal}. We have
	\begin{align*}
		[\widehat{\pi}[x,L](\phi)](\chi_y) &= |\det{L}|\int_{\mathbb{R}^n} f(u-x) e^{2\pi i y\cdot(Lu)}\ud u\\
		&= \int_{\mathbb{R}^n} f(u-x) e^{2\pi i (L^T y)\cdot u}\ud u \\
		&= \int_{\mathbb{R}^n} f(u) e^{2\pi i (L^{T}y)\cdot(u+x)}\ud u,
	\end{align*} 
	by translation-invariance of the Lebesgue integral (Theorem 2.42 in \cite{FollandReal}); here $L^T$ denotes the transpose of the matrix $L$. Since $L \in \mathrm{O}(n)$ we have $L^T = L^{-1}$, we have
	\[
	[\widehat{\pi}[x,L](\phi)](\chi_y) = e^{2\pi i (L^{-1}y)\cdot x }\int_{\mathbb{R}^n} f(u) e^{2\pi i (L^{-1}y)\cdot u} \ud u 
		= e^{2\pi i (L^{-1}y)\cdot x }\phi(\chi_{L^{-1}y}).
\]
	Thus the unitary representation $\widehat{\pi}$ of $\Gamma$ on $L^2(\widehat{\mathbb{R}^n})$ is given by
	\begin{equation}\label{FT}
		\boxed{[\widehat{\pi}[x,L](\phi)](\chi_y) = \chi_{L^{-1}y}(x)\phi(\chi_{L^{-1}y})}\, 
	\end{equation}
	for almost every $\chi_y \in\widehat{\mathbb{R}^n}$.

	\section{Induced Representations}\label{inducedRep}
	Later in our work, a special type of unitary representation will play an important role. Given a locally compact group $G$ with a closed subgroup $H$, and given a unitary representation $(\sigma, \mathscr{H}_\sigma)$ of $H$, one can make a unitary representation $\mathrm{ind}_H^G \sigma$\label{induced} of $G$, called the \textit{representation of $G$ induced from $\sigma$}, or simply the \textit{induced representation}. When $H$ is an open subgroup of $G$---as is the case with $\mathrm{T}$ in the discrete group $\Gamma$---the theory of the induced representation is simpler; see \S\S 2.1 -- 2.2 of \cite{TaylorKaniuth} for more details.
	
	Let $\Gamma\subset \mathrm{I}(\mathbb{R}^n)$ be a crystal group with translation subgroup $\mathrm{T}$ and Bravais lattice $\mathcal{L}$. As in definition \ref{Rdomain} we take $R$ to be a fundamental domain for $\Gamma^* \cong \mathcal{L}^*\rtimes \Pi \subset \mathrm{I}(\mathbb{R}^n)$. We also let $\Omega_{\mathrm{T}^*}$ be a fundamental set for $\mathrm{T}^*$ with $\Pi R\subset \Omega_{\mathrm{T}^*} \subset \Pi\overline{R}$, as in definition \ref{paramdomaindef}. We saw previously that $\Omega_{\mathrm{T}^*}$ parametrizes the characters $\chi_\nu^{\mathrm{T}}$ of $\mathrm{T}$. For $\nu \in \Omega_{\mathrm{T}^*}$, we wish to describe $\mathrm{ind}_{\mathrm{T}}^\Gamma \chi_\nu^\mathrm{T}$. In the present case $\Gamma/\mathrm{T}$ is finite, which will simplify matters. Let $\gamma$ be a fixed cross-section for $\Gamma$, as in \ref{crosssection}. For each  $\nu \in \Omega_{\mathrm{T}^*}$, we define the Hilbert space
	\begin{multline*}
		\mathscr{H}(\Gamma,\chi_\nu^{\mathrm{T}}) = \Big\{f:\Gamma\to \mathbb{C} : f([x,L][y,\text{id}]) =  \chi_\nu^{\mathrm{T}}([y,\text{id}]^{-1})f([x,L]),\\ 
		\text{for all}\ [x,L] \in\Gamma\ \text{and all}\ [y,\text{id}]\in\mathrm{T} \Big\},
	\end{multline*}
	with inner product
	\[
	\langle f,g\rangle = \sum_{L\in \Pi} \langle f(\gamma(L)),g(\gamma(L)) \rangle.
	\]
	We define the induced representation by
	\[
	\mathrm{ind}_{\mathrm{T}}^\Gamma \chi_\nu^\mathrm{T}([x,L]) f([y,M]) = f([x,L]^{-1}[y,M]).
	\]
	for $f\in\mathscr{H}(\Gamma,\chi_\nu^{\mathrm{T}})$ and $[x,L],[y,M]\in \Gamma$. We presently give a more explicit formula for $\mathrm{ind}_{\mathrm{T}}^\Gamma\chi_\nu^\mathrm{T}$. By definition, for $f\in\mathscr{H}(\Gamma,\chi_\nu^{\mathrm{T}})$, $L\in\Pi$, and $[y,\mathrm{id}]\in\mathrm{T}$, 
	\begin{equation}\label{first}
		f\big(\gamma(L)[y,\text{id}]\big) = \chi_\nu^{\mathrm{T}}([y,\text{id}]^{-1})f\big(\gamma(L)\big).
	\end{equation}
	Define $W:\mathscr{H}(\Gamma,\chi_\nu^{\mathrm{T}}) \to l^2(\Pi)$ by 
	\begin{equation*}\label{second}
		(Wf)(L) = f\big(\gamma(L)\big)\ \ \text{for}\ L\in\ \Pi.	
	\end{equation*}
	Then for $\zeta\in l^2(\Pi)$, we have
	\begin{equation*}
		(W^{-1}\zeta)\big(\gamma(L)[y,\text{id}]\big)=\chi_\nu^{\mathrm{T}}([-y,\text{id}])\zeta(L)
	\end{equation*}
	for all $[y,\mathrm{id}]\in \mathrm{T}$. In particular, 
	\begin{equation}\label{third}
		(W^{-1}\zeta)\big(\gamma(L)\big) = \zeta(L).
	\end{equation} 
	It is easily verified that $W$ is a unitary transformation from $\mathscr{H}(\Gamma,\chi_\nu^{\mathrm{T}})$ onto $ l^2(\Pi)$\label{l2pi}. Of course, since $\Pi$ is finite, $l^2(\Pi)$ is just $\mathbb{C}^\Pi$, but we will continue to write this as $l^2(\Pi)$. For each $[x,L]\in \Gamma$, define $U^{\chi_\nu^\mathrm{T}}[x,L] \in \mathrm{U}(l^2(\Pi))$ by\label{uchinotation0}
	\[
	\boxed{U^{\chi_\nu^\mathrm{T}}[x,L] = W \mathrm{ind}_{\mathrm{T}}^\Gamma \chi_\nu^\mathrm{T}([x,L])  W^{-1}}\,.
	\]
	So $U^{\chi_\nu^\mathrm{T}}$ is a unitary representation of $\Gamma$ on $l^2(\Pi)$, and 
	\[
	\boxed{U^{\chi_\nu^\mathrm{T}} \sim \mathrm{ind}_{\mathrm{T}}^\Gamma \chi_\nu^\mathrm{T}}\,.
	\]
	\begin{tcolorbox}
		\begin{notation}\label{uchinotation}
In the sequel, when there is no loss of clarity, we shall write $U^\nu$ in place of $U^{\chi_\nu^\mathrm{T}}$.
		\end{notation}
	\end{tcolorbox}
	For $\zeta\in l^2(\Pi)$, $f = W^{-1}\zeta$, and $M \in \Pi$, we have
		\begin{equation}\label{fourth}
			\begin{aligned}
				U^\nu[x,L] \zeta(M) &= W \mathrm{ind}_{\mathrm{T}}^\Gamma \chi_\nu^\mathrm{T} ([x,L]) W^{-1} \zeta (M)  \\
				&= \mathrm{ind}_{\mathrm{T}}^\Gamma \chi_\nu^\mathrm{T} ([x,L]) f\big(\gamma(M)\big)  \\
				&= f\big([x,L]^{-1} \gamma(M)\big) .
			\end{aligned}
		\end{equation}
		Recall that the group homomorphism $Q\vert_\Gamma: \Gamma \to \Pi$ was given by $Q\vert_\Gamma([x,L]) = L$. We have
		\[
		Q\vert_\Gamma\big([x,L]^{-1}\gamma(M) \big) = \big(Q\vert_\Gamma([x,L])\big)^{-1} Q\vert_\Gamma\big(\gamma(M)\big) = L^{-1}M,
		\]
		so $[x,L]^{-1}\gamma(M)$ is in the fiber of $Q$ over $L^{-1}M$. Since $\gamma(L^{-1}M)$ is also in this fiber and $\ker(Q) =\mathrm{T}$, we know from algebra that this fiber is $\gamma(L^{-1}M)\mathrm{T}$, and hence
		\[
		[x,L]^{-1}\gamma(M) \in \gamma(L^{-1}M)\mathrm{T}.
		\]
		Thus
	\begin{equation}\label{elementinT}
		\gamma(L^{-1}M)^{-1}[x,L]^{-1}\gamma(M)\in\mathrm{T},
	\end{equation}
	and hence by equation \eqref{fourth},
	\begin{align*}
		U^\nu[x,L]\zeta(M) &= f\big([x,L]^{-1} \gamma(M)\big) \\ &=f\big(\gamma(L^{-1}M)\big(\gamma(L^{-1}M)^{-1}[x,L]^{-1}\gamma(M)\big)\big) \\
		&= \chi_\nu^{\mathrm{T}}\big(\gamma(M)^{-1}[x,L]\gamma(L^{-1}M)\big) f\big(\gamma(L^{-1}M)\big)\ \ \text{by}\ \eqref{first}.
	\end{align*}
	Therefore, by \eqref{third}, we have
	\begin{equation}\label{Uchi}
		\boxed{[U^\nu[x,L]\zeta](M) = \chi_\nu^{\mathrm{T}}\big(\gamma(M)^{-1}[x,L]\gamma(L^{-1}M)\big) \zeta(L^{-1}M)}\, 
	\end{equation}
for all $\nu\in \Omega_{\mathrm{T}^*}$ and $\zeta \in l^2(\Pi)$. It follows that we can define $\chi_\nu^{\mathrm{T}}$ on point masses by
\begin{equation}\label{uchispointmass} 
U^\nu[x,L] \delta_M =\chi_\nu^{\mathrm{T}}\big(\gamma(LM)^{-1}[x,L]\gamma(M)\big) \delta_{LM}
\end{equation}
for all $\nu\in \Omega_{\mathrm{T}^*}$ and $M\in \Pi$. Recall that $\Gamma$ has cross-section $\gamma: \Pi \to \Gamma$, with $\Gamma(\mathrm{id}_\Pi) = [0,\mathrm{id}]$; we write $\gamma(L) = [x_L,L]$ for every $L\in \Pi$, as in \ref{crosssection}. 
\begin{lemma}\label{simplifiedchi}
We have
\begin{equation}\label{simplifiedchieq}
	\chi_y^{\mathrm{T}}\big(\gamma(M)^{-1}[x,L]\gamma(L^{-1}M)\big) = \chi_y(x_{L^{-1}M} - x_M + M^{-1}Lx)
\end{equation}
for all $y\in \mathbb{R}^n$.
\end{lemma}
\begin{proof}
We have
\begin{align*}
		\gamma(M)^{-1}[x,L]\gamma(L^{-1}M)&=[-Mx_M,M^{-1}][x,L][x_{L^{-1}M}, L^{-1}M]\\
		 	&=  [-Mx_M,M^{-1}][(M^{-1}Lx+x_{L^{-1}M} ,M]\\
		 &=[-x_M + M^{-1}Lx + x_{L^{-1}M},\text{id}];
\end{align*}
whence
	\begin{align*}
		\chi_y^{\mathrm{T}}\big(\gamma(M)^{-1}[x,L]\gamma(L^{-1}M) \big)
		&= \chi_y^{\mathrm{T}}([x_{L^{-1}M} - x_M + M^{-1}Lx,\text{id}])\\
&= \chi_y(x_{L^{-1}M} - x_M + M^{-1}Lx).
	\end{align*}
\end{proof}
\begin{cor}\label{chicor}
\[
	\chi_z^{\mathrm{T}}\big(\gamma(M)^{-1}[x,L]\gamma(L^{-1}M)\big) =1 \quad \text{for all}\ z\in\mathcal{L}^*
\]
\end{cor}
\begin{proof}
This is immediate from the preceding lemma and \eqref{elementinT}.
\end{proof}
The above lemma allows us to express $U^\nu[x,L]$ in the following simplified form:
\begin{equation}\label{Uchi1}
		\boxed{[U^\nu[x,L]\zeta](M) = \chi_\nu(x_{L^{-1}M} - x_M + M^{-1}Lx) \zeta(L^{-1}M)}\, .
	\end{equation}
		Recall again that $\Omega_{\mathrm{T}^*}$ is a fundamental set for $\mathrm{T}^*$, and that the elements of $\widehat{\mathrm{T}}$ are parametrized by  $\Omega_{\mathrm{T}^*}$. That is, every element of $\widehat{\mathrm{T}}$ is of the form $\chi_\nu^{\mathrm{T}}$ for a unique $\nu \in \Omega_{\mathrm{T}^*}$. Also recall that $\Omega$ is a transversal for the action $\cdot_{\idref{12}}$ of $\Pi$ on $\Omega_{\mathrm{T}^*}$, with $R\subset \Omega \subset\overline{R}$. 
		\begin{prop} $\mathrm{ind}_{\mathrm{T}}^\Gamma (N\cdot_{\idref{6}}\chi_\nu^\mathrm{T})\sim \mathrm{ind}_{\mathrm{T}}^\Gamma \chi_\nu^\mathrm{T}$ for all $N\in \Pi$.
		\end{prop}
		\begin{proof}
		Recall that $\Gamma$ acts on $\widehat{\mathrm{T}}$ by
			\[
			([x,L]\cdot_{\idref{5}} \chi_\nu^{\mathrm{T}})([y,\text{id}]) = \chi_\nu^{\mathrm{T}}([x,L]^{-1}\cdot_{\idref{3}} [y,\text{id}])= \chi_\nu^{\mathrm{T}}([x,L]^{-1}[y,\mathrm{id}][x,L]),
			\]
	and that $\Pi$ acts on $\widehat{\mathrm{T}}$, via the cross-section $\gamma$:
		\[
		N\cdot_{\idref{6}} \chi^{\mathrm{T}} = \gamma(N)\cdot_{\idref{5}} \chi^{\mathrm{T}}.
		\]	
		For any $[x,L]\in\Gamma$, the induced representations $\mathrm{ind}_{\mathrm{T}}^\Gamma ([x,L]\cdot_{\idref{5}}\chi_\nu^\mathrm{T})$ and $\mathrm{ind}_{\mathrm{T}}^\Gamma \chi_\nu^\mathrm{T}$ are equivalent, by Theorem 2.39 in \cite{TaylorKaniuth}, and hence the result follows.
		\end{proof}
		\begin{cor}\label{identificationinduced}
			We have
			\begin{equation}
				\{[\mathrm{ind}_{\mathrm{T}}^\Gamma \chi_\nu^\mathrm{T}]: \nu \in \Omega_{\mathrm{T}^*}\} = 	\{[\mathrm{ind}_{\mathrm{T}}^\Gamma \chi_\omega^\mathrm{T}]: \omega \in \Omega\},
			\end{equation}
			where the square brackets denote equivalence classes.
		\end{cor}	
		\begin{proof}
It was shown in \S\ref{PionOmegaTStar} that $\{\chi_\omega^{\mathrm{T}}: \omega\in \Omega\}$ is a transversal for the action of $\Pi$ on $\widehat{\mathrm{T}}$. Thus given any $\nu \in \Omega_{\mathrm{T}^*}$, there exists $\omega \in \Omega$ and $N\in \Pi$ such that $\chi_\nu^\mathrm{T} = N\cdot_{\idref{6}}\chi_\omega^\mathrm{T}$, whence 
\[
\mathrm{ind}_{\mathrm{T}}^\Gamma \chi_\nu^\mathrm{T} = \mathrm{ind}_{\mathrm{T}}^\Gamma (N\cdot_{\idref{6}}\chi_\omega^\mathrm{T})\sim \mathrm{ind}_{\mathrm{T}}^\Gamma \chi_\omega^\mathrm{T}.
\]
\end{proof}
	
	\begin{rem}
		The intertwining operator establishing the equivalence $U^{N\cdot_{\idref{6}}\chi_\omega^\mathrm{T}} \sim U^\omega$ for $N\in \Pi$ and $\omega\in \Omega$ can be given explicitly. Indeed, define the unitary operator $\Psi_N$ on $l^2(\Pi)$ by 
		\[
		\boxed{\Psi_N \zeta(M)= \chi_\omega^{\mathrm{T}}(\alpha(M,N)^{-1})\zeta(MN)}\,,\ \ \ \text{for}\ \zeta \in l^2(\Pi) ,
		\]
		where $\alpha: \Pi\times \Pi \to \mathrm{T}$\label{cocycle0} is the \textit{cocycle} defined in \ref{cocycle}. Recall that $\alpha$ satisfies equation \eqref{cocycledef} for all $L,M\in \Pi$. We verify directly that $\Psi_N$ satisfies  
		\begin{equation}\label{PsiConj}
			U^{N\cdot_{\idref{6}}\chi_\omega^\mathrm{T}}[x,L] =\Psi_N U^\omega[x,L]\Psi_N^{-1}\ \ \ \text{for all}\ [x,L]\in \Gamma,
		\end{equation}
		so that $U^{N\cdot_{\idref{6}}\chi_\omega^\mathrm{T}} \sim U^\omega$ via the intertwining operator $\Psi_N$. First note that, for any $\zeta\in l^2(\Pi)$ and $M\in\Pi$
		\[
		\Psi_N^{-1}\zeta(M) = \chi_\omega^{\mathrm{T}}(\alpha(MN^{-1},N))\zeta(MN^{-1}).
		\]
		Thus we have
		\begin{align}\label{psiequiv}
			&[\Psi_N\big(U^\omega[x,L]\Psi_N^{-1}\zeta\big)](M)\nonumber \\
			&\qquad = \chi_\omega^{\mathrm{T}}(\alpha(M,N)^{-1})[U^\omega[x,L]\Psi_N^{-1}\zeta](MN)\nonumber\\
			&\qquad=\chi_\omega^{\mathrm{T}}(\alpha(M,N)^{-1})\chi_\omega^{\mathrm{T}}\big(\gamma(MN)^{-1}[x,L]\gamma(L^{-1}MN)\big) [\Psi_N^{-1}\zeta](L^{-1}MN)\nonumber\\
			&\qquad = \chi_\omega^{\mathrm{T}}(\alpha(M,N)^{-1})\chi_\omega^{\mathrm{T}}\big(\gamma(MN)^{-1}[x,L]\gamma(L^{-1}MN)\big) \chi_\omega^{\mathrm{T}}(\alpha(L^{-1}M,N))\zeta(L^{-1}M)\nonumber\\
			&\qquad = \chi_\omega^{\mathrm{T}}\big(\alpha(M,N)^{-1}\gamma(MN)^{-1}[x,L]\gamma(L^{-1}MN)\alpha(L^{-1}M,N)\big)\zeta(L^{-1}M).
		\end{align}
		Now it follows from \eqref{cocycledef} that
		\[
		\alpha(M,N)^{-1}\gamma(MN)^{-1} = \gamma(N)^{-1}\gamma(M)^{-1},
		\]
		and
		\[
		\gamma(L^{-1}MN)\alpha(L^{-1}M,N) = \gamma(L^{-1}M)\gamma(N).
		\]
		Therefore \eqref{psiequiv} becomes
		\begin{align*}
			&= \chi_\omega^{\mathrm{T}}\big(\gamma(N)^{-1}\gamma(M)^{-1}[x,L] \gamma(L^{-1}M)\gamma(N)\big)\zeta(L^{-1}M)\\
			&\qquad= (N\cdot_{\idref{6}} \chi_\omega^{\mathrm{T}})\big(\gamma(M)^{-1}[x,L] \gamma(L^{-1}M)\big)\zeta(L^{-1}M)\\
			&\qquad=	U^{N\cdot_{\idref{6}}\chi_\omega^\mathrm{T}}[x,L]\zeta(M).
		\end{align*}
		Thus we have verified \eqref{PsiConj}. Moreover, it follows from Proposition \ref{actionsrel} that $U^{N\cdot_{\idref{12}}\omega} \sim U^\omega$ for all $N \in \Pi$.
	\end{rem}
	\begin{defn}\label{stab} Given a group $G$ acting on a set $X$, and given $x\in X$, we define $G_x = \{g\in G: g\cdot x = x\}$. This is called the \textit{stabilizer of $x$}.
	\end{defn}
	
	Given a unitary representation $\sigma$ of a locally compact group $G$ on $\mathscr{H}$, recall that $\rho$ is a subrepresentation of $\sigma$ if there exists a closed subspace $K$ of $\mathscr{H}$ such that $\sigma|_K = \rho$. Recall also that if $\sigma= \bigoplus_{i\in I} \sigma_i$ with $\sigma_i = \rho$ for some fixed representation $\rho$, then we say that $\sigma$ is a multiple of $\rho$. For a proof of the following lemma, see pp.~144--145 in \cite{TaylorKaniuth} or Proposition 6.40 in \cite{FollandAbstract}.
	\begin{lemma}\label{InducedIrreducible}
		Let $G$ be a locally compact group, and let $N$ be an abelian closed normal subgroup of finite index in $G$. Let $\chi \in \widehat{N}$. If $\sigma$ is an irreducible representation of the stabilizer $G_{\chi}$ whose restriction to $N$ is a multiple of $\chi$, then $\mathrm{ind}_{G_{\chi}}^G \sigma$ is an irreducible representation of $G$.
	\end{lemma}

	\begin{theorem}[Theorem 4.3 in \cite{TaylorKaniuth}]\label{bigthm}
		Let $G$ be a locally compact group and $N$ an abelian closed normal subgroup of $G$ such that $[G:N]<\infty$. Let $S \subset \widehat{N}$ be such that the intersection of $S$ with each $G$-orbit in $\widehat{N}$ is a singleton. For $\chi \in S$, let $\widehat{G_\chi}^\chi$ denote the equivalence classes of irreducible representations of $G_\chi$ which restrict to a multiple of $\chi$ on $N$. Then 
		\[
		\widehat{G} = \bigcup_{\chi\in S} \big\{[\mathrm{ind}_{G_\chi}^G \sigma]: [\sigma] \in\widehat{G_\chi}^\chi \big\}.
		\]
		Here, $[\mathrm{ind}_{G_\chi}^G \sigma]$ refers to the equivalence class containing $\mathrm{ind}_{G_\chi}^G \sigma$.
	\end{theorem}
	This theorem applies in our situation with $\Gamma$ playing the part of $G$ and with $\mathrm{T}$ in place of $N$. Note that the role of $S$ is satisfied by $\upsilon(\Omega) \subset \widehat{\mathrm{T}}$, where $\upsilon$ is defined in \ref{upsilon}. Indeed, we saw in \eqref{transveralpiThat} that $\upsilon(\Omega)$ is a transversal for the action of $\Pi$ on $\widehat{\mathrm{T}}$. Let $\Gamma_{\chi_\nu^{\mathrm{T}}}$ denote the stabilizer of $\chi_\nu^{\mathrm{T}}$ under the action $\idref{5}$, when $\nu\in\Omega_{\mathrm{T}^*}$, and let $\Pi_{\chi_\nu^\mathrm{T}}$ denote the stabilizer of $\chi_\nu^{\mathrm{T}}$ under the action $\idref{6}$. We make the following definition to reduce notational clutter:
	\begin{tcolorbox}
		\begin{notation}\label{GammaandPiStab}
			In what follows, we write $\Gamma_{\nu}$ in place of $\Gamma_{\chi_\nu^\mathrm{T}}$, and $\Pi_\nu$ in place of $\Pi_{\chi_\nu^\mathrm{T}}$. We write $\widehat{\Gamma_{\nu}}^\nu$ for the equivalence classes of irreducible representations of $\Gamma_{\chi_\nu^\mathrm{T}}$ which restrict to a multiple of $\chi_\nu^{\mathrm{T}}$ on $\mathrm{T}$.
		\end{notation}
	\end{tcolorbox}
Theorem \ref{bigthm} applies to give:
	\begin{tcolorbox}
		\begin{theorem}\label{bigtheoremapp}  Let $\Gamma$ be a crystallographic group, and let $\Omega$ be the fundamental set defined in \ref{Omega}. With the notation defined above,
	\begin{equation}\label{GammaHatParam}
		\widehat{\Gamma}=\bigcup_{\omega\in\Omega} \big\{[\mathrm{ind}_{\Gamma_{\omega}}^\Gamma \sigma]: [\sigma] \in \widehat{\Gamma_{\omega}}^\omega\big\},
	\end{equation}
	where the brackets denote the equivalence classes containing these representations.
		\end{theorem}
	\end{tcolorbox}
 We show that the action of $\Gamma$ on $\widehat{\mathrm{T}}$ reduces to the action of $\Pi$. Indeed, for $\nu\in\Omega_{\mathrm{T}^*}$, we have by Proposition \ref{actionsrel} that
 \begin{equation}\label{action5anygamma}
	[x,L]\cdot_{\idref{5}} \chi_\nu^{\mathrm{T}} = L\cdot_{\idref{6}} \chi_\nu^{\mathrm{T}}	= Q\vert_\Gamma([x,L])\cdot_{\idref{6}} \chi_\nu^{\mathrm{T}}.
	\end{equation}
	Thus, for $\nu \in \Omega_{\mathrm{T}^*}$, we have
	\begin{equation}\label{stabilizerofchiomega}
		\begin{aligned}
			\Gamma_{\nu}&= \{[x,L]\in\Gamma: [x,L]\cdot_{\idref{5}} \chi_\nu^{\mathrm{T}}= \chi_\nu^{\mathrm{T}}\}\\
			&=\{[x,L]\in\Gamma:  Q\vert_\Gamma([x,L])\cdot_{\idref{6}}\chi_\nu^{\mathrm{T}} = \chi_\nu^{\mathrm{T}}\}\\
			&=Q\vert_\Gamma^{-1}(\Pi_{\nu}).
		\end{aligned}
	\end{equation}
		The following set will play a fundamental role in the sequel.
		\begin{tcolorbox}
		\begin{defn}\label{Omega0}
			We define the set of \textit{free points} in $\Omega$ to be
			\[
			\Omega_0 = \{\omega\in \Omega : \Pi_\omega = \{\mathrm{id}\}\}.
			\]
		\end{defn}
	\end{tcolorbox}
	It follows from \eqref{stabilizerofchiomega} that
		\begin{equation}\label{stabilizer}
		\Gamma_\omega = (Q\vert_\Gamma)^{-1}(\{\mathrm{id}\}) = \mathrm{T}\ \ \text{for}\ \omega\in \Omega_0.
	\end{equation}
	\begin{lemma}
	Let $R$ be the fixed fundamental domain for $\Gamma^*$, as defined in \ref{Rdomain}. Then $R \subset \Omega_0$.
	\end{lemma}
	\begin{proof}
	Recall from \eqref{pionthatreflectparamdomain} that when $\omega\in R$, 
	\[
	L\cdot_{\idref{6}}\chi_\omega^{\mathrm{T}} = \chi_{L\omega}^{\mathrm{T}}.
	\]
	Moreover, it follows from the definition of $R$ that the sets $\{[z,M]R:[z,M] \in\Gamma^*\}$ are mutually disjoint, and hence
	\[
	[z,M]R\cap R = \emptyset\ \ \text{when}\ [z,M]\neq [0,\mathrm{id}].
	\]
	In particular, 
	\[
	[0,L]R\cap R = \emptyset\ \ \text{when}\ L\neq \mathrm{id}.
	\]
	Hence, $L\omega \notin R$, when $\omega\in R$ and $L\neq\mathrm{id}$. In particular, $L\omega \neq \omega$ for such $L$ and $\omega$. It follows that $\Pi_\omega = \{\mathrm{id}\}$ whenever $\omega \in R$. 
	\end{proof}

	\begin{lemma}\label{Omega0onbase}
\begin{equation}
	\boxed{\Omega_0 = \{\omega \in\Omega :  L\cdot_{\idref{12}}\omega =\omega \implies L=\mathrm{id} \} }\, .
	\end{equation}
	\end{lemma}
	\begin{proof}
 Let $\omega \in\Omega$. Recall from \eqref{actionsixtwelverel} that
\[
L\cdot_{\idref{6}} \chi_\omega^{\mathrm{T}}=\chi_{L\cdot_{\idref{12}}\omega}^{\mathrm{T}}.
\] 
Thus, 
\begin{equation*}
	\Pi_\omega = \{L\in \Pi: 	L\cdot_{\idref{6}}\chi_\omega^{\mathrm{T}}=\chi_\omega^{\mathrm{T}} \} = \{L\in \Pi: \chi_{L\cdot_{\idref{12}}\omega}^{\mathrm{T}} = \chi_\omega^{\mathrm{T}} \}.
\end{equation*}
Now it follows from \eqref{ParamOfThat} and the above equation that
\begin{equation*}
	\Pi_\omega= \{L\in \Pi: L\cdot_{\idref{12}}\omega = \omega \}.
\end{equation*}	
The result now follows.
	\end{proof}
If we let $\Pi^{\cdot_{\idref{12}}}_{\omega}$ denote the stabilizer of $\omega$ under the action $\cdot_{\idref{12}}$, then Lemma \ref{Omega0onbase} can be expressed as 
\begin{equation}
\Omega_0 = \{\omega\in\Omega: \Pi^{\cdot_{\idref{12}}}_{\omega} = \{\mathrm{id}\}\}.
\end{equation}

	\begin{tcolorbox}
		\begin{lemma}\label{uchiirred}
			$U^\omega$ is irreducible for each $\omega \in \Omega_0$.
		\end{lemma}
	\end{tcolorbox}
	\begin{proof}
Let $\omega \in \Omega_0$. By \eqref{stabilizer} we have $\Gamma_{\omega} = \mathrm{T}$, whence
\begin{equation}
	U^\omega \sim \mathrm{ind}_{\mathrm{T}}^\Gamma \chi_\omega^\mathrm{T} =\mathrm{ind}_{\Gamma_{\omega}}^\Gamma \chi_\omega^\mathrm{T} \ \ \ \text{for}\ \omega \in \Omega_0.
\end{equation}
The result follows immediately from Lemma \ref{InducedIrreducible}.
	\end{proof}
	Note that by \eqref{stabilizer}, $\Gamma_{\omega}=\mathrm{T}$ precisely when $\omega \in\Omega_0$. Hence for $\omega\in\Omega_0$, we have that $\widehat{\Gamma_{\omega}}^\omega$ is simply the equivalence classes of irreducible representations of $\mathrm{T}$ that are multiples of $\chi_\omega^{\mathrm{T}}$. Irreducible representations of $\mathrm{T}$ are one-dimensional, so $\widehat{\Gamma_{\omega}}^\omega$ is just the equivalence classes of irreducible representations of $\mathrm{T}$ that coincide with the equivalence class of $\chi_\omega^{\mathrm{T}}$.  But from our parametrization result, \eqref{ParamOfThat}, $\chi_\omega^{\mathrm{T}} \nsim \chi_\nu^{\mathrm{T}}$ for any other $\nu \in 	\Omega_{\mathrm{T}^*}$. Thus, 
	\begin{equation}\label{stabmult}
		\widehat{\Gamma_{\omega}}^\omega= [\{\chi_\omega^\mathrm{T}\}]\ \ \ \text{for all}\ \omega \in \Omega_0.
	\end{equation}
Now Theorem \ref{bigtheoremapp} gives
	\begin{equation*}
		\widehat{\Gamma}=  \big\{[\mathrm{ind}_{\mathrm{T}}^\Gamma \chi_\omega^{\mathrm{T}}]: \omega\in \Omega_0 \big\} \cup  \bigcup_{\omega\in\Omega - \Omega_0} \big\{[\mathrm{ind}_{\Gamma_{\omega}}^\Gamma \sigma]: [\sigma] \in \widehat{\Gamma_{\omega}}^\omega\big\}.
	\end{equation*}
	In other words, if we can show that $\Omega - \Omega_0$ is of Lebesgue measure zero, then almost all the elements of $\widehat{\Gamma}$ are of the form $[\mathrm{ind}_{\mathrm{T}}^\Gamma \chi_\omega^{\mathrm{T}}]$, for $\omega \in \Omega$. Moreover, since $\mathrm{ind}_{\mathrm{T}}^\Gamma \chi_\omega^{\mathrm{T}} \sim U^\omega$, it follows that almost all the elements of $\widehat{\Gamma}$ are of the form $[U^\omega]$, for $\omega \in \Omega$. We invoke a result in Ratcliffe \cite{Ratcliffe} to show that $\Omega - \Omega_0$ is null. For this we need another piece of terminology: a fundamental region $\mathcal{R}$ for a discrete group of isometries is called \textit{proper} if $\partial\mathcal{R}$ is of measure zero. Recall that our fixed fundamental domain $R$ for $\Gamma^*$ is convex and locally finite (Remarks \ref{fundDomainFact1} and  \ref{fundDomainFact2}). Thus by Corollary 1 on p.~244 of \cite{Ratcliffe}, $R$ is proper, whence $\partial R = \overline{R}-R$ is of measure zero. Now since $\Omega - \Omega_0\subset \partial R$, it follows that this too is of measure zero. Thus we have may modify Theorem \ref{bigtheoremapp} to give:
	\begin{tcolorbox}
	\begin{theorem}\label{bigtheoremcor}
		Let $\Gamma$ be a crystal group with fundamental set $\Omega$ defined in \ref{Omega} and $\Omega_0$ as defined in \ref{Omega0}. Then $\Omega -\Omega_0$ is of measure zero, and
		\begin{equation}\label{Ghatform}
			\widehat{\Gamma}= \big\{[U^\omega]: \omega\in \Omega_0 \big\} \cup   \bigcup_{\omega\in\Omega - \Omega_0} \big\{[\mathrm{ind}_{\Gamma_{\omega}}^\Gamma \sigma]: [\sigma] \in \widehat{\Gamma_{\omega}}^\omega\big\},
		\end{equation}
		where $[U^\omega]$ and $[\mathrm{ind}_{\Gamma_{\omega}}^\Gamma \sigma]$ denote the equivalence classes of these representations.
	\end{theorem}
	\end{tcolorbox}
	\begin{rem}
	Recall that we have the explicit formula \eqref{Uchi1} for $U^\omega$. When $\Gamma$ is symmorphic, we can take $\gamma(M) = [0,M]$ and $\gamma(L^{-1}M) = [0,L^{-1}M]$; that is, $x_M = x_{L^{-1}M} = 0$. It follows that
\[	
		[U^\omega[x,L]\zeta](M) = \chi_\omega(M^{-1}Lx)\zeta(L^{-1}M)= e^{2\pi i L^{-1}M\omega\cdot x}\zeta(L^{-1}M).
	\]
	\end{rem}
	The decomposition \eqref{Ghatform}, together with the existence of a closed formula for $U^\omega$ when $\omega\in \Omega_0$, as well as Theorem \ref{centraldecomptheorem}, motivates us to look for decomposition of the natural representation $\widehat{\pi}$ as a direct integral of induced representations over $\Omega_0$. We shall see that it will be convenient to consider a direct integral over $R$, which is conull in $\Omega_0$ and has the advantage of being clearly Borel.

	\newpage
	
	\chapter{Direct Integral Decomposition of \texorpdfstring{${\boldsymbol{\widehat\pi}}$}{Pi-Hat}}\label{DecompChapter}
	
	\section[An identification of \texorpdfstring{$L^2$}{L2} spaces]{An identification of \texorpdfstring{$\boldsymbol{L^2}$}{L2} spaces}
In \S\ref{EssentialFourier} we introduced the Plancherel transform, which is an isometric isomorphism from $L^2(G)$ to $L^2(\widehat{G})$, where $G$ is a locally compact abelian group. In this section, we introduce another identification between $L^2(\mathbb{R}^n)$ and $L^2(\widehat{\mathbb{R}^n})$.
\begin{lemma}\label{L2identification}
	Let $X$ and $Y$ be topological spaces, and let $\mu$ be a Radon measure $\mu$ on $(X,\mathcal{B}_X)$. Given $f:X\to Y$ a continuous map, $I(g) := \int_X g\circ f \ud\mu(x)$ is a positive linear functional on $C_c(Y)$---the space of compactly supported continuous functions on $Y$. Let $\nu$ be the unique Radon measure on $(Y,\mathcal{B}_Y)$---whose existence is guaranteed by the Riesz Representation Theorem (Theorem 7.2 in \cite{FollandReal})---such that $I(g) = \int_Y g \ud\nu(y)$. If $f$ is a homeomorphism, then the map $\Psi_f: C_c(Y) \to C_c(X)$ defined by 
	\begin{equation}\label{psif}
		\Psi_f(g) = g\circ f
	\end{equation}
	is a bijection, and extends continuously to an isometric isomorphism from $L^2(Y,\mathcal{B}_Y,\nu)$ to $L^2(X,\mathcal{B}_X,\mu)$ in such a way that $\Psi_f$ is given by \eqref{psif} on all of $L^2(Y,\mathcal{B}_Y,\nu)$.
\end{lemma}	 	
\begin{proof}
	Let $f: X\to Y$ be a homeomorphism.	We use that fact that if $\mu$ is a Radon measure on $X$, then $C_c(X)$ is dense in $L^p(X,\mu)$ for $1\leq p <\infty$ (Proposition 7.9 in \cite{FollandReal}). To see exactly how $\Psi_f$ extends, suppose that $g\in L^2(Y,\mathcal{B}_Y,\nu)$ and that $g_n\in C_c(Y)$ is such that $g_n\to g$ in $L^2(Y)$. Then
	\begin{align*}
		\|\Psi_f(g_n) - \Psi_f(g_m)\|_X^2&= \int_X |(g_n-g_m)\circ f|^2\ud\mu\\
		&= \int_X |g_n-g_m|^2\circ f\ud\mu\\
		&= I(|g_n-g_m|^2)\\
		&= \int_Y |g_n-g_m|^2\ud\nu\\
		&= \|g_n-g_m\|_Y^2 \to 0
	\end{align*}
	as $m,n\to\infty$. Since $L^2(X,\mathcal{B}_X,\mu)$ is complete, it follows that $\Psi_f(g_n)$ converges; we define $\Psi_f(g) \in L^2(X,\mathcal{B}_X,\mu)$ to be its limit. Moreover, it is easy to show that $\Psi_f(g) = g\circ f$; for there exists a subsequence $\{g_{n_k}\}$ such that $\Psi_f(g_{n_k})$ converges to $\Psi_f(g)$ almost everywhere. And $\{g_{n_k}\}$ converges to $g$ almost everywhere, so that $g_{n_k} \circ f\to g\circ f$ almost everywhere, and hence $\Psi_f(g) = g\circ f$ almost everywhere. 
	
	Lastly, we check that $\Psi_f$ is a unitary map. Linearity is clear. That $\Psi_f$ is bijective follows from \eqref{psif} and the fact that $f$ is assumed to be a bijection. And $\Psi_f$ is an isometry, because
	\[
	\|\Psi_f(g)\|^2 = \int_X |g\circ f|^2\ud \mu  =\int_X |g|^2 \circ f \ud \mu = I(|g|^2) = \int_Y |g|^2 \ud \nu = \|g\|^2,
	\]
	for all $g\in L^2(Y,\mathcal{B}_Y,\nu)$.
\end{proof}
\begin{rem}
	It can be shown that the measure $\nu$ in the above lemma is just $\mu\circ f^{-1}$, the \textit{image of the measure $\mu$ under $f$}, defined by $\mu\circ f^{-1}(E) := \mu\big(f^{-1}(E)\big)$. See Theorem \ref{imagemeasurethm} and Corollary \ref{L2isom}.
\end{rem}
Recall that we sometimes identify $\widehat{\mathbb{R}^n}$ with $\mathbb{R}^n$ via $\chi_y \leftrightarrow y$, where $\chi_y\in  \widehat{\mathbb{R}^n}$ is the character $\chi_y(u) = e^{2\pi i y\cdot u}$. It will be convenient to give this identification a name.
\begin{tcolorbox}
	\begin{defn}\label{PhiMap}
		We define $\Phi: \mathbb{R}^n \to \widehat{\mathbb{R}^n}$ by $\Phi(y)= \chi_y$, where $\chi_y(u) = e^{2\pi i y\cdot u}$.
	\end{defn}
\end{tcolorbox} 
\begin{ex}\label{Phimap}
	By Example 23.27(f) in \cite{HewittRoss}, $\Phi$ is a topological group isomorphism. Then it follows from Theorem \ref{imagemeasurethm} that
	\begin{equation}\label{Phimapimage0}
		\int_{\widehat{\mathbb{R}^n}} g(\chi) \ud(m\circ \Phi^{-1})(\chi) = \int_{\mathbb{R}^n} g(\Phi(y))\ud m(y)
	\end{equation}
	for all measurable $g$ on $\widehat{\mathbb{R}^n}$. It then follows from Lemmas \ref{L2identification} and \ref{imageofhomeoradon} that the unique Radon measure on $\widehat{\mathbb{R}^n}$ is $m\circ \Phi^{-1}$. Let $\widehat{\mathscr{L}}$\label{Lebesguesigalg0} denote the completion of $\mathcal{B}_{\widehat{\mathbb{R}^n}}$ with respect to $m\circ \Phi^{-1}$, and let $\M$ denote the completion of $m\circ \Phi^{-1}$. Then by Theorem \ref{imagemeasureextended}, we have
	\begin{equation}\label{Phimapimage}
		\int_{\widehat{\mathbb{R}^n}} g(\chi) \ud\M(\chi) = \int_{\mathbb{R}^n} g(\chi_y)\ud m(y).
	\end{equation}
	Moreover, by Corollary \ref{imagemeasExcor} that $g\mapsto g\circ \Phi$ is an isometric isomorphism from $L^2(\widehat{\mathbb{R}^n},\widehat{\mathscr{L}},\M)$ onto $L^2(\mathbb{R}^n, \mathscr{L}, m)$.\label{Lebesguesigalg}
\end{ex}

\section[Decomposition of the unitary representation \texorpdfstring{$\widehat{\pi}$}{Pi-Hat}]{Decomposition of the unitary representation \texorpdfstring{$\boldsymbol{\widehat{\pi}}$}{Pi-Hat}}\label{Decomposition of unitary}

Recall our formula for the unitary representation $\widehat{\pi}$:
\begin{equation}
	[\widehat{\pi}[x,L](\phi)](\chi_r) = \chi_{L^{-1}r}(x)\phi(\chi_{L^{-1}r}),
\end{equation}
	with the usual understanding that we are dealing with a.e. equivalence classes of functions. As before, let $\Gamma^* = \mathcal{L}^*\rtimes \Pi\subset\mathrm{I}(\mathbb{R}^n)$ and let $R$ be as defined in \ref{Rdomain}. 
\begin{tcolorbox}
\begin{defn}\label{Upsilondef}
We define a map
	\[
	L^2(\widehat{\mathbb{R}^n}) \stackrel{\Upsilon }{\to} L^2\big(R; l^2(\Pi) \otimes l^2(\mathcal{L}^*)\big) 
	\]
	by
\begin{equation*}
(\Upsilon\phi)(r) = \sum_{[z,M]\in \Gamma^*} \chi_{-r-z} (x_M) \phi(\chi_{M(r+z)})\delta_M \otimes \delta_z
\end{equation*}
\end{defn}
\end{tcolorbox}
\label{l2lstar}This can also be expressed as
\[
[(\Upsilon\phi)(r)](M,z) = \chi_{-r-z} (x_M) \phi(\chi_{M(r+z)}).
\]
It is clear that $\Upsilon$ is linear. We show that $\Upsilon$ is an isometry. 	First note that 
\[
\bigcup_{[z,M]\in \Gamma^*} [z,M] \cdot_{\idref{1}}\overline{R} = \mathbb{R}^n,
\] 
since $R$ is a fundamental domain for $\Gamma^*$. Since $R$ is proper, we have that $\partial R = \overline{R}-R$ is null. It follows from this and the fact that $\Gamma^*$ is countable, that $\bigcup_{[z,M]\in \Gamma^*} [z,M] \cdot_{\idref{1}}R$ is almost all of $\mathbb{R}^n$. Hence $\sum_{[z,M]\in\Gamma^*}\mathbbm{1}_{[z,M]\cdot_{\idref{1}}{R}}(r) = 1$ for a.e. $r \in \mathbb{R}^n$, because the tiles $[z,M]\cdot_{\idref{1}}{R}$ and $[w,N]\cdot_{\idref{1}}{R}$ are disjoint whenever $[w,N]\neq [z,M]$.  Let $\phi\in L^2(\widehat{\mathbb{R}^n})$. By Tonelli's theorem, we have
\begin{equation}\label{upsisom}
\begin{aligned}
\|\Upsilon\phi\|^2 &= \int_R \sum_{[z,M]\in \Gamma^*} \big\|\chi_{-r-z} (x_M) \phi(\chi_{M(r+z)})\delta_M \otimes \delta_z \big\|^2 \ud r\\
&=  \int_R \sum_{[z,M]\in \Gamma^*} |\chi_{-r-z} (x_M) \phi(\chi_{[z,M]\cdot_{\idref{1}} r})|^2\|\delta_M \otimes \delta_z \|^2 \ud r\\
&= \sum_{[z,M]\in \Gamma^*}  \int_R |\phi(\chi_{[z,M]\cdot_{\idref{1}} r})|^2\|\delta_M \|^2\|\delta_z \|^2 \ud r\\
&= \sum_{[z,M]\in \Gamma^*} \int_R  |\phi(\chi_{[z,M]\cdot_{\idref{1}} r})|^2 \mathbbm{1}_R(r)\ud r\\
&= \sum_{[z,M]\in \Gamma^*} |\det M^{-1}| \int_R |\phi(\chi_r)|^2  \mathbbm{1}_R([z,M]^{-1}\cdot_{\idref{1}} r) \ud r\\
&= \int_R |\phi(\chi_r)|^2  \sum_{[z,M]\in \Gamma^*}  \mathbbm{1}_R([z,M]^{-1}\cdot_{\idref{1}} r) \ud r\\
&= \int_R |\phi(\chi_r)|^2  \sum_{[z,M]\in \Gamma^*}  \mathbbm{1}_{[z,M]\cdot_{\idref{1}} R}(r) \ud r\\
&= \int_R |\phi(\chi_r)|^2  \ud r\\
&=\int_{\widehat{\mathbb{R}^n}}\big|\phi(\chi) \big|^2\ud \M(\chi)\\
&= \|\phi\|^2,
\end{aligned}
\end{equation}
	where we have used the isometric isomorphism from $L^2(\widehat{\mathbb{R}^n},\widehat{\mathscr{L}},\M)$ onto $L^2(\mathbb{R}^n, \mathscr{L}, m)$ described in Example \ref{Phimap}. It is also not hard to see
that $\Upsilon$ is surjective; indeed, given $F\in L^2\big(R;l^2(\Pi)\otimes l^2(\mathcal{L}^*)\big)$, if we define $\phi$ almost everywhere by
	\[
\phi(\chi_{M(r + z)}) = \chi_{r + z}(x_M)\langle F(r),\delta_M \otimes \delta_{z}\rangle,
	\]
then 
\begin{align*}
(\Upsilon \phi)(r) &= \sum_{[z,M]\in \Gamma^*} \chi_{-r-z} (x_M)  \chi_{r + z}(x_M)\langle F(r),\delta_M \otimes \delta_{z}\rangle\delta_M \otimes \delta_z\\
&= \sum_{[z,M]\in \Gamma^*} \langle F(r),\delta_M \otimes \delta_{z}\rangle\delta_M \otimes \delta_z\\
&= F(r),
\end{align*}
and $\phi \in L^2(\widehat{\mathbb{R}^n})$ by equation \eqref{upsisom}. Thus we have shown that 
\begin{equation*}
(\Upsilon^{-1}F)(\chi_{M(r + z)}) = \chi_{r + z}(x_M)\langle F(r),\delta_M \otimes \delta_{z}\rangle ,
\end{equation*}
That is, 
\begin{equation}\label{upsinv1}
\boxed{(\Upsilon^{-1}F)(\chi_{[z,M]\cdot_{\idref{1}} r}) = \chi_{r + z}(x_M)\langle F(r),\delta_M \otimes \delta_{z}\rangle}\, ,
\end{equation}
for $r\in R$ and $[z,M]\in\Gamma^*$. Recall that 
\[
L^2\big(R; l^2(\Pi) \otimes l^2(\mathcal{L}^*)\big) = \int_R^\oplus \mathscr{H}_r \ud r,
\]
where $\mathscr{H}_r \equiv l^2(\Pi) \otimes l^2(\mathcal{L}^*)$ is a constant field.
	\begin{lemma}\label{Uchimeas1}
		 $\{U^r[x,L]\otimes \rho\}_{r\in R}$ is a measurable field of operators on the constant field $\{\mathscr{K}_ r \equiv l^2(\Pi)\otimes \mathscr{H}_\rho\}$, where $U^r[x,L]\otimes \rho$ denotes the inner tensor product. 
		 	\end{lemma}
		 	\begin{proof}
This is equivalent to showing that, for any unitary representation $\rho$ of $\Gamma$ on $\mathscr{H}_\rho$, 
	\begin{equation}\label{UchiMeasExt1}
		 r \mapsto \langle (U^r[x,L]\otimes \rho) g,h\rangle
	\end{equation}
	is measurable for all $g, h \in l^2(\Pi)\otimes \mathscr{H}_\rho$. First note that, by Proposition 4.11 in \cite{FollandReal}, for each $\zeta\in l^2(\Pi)$, 
			\begin{equation}\label{uCHiCont1}
				 r\mapsto U^r[x,L]\zeta
			\end{equation}
			is a continuous mapping from $R$ to $l^2(\Pi)$. 
Now note that, for $\zeta\in l^2(\Pi)$ and $\beta\in \mathscr{H}_\rho$, 
			\[
			 r\mapsto U^r[x,L](\zeta)\otimes \rho(\beta) = (U^r[x,L]\otimes \rho)(\zeta\otimes \beta)
			\]
			is continuous from ${R}$ to $l^2(\Pi)\otimes \mathscr{H}_\rho$, as a consequence of \eqref{uCHiCont1} being continuous. Since the family $\mathscr{H}_0$ of finite sums of simple tensors is dense in $l^2(\Pi)\otimes \mathscr{H}_\rho$, we can take $g_n \in \mathscr{H}_0$ with $g_n\to g$. It follows from the above equation that $ r\mapsto (U^r[x,L]\otimes \rho)g_n$ is continuous for each $n$. Let $f_n( r) = (U^r[x,L]\otimes \rho)g_n$, and let $f(r) = (U^r[x,L]\otimes \rho)g$. Then $f_n \to f$ uniformly on ${R}$, since
			\begin{align*}
				\|f_n( r)-f( r) \| &= \|(U^r[x,L]\otimes \rho)(g_n-g)\|\\
				&\leq \|U^r[x,L] \otimes \rho\| \|g_n -g\|\\
				&\leq \|U^r[x,L]\|\cdot\|\rho\|\cdot\|g_n -g\|\\
				&= \|g_n -g\|.
			\end{align*}
			It follows by the uniform limit theorem (Theorem 21.6 in\cite{Munkres}) that $f$ is continuous. Then, by the continuity of the inner product, we have that
			\[
			 r \mapsto \langle f( r),h\rangle
			\]
			is continuous; i.e.: \eqref{UchiMeasExt1} is continuous (and hence measurable), where $g, h \in l^2(\Pi)\otimes \mathscr{H}_\rho$ are arbitrary. It follows from Remark \ref{measopconstantfield} that $\{U^r[x,L]\otimes \rho\}_{ r\in{R}}$ is a measurable field of operators on the constant field $\{\mathscr{K}_ r \equiv l^2(\Pi)\otimes \mathscr{H}_\rho\}$. 
\end{proof}
\begin{tcolorbox}
\begin{defn}\label{taudef}
 We define a unitary operator $\tau$ on $L^2\big(R; l^2(\Pi) \otimes l^2(\mathcal{L}^*)\big)$ by
\begin{equation*}
	\tau = \int_{R}^\oplus U^r \otimes I \ud r
	\end{equation*}
\end{defn}
\end{tcolorbox}
\begin{prop}
The map $\Upsilon$ intertwines $\hat{\pi}$ with $\tau$.
\end{prop}
\begin{proof}
We have
\begin{align}\label{mainequiv1}
&\Big[ \Big( \Big(\int_R^\oplus U^r \otimes I \ud r \Big)[x,L] (\Upsilon\phi) \Big)(r_0) \Big](M,z) \nonumber \\
&=\Big[ \Big( \big( U^{r_0}[x,L] \otimes I  \big) (\Upsilon\phi) \Big)(r_0) \Big](M,z) \nonumber\\
&=\Big[ \Big( \big( U^{r_0}[x,L] \otimes I  \big) (\Upsilon\phi) \Big)(r_0) \Big](M,z) \nonumber\\
&=\big( U^{r_0}[x,L] \otimes I \big) \Big(\sum_{[z,M]\in \Gamma^*} \chi_{-r_0-z} (x_M) \phi(\chi_{M(r_0+z)})\delta_M \otimes \delta_z\Big) \nonumber\\
&= \sum_{[z,M]\in \Gamma^*}  U^{r_0}[x,L] \big(\chi_{-r_0-z} (x_M) \phi(\chi_{M(r_0+z)})\delta_M \big)\otimes \delta_z \nonumber\\
&=\sum_{[z,M]\in \Gamma^*}\nonumber  \chi_{r_0}^{\mathrm{T}}\big(\gamma(LM)^{-1}[x,L]\gamma(M)\big)\chi_{-r_0-z}(x_M) \phi(\chi_{M(r_0+z)}) \delta_{LM} \otimes \delta_z\nonumber \\
&= \sum_{[z,M]\in \Gamma^*}\chi_{r_0}^{\mathrm{T}}\big(\gamma(M)^{-1}[x,L]\gamma(L^{-1}M)\big) \chi_{-r_0-z}(x_{L^{-1}M}) \phi(\chi_{L^{-1}M(r_0+z)}) \delta_{M} \otimes \delta_z,
\end{align}
where we have used \eqref{uchispointmass} in the penultimate equality. On the other hand, we have 
\begin{align}\label{mainequiv2}
&\big[\Upsilon(\hat{\pi}[x,L]\phi)\big](r_0) =  \sum_{[z,M]\in \Gamma^*} \chi_{-r_0-z} (x_M) (\hat{\pi}[x,L]\phi)(\chi_{M(r_0+z)})\delta_M \otimes \delta_z\nonumber\\
&=  \sum_{[z,M]\in \Gamma^*} \chi_{-r_0-z} (x_M) \chi_{L^{-1}M(r_0+z)} \phi(\chi_{L^{-1}M(r_0+z)})(x) \delta_M \otimes \delta_z. 
\end{align}
To show that \eqref{mainequiv1} and \eqref{mainequiv2} are equal, it suffices to show that 
\begin{equation}\label{desiredeq}
\chi_{r_0}^{\mathrm{T}}\big(\gamma(M)^{-1}[x,L]\gamma(L^{-1}M)\big) \chi_{-r_0-z}(x_{L^{-1}M})  = \chi_{-r_0-z} (x_M) \chi_{L^{-1}M(r_0+z)}(x).
\end{equation}
Recall from Corollory \ref{chicor} that we have
\[
\chi_z(x_{L^{-1}M} -x_M +M^{-1}Lx)=1 
\]
for all $z\in\mathcal{L}^*$; it follows that \eqref{desiredeq} reduces to 
\[
\chi_{r_0}\big(\gamma(M)^{-1}[x,L]\gamma(L^{-1}M) \big) 
= \chi_{r_0}(x_{L^{-1}M}- x_M + M^{-1}Lx),
\]
which holds by Lemma \ref{simplifiedchi}.
\end{proof}

	\section{Summary}\label{Ch5summary}
	
	We showed in this chapter that $\widehat{\pi}$ was equivalent to  $\tau = \int_{{R}}^{\oplus} U^r \otimes I \ud r$ via the intertwining operator $\Upsilon$. Recall from Section \ref{pihatpi} that $\widehat{\pi}$ was equivalent to $\pi$ via the Plancherel transform $\mathcal{F}$. Let
	\begin{equation}\label{Udef}
	\boxed{	U =  \Upsilon \circ \mathcal{F}}\, ,
	\end{equation}
	Thus we have proved one of the main results in this thesis:
	\begin{tcolorbox}[]
		\begin{theorem}\label{DecompTheorem}
			There is a unitary transformation 
			\[
			U: L^2(\mathbb{R}^n) \to L^2\big(R;l^2(\Pi)\otimes l^2(\mathcal{L}^*)\big)
			\]
			under which the natural shift representation $\pi$ is equivalent to 
			\[
			\tau = \int_{R}^\oplus U^r \otimes I \ud r.
			\]
		\end{theorem}
	\end{tcolorbox}
	\begin{rem}
	Recall from Lemma \ref{uchiirred} that $U^r$ is irreducible for each $ r\in R$. Thus, by Corollary \ref{tensorrepfactor}, $\tau$ is a direct integral of factor representations.
	\end{rem}
\begin{rem}
The map $\Upsilon$ was originally found by composing a sequence of five intertwining operators. Having found its form, we spare the reader these details.
\end{rem}

	\chapter[The Orthogonality Condition for \texorpdfstring{$\boldsymbol{\tau}$}{Tau}-Invariant \\Subspaces]{The Orthogonality Condition for \texorpdfstring{$\boldsymbol{\tau}$}{Tau}-Invariant Subspaces}\label{OrthConditionChapter}
	
	\section[Revisiting the parametrization of  \texorpdfstring{$\widehat{\mathrm{T}}$}{T-Hat}]{Revisiting the parametrization of  \texorpdfstring{$\boldsymbol{\widehat{\mathrm{T}}}$}{T-Hat}}

	Let $\mathrm{T}\subset \Gamma$ be the translation subgroup of $\Gamma$, as usual. Recall from \eqref{ParamOfThat} that the elements of $\widehat{\mathrm{T}}$ are parametrized by $\Omega_{\mathrm{T}^*}$. This was accomplished by the composition of maps
	\[
	\begin{tikzcd} 
		\Omega_{\mathrm{T}^*}\arrow[to=1-3, bend right=30, "\upsilon"{name=D}] \arrow[r, "X"{name=U}]
		& \widehat{\mathcal{L}} \arrow[r, "\vartheta^{-1}"{name=vartheta}]
		& \widehat{\mathrm{T}} \, ,
	\end{tikzcd}
	\]
	where $X(\nu) = \chi_\nu^{\mathcal{L}}$ and  $\vartheta^{-1}(\chi^{\mathcal{L}})= \chi^{\mathcal{L}} \circ \pi_1\vert_{\mathrm{T}}$; hence $\upsilon(\nu) = \chi_\nu^\mathrm{T}$ (recall definitions \ref{chiT} and \ref{upsilon}). Here $\pi_1\vert_{\mathrm{T}}$ is the projection which gives a bijection from $\mathrm{T}$ to $\mathcal{L}$, namely $\pi_1\vert_{\mathrm{T}}([x,\mathrm{id}]) = x$. Alternatively, we may view $\upsilon$ as the composition of mappings
	\[
	\begin{tikzcd} 
		\Omega_{\mathrm{T}^*}\arrow[to=1-3, bend right=30, "\upsilon"{name=D}] \arrow[r, "X_1"{name=X1}]
		& \mathbb{R}^n/\mathcal{L}^* \arrow[r, "\mathcal{C}"{name=C}]
		& \widehat{\mathrm{T}} \, ,
	\end{tikzcd}
	\]
	where
	\[
	X_1(\nu) = \nu+\mathcal{L}^*,
	\]
	and \label{C}
	\[
	\mathcal{C}(\nu+\mathcal{L}^*)= \chi_\nu^\mathrm{T},
	\]
	where $\chi_\nu^\mathrm{T}([x,\mathrm{id}]) = e^{2\pi i\nu\cdot x}$. We note that $X_1$ is just the restriction of $q$ to $\Omega_{\mathrm{T}^*}$, where $q$ is the natural quotient map
	\[
	q: \mathbb{R}^n \to \mathbb{R}^n/\mathcal{L}^*
	\]
	given by $q(y) = y+\mathcal{L}^*$. Recall as well that defined  a bijection $X_2: \mathbb{R}^n/\mathcal{L}^* \to \widehat{\mathbb{R}^n}/\mathcal{L}^\perp$ by 
	\[
	X_2(y+\mathcal{L}^*) = \chi_y\mathcal{L}^\perp.
	\]
Thus, with $\Phi$ as defined in \ref{PhiMap}, we can write
	\[
X_2(y+\mathcal{L}^*) = \Phi(y)\mathcal{L}^\perp.
\]
 Moreover, $X_2^{-1}: \widehat{\mathbb{R}^n}/\mathcal{L}^\perp \to \mathbb{R}^n/\mathcal{L}^*$ is given by 
\[
X_2^{-1}(\chi_y \mathcal{L}^\perp)= \Phi^{-1}(\chi_y) + \mathcal{L}^*  = y+ \mathcal{L}^*\ \ \text{for}\ \chi_y\in\widehat{\mathbb{R}^n}.
\]
Now note that
		\begin{align*}
	\mathcal{L}^\perp &= \{\chi_y \in \widehat{\mathbb{R}^n}: \chi_y(x) = 1\ \ \text{for all}\ x\in\mathcal{L}\}\\
	&= \{\Phi(y) \in \widehat{\mathbb{R}^n}: e^{2\pi i y\cdot x} = 1\ \ \text{for all}\ x\in\mathcal{L}\}\\
	&= \{\Phi(y) \in \widehat{\mathbb{R}^n}: y\cdot x \in \mathbb{Z}\ \ \text{for all}\ x\in\mathcal{L}\}\\
	&= \Phi(\mathcal{L}^*).
	\end{align*}
	Thus $X_2: \mathbb{R}^n/\mathcal{L}^* \to \Phi(\mathbb{R}^n)/\Phi(\mathcal{L}^*)$.
\begin{theorem}\label{X2Homeo}
The map $X_2$ described above is a homeomorphism.
\end{theorem}
\begin{proof}\footnote{Credit for the essence of this argument goes to \cite{QuotientHomeo}} We first note that the map $\widehat{q}\circ \Phi$ in the diagram below is constant on each set $q^{-1}(y+\mathcal{L}^*)$. \label{qhat}
		\begin{equation*}
		\begin{tikzcd}
			\mathbb{R}^n \arrow[d,"q"] \arrow[r, "\Phi"] & \widehat{\mathbb{R}^n} \arrow[d, "\widehat{q}"]\\
			\mathbb{R}^n/\mathcal{L}^* \arrow[r, "\tilde{\Phi}"] & \widehat{\mathbb{R}^n}/\mathcal{L}^\perp
		\end{tikzcd}
	\end{equation*} 
Indeed,
	\[
	q^{-1}(y+\mathcal{L}^*) = \{y + z:z\in\mathcal{L}^*\};
	\]
and for $z\in\mathcal{L}^*$, we have $\chi_z =\Phi(z) \in \mathcal{L}^\perp$, so that
\[
(\widehat{q}\circ \Phi)(y+z) = \widehat{q}(\chi_{y+z}) = \chi_{y+z}\mathcal{L}^\perp = \chi_{y}\chi_z\mathcal{L}^\perp = \chi_y\mathcal{L}^\perp
\]
for all $y + z\in q^{-1}(y+\mathcal{L}^*)$. That is, $\widehat{q}\circ \Phi \equiv \chi_y\mathcal{L}^\perp$ on $q^{-1}(y+\mathcal{L}^*)$. It follows from the universal property of the quotient map (Theorem 22.2 in \cite{Munkres}) that $\widehat{q}\circ \Phi$ induces a map $\tilde{\Phi}:  \mathbb{R}^n/\mathcal{L}^* \to \widehat{\mathbb{R}^n}/\mathcal{L}^\perp$ such that $\tilde{\Phi}\circ q = \widehat{q}\circ \Phi$; and moreover, $\tilde{\Phi}$ is continuous if and only if $\widehat{q}\circ \Phi$ is continuous. Since $\Phi$ is a homeomorphism and $\widehat{q}$ is a quotient map, this proves that $\tilde{\Phi}$ is continuous. Now consider the diagram
\begin{equation*}
\begin{tikzcd}
	\mathbb{R}^n \arrow[d,"q"]  & \widehat{\mathbb{R}^n} \arrow[l, "\Phi^{-1}"] \arrow[d, "\widehat{q}"]\\
	\mathbb{R}^n/\mathcal{L}^*  & \widehat{\mathbb{R}^n}/\mathcal{L}^\perp \arrow[l, "\tilde{\Psi}"]
\end{tikzcd}
\end{equation*} 
Then $q\circ \Phi^{-1}$ is constant on each set $\widehat{q}^{-1}(\chi_y\mathcal{L}^\perp)$. Indeed,
\[
\widehat{q}^{-1}(\chi_y\mathcal{L}^\perp) = \{\chi_y\chi_z: \chi _z\in \mathcal{L}^\perp\};
\]
and for $\chi_z \in\mathcal{L}^\perp$, we have $\chi_z = \Phi(z)$ with $z\in\mathcal{L}^*$, whence
\[
(q\circ \Phi^{-1})(\chi_y \chi_z) =  (q\circ \Phi^{-1})(\chi_{y+z}) = q(y+z) = y + z + \mathcal{L}^* = y +  \mathcal{L}^*.
\]
for all $\chi_y\chi_z \in \widehat{q}^{-1}(\chi_y\mathcal{L}^\perp)$. That is, $q\circ \Phi^{-1} \equiv y +  \mathcal{L}^*$ on $\widehat{q}^{-1}(\chi_y\mathcal{L}^\perp)$. It follows from Theorem 22.2 in \cite{Munkres} that $q\circ \Phi^{-1}$ induces a map $\tilde{\Psi}: \widehat{\mathbb{R}^n}/\mathcal{L}^\perp \to \mathbb{R}^n/\mathcal{L}^*$ such that $\tilde{\Psi} \circ \widehat{q} = q \circ \Phi^{-1}$; and moreover, $\tilde{\Psi}$ is continuous precisely when $q \circ \Phi^{-1}$ is. This proves that $\tilde{\Psi}$ is continuous. Now we show that $\tilde{\Psi} = \tilde{\Phi}^{-1}$. We have
\[
(\tilde{\Psi}\circ\tilde{\Phi})(y + \mathcal{L}^*) = \tilde{\Psi}\big(\tilde{\Phi}(q(y))\big) = \tilde{\Psi}\big((\widehat{q}\circ\Phi)(y)\big) = (q\circ \Phi^{-1})\big(\Phi(y)\big) = q(y) = y + \mathcal{L}^*,
\]
so that $\tilde{\Psi}\circ\tilde{\Phi}$ is the identity map on $\mathbb{R}^n/\mathcal{L}^*$; and 
\[
(\tilde{\Phi}\circ\tilde{\Psi})(\chi_z\mathcal{L}^\perp) = \tilde{\Phi}\big(\tilde{\Psi}(\widehat{q}(\chi_z))\big) = \tilde{\Phi}\big((q\circ \Phi^{-1})(\chi_z)\big) = (\widehat{q}\circ\Phi)\big(\Phi^{-1}(\chi_z)\big) = \widehat{q}(\chi_z) = \chi_z\mathcal{L}^\perp,
\]
so that  $\tilde{\Phi}\circ\tilde{\Psi}$ is the identity map on $\widehat{\mathbb{R}^n}/\mathcal{L}^\perp$. Therefore, we have shown that $\tilde{\Phi} = \tilde{\Psi}^{-1}$ is a homeomorphism. Lastly, we note that $\tilde{\Phi} = X_2$: given $y \in \mathbb{R}^n$, we have
\[
\tilde{\Phi}(y + \mathcal{L}^*) = \tilde{\Phi}(q(y)) = \widehat{q}(\Phi(y)) = \widehat{q}(\chi_y) = \chi_y \mathcal{L}^\perp = X_2(y + \mathcal{L}^*).
\]
Therefore, $X_2 =\tilde{\Phi}$ is a homeomorphism.
	\end{proof}
	\begin{cor}
		The map $\mathcal{C}$ described above is a homeomorphism.
	\end{cor}
	\begin{proof} Since $\mathcal{C}(\nu+\mathcal{L}^*) := \chi_\nu^\mathrm{T}$, it is immediate that $\mathcal{C}$ is the composition $\vartheta^{-1}\circ X_3 \circ X_2$, as illustrated below.
			\[
		\begin{tikzcd} 
		\mathbb{R}^n/\mathcal{L}^* 	\arrow[to=1-4, bend right=30, "\mathcal{C}"{name=C}] \arrow[r, "X_2"{name=X2}]
			& \widehat{\mathbb{R}^n}/\mathcal{L}^\perp \arrow[r, "X_3"{name=X3}]
			& \widehat{\mathcal{L}} \arrow[r, "\vartheta^{-1}"{name=vartheta}]
			& \widehat{\mathrm{T}}\, ,
		\end{tikzcd}
		\]
		Proposition 1.86 in \cite{TaylorKaniuth} states that $X_3$ is a topological group isomorphism; so in particular, it is a homeomorphism. By Lemma \ref{dualiso}, we have that $\vartheta^{-1}$ is an isomorphism of topological groups; and we proved in the above theorem that $X_2$ is a homeomorphism. It follows that $\mathcal{C}=\vartheta^{-1}\circ X_3 \circ X_2$ is a homeomorphism.
	\end{proof}
Let \label{upsilon0}
\[
\boxed{\upsilon_0 = \mathcal{C} \circ q\vert_{\Pi R}}\,,
\]
so that $\upsilon_0 = \upsilon\vert_{\Pi R}$. Define \label{T0}
\[
\boxed{\widehat{\mathrm{T}}_0 = \upsilon(\Pi R)}\,;
\]
 then $\widehat{\mathrm{T}}_0 =\upsilon_0(\Pi R)$. Since, by Remark \ref{qpir}, $q(\Pi R)$ is open in $\mathbb{R}^n/\mathcal{L}^*$, it follows that $\widehat{\mathrm{T}}_0$ is open in $\widehat{\mathrm{T}}$. Since, by \ref{restrictedqBoreliso}, $q\vert_{\Pi R}$ is a Borel isomorphism onto its image, we have the following corollary:
	\begin{cor}
		$\upsilon_0:\Pi R \to \widehat{\mathrm{T}}_0$ is a Borel isomorphism.
	\end{cor} 

	\section[The \texorpdfstring{$\kappa$}{Kappa} intertwiner]{The \texorpdfstring{$\boldsymbol{\kappa}$}{Kappa} intertwiner}
	
	\subsection{Preliminaries}
	
	In this section we introduce an isometric isomorphism 
	\[
	\kappa: L^2\big(R;l^2(\Pi)\otimes l^2(\mathcal{L}^*)\big) \to L^2\big(\widehat{\mathrm{T}}; l^2(\mathcal{L}^*)\big),
	\]
which we later use to prove that $\tau$-invariant subspaces decompose as direct integrals of $U^r\otimes I$-invariant subspaces (Theorem \ref{MainDecompTheorem}). This map is composed from the following sequence of maps, which will be described in more detail below:
	\begin{equation}
		\begin{tikzcd}[column sep=small, row sep=small, cramped]       
			L^2\big(R;l^2(\Pi)\otimes l^2(\mathcal{L}^*)\big)\arrow[r,"\kappa_1"]
			& L^2(\Pi\times R\times\mathcal{L}^*) \arrow[r,"\kappa_2"]
			& L^2(\Pi R \times\mathcal{L}^*) \\
			\arrow[r,"\kappa_3"]
			&  L^2(\widehat{\mathrm{T}}_0\times\mathcal{L}^*) \arrow[r,"\kappa_4"]
			&L^2\big(\widehat{\mathrm{T}}_0; l^2(\mathcal{L}^*)\big)  \arrow[r,"\kappa_5"]	&L^2\big(\widehat{\mathrm{T}}; l^2(\mathcal{L}^*)\big) 
		\end{tikzcd}
	\end{equation}
	Before describing these maps, it will be helpful to recall from \eqref{ParamOfThat} that every element of $\widehat{\mathrm{T}}$ is of the form $\chi_\nu^{\mathrm{T}}$ for a unique $\nu\in \Omega_{\mathrm{T}^*}$. Here, $\Omega_{\mathrm{T}^*}$ is our chosen fundamental set for $\mathrm{T}^*$ such that $\Pi R \subset \Omega_{\mathrm{T}^*}\subset \overline{\Pi R} = \Pi \overline{R}$. In particular, we note that since $\Pi R$ is open, $\partial(\Pi R) = \overline{\Pi R} - \Pi R$. It will be useful in the sequel to note that $\Pi R$ is proper, meaning that its boundary has measure zero.
	\begin{lemma}\label{paramdomainminusPiRNull}
		Let $m$ denote the Lebesgue measure on $\mathbb{R}^n$. Then $m(\overline{\Pi R} - \Pi R) = 0$. In particular, $m(\Omega_{\mathrm{T}^*} - \Pi R) = 0$.
	\end{lemma}  
	\begin{proof}
		First, note that for any subset $X\subset \mathbb{R}^n$,
		\[
		\Pi X = \bigcup_{L\in \Pi} LX = \bigcup_{L\in \Pi} \{Lu: u\in X\}.
		\]
		Thus, if $u\in \Pi\overline{R}- \Pi R$, then $u\in L\overline{R}$ for some $L \in\Pi$, and $u\notin M R$ for any $M\in \Pi$. In particular, 
		\[
		u \in L\overline{R}- LR = L(\overline{R}-R) \subset \Pi(\overline{R}-R).
		\]
		Thus we have 
		\[
		\Pi\overline{R}- \Pi R\subset  \Pi(\overline{R}-R).
		\]
		But now, by Theorem 2.44(b) in \cite{FollandReal} and the fact that $R$ is proper, we have
		\[
		m\big(L(\overline{R}-R)\big) = |\det L|m(\overline{R}-R) = m(\overline{R}-R) = m(\partial R) = 0,
		\]
		for each $L\in \Pi$. Thus
		\[
		m\big(\Pi (\overline{R}-R)\big) \leq\sum_{L\in\Pi} m\big(L(\overline{R}-R)\big) = 0,
		\]
		and hence $m(\Pi\overline{R}- \Pi R) = 0$. Since $\overline{\Pi R} =\Pi\overline{R}$, it follows that 
		\[
		m\big( \partial(\Pi R)\big) = m(\overline{\Pi R} - \Pi R)=0.
		\]\end{proof}
	\begin{cor}\label{OmegaTstarmeas}
		Let $|\Pi|$\label{order} denote the (finite) order of $\Pi$. Then $m(\Omega_{\mathrm{T}^*}) = |\Pi|m(R)$.
	\end{cor}
	\begin{proof}
		Since $\Pi R \subset \Omega_{\mathrm{T}^*} \subset  \overline{\Pi R}$, we have
		\[
		m(\Pi R) \leq m(\Omega_{\mathrm{T}^*}) \leq m(\overline{\Pi R}).
		\]
		But $\Pi R$ is open, and hence $\overline{\Pi R} = \Pi R \,\dot\cup\, (\overline{\Pi R} - \Pi R)$. Then
		\[
		m(\overline{\Pi R}) = m(\Pi R) + m(\overline{\Pi R} - \Pi R) = m(\Pi R),
		\]
		whence $m(\Omega_{\mathrm{T}^*}) = m(\Pi R)$. But $\Pi R = \dot\cup_{L\in \Pi} LR$, and therefore
		\[
		m(\Omega_{\mathrm{T}^*}) = \sum_{L\in\Pi} m(LR) = \sum_{L\in\Pi} |\det L|m(R) = |\Pi|m(R).
		\]
	\end{proof}
	
\subsection[The \texorpdfstring{$\kappa_1$}{Kappa1} map]{The \texorpdfstring{$\boldsymbol{\kappa_1}$}{Kappa1} map}
	
	To define the $\kappa_1$ map, we first recall a basic fact about tensor products of $l^2$ spaces. 
	\begin{lemma}
		Let $A$ and $B$ be arbitrary sets. Given $x\in l^2(A)$ and $y\in l^2(B)$, we define $p_{x,y} \in l^2(A\times B)$ by
		\[
		p_{x,y}(a,b) = x(a)y(b) \ \ \ \text{for all}\ a\in A, b\in B.
		\]	
		There exists a unique unitary transformation $p$ from $l^2(A)\otimes l^2(B)$ onto $l^2(A\times B)$ such that
		\[
		p(x\otimes y) = p_{x,y} \ \ \ \ \text{for all simple tensors}\ x\otimes y \in l^2(A)\otimes l^2(B).
		\]
	\end{lemma}
	\begin{proof}
		This is Example 2.6.10 in \cite{KadRing}.
	\end{proof}
	
	In what follows we take $p$\label{p} to be the unique unitary transformation from $l^2(\Pi)\otimes l^2(\mathcal{L}^*)$ to $l^2(\Pi\times \mathcal{L}^*)$ guaranteed by the above lemma. Note that $p$ naturally gives rise to a unitary transformation\label{P}
	\[
	\mathrm{P}: L^2\big(R;l^2(\Pi)\otimes l^2(\mathcal{L}^*)\big) \to L^2\big(R;l^2(\Pi\times\mathcal{L}^*)\big),
	\]
	which is given by 
	\[
\boxed{	(\mathrm{P}F)(r) = p\big(F(r)\big)}\,.
	\]
	The inverse of this map is given by $[\mathrm{P}^{-1}\ell](r) = p^{-1}\big(\ell(r)\big)$. We now define the map \label{kappa1}
	\[
	\kappa_1: L^2\big(R;l^2(\Pi)\otimes l^2(\mathcal{L}^*)\big) \to L^2(\Pi\times R \times \mathcal{L}^*)
	\]
	by
	\[
\boxed{	(\kappa_1F)(M,r,z) = [p(F(r))](M,z)= [(\mathrm{P}F)(r)](M,z)}\,,
	\]
	for $F \in L^2\big(R;l^2(\Pi)\otimes l^2(\mathcal{L}^*)\big)$, $r \in R$, and $z\in \mathcal{L}^*$, and $M\in \Pi$. We show that $\kappa_1$ is an isometric isomorphism. To see that $\kappa_1$ is surjective, let $J \in L^2(\Pi\times R\times \mathcal{L}^*)$, and define $\ell \in L^2\big(R;l^2(\Pi\times \mathcal{L}^*)\big)$ by
	\[
	[\ell(r)](M,z) = J(M,r,z).
	\]
	Now define $F = \mathrm{P}^{-1}\ell$. Then it is immediate that $\kappa_1(F) = J$. We now show that $\kappa_1$ is an isometry. Give $R$ the relative Lebesgue $\sigma$-algebra inherited from $\mathbb{R}^n$, which is the completion of the Borel $\sigma$-algebra $\mathcal{B}_R$ with respect to the Lebesgue measure, by \ref{BorelClosure}. Then $\Pi\times R \times \mathcal{L}^*$ has the product $\sigma$-algebra $\mathcal{P}(\Pi) \times \mathscr{L}\vert_R \times \mathcal{P}(\mathcal{L}^*)$ and product measure $\mu_{\Pi}\times m \times \mu_{\mathcal{L}^*}$. We have
	\begin{align*}
		\|\kappa_1F\|^2 &= \int_{\Pi\times R\times\mathcal{L}^*} |(\kappa_1F)(M,r,z)|^2 \ud(\mu_{\Pi}\times m \times \mu_{\mathcal{L}^*})(M,r,z)\\
		&=\int_{\Pi\times R \times\mathcal{L}^*} |[(\mathrm{P}F)(r)](M,z)|^2 \ud(\mu_{\Pi}\times m \times \mu_{\mathcal{L}^*})(M,r,z)\\
		&= \int_R \int_{\Pi\times\mathcal{L}^*}  |[(\mathrm{P}F)(r)](M,z)|^2 \ud(\mu_{\Pi}\times \mu_{\mathcal{L}^*})(M,z) \ud r\\
		&= \int_R \|(\mathrm{P}F)(r)\|^2 \ud r\\
		&= \|\mathrm{P}F\|^2\\
		&=\|F\|^2.
	\end{align*}
	where we have used Tonelli's theorem in the third equality.
	
\subsection[The \texorpdfstring{$\kappa_2$}{Kappa2} map]{The \texorpdfstring{$\boldsymbol{\kappa_2}$}{Kappa2} map}
		
	Recall that $LR = \{Lr: r\in R\}$ for each $L \in \Pi$. Since $R$ is a fundamental region for $\Gamma^*$, it follows that $LR\hspace{.02cm} \cap\hspace{.02cm} MR = \emptyset$ for $L,M \in \Pi$ with $L \neq M$, and hence $\Pi R = \dot{\bigcup}_{L\in\Pi} LR$. We define a bijection\label{Lambda}
	\[
	\Lambda: \Pi\times R \to \Pi R,
	\]
	by $\Lambda(L,r) = L\cdot_{\idref{12}}r = Lr$. Since each element of $\Pi R$ is of the form $Lr$ for a unique $L\in \Pi$ and unique $r\in R$, we can define $\Lambda^{-1}(Lr) = (L,r)$. We note that $\Lambda$ is a homeomorphism when $\Pi R$ is given the relative topology inherited from the usual topology on $\mathbb{R}^n$, and $\Pi\times R$ has the product topology from the discrete topology on $\Pi$ and the relative Euclidean metric topology on $R$. To see this, suppose that $(L_n,r_n) \to (L,r)$, so that $L_n \to L$ and $r_n \to r$. Let $\epsilon>0$. Since $\Pi$ is discrete, there exists $N_1\in\mathbb{N}$ such that $n\geq N_1 \implies L_n = L$. Since $r_n \to r$ we can find $N_2\in\mathbb{N}$ such that $n\geq N_1 \implies |r_n-r|<\epsilon/\|L\|$. Thus for $n\geq\max(N_1,N_2)$,
	\[
	|L_nr_n - Lr| = |L(r_n-r)| \leq \|L\||r_n-r|<\epsilon.
	\]
	Thus $\Lambda$ is continuous. To show that $\Lambda^{-1}$ is continuous, it suffices (Proposition 4.11 in \cite{FollandReal}) to show that $\pi_1\circ \Lambda^{-1}$ and $\pi_2\circ \Lambda^{-1}$ are continuous, where $\pi_1$ and $\pi_2$ are the canonical projection maps. Here
	\[
	\pi_1\circ \Lambda^{-1}(Lr) = \pi_1(L,r) = L,
	\]
	and 
	\[
	\pi_2\circ \Lambda^{-1}(Lr) = \pi_2(L,r) = r.
	\]
	If $\mathcal{U}\subset \Pi$, then 
	\[
	(\pi_1\circ \Lambda^{-1})^{-1}(\mathcal{U}) =\bigcup_{L\in\mathcal{U}}LR=: \mathcal{U}R,
	\]
	which is open in $\Pi R$, since $\mathcal{U} R$ is open in $\mathbb{R}^n$. Likewise, if $\mathcal{V}\subset R$ is open, then 
	\[
	(\pi_2\circ \Lambda^{-1})^{-1}(\mathcal{V}) = \bigcup_{L\in\Pi} L\mathcal{V} =:\Pi\mathcal{V},
	\]
	which is open in $\Pi R$, since $L\mathcal{V}$ is open in $\mathbb{R}^n$. Thus we have shown that $\pi_1\circ \Lambda^{-1}$ and $\pi_2\circ \Lambda^{-1}$ are continuous, and hence $\Lambda^{-1}$ is continuous. The continuity of $\Lambda^{-1}$ can also be seen by the observing that if $L_nr_n \to Lr$, where $r\in R$, then by the mutual disjointness of the open sets $\{LR:L\in \Pi\}$, it follows that $L_n = L$ for $n$ sufficiently large, and $r_n \to r$. 
	
	We now define $\widetilde{\Lambda}: L^2(\Pi R\times \mathcal{L}^*)\to  L^2(\Pi \times R\times \mathcal{L}^*)$ by \label{tildeLambda}
	\[
	\widetilde{\Lambda}(\eta)(L,r,z) = \eta\big(\Lambda(L,r),z\big) = \eta(Lr,z).
	\]
	Thus $\widetilde{\Lambda}(\eta) = \eta\circ (\Lambda\times I)$. By definition, $\widetilde{\Lambda}$ is a linear transformation. It is also an isometric isomorphism, as we argue presently. Let $\mu_{\mathcal{L}^*}$ denote the counting measure on $\mathcal{L}^*$, and let $\mu_{\Pi}$ denote the counting measure on $\Pi$. The spaces $\Pi R \times \mathcal{L}^*$ and $\Pi \times R\times \mathcal{L}^*$ have Borel $\sigma$-algebras. Note that $\mathcal{B}_{\Pi R \times\mathcal{L}^*}$ coincides with the product $\sigma$-algebra $\mathcal{B}_{\Pi R} \times \mathcal{B}_{\mathcal{L}^*}$, by Proposition \ref{SepBorelAlg}. Likewise, $\mathcal{B}_{\Pi \times R\times \mathcal{L}^*} = \mathcal{B}_\Pi \times \mathcal{B}_R \times \mathcal{B}_{\mathcal{L}^*}$. Taking the completions of these $\sigma$-algebras, we have by Corollary \ref{completionofproduct} that
	\begin{equation}\label{BorelClosure1}
		\begin{aligned}
			\overline{\mathcal{B}_{\Pi R \times\mathcal{L}^*}}&=\overline{\mathcal{B}_{\Pi R} \times \mathcal{B}_{\mathcal{L}^*}}\\
			&=\overline{\mathcal{B}_{\Pi R} \times \mathcal{P}(\mathcal{L}^*)}\\
			&=  \overline{\mathcal{B}_{\Pi R}}  \times \mathcal{P}(\mathcal{L}^*)\\
			&=\mathscr{L}\vert_{\Pi R}  \times \mathcal{P}(\mathcal{L}^*),
		\end{aligned}
	\end{equation}
	and likewise
	\begin{equation}\label{BorelClosure2}
		\begin{aligned}
			\overline{\mathcal{B}_{\Pi \times R\times \mathcal{L}^*}} &= \overline{\mathcal{P}(\Pi) \times \mathcal{B}_R \times \mathcal{P}(\mathcal{L}^*) }\\
			&= \mathcal{P}(\Pi) \times \overline{\mathcal{B}_R} \times \mathcal{P}(\mathcal{L}^*) \\
			&=  \mathcal{P}(\Pi) \times \mathscr{L}\vert_R \times \mathcal{P}(\mathcal{L}^*) .
		\end{aligned}
	\end{equation}
	Since $\Lambda\times I$ is a homeomorphism it is, \textit{a fortiori}, a Borel isomorphism, that is, a Borel bijection whose inverse is Borel also. Thus, if we show that the product measure $m\times \mu_{\mathcal{L}^*}$ coincides with the pushfoward measure $(\mu_{\Pi}\times m \times \mu_{\mathcal{L}^*})\circ (\Lambda\times I)^{-1}$ on Borel sets, then it follows from Corollary \ref{imagemeasExcor} that $\widetilde{\Lambda}$ is an isometric isomorphism between the $L^2$ spaces with the completed measures.
	\begin{lemma}
		Let $\Pi R$ have the relative Lebesgue $\sigma$-algebra inherited from $\mathbb{R}^n$, and Lebesgue measure $m$. Then 
		\[
		m = (\mu_\Pi \times m)\circ \Lambda^{-1}.
		\]
	\end{lemma}
	\begin{proof}
		First note that for any Borel set $U \in  \mathcal{B}_{\Pi R}$, since $L R \cap M R = \emptyset$ whenever $L\neq M$, we have 
		\begin{equation}\label{breakupU}
			U = \dot\bigcup_{L\in \Pi} LR\cap U,
		\end{equation}
		with $LR\cap U\in \mathcal{B}_{\Pi R}$ for each $L\in \Pi$. Hence
		\[
		m(U) = \sum_{L\in\Pi} m(LR\cap U).
		\] 
		Now note that
		\[
		\Lambda^{-1} \big( \dot\bigcup_{L\in \Pi} LR\cap U\big) = \dot{\bigcup}_{L\in \Pi}\{L\}\times L^{-1}(LR\cap U),
		\]
		and hence
		\begin{align*}
			\big((\mu_\Pi \times m)\circ \Lambda^{-1}\big)(U) &= (\mu_{\Pi}\times m)\Big(\dot{\bigcup}_{L\in \Pi}\{L\}\times L^{-1}(LR\cap  U)\Big)\\
			&=\sum_{L\in\Pi} m\big(\big(L^{-1}(LR\cap U)\big)\\
			&= \sum_{L\in\Pi} |\det{L^{-1}}|m(LR\cap U)\\
			&=\sum_{L\in\Pi} m(LR\cap U)\\
			&=m(U).
		\end{align*}
		Since the measures agree on $\mathcal{B}_{\Pi R}$, it follows by Proposition \ref{BorelClosure} and the uniqueness of the complete measure (Theorem 1.9 in \cite{FollandReal}) that they agree on the relative Lebesgue $\sigma$-algebra $\mathscr{L}\vert_{\Pi R}$.
	\end{proof}
	\begin{cor}\label{prodagreeBorel}
		The product measure $m\times \mu_{\mathcal{L}^*}$ coincides with the pushfoward measure $(\mu_{\Pi}\times m \times \mu_{\mathcal{L}^*})\circ (\Lambda\times I)^{-1}$ on Borel sets.
	\end{cor}
	\begin{proof}
		Let $E\in\mathcal{B}_{\Pi R \times\mathcal{L}^*}=  \mathcal{B}_{\Pi R}\times \mathcal{P}(\mathcal{L}^*)$; then $E=U\times K$, where $U$ is Borel in $\Pi R$, and $F \subset \mathcal{L}^*$. In particular, $E$ is a measurable rectangle, and 
		\[
		(m\times\mu_{\mathcal{L}^*})(E) = m(U)\mu_{\mathcal{L}^*}(K).
		\]
		We have
		\begin{align*}
			&\big((\mu_{\Pi}\times m \times \mu_{\mathcal{L}^*})\circ (\Lambda\times I)^{-1}\big)(E) \\
			&= (\mu_{\Pi}\times m \times \mu_{\mathcal{L}^*})\big((\Lambda\times I)^{-1}(U\times K)\big)\\
			&= (\mu_{\Pi}\times m \times \mu_{\mathcal{L}^*})\big(\Lambda^{-1}(U) \times K \big)\\
			&=(\mu_{\Pi}\times m)\big(\Lambda^{-1}(U)\big)\mu_{\mathcal{L}^*}(K) \\
			&= m(U)\mu_{\mathcal{L}^*}(K) \\
			&= (m\times\mu_{\mathcal{L}^*})(E).
		\end{align*}
	\end{proof}
	\begin{cor}
		The map
		\begin{multline*}
			\widetilde{\Lambda}: 	L^2\big(\Pi R\times \mathcal{L}^*,\mathscr{L}\vert_{\Pi R}  \times \mathcal{P}(\mathcal{L}^*), m\times \mu_{\mathcal{L}^*} \big)\\
			\to  L^2(\Pi \times R\times \mathcal{L}^*,\mathcal{P}(\Pi) \times \mathscr{L}\vert_R \times \mathcal{P}(\mathcal{L}^*) ,\mu_{\Pi}\times m \times \mu_{\mathcal{L}^*}).
		\end{multline*}
		is a unitary operator.
	\end{cor}
	\begin{proof}
		As mentioned above, $\widetilde{\Lambda}$ is linear, and it follows from Corollary \ref{imagemeasExcor} that 
		\begin{multline*}
			\widetilde{\Lambda}:  L^2\big(\Pi R\times \mathcal{L}^*, \overline{\mathcal{B}_{\Pi R \times\mathcal{L}^*}},\overline{(\mu_{\Pi}\times m \times \mu_{\mathcal{L}^*})\circ (\Lambda\times I)^{-1} }\big)\\
			\to  L^2(\Pi \times R\times \mathcal{L}^*,\overline{\mathcal{B}_{\Pi \times R\times \mathcal{L}^*}},\mu_{\Pi}\times m \times \mu_{\mathcal{L}^*}).
		\end{multline*}
		is an isometric isomorphism. Note that $\overline{m\times \mu_{\mathcal{L}^*}} = m\times \mu_{\mathcal{L}^*}$ by \ref{completionofproduct}. Now by Corollary \ref{prodagreeBorel} above, together with \eqref{BorelClosure1} and \eqref{BorelClosure2}, we get that
		\begin{multline*}
			\widetilde{\Lambda}: 	L^2\big(\Pi R\times \mathcal{L}^*,\mathscr{L}\vert_{\Pi R}  \times \mathcal{P}(\mathcal{L}^*), m\times \mu_{\mathcal{L}^*} \big)\\
			\to  L^2(\Pi \times R\times \mathcal{L}^*,\mathcal{P}(\Pi) \times \mathscr{L}\vert_R \times \mathcal{P}(\mathcal{L}^*) ,\mu_{\Pi}\times m \times \mu_{\mathcal{L}^*}).
		\end{multline*}
		is an isometric isomorphism.
	\end{proof}
	We define\label{kappa2}
	\[
	\kappa_2: L^2(\Pi\times R\times \mathcal{L}^*) \to L^2(\Pi R\times \mathcal{L}^*)
	\]
	by
	\[
	\boxed{\kappa_2 := \widetilde{\Lambda}^{-1} }\,,
	\]
so that
	\begin{equation}\label{kappa2quality}
		(\kappa_2J)(Mr, z) = \big(J\circ (\Lambda^{-1}\times I)\big)(Mr,z) = J(M,r,z) 
	\end{equation}
	for all $(Mr,z) \in \Pi R\times\mathcal{L}^*$.
	
\subsection[The \texorpdfstring{$\kappa_3$}{Kappa3} map]{The \texorpdfstring{$\boldsymbol{\kappa_3}$}{Kappa3} map}
	
	It follows from the Appendix section \ref{MOUHaar} that the image measure $m\circ \upsilon^{-1}$ is a Haar measure on $\widehat{\mathrm{T}}$. Furthermore, it is shown in that section that if $\mu_{\mathrm{T}}$ is the counting measure on $\mathrm{T}$, so that its dual measure $\mu_{\widehat{\mathrm{T}}}$ satisfies $\mu_{\widehat{\mathrm{T}}}(\widehat{\mathrm{T}}) =1$, then 
	\begin{equation}\label{measurescaling}
		\overline{m\circ \upsilon^{-1}} = |\Pi|m(R) \mu_{\widehat{\mathrm{T}}},
	\end{equation}
where the overline denotes the completed measure. From this, we can show that the image of $\Pi R$ under $\upsilon$ is $\mu_{\widehat{\mathrm{T}}}$-almost all of $\widehat{\mathrm{T}}$:
	\begin{lemma}
As before, let $\widehat{\mathrm{T}}_0 = \upsilon(\Pi R)$. Then 
		\[
		\mu_{\widehat{\mathrm{T}}}(\widehat{\mathrm{T}} - \widehat{\mathrm{T}}_0) =0. 
		\]
	\end{lemma}
	\begin{proof}
		Since $\upsilon(\Omega_{\mathrm{T}^*}) = \widehat{\mathrm{T}}$, we have
		\begin{align*}
			m(\Omega_{\mathrm{T}^*}) &= |\Pi|m(R) \mu_{\widehat{\mathrm{T}}}(\widehat{\mathrm{T}})\\
			&= |\Pi|m(R)\big[\mu_{\widehat{\mathrm{T}}}(\widehat{\mathrm{T}}_0)+ \mu_{\widehat{\mathrm{T}}}(\widehat{\mathrm{T}}-\widehat{\mathrm{T}}_0)\big]\\
			&=  |\Pi|m(R)\Big[\frac{1}{|\Pi|m(R)} (m\circ \upsilon^{-1})(\widehat{\mathrm{T}}_0) + \mu_{\widehat{\mathrm{T}}}(\widehat{\mathrm{T}}-\widehat{\mathrm{T}}_0)\Big]\\
			&= m(\Pi R) +  |\Pi|m(R)\mu_{\widehat{\mathrm{T}}}(\widehat{\mathrm{T}}-\widehat{\mathrm{T}}_0).
		\end{align*}
		Thus 
		\begin{equation*}
			\mu_{\widehat{\mathrm{T}}}(\widehat{\mathrm{T}}-\widehat{\mathrm{T}}_0) = \frac{1}{|\Pi|m(R)}\big[m(\Omega_{\mathrm{T}^*}) - m(\Pi R)\big] = 0,
		\end{equation*}
		since $m(\Omega_{\mathrm{T}^*}) = m(\Pi R)$, as shown in the proof of Corollary \ref{OmegaTstarmeas}.
	\end{proof}
	\begin{rem}\label{T0form}
	Recall that $\upsilon_0 = \upsilon\vert_{\Pi R}$. Each element of $\widehat{\mathrm{T}}_0$ is of the form $\chi_{Mr}^{\mathrm{T}}$ for some $M\in\Pi$ and $r\in R$. For if $\chi \in \widehat{\mathrm{T}}_0$, then $\chi = \upsilon_0(M r)$ for some $M\in\Pi$ and $r\in R$, and hence
		\[
		\chi = \upsilon_0(Mr) =  \chi_{Mr}^{\mathrm{T}}.
		\]
		Moreover, the $M \in \Pi$ and $r\in R$ for which $\chi =  \chi_{Mr}^{\mathrm{T}}$ are uniquely determined, since $\Pi R\cong \Pi \times R$.
	\end{rem}
	Note that
	\[
	\upsilon_0 \times I: \Pi R \times \mathcal{L}^* \to \widehat{\mathrm{T}}_0\times \mathcal{L}^* 
	\]
is a Borel isomorphism. We define a linear operator\label{upsilon0I}
	\begin{multline*}
	\widetilde{\upsilon_0\times I}: L^2\big(\widehat{\mathrm{T}}_0\times \mathcal{L}^*, \overline{\mathcal{B}_{\widehat{\mathrm{T}}_0} \times\mathcal{B}_{\mathcal{L}^*}}, \overline{(m\times\mu_{\mathcal{L}^*})\circ(\upsilon_0\times I)^{-1}}\big) \\
	\to  L^2(\Pi R \times\mathcal{L}^*, \overline{\mathcal{B}_{\Pi R} \times\mathcal{L}^*}, m\times\mu_{\mathcal{L}^*})
	\end{multline*}
	by 
	\[
\boxed{[(\widetilde{\upsilon_0\times I})\phi](\nu,z) = \phi\big(\upsilon_0(\nu),z\big)}\,.
	\]
	That is, $(\widetilde{\upsilon_0\times I})\phi = \phi\circ (\upsilon_0 \times I)$. By Corollary \ref{imagemeasExcor}, $\widetilde{\upsilon_0\times I}$ is an isometric isomorphism.  We note that 
	\[
	(m\times\mu_{\mathcal{L}^*})\circ(\upsilon_0\times I)^{-1} = (m\circ \upsilon_0^{-1})\times \mu_{\mathcal{L}^*},
	\]
	so that by Corollaries \ref{completionofproduct} and \ref{mofupsiloninverse1},
	\begin{align*}
	\overline{(m\times\mu_{\mathcal{L}^*})\circ(\upsilon_0\times I)^{-1}} &=  \overline{(m\circ \upsilon_0^{-1})\times \mu_{\mathcal{L}^*}}\\
	&=  \overline{(m\circ \upsilon_0^{-1})}\times \mu_{\mathcal{L}^*}\\
	&= |\Pi|m(R)\mu_{\widehat{\mathrm{T}}}\times \mu_{\mathcal{L}^*}.
	\end{align*}
 Moreover, it follows from Corollary \ref{completionofproduct} that
	\[
	\overline{\mathcal{B}_{\widehat{\mathrm{T}}_0} \times\mathcal{B}_{\mathcal{L}^*}} = \overline{\mathcal{B}_{\widehat{\mathrm{T}}_0}} \times \mathcal{B}_{\mathcal{L}^*} = \overline{\mathcal{B}_{\widehat{\mathrm{T}}_0}} \times \mathcal{P}(\mathcal{L}^*)
	\]
	and
	\[
	\overline{\mathcal{B}_{\Pi R\times\mathcal{L}^*}} = 	\overline{\mathcal{B}_{\Pi R}} \times \mathcal{B}_{\mathcal{L}^*}
	= \mathscr{L}\vert_{\Pi R}\times\mathcal{P}(\mathcal{L}^*).
	\]
	The following definition will streamline notation:
	\begin{tcolorbox}
		\begin{defn}\label{tildemuthat}
			We define
			\[
			\tilde\mu_{\widehat{\mathrm{T}}} = |\Pi|m(R)\mu_{\widehat{\mathrm{T}}}.
			\]
		\end{defn}
	\end{tcolorbox}
	In particular, we have 
	\[
	\boxed{\overline{m\circ \upsilon^{-1}} =\overline{m\circ \upsilon_0^{-1}} = \tilde\mu_{\widehat{\mathrm{T}}}}\,.
	\]
 In what follows, we shall assume $\tilde\mu_{\widehat{\mathrm{T}}}$ to be the measure on $\widehat{\mathrm{T}}$ and $\widehat{\mathrm{T}}_0$. It follows from the above identifications that
	\begin{multline*}
		\widetilde{\upsilon_0\times I}: L^2\big(\widehat{\mathrm{T}}_0\times \mathcal{L}^*, \overline{\mathcal{B}_{\widehat{\mathrm{T}}_0}} \times \mathcal{P}(\mathcal{L}^*), \tilde\mu_{\widehat{\mathrm{T}}}\times \mu_{\mathcal{L}^*}\big) \\ \to L^2(\Pi R \times\mathcal{L}^*, \mathscr{L}\vert_{\Pi R}\times\mathcal{P}(\mathcal{L}^*), m\times\mu_{\mathcal{L}^*})
	\end{multline*}
	is a unitary operator. We define \label{kappa3}
	\[
\boxed{	\kappa_3:= (\widetilde{\upsilon_0\times I})^{-1}}\,;
	\]
	that is, $\kappa_3\psi = \psi\circ (\upsilon_0\times I)^{-1}$. In particular, for $\chi_{Mr}^{\mathrm{T}} \in \widehat{\mathrm{T}}_0$, we have
	\[
	(\kappa_3\psi)(\chi_{Mr}^{\mathrm{T}},z) = \psi\big(Mr,z\big).
	\]
	
\subsection[The \texorpdfstring{$\kappa_4$}{Kappa4} map]{The \texorpdfstring{$\boldsymbol{\kappa_4}$}{Kappa4} map}
	We now define $\kappa_4$. By Theorem II.10(c) in \cite{ReedSimon} (p.~52), if $(\mathcal{M}_1,\mu_1)$ and $(\mathcal{M}_2,\mu_2)$ are measure spaces so that $L^2(\mathcal{M}_1,\mu_1)$ and $L^2(\mathcal{M}_2, \mu_2)$ are separable, then there is a unique isomorphism $W$ from $L^2(\mathcal{M}_1\times \mathcal{M}_2, \mu_1 \times \mu_2)$ to $L^2\big(\mathcal{M}_1,\mu_1; L^2(\mathcal{M}_2, \mu_2)\big)$ such that $[Wf(x)](y) = f(x,y)$ for $f\in L^2(\mathcal{M}_1\times \mathcal{M}_2, \mu_1 \times \mu_2)$, where $(x,y) \in \mathcal{M}_1\times \mathcal{M}_2$. Applying this theorem gives the unique isomorphism\label{kappa4}
	\[
	\kappa_4: L^2(\widehat{\mathrm{T}}_0\times\mathcal{L}^*, \tilde\mu_{\widehat{\mathrm{T}}}\times \mu_{\mathcal{L}^*}) \to
	L^2\big(\widehat{\mathrm{T}}_0, \tilde\mu_{\widehat{\mathrm{T}}}; l^2(\mathcal{L}^*)\big)
	\]
	satisfying 
	\[
\boxed{	[\kappa_4 \phi(\chi)](z) =\phi(\chi,z)}
	\]
	 for $\phi\in  L^2(\widehat{\mathrm{T}}_0\times\mathcal{L}^*)$, $\chi\in \widehat{\mathrm{T}}_0$, and $z\in \mathcal{L}^*$.
	
\subsection[The \texorpdfstring{$\kappa_5$}{Kappa5} map]{The \texorpdfstring{$\boldsymbol{\kappa_5}$}{Kappa5} map}
	
Since $\tilde\mu_{\widehat{\mathrm{T}}}(\widehat{\mathrm{T}}-\widehat{\mathrm{T}}_0) =0$, it follows from Remark \ref{L2aevector} that we can embed $L^2\big(\widehat{\mathrm{T}}_0; l^2(\mathcal{L}^*)\big)$ isomorphically in $L^2\big(\widehat{\mathrm{T}}; l^2(\mathcal{L}^*)\big)$. We define $\kappa_5$ to be this embedding. Explicitly, for $F \in L^2\big(\widehat{\mathrm{T}}_0; l^2(\mathcal{L}^*)\big)$, we define 
	\begin{align*}
		\tilde{F}(\chi) = \left\{ \begin{array}{ll}
			F(\chi) &\text{if } \chi \in \widehat{\mathrm{T}}_0 \\
			0 &\text{if } \chi\in \widehat{\mathrm{T}}-\widehat{\mathrm{T}}_0
		\end{array} \right. ;
	\end{align*}
and we define \label{kappa5}
\[
\kappa_5:	L^2\big(\widehat{\mathrm{T}}_0; l^2(\mathcal{L}^*)\big) \to 	L^2\big(\widehat{\mathrm{T}}; l^2(\mathcal{L}^*)\big)
\] 
by
\[
\kappa_5([F]_{L^2(\widehat{\mathrm{T}}_0; l^2(\mathcal{L}^*))} = [\tilde{F}]_{L^2(\widehat{\mathrm{T}}; l^2(\mathcal{L}^*))}.
\] 

\section[The \texorpdfstring{$\mathrm{M}$}{M} Representation]{The \texorpdfstring{$\boldsymbol{\mathrm{M}}$}{M} Representation}
	
	We define $\kappa: L^2\big({R};l^2(\Pi)\otimes l^2(\mathcal{L}^*)\big) \to L^2\big(\widehat{\mathrm{T}}, \tilde\mu_{\widehat{\mathrm{T}}}; l^2(\mathcal{L}^*)\big)$ to be the composition \label{kappadef}
	\[
\boxed{	\kappa = \kappa_5\circ \kappa_4\circ\kappa_3\circ \kappa_2 \circ\kappa_1}\,.
	\]
	We note that $\kappa$ is an isometric isomorphism, being a composition of such maps. Moreover, it follows from the definitions that, for $\chi_{Mr}^{\mathrm{T}} \in \widehat{\mathrm{T}}_0$, 
	\[
	\boxed{[(\kappa F)(\chi_{Mr}^{\mathrm{T}})](z) = [(\mathrm{P}F)(r)](M,z)}\,.
	\]
	Indeed, given $\chi_{Mr}^{\mathrm{T}} \in \widehat{\mathrm{T}}_0$ we have
	\begin{align*}
		[(\kappa F)(\chi_{Mr}^{\mathrm{T}})](z)  &= [\kappa_5((\kappa_4\circ \kappa_3\circ\kappa_2\circ \kappa_1)F)(\chi_{Mr}^{\mathrm{T}})](z)\\
		&= [\kappa_4((\kappa_3\circ\kappa_2\circ \kappa_1)F)(\chi_{Mr}^{\mathrm{T}})](z)\\
		&= [\kappa_3((\kappa_2\circ \kappa_1)F)](\chi_{Mr}^{\mathrm{T}},z)\\
		&= ( (\kappa_2\circ \kappa_1)F)(Mr,z)\\
		&= (\kappa_2(\kappa_1F))(Mr,z)\\
		&= (\kappa_1F)(M,r,z)\\
		&= [(\mathrm{P}F)(r)](M,z)
	\end{align*}
	Recall from Theorem \ref{PontDual} that for each $[y,\mathrm{id}]\in \mathrm{T}$, we can define a character $E_{[y,\mathrm{id}]}$ on $\widehat{\mathrm{T}}$ by $E_{[y,\mathrm{id}]}(\chi) = \chi([y,\mathrm{id}])$\label{evalcharonThat}, and moreover the map $[y,\mathrm{id}] \mapsto E_{[y,\mathrm{id}]}$ gives a topological isomorphism $\mathrm{T} \cong \widehat{\widehat{\mathrm{T}}}$. In particular, 
	\begin{equation}\label{That}
		\widehat{\widehat{\mathrm{T}}} = \{E_{[y,\mathrm{id}]}: [y,\mathrm{id}] \in \mathrm{T}\}.
	\end{equation}
	\begin{defn}
		For each $[y,\mathrm{id}] \in \mathrm{T}$, we define the operator $\mathrm{M}_{[y,\mathrm{id}]}\in\mathscr{B}\big( L^2\big(\widehat{\mathrm{T}}; l^2(\mathcal{L}^*)\big)\big)$ by 
		\[
	\boxed{	[(\mathrm{M}_{[y,\mathrm{id}]}\upphi)(\chi)](z) = E_{[y,\mathrm{id}]}(\chi)[\upphi(\chi)](z) \ \ \ \text{for}\ \chi \in \widehat{\mathrm{T}}, z\in\mathcal{L}^*}\,.
		\]
	\end{defn}
	Thus $(\mathrm{M}_{[y,\mathrm{id}]}\upphi)(\chi) = \chi([y,\mathrm{id}])\upphi(\chi)$. By definition, $\mathrm{M}_{[y,\mathrm{id}]}$ is linear, and since $\chi([y,\mathrm{id}])\in\mathbb{T}$ for all $\chi \in \widehat{\mathrm{T}}$, it follows immediately that $\mathrm{M}_{[y,\mathrm{id}]}$ is an isometry. Moreover, $\mathrm{M}_{[y,\mathrm{id}]}$ is clearly surjective, for given $\upphi \in L^2\big(\widehat{\mathrm{T}}; l^2(\mathcal{L}^*)\big)$, we have 
	\[
	\upphi = \mathrm{M}_{[y,\mathrm{id}]}\big(\mathrm{M}_{[-y,\mathrm{id}]}\upphi\big).
	\]
	Hence $\mathrm{M}_{[y,\mathrm{id}]}$ is a unitary operator on $L^2\big(\widehat{\mathrm{T}}; l^2(\mathcal{L}^*)\big)$. We define a unitary representation of $\mathrm{T}$ on $L^2\big(\widehat{\mathrm{T}}; l^2(\mathcal{L}^*)\big)$ by \label{Mrep}
	\[
	\mathrm{M}([y,\mathrm{id}]) = \mathrm{M}_{[y,\mathrm{id}]}.
	\]
	Indeed, it only remains to check that this is a group homomorphism, which follows from the fact that $E_{[y,\mathrm{id}]}$ is a character on $\widehat{\mathrm{T}}$. We denote this representation by $\mathrm{M}$.
	
We now show that $\kappa$ intertwines $\mathrm{M}$ with $\tau\vert_\mathrm{T}$, where $\tau[x,L]= \int_R^{\oplus} U^r[x,L]\otimes I \ud r $. Let $F\in L^2\big(R;l^2(\Pi)\otimes l^2(\mathcal{L}^*)\big)$. To this end, we define, for each $[z,M]\in \Gamma^*$,
\begin{equation}\label{Szmdef}
	\boxed{	\mathscr{S}_{[z,M]} =\{F\in L^2\big({R};l^2(\Pi)\otimes l^2(\mathcal{L}^*)\big): F( r) \in \langle \delta_M \otimes \delta_z\rangle\ \text{for a.e.}\  r \in {R}\}}\,,
\end{equation}
where $\langle \delta_M \otimes \delta_z\rangle = \{k(\delta_M \otimes \delta_z): k\in \mathbb{C}\}$. That $\mathscr{S}_{[z,M]}$ is a closed subspace follows from \ref{Mjsubspace} or by checking directly. Moreover, these spaces are clearly mutually orthogonal. We show that
\begin{equation}\label{directSumSzm}
	L^2\big({R};l^2(\Pi)\otimes l^2(\mathcal{L}^*)\big) = \sideset{}{\oplus}\sum_{[z,M]\in\Gamma^*} \mathscr{S}_{[z,M]}
\end{equation}
It follows from Theorem II.10 on p.~52 of \cite{ReedSimon} that there is a unique unitary transformation\label{udef0}
\[
u: L^2({R}) \otimes l^2(\Pi)\otimes l^2(\mathcal{L}^*) \to L^2\big({R};l^2(\Pi)\otimes l^2(\mathcal{L}^*)\big)
\]
satisfying 
\begin{equation}\label{udeF}
	u(a\otimes \zeta\otimes \beta) = u_{a,\zeta\otimes \beta},
\end{equation}
where \label{udef2}
\[
\boxed{	u_{a,g}( r) := a( r)g}
\]
for each $g\in l^2(\Pi)\otimes l^2(\mathcal{L}^*)$. We define, for $[z,M]\in \Gamma^*$,
\begin{equation}\label{Lzmdef}
	\boxed{	\mathscr{L}_{[z,M]} = \{u_{a,\delta_M\otimes \delta_z}: a \in L^2({R})\}\subset \mathscr{S}_{[z,M]}}\,.
\end{equation}
Thus,
\[
\mathscr{L}_{[z,M]}  = u\big(L^2({R})\otimes \{\delta_M\otimes \delta_z\} \big).
\]
We also define, for $[w,N]\in \Gamma^*$,
\[
I_{[w,N]}: L^2\big({R};l^2(\Pi)\otimes l^2(\mathcal{L}^*)\big) \to L^2({R})
\] 
by 
\begin{equation}\label{Izm}
	\boxed{	\big(I_{[w,N]}(F)\big)( r) = \langle F( r),\delta_N\otimes \delta_w\rangle}\,.
\end{equation}
That $I_{[w,N]} \in L^2({R})$ follows from the Cauchy-Schwarz inequality:
\[
\|I_{[w,N]}(F)\|^2 = \int_R |\langle F( r),\delta_N\otimes \delta_w\rangle|^2\ud  r \leq \int_{R} \|F(r)\|^2\|\delta_N\otimes \delta_w\|^2 \ud  r = \|F\|^2 <\infty.
\]
This also shows that $I_{[w,N]}$ is bounded. Since 
\[
\{\delta_M\otimes \delta_z: M\in \Pi, z\in\mathcal{L}^*\}
\] 
is an orthonormal basis for $l^2(\Pi)\otimes l^2(\mathcal{L}^*)$, it follows that, for an everywhere-defined representative $F\in L^2\big({R};l^2(\Pi)\otimes l^2(\mathcal{L}^*)\big)$,
\begin{equation}\label{Xipointwisedecomp}
	F( r)= \sum_{[z,M]\in \Gamma^*} u_{ I_{[z,M]}(F),\delta_M\otimes\delta_z}( r),
\end{equation}
for each $ r\in {R}$. From this it follows that
\begin{equation}\label{Szmint}
	\mathscr{S}_{[z,M]} = \bigcap\{\ker I_{[w,N]}: [w,N]\neq [z,M]\},
\end{equation}
from which we deduce that $\mathscr{S}_{[z,M]} = \mathscr{L}_{[z,M]}$ for all $[z,M]\in \Gamma^*$. Indeed, if $u_{a,\delta_M\otimes \delta_z} \in \mathscr{L}_{[z,M]}$, then
\[
\big(I_{[w,N]}(u_{a,\delta_M\otimes \delta_z})\big)(r) = \langle u_{a,\delta_M\otimes \delta_z}(r),\delta_N\otimes \delta_w\rangle = \langle a(r)\delta_M\otimes \delta_z ,\delta_N\otimes \delta_w\rangle = 0
\]
for a.e. $r\in R$, whenever $[w,N]\neq [z,M]$. Thus $u_{a,\delta_M\otimes \delta_z} \in \ker I_{[w,N]}$ whenever $[w,N]\neq [z,M]$, whence $u_{a,\delta_M\otimes \delta_z}  \in \mathscr{S}_{[z,M]}$.

 We show that the sum \eqref{Xipointwisedecomp} converges in the $L^2$ sense. By the Pythagorean theorem, we have
\[
\|F( r)\|^2 = \sum_{[z,M]\in \Gamma^*} \| u_{ I_{[z,M]}(F),\delta_M\otimes\delta_z}( r)\|^2,
\]
and using Tonelli's theorem we get
\begin{equation}\label{L2summableXi}
	\|F\|^2 =  \sum_{[z,M]\in \Gamma^*} \|u_{ I_{[z,M]}(F),\delta_M\otimes\delta_z}\|^2.
\end{equation}
It follows that the orthogonal family of vectors $\{u_{ I_{[z,M]}(F),\delta_M\otimes\delta_z}\}_{[z,M]\in \Gamma^*}$ is $l^2$-summable, and therefore converges to $F$ in the $L^2$ sense. Thus \eqref{directSumSzm} holds, and for each $F\in L^2\big({R};l^2(\Pi)\otimes l^2(\mathcal{L}^*)\big)$, we have
	\begin{equation}\label{Fdecomp}
		\boxed{F = \sideset{}{\oplus}\sum_{[z,M]\in\Gamma^*}F_{[z,M]}}\,,
	\end{equation}
	where
	\begin{equation}
		\boxed{	F_{[z,M]} =u_{ I_{[z,M]}(F),\delta_M\otimes\delta_z} \in \mathscr{S}_{[z,M]}}\,.
	\end{equation}
 For $N\in \Pi$, $r_0\in R$, and $w\in\mathcal{L}^*$, we have
\begin{align*}
		&[\big((\kappa\circ \tau[y,\mathrm{id}])F_{[z,M]}\big)(\chi_{Nr_0}^{\mathrm{T}})](w)\\
		&=  [\mathrm{P}(\tau[y,\mathrm{id}]F_{[z,M]})(r_0)\big)](N,w)\\
		&=  [\mathrm{P}\big(\Big(\int_{{R}}^{\oplus} U^r[y,\mathrm{id}]\otimes I \ud r \Big)F_{[z,M]}\big)(r_0)](N,w)\\ 
		&=  [\mathrm{P}\big((U^{r_0}[y,\mathrm{id}]\otimes I)u_{ I_{[z,M]}(F),\delta_M\otimes\delta_z} \big)({r_0})](N,w)\\
		&=  [p\big((U^{r_0}[y,\mathrm{id}]\otimes I)u_{ I_{[z,M]}(F),\delta_M\otimes\delta_z}({r_0})\big)](N,w)\\
		&= [p\big((U^{r_0}[y,\mathrm{id}]\otimes I)(I_{[z,M]}(F)({r_0})(\delta_M\otimes \delta_z)  \big)](N,w)\\
		&= [p\big(U^{r_0}[y,\mathrm{id}]\big(I_{[z,M]}(F)({r_0})\delta_M\big) \otimes \delta_z\big)](N,w)\\
		&= [U^{r_0}[y,\mathrm{id}]\big(I_{[z,M]}(F)({r_0})\delta_M\big)](N)\delta_z(w).
	\end{align*}
	Now 
	\begin{align*}
		&[U^{r_0}[y,\mathrm{id}]\big(I_{[z,M]}(F)({r_0})\delta_M\big)](N) \\
		&= \chi_{r_0}^\mathrm{T}\big(\gamma(N)^{-1}[y,\mathrm{id}]\gamma(N)\big)[I_{[z,M]}(F)({r_0})\delta_M](N)\\
		&= (N\cdot_{\idref{6}}\chi_{r_0}^\mathrm{T})([y,\mathrm{id}])[I_{[z,M]}(F)({r_0})\delta_M](N)\\
		&= E_{[y,\mathrm{id}]}(N\cdot_{\idref{6}}\chi_{r_0}^\mathrm{T})[I_{[z,M]}(F)({r_0})\delta_M](N)\\
		&=  E_{[y,\mathrm{id}]}(\chi_{N{r_0}}^\mathrm{T})[I_{[z,M]}(F)({r_0})\delta_M](N),
	\end{align*}
	where we have used \eqref{pionthatreflectparamdomain} in the last equality. Thus,
	\begin{align*}
		&[\big((\kappa\circ \tau[y,\mathrm{id}])F_{[z,M]}\big)(\chi_{N{r_0}}^{\mathrm{T}})](w)\\
		&= E_{[y,\mathrm{id}]}(\chi_{N{r_0}}^\mathrm{T})[I_{[z,M]}(F)({r_0})\delta_M](N)\delta_z(w)\\
		&= E_{[y,\mathrm{id}]}(\chi_{N{r_0}}^\mathrm{T})[p\big(I_{[z,M]}(F)({r_0})\delta_M\otimes \delta_z\big)](N,w)\\
		&= E_{[y,\mathrm{id}]}(\chi_{N{r_0}}^\mathrm{T})[p\big(u_{ I_{[z,M]}(F),\delta_M\otimes\delta_z}({r_0})\big)](N,w)\\
		&= E_{[y,\mathrm{id}]}(\chi_{N{r_0}}^\mathrm{T})[\mathrm{P}(u_{ I_{[z,M]}(F),\delta_M\otimes\delta_z})({r_0})](N,w)\\
		&= 	E_{[y,\mathrm{id}]}(\chi_{N{r_0}}^\mathrm{T})[(\kappa u_{ I_{[z,M]}(F),\delta_M\otimes\delta_z})(\chi_{N{r_0}}^\mathrm{T})](w)\\
		&= [\big(\mathrm{M}[y,\mathrm{id}] (\kappa F_{[z,M]})\big)(\chi_{N{r_0}}^\mathrm{T})](w)\\
		&= [\big((\mathrm{M}[y,\mathrm{id}]\circ \kappa)F_{[z,M]}\big)(\chi_{N{r_0}}^\mathrm{T})](w).
	\end{align*}
	Since $w\in\mathcal{L}^*$ is arbitrary, we have
	\[
	\big((\kappa\circ \tau[y,\mathrm{id}])F_{[z,M]}\big)(\chi_{N{r_0}}^{\mathrm{T}}) = \big((\mathrm{M}[y,\mathrm{id}]\circ \kappa)F_{[z,M]}\big)(\chi_{N{r_0}}^\mathrm{T}).
	\]
	Moreover, since $\chi_{N{r_0}}^\mathrm{T} \in \widehat{\mathrm{T}}_0$ was arbitrary, and $\widehat{\mathrm{T}} -  \widehat{\mathrm{T}}_0$ is of measure zero, we have
	\[
	(\kappa\circ \tau[y,\mathrm{id}])F_{[z,M]} =(\mathrm{M}[y,\mathrm{id}]\circ \kappa)F_{[z,M]}
	\]
	as elements of $L^2\big(\widehat{\mathrm{T}}; l^2(\mathcal{L}^*)\big)$. Since $\kappa\circ \tau[y,\mathrm{id}]$ and $\mathrm{M}[y,\mathrm{id}]\circ \kappa$ are unitary operators on $L^2\big(\widehat{\mathrm{T}}; l^2(\mathcal{L}^*)\big)$, it follows from \eqref{Fdecomp} that 
	\[
	(\kappa\circ \tau[y,\mathrm{id}])F =(\mathrm{M}[y,\mathrm{id}]\circ \kappa)F.
	\]
	Lastly, since $F\in L^2\big(R;l^2(\Pi)\otimes l^2(\mathcal{L}^*)\big)$ was arbitrary, we have
	\[
	\kappa\circ \tau[y,\mathrm{id}] = \mathrm{M}[y,\mathrm{id}]\circ \kappa,
	\]
	and this holds for all $[y,\mathrm{id}] \in \mathrm{T}$, whence $\mathrm{M} \sim \tau\vert_{\mathrm{T}}$, as claimed.
	
	\begin{defn}\label{detset}
		Given a measure space $(X, \mathcal{M},\mu)$, a \textit{determining set} for $L^1(X)$ is defined as a subset $\mathcal{D}$ of $L^\infty(X)$ such that, for each $f\in L^1(X)$, 
		\[
		\text{if}\ \int_X fg \ud \mu = 0 \ \ \text{for all}\ g\in \mathcal{D}, \ \text{then}\ f=0.
		\]
		This is definition 2.2 in \cite{BownikRoss}. We note that some authors use a slightly different definition of determining set (\textit{cf.} Definition 2.8.1, p.~34, \cite{HillePhillips}).
	\end{defn}
	\begin{lemma}
		Given $\upphi \in L^1(\widehat{\mathrm{T}})$, we have
		\begin{equation}\label{lemmadetset}
			\sum_{[y,\mathrm{id}] \in \mathrm{T}}\Big|\int_{\widehat{\mathrm{T}}}\upphi(\chi) E_{[y,\mathrm{id}]}(\chi)\ud \chi \Big|^2 = \|\upphi\|^2.
		\end{equation}
	\end{lemma}
	\begin{proof}
		Given $\upphi \in L^1(\widehat{\mathrm{T}})$ and $E_{[y,\mathrm{id}]}\in \widehat{\widehat{\mathrm{T}}}$, we have 
		\[
		\widehat{\upphi}(E_{[y,\mathrm{id}]}) = \int_{\widehat{\mathrm{T}}} \upphi(\chi) E_{[y,\mathrm{id}]}(\chi) \ud\chi.
		\]
		Hence
		\[
		\sum_{[y,\mathrm{id}] \in \mathrm{T}}\Big|\int_{\widehat{\mathrm{T}}}\upphi(\chi) E_{[y,\mathrm{id}]}(\chi)\ud \chi \Big|^2 = 	\sum_{[y,\mathrm{id}] \in \mathrm{T}}|\widehat{\upphi}(E_{[y,\mathrm{id}]})|^2 = \|\widehat{\upphi}\|^2.
		\]
		If $\upphi \in L^2(\widehat{\mathrm{T}})$, then by $\|\widehat{\upphi}\|^2 = \|\upphi\|^2 < \infty$, by Plancherel's Theorem (\ref{Plancherel}), and \eqref{lemmadetset} follows. On the other hand, if $\upphi \notin L^2(\widehat{\mathrm{T}})$, then $\|\widehat{\upphi}\|^2 = \|\upphi\|^2 = \infty$, and both sides of \eqref{lemmadetset} equal $\infty$.
	\end{proof}
	\begin{cor}\label{detsetcor}
		$\widehat{\widehat{\mathrm{T}}}$ is a determining set for $L^1(\widehat{\mathrm{T}})$.
	\end{cor}
	\begin{proof}
		If $\upphi \in L^1(\widehat{\mathrm{T}})$ and $\int_{\widehat{\mathrm{T}}} \upphi E_{[y,\mathrm{id}]} \ud\tilde\mu_{\widehat{\mathrm{T}}}=0$ for all  $E_{[y,\mathrm{id}]}\in \widehat{\widehat{\mathrm{T}}}$, then $\|\upphi\|= 0$, by \eqref{lemmadetset}, whence $\upphi = 0$ a.e. in $\widehat{\mathrm{T}}$. Thus $\upphi = 0$ in $L^1(\widehat{\mathrm{T}})$.
	\end{proof}
	
	\begin{prop}\label{kappaOrth}
		Let $W$ be a closed, $\tau$-invariant subspace of $L^2\big(R;l^2(\Pi)\otimes l^2(\mathcal{L}^*)\big)$ and let $F \in W$ and let $G\in L^2\big({R};l^2(\Pi)\otimes l^2(\mathcal{L}^*)\big)$. Then $F\perp G$ if and only if $F(r) \perp G(r)$ for almost every $r\in R$. 
	\end{prop}
	\begin{proof}
		First note that if $F(r) \perp G(r)$ for almost every $r \in {R}$, then 
		\[
		\langle F,G\rangle = \int_R \langle F(r),G(r)\rangle  \ud r =0,
		\]
		so that $F\perp G$. For the converse, suppose that $\xi \in W^\perp$. Recall from Definition \ref{invarDef} that $W$ being $\tau$-invariant means that $\tau[x,L]F \in W$ for any $F \in W$ and $[x,L]\in\Gamma$. In particular, if $W$ is $\tau$-invariant and $F \in W$, we have $\tau[y,\mathrm{id}]F \in W$ for all $[y,\mathrm{id}]\in\mathrm{T}$. Hence
		\begin{align*}
			0 &= \langle \tau[y,\mathrm{id}]F, G\rangle\\
			&= \langle \kappa \big(\tau[y,\mathrm{id}]F\big), 
			\kappa G\rangle	 \\
			&=  \langle \mathrm{M}[y,\mathrm{id}](\kappa F	), \kappa G\rangle	 \\
			&= \int_{\widehat{\mathrm{T}}} \langle \mathrm{M}_{[y,\mathrm{id}]} (\kappa F)(\chi), (\kappa G)(\chi)\rangle_{l^2(\mathcal{L}^*)} \ud \tilde\mu_{\widehat{\mathrm{T}}}(\chi)\\
			&= \int_{\widehat{\mathrm{T}}} \sum_{z\in \mathcal{L}^*} [\mathrm{M}_{[y,\mathrm{id}]} (\kappa F)(\chi)](z) \overline{[(\kappa G)(\chi)](z)}\ud \tilde\mu_{\widehat{\mathrm{T}}}(\chi)\\
			&= \int_{\widehat{\mathrm{T}}} \sum_{z\in \mathcal{L}^*} E_{[y,\mathrm{id}]}(\chi)[(\kappa F)(\chi)](z) \overline{[(\kappa G)(\chi)](z)}\ud \tilde\mu_{\widehat{\mathrm{T}}}(\chi)\\
			&= \int_{\widehat{\mathrm{T}}} E_{[y,\mathrm{id}]}(\chi) \sum_{z\in \mathcal{L}^*} [(\kappa F)(\chi)](z) \overline{[(\kappa G)(\chi)](z)}\ud \tilde\mu_{\widehat{\mathrm{T}}}(\chi)\\
			&= \int_{\widehat{\mathrm{T}}} E_{[y,\mathrm{id}]}(\chi) \sum_{z\in \mathcal{L}^*} [(\kappa F)(\chi)](z) \overline{[(\kappa G)(\chi)](z)}\ud \tilde\mu_{\widehat{\mathrm{T}}}(\chi)\\
			&= \int_{\widehat{\mathrm{T}}} \langle (\kappa F)(\chi) , (\kappa G)(\chi)\rangle_{l^2(\mathcal{L}^*)} E_{[y,\mathrm{id}]}(\chi) \ud \tilde\mu_{\widehat{\mathrm{T}}}(\chi).
		\end{align*}
		Note that $\chi\mapsto \langle (\kappa F)(\chi) , (\kappa G)(\chi)\rangle_{l^2(\mathcal{L}^*)}$ is in $L^1\big(\widehat{\mathrm{T}}\big)$, by the Cauchy-Schwarz and H\"older inequalities. It follows from Corollary \ref{detsetcor} that $\langle (\kappa F)(\chi) , (\kappa G)(\chi)\rangle_{l^2(\mathcal{L}^*)}=0$ for almost every $\chi \in \widehat{\mathrm{T}}$.
		By two applications of \ref{imagemeasureextended}, we have
		\small{
			\begin{align*}
	&0 = \int_{\widehat{\mathrm{T}}_0} | \langle (\kappa F)(\chi_{Mr}^{\mathrm{T}}) , (\kappa G)(\chi_{Mr}^{\mathrm{T}})\rangle_{l^2(\mathcal{L}^*)}|\ud \tilde\mu_{\widehat{\mathrm{T}}}(\chi_{Mr}^{\mathrm{T}})\\
	&= \int_{\widehat{\mathrm{T}}_0} \Big|\sum_{z\in\mathcal{L}^*}[(\kappa F)(\chi_{Mr}^{\mathrm{T}})](z) \overline{[(\kappa G)(\chi_{Mr}^{\mathrm{T}})](z)} \Big|\ud \tilde\mu_{\widehat{\mathrm{T}}}(\chi_{Mr}^{\mathrm{T}})\\
	&= \int_{\widehat{\mathrm{T}}_0} \Big|\sum_{z\in\mathcal{L}^*}(\kappa_3\circ \kappa_2 \circ \kappa_1)F(\chi_{Mr}^{\mathrm{T}},z)\overline{(\kappa_3\circ \kappa_2 \circ \kappa_1)G(\chi_{Mr}^{\mathrm{T}},z)} \Big|\ud \tilde\mu_{\widehat{\mathrm{T}}}(\chi_{Mr}^{\mathrm{T}})\\
	&= \int_{\widehat{\mathrm{T}}_0} \Big|\sum_{z\in\mathcal{L}^*}(\kappa_3\circ \kappa_2 \circ \kappa_1)F(\chi_{Mr}^{\mathrm{T}},z)\overline{(\kappa_3\circ \kappa_2 \circ \kappa_1)G(\chi_{Mr}^{\mathrm{T}},z)} \Big|\ud (\overline{m\circ \upsilon_0^{-1}})(\chi_{Mr}^{\mathrm{T}})\\
	&= \int_{\Pi R} \Big|\sum_{z\in\mathcal{L}^*}(\kappa_3\circ \kappa_2 \circ \kappa_1)F(\upsilon_0(Mr),z)\overline{(\kappa_3\circ \kappa_2 \circ \kappa_1)G(\upsilon_0(Mr),z)} \Big|\ud m(Mr)\\
	&= \int_{\Pi R} \Big|\sum_{z\in\mathcal{L}^*}(\kappa_2 \circ \kappa_1)F(Mr,z)\overline{(\kappa_2 \circ \kappa_1)G(Mr,z)} \Big|\ud [(\mu_\Pi \times m)\circ \Lambda^{-1}](Mr)\\
	&= \int_{\Pi \times R} \Big|\sum_{z\in\mathcal{L}^*}(\kappa_2 \circ \kappa_1)F\big(\Lambda(M,r),z\big)\overline{(\kappa_2 \circ \kappa_1)G\big(\Lambda(M,r),z\big)}\Big|\ud(\mu_\Pi \times m)(M,r)\\
	&= \int_{\Pi \times R} \Big|\sum_{z\in\mathcal{L}^*}(\kappa_1F)(M,r,z)\overline{(\kappa_1G)(M,r,z)}\Big|\ud(\mu_\Pi \times m)(M,r)\\
	&= \int_{\Pi \times R} \Big|\sum_{z\in\mathcal{L}^*}[(\mathrm{P}F)(r)](M,z) \overline{[(\mathrm{P}G)(r)](M,z)}\Big|\ud(\mu_\Pi \times m)(M,r)\\
	&=\int_R \sum_{M\in \Pi} \Big|\sum_{z\in\mathcal{L}^*}[(\mathrm{P}F)(r)](M,z) \overline{[(\mathrm{P}G)(r)](M,z)}\Big|\ud r\\
	&\geq  \int_R \Big|\sum_{M\in \Pi}\sum_{z\in\mathcal{L}^*}[(\mathrm{P}F)(r)](M,z) \overline{[(\mathrm{P}G)(r)](M,z)}\Big|\ud r\\
	&=  \int_R \big| \langle (\mathrm{P}F)(r) ,(\mathrm{P}G)(r) \rangle_{l^2(\Pi\times\mathcal{L}^*)} \big|\ud r\\
	&=  \int_R \big| \langle p(F(r)),p(G(r)) \rangle_{l^2(\Pi\times\mathcal{L}^*)} \big|\ud r\\
	&=  \int_R \big| \langle F(r),G(r) \rangle_{l^2(\Pi) \otimes l^2(\mathcal{L}^*)} \big|\ud r.
	\end{align*}
		}
		Thus $F(r) \perp G(r)$ for almost every $r\in R$, as claimed.
	\end{proof}
	Note that $L^2(R)$ is separable, since $L^2(\mathbb{R}^n)$ is---see \cite{FollandReal}, p.~187, Exercise 13. Moreover, it is clear that $l^2(\Pi)\otimes l^2(\mathcal{L}^*)$ is separable, with countable orthonormal basis $\{\delta_M \otimes \delta_z: M\in \Pi, z\in \mathcal{L}^*\}$. Thus, $L^2\big(R;l^2(\Pi)\otimes l^2(\mathcal{L}^*)\big)$ is separable, by Lemma \ref{vectorsep}. Thus, given a closed subspace $W \subset L^2\big(R;l^2(\Pi)\otimes l^2(\mathcal{L}^*)\big)$ we may, as in Theorem \ref{mainpointwisedecompthm}, select a countable dense subset $\mathcal{A}$ of $L^2\big(R;l^2(\Pi)\otimes l^2(\mathcal{L}^*)\big)$, and define a range function $J_W^\mathcal{A}:R \to \{\text{closed subspaces of}\ l^2(\Pi)\otimes l^2(\mathcal{L}^*)\}$ associated with~$W$. More precisely, $J_W^\mathcal{A}$ is a representative from an equivalence class of a.e.-equal range functions.
	
	\begin{theorem}\label{RangeFunctionOrthTheorem}
		Let $W$ be a closed, $\tau$-invariant subspace of $L^2\big(R;l^2(\Pi)\otimes l^2(\mathcal{L}^*)\big)$. Let $\mathcal{A}$ be a countable dense subset of $L^2\big(R;l^2(\Pi)\otimes l^2(\mathcal{L}^*)\big)$, and $J_W^\mathcal{A}$ a range function as described in the preceding paragraph. Let $G$ be an element of $L^2\big(R;l^2(\Pi)\otimes l^2(\mathcal{L}^*)\big)$. Then $G \in W^\perp$ if and only if $G(r) \perp J_W^\mathcal{A}(r)$ for almost every $r \in R$. In particular, $W$ satisfies the orthogonality condition \eqref{orthCond}.
	\end{theorem}
	\begin{proof}
		Recall that $J_W^\mathcal{A}$ is defined as in \eqref{rangefromsubspace}: namely, for each $\Phi \in \mathcal{A}$, we make a choice $f(\Phi)\in [P_W\Phi]_{L^2}$, and define
		\[
		J_W^\mathcal{A}[f](r) = \bigvee\{f(\Phi)(r): \Phi\in\mathcal{A}\}.
		\]
		It was shown that this range function is independent (up to a.e.-equivalence) of the choice of representatives $f(\Phi)\in [P_W\Phi]_{L^2}$, which justifies writing $J_W^\mathcal{A}$ without specifying the choice $f$. Suppose $G(r)\perp J_W^\mathcal{A}(r)$ for a.e. $r\in R$. Then $G(r)\perp f(\Phi)(r)$ a.e., for all $\Phi \in \mathcal{A}$. Thus $G \perp f(\Phi)= P_W \Phi$ in $L^2\big(R;l^2(\Pi)\otimes l^2(\mathcal{L}^*)\big)$, for all $\Phi \in \mathcal{A}$. But $\{P_W \Phi: \Phi \in \mathcal{A}\}$ is dense in $W$, by the argument in the proof of \ref{mainpointwisedecompthm}(1). Thus $G \perp W$.

		Conversely, let $G  \in W^\perp$. Since $f(\Phi) = P_W \Phi$ a.e., and $P_W \Phi \in W$, it follows from Proposition \ref{kappaOrth} that $G(r) \perp f(\Phi)(r)$ for almost all $r\in R$; i.e., there exists a null set $N_\Phi \subset {R}$ such that  $G(r) \perp f(\Phi)(r)$ on $N_\Phi^c$. Letting $N= \bigcup_{\Phi\in\mathcal{A}} N_\Phi$, we have $G(r) \perp J_W^\mathcal{A}[f](r)$ on $N^c$. Thus $G(r) \perp J_W^\mathcal{A}(r)$ almost everywhere, as claimed.
	\end{proof}
	It follows from the preceding theorem and Definition \ref{mblerangefun} that each $\tau$-invariant subspace $W\subset L^2\big(R;l^2(\Pi)\otimes l^2(\mathcal{L}^*)\big)$ has associated to it a unique (up to a.e.-equivalence) measurable range function $J_W$. The following theorem will be essential for characterizing $\tau$-invariant subspaces.
	\begin{theorem}\label{MainDecompTheorem}\hspace{1cm}
		\begin{list}{(\arabic{ctr})}{\usecounter{ctr}}
			\item	Let $W$ be a closed, $\tau$-invariant subspace of $L^2\big({R};l^2(\Pi) \otimes l^2(\mathcal{L}^*)\big)$, and let $J_W$ denote the measurable range function associated with $W$. Let $P_W$ be the orthogonal projection of $L^2\big(R;l^2(\Pi) \otimes l^2(\mathcal{L}^*)\big)$ onto $W$, and $P_J(r)$ denote the orthogonal projection from $l^2(\Pi) \otimes l^2(\mathcal{L}^*)$ onto $J_W(r)$. Then
			\begin{list}{(\alph{ctr1})}{\usecounter{ctr1}}
				\setcounter{ctr1}{0}
				\item Regarding $\int_R^\oplus J_W(r)\ud r$ as a subspace of $L^2\big(R;l^2(\Pi) \otimes l^2(\mathcal{L}^*)\big)$, we have
				\[
				W = \int_R^\oplus J_W(r)\ud r\ \ \ \text{and}\ \ \ 	P_W =\int_R^\oplus P_J(r)\ud r.
				\]
				\item $J_W(r)$ is invariant under $U^r\otimes I$ for almost every $r \in R$.
			\end{list}
			\item Conversely, if $r\mapsto W_r$ is a measurable range function from $R$ into subspaces of $l^2(\Pi) \otimes l^2(\mathcal{L}^*)$, with $W_r$ invariant under $U^r\otimes I$ for a.e. $r$, then setting
			\[
			W=\int_R^\oplus W_r \ud r
			\]
		and regarding $W$ as a subspace of $L^2\big(R;l^2(\Pi) \otimes l^2(\mathcal{L}^*)\big)$, we have that $W$ is $\tau$-invariant. Moreover, $W_r = J_W(r)$ for a.e. $r$, and $W = \int_{R}^\oplus J_W(r) \ud r$. 
		\end{list}
	\end{theorem} 
	\begin{proof}
		(a) is an application of Theorem \ref{pointwisedecompcor1}, since any $\tau$-invariant subspace satisfies the orthogonality condition \eqref{orthCond}, by Theorem \ref{RangeFunctionOrthTheorem}. For (b), note that $W$ being $\tau$-invariant is equivalent to $P_W \in \tau(\Gamma)^\prime$, by Proposition \ref{invorth}. Thus, if $W$ is $\tau$-invariant, then 
		\[
		\int_R^\oplus U^r[x,L] \otimes I\ud r \int_R^\oplus P_{J}(r)\ud r = \int_{{R}}^\oplus P_{J}(r)\ud r \int_R^\oplus U^r[x,L] \otimes I\ud r
		\]
		for all $[x,L]\in \Gamma$. By Lemma \ref{directintegralprops}, this is equivalent to 
		\[
		\int_R^\oplus (U^r[x,L] \otimes I)P_{J}(r)\ud r = \int_R^\oplus P_{J}(r)(U^r[x,L] \otimes I)\ud r. 
		\]
		By Lemma \ref{aeequalops} it follows that 
		\[
		(U^r[x,L] \otimes I)P_{J}(r) = P_{J}(r) (U^r[x,L] \otimes I)
		\]
		for a.e. $r\in R$, and all $[x,L]\in \Gamma$. Thus $P_{J}(r)\in \big((U^r\otimes I)(\Gamma)\big)^\prime$ for a.e. $r$, and by Proposition \ref{invorth} again, we have that $J_W(r)$ is $(U^{\chi_r^\mathrm{T}}\otimes I)$-invariant, for a.e. $r\in {R}$.
		
		To prove (2), we note that for $G \in \int_R^\oplus W_r \ud r$, we have $G(r)\in W_r$ for a.e. $r\in R$, whence $(U^r[x,L]\otimes I)G(r) \in W_r$ for a.e. $r\in R$ and all $[x,L]\in \Gamma$. Since 
		\[
		\{U^r[x,L]\otimes I\}_{r\in{R}}
		\]
		is a measurable field of operators on the constant field $\mathscr{H}_r \equiv l^2(\Pi)\otimes l^2(\mathcal{L}^*)$ (as shown in Lemma \ref{Uchimeas1}), it follows that $r\mapsto (U^r[x,L]\otimes I)G(r)$ is a measurable vector field with respect to this constant field. Moreover, 
		\[
		\|U^r[x,L]\otimes I\|_\infty = 1 < \infty,
		\] 
		and therefore the measurable operator field $\{U^r[x,L]\otimes I\}_{r\in R}$ gives rise to a bounded linear operator $\tau[x,L] = \int_R^\oplus U^r[x,L]\otimes I \ud r$ on $\int_{{R}}^\oplus \mathscr{H}_r \ud r  =  L^2\big({R};l^2(\Pi) \otimes l^2(\mathcal{L}^*)\big)$, defined by
		\[
		\tau[x,L] G = \int_R^\oplus(U^r[x,L]\otimes I)G(r)\ud r.
		\]
		Now by the criteria described in Corollary \ref{criteriaForDirectIntSubspace}, we have 
		\[
		\tau[x,L] G  \in \int_R^\oplus W_r \ud r,
		\]
		which proves that $W$ is $\tau$-invariant. The statement that $W_r = J_W(r)$ for a.e. $r$ follows from Lemma \ref{aeequalops}. The assertion $W = \int_R^\oplus J_W(r) \ud r$ follows immediately from (1).
	\end{proof}
	\begin{cor}\label{taucor}
 Every subrepresentation of $\tau$ decomposes as a direct integral of subrepresentations of the $U^r \otimes I$. 
	\end{cor}
	\begin{rem}
		We note that the argument in this section rests on the idea of using a determining set to establish an orthogonality condition \eqref{orthCond}. This idea came from \cite{BownikRoss} (see in particular Example 2.2.(b), the proof of Theorem 2.4, and Lemma 3.5). However, in our setting, the structure is different---$\Gamma$ is not a subgroup of $\mathbb{R}^n$, but rather of $\mathrm{I}(\mathbb{R}^n)$. And we do not have a powerful structure theorem such as Theorem 3.1 from \cite{BownikRoss} available. Hence we introduce the $\kappa$ map and utilize Theorem \ref{imagemeasureextended}.
	\end{rem}

	\begin{rem}
		The use of complete measures in this section not necessary. It is included merely to show that there is nothing to stop us from endowing the subsets of $\mathbb{R}^n$ with the usual Lebesgue $\sigma$-algebra and measure.
	\end{rem}
	
	\chapter{Invariant Subspaces}\label{InvSubspaceChapter}
	
	\section{Preimages of Subspaces}\label{InvSubspaceSection}
	
	We showed in Section \ref{Decomposition of unitary} that $\widehat{\pi}$ was equivalent to $\tau =\int_{R}^{\oplus} U^r \otimes I \ud r$ via the unitary operator $\Upsilon$. We also defined $U = \mathcal{F}\circ \Upsilon$, where $\mathcal{F}: L^2(\mathbb{R}^n) \to L^2(\widehat{\mathbb{R}^n})$ is the Plancherel transform. We here rephrase the essential content of Theorem \ref{MainDecompTheorem}.
	\begin{lemma}\label{invChar}
		The closed subspace $W\subset L^2\big(R;l^2(\Pi) \otimes l^2(\mathcal{L}^*)\big)$ is $\tau$-invariant if and only if there is a measurable range function $r \mapsto W_r$ from $R$ into the closed subspaces of $l^2(\Pi)\otimes l^2(\mathcal{L}^*)$, with $W_r$ invariant under $U^r\otimes I$ for almost every $r \in R$, such that 
		\[
		W = \int_{R}^\oplus W_r \ud r,
		\]
		where this direct integral is understood as a subspace of  $L^2\big({R};l^2(\Pi)\otimes l^2(\mathcal{L}^*)\big)$. When this holds,  $W_r = J_W(r)$ for a.e. $r$, and $W = \int_{R}^\oplus J_W(r) \ud r$, where $J_W$ is the measurable range function associated with $W$. Moreover, every $\tau$-invariant subspace of $L^2\big(R;l^2(\Pi) \otimes l^2(\mathcal{L}^*)\big)$ arises in this way, and the correspondence $W \leftrightarrow J_W$ is one-to-one.
	\end{lemma}
	Since the $U^r$ are irreducible for $r \in R$, we may apply Corollary \ref{CorforTensorSubspaces} and Proposition \ref{measrangefntensorprod} to get:
	\begin{tcolorbox}	
		\begin{cor}\label{invChar1}
			The closed subspace $W\subset L^2\big(R;l^2(\Pi) \otimes l^2(\mathcal{L}^*)\big)$ is $\tau$-invariant if and only if there is a measurable range function $r\mapsto K_r$ from $R$ into the closed subspaces of $l^2(\mathcal{L}^*)$ such that
			\[
			W = \int_{R}^\oplus l^2(\Pi)\otimes K_r\ud r,
			\]
			where this direct integral is understood as a subspace of  $L^2\big({R};l^2(\Pi)\otimes l^2(\mathcal{L}^*)\big)$. Moreover, every $\tau$-invariant subspace of $L^2\big(R;l^2(\Pi) \otimes l^2(\mathcal{L}^*)\big)$ arises in this way, and the correspondence $W \leftrightarrow \{K_r\}_{r\in R}$ is one-to-one, up to modification on a null subset of $R$.
		\end{cor}
	\end{tcolorbox}
Recall from \eqref{upsinv1} that
	\begin{equation}\label{FTdescription}
		\Upsilon^{-1}(F)(\chi_{[z,M]\cdot_{\idref{1}} r}) = \chi_{r+z}(x_M)\langle F(r),\delta_M \otimes \delta_{z}\rangle; \ \ r\in R,\ [z,M]\in\Gamma^*.
	\end{equation}
	We note that
	\[
	\widehat{\mathbb{R}^n} = \{\chi_{[z,M]\cdot_{\idref{1}} \omega}: [z,M]\in \Gamma^*, \omega\in \Omega\},
	\]
	so that \eqref{FTdescription} defines $\Upsilon^{-1}(F)$ almost everywhere. Given a closed subspace $W \subset L^2\big(R;l^2(\Pi)\otimes l^2(\mathcal{L}^*)\big)$, we have $\Upsilon^{-1}(W) = \{\Upsilon^{-1}(F): F\in W\}$. Thus, from \eqref{FTdescription} and the preceding corollary, we get the following characterization of $\pi$-invariant subspaces:
	\begin{tcolorbox}
\begin{prop}\label{almostTheorem}
			A closed subspace  $V \subset L^2(\mathbb{R}^n)$ is $\pi$-invariant if and only if there exists a measurable range function $r \mapsto K_r$ from $R$ into the closed subspaces of $l^2(\mathcal{L}^*)$, uniquely determined almost everywhere by $V$, such that $V$ consists of all $f \in L^2(\mathbb{R}^n)$ satisfying
			\[
			\widehat{f}(\chi_{[z,M]\cdot_{\idref{1}} r}) = \chi_{r+z}(x_M)\langle F(r),\delta_M \otimes \delta_{z}\rangle
			\]
			for a.e. $r\in R$ and all $[z,M]\in \Gamma^*$, for some $F\in \int_R^\oplus l^2(\Pi)\otimes K_r \ud r$, where $x_M$ comes from the cross-section $\gamma:\Pi \to \Gamma$ with $\gamma(M) = [x_M, M]$.  Moreover, every $\pi$-invariant subspace of $L^2\big(R;l^2(\Pi) \otimes l^2(\mathcal{L}^*)\big)$ arises in this way, and the correspondence $V \leftrightarrow \{K_r\}_{r\in R}$ is one-to-one, up to modification on a null subset of $R$.
\end{prop}
	\end{tcolorbox}

\begin{lemma}
Let $F\in \int_{R}^\oplus l^2(\Pi)\otimes K_r \ud r$, where $r\mapsto K_r$ is a measurable range function. Let $\{\xi_n\}_{n\in\mathbb{N}}$ be a fundamental sequence of measurable vector fields for $\{K_r\}_{r\in R}$ such that\footnote{Such a sequence is guaranteed by Lemma \ref{basissequence}.}
\begin{list}{(\roman{ctr})}{\usecounter{ctr}}
	\item For each $r\in R$, $\{\xi_1(r),\xi_2(r),\ldots,\xi_{d(r)}(r)\}$ is an orthonormal basis for $K_r$, where $d(r) =\dim K_r$.
	\item $\xi_n(r) = 0$ for $n>d(r)$, if $d(r)<\infty$.
\end{list} 
Define $\bar{F}\in \int_{R}^\oplus l^2(\Pi)\otimes K_r \ud r$ by
\[
\bar{F}(r) = \sum_{(N,n) \in \Pi\times \mathbb{N}} \big(\chi_r(x_N) \langle F(r),\delta_N\otimes \xi_n(r)\rangle \big) \delta_N\otimes \xi_n(r).
\]
Then for each $M\in \Pi$ and $z\in \mathcal{L}^*$, we have
\[
\langle \bar{F}(r),\delta_M\otimes \delta_z \rangle= \chi_r(x_M)\langle F(r),\delta_M\otimes \delta_z\rangle. 
\]
\end{lemma}
\begin{proof}
Note that, for each $r\in R$ and $z\in\mathcal{L}^*$, we can write 
\[
\delta_z = P_K(r)\delta_z + (I_r-P_K(r))\delta_z,
\]
where $P_K(r)$ denotes the orthogonal projection of $l^2(\mathcal{L}^*)$ onto $K_r$, and $I_r$ denotes the identity operator on $K_r$. In particular, we have 
\[
\langle \delta_z,\xi_n(r) \rangle = \langle P_K(r)\delta_z,\xi_n(r)\rangle \ \ \text{for all}\ z\in\mathcal{L}^*, n\in\mathbb{N}, \ \text{and}\ r \in R.
\]
Also, since $\{\xi_n(r)\}$ is an orthonormal basis for $K_r$, we have
\[
P_K(r)\delta_z=\sum_{n\in\mathbb{N}} \langle P_K(r)\delta_z,\xi_n(r)\rangle \xi_n(r)\ \ \text{for each}\ z\in \mathcal{L}^*\ \text{and}\ r\in R.
\]
Therefore,
\begin{align*}
\langle \bar{F}(r),\delta_M\otimes \delta_z \rangle&= \Big\langle \sum_{(N,n) \in \Pi\times \mathbb{N}} \big(\chi_r(x_N) \langle F(r),\delta_N\otimes \xi_n(r)\rangle \big)\delta_N\otimes \xi_n(r), \delta_M\otimes \delta_z \Big\rangle \\
&= \sum_{n\in\mathbb{N}} \chi_r(x_M) \langle F(r),\delta_M\otimes \xi_n(r)\rangle \langle \xi_n(r),\delta_z\rangle \\
&= \chi_r(x_M) \sum_{n\in\mathbb{N}} \big\langle F(r), \overline{\langle \xi_n(r),\delta_z\rangle} \big(\delta_M\otimes \xi_n(r)\big) \big\rangle \\
&= \chi_r(x_M) \Big\langle F(r), \sum_{n\in\mathbb{N}}\langle \delta_z,\xi_n(r)\rangle \big(\delta_M\otimes \xi_n(r)\big) \Big\rangle \\
&= \chi_r(x_M) \Big\langle F(r),  \delta_M \otimes \sum_{n\in\mathbb{N}}\langle P_K(r)\delta_z,\xi_n(r)\rangle \xi_n(r) \Big\rangle \\
&= \chi_r(x_M) \langle F(r), \delta_M\otimes P_K(r)\delta_z \rangle\\
&= \chi_r(x_M) \langle F(r), \delta_M\otimes \delta_z \rangle,
	\end{align*}
	as desired.
\end{proof}

	\begin{cor}\label{absorbchiomega}
In the description of $\widehat{f}$ in the preceding characterization, the $\chi_{r+z}(x_M)$ factor can be replaced with $\chi_z(x_M)$---that is, the $\chi_r(x_M)$ term gets absorbed into $\int_{R}^\oplus l^2(\Pi)\otimes K_r \ud r$.
\end{cor}

\begin{rem}
	The $\pi$-invariant subspaces should, of course, be independent of the choice of cross-section $\gamma$. To verify this, recall from Proposition \ref{differencecrosssec} that if $\gamma, \gamma^\prime$ are two cross-sections with $\gamma(M) = [x_M,M]$ and $\gamma^\prime(M) = [x_M^\prime,M]$ for all $M\in\Pi$, then $x_M - x_M^\prime \in\mathcal{L}$. Moreover, by the above corollary, the factor involving the cross-section can be replaced with $\chi_z(x_M)$. Now
	\[
 \chi_z(x_M^\prime) = \chi_z\big(x_M + (x_M^\prime - x_M)\big) = \chi_z(x_M)\chi_z(x_M^\prime - x_M) =\chi_z(x_M).
	\]
Thus the choice of cross-section makes no difference in our characterization.
		\end{rem}
With the above Corollary and Remark, we can strengthen the characterization in Proposition \ref{almostTheorem}:
	\begin{tcolorbox}
	\begin{theorem}[Invariant Subspaces Characterization]\label{BigInvTheorem2} 
		A closed subspace  $V \subset L^2(\mathbb{R}^n)$ is $\pi$-invariant if and only if there exists a measurable range function $r \mapsto K_r$ from $R$ into the closed subspaces of $l^2(\mathcal{L}^*)$, uniquely determined almost everywhere by $V$, such that $V$ consists of all $f \in L^2(\mathbb{R}^n)$ satisfying
		\[
		\widehat{f}(\chi_{[z,M]\cdot_{\idref{1}}r}) = \chi_z(x_M)\langle F(r),\delta_M \otimes \delta_{z}\rangle
		\]
		for a.e. $r\in R$ and all $[z,M]\in \Gamma^*$, for some $F\in \int_R^\oplus l^2(\Pi)\otimes K_r \ud r$, where $x_M$ comes from any cross-section $\gamma:\Pi \to \Gamma$ with $\gamma(M) = [x_M, M]$.  Moreover, every $\pi$-invariant subspace of $L^2\big(R;l^2(\Pi) \otimes l^2(\mathcal{L}^*)\big)$ arises in this way, and the correspondence $V \leftrightarrow \{K_r\}_{r\in R}$ is one-to-one, up to modification on a null subset of $R$.
	\end{theorem}
\end{tcolorbox}
		\section[Examples with \texorpdfstring{$\Gamma = pg$}{Gamma=pg}]{Illustrative Examples with \texorpdfstring{$\boldsymbol{\Gamma = pg}$}{Gamma=pg}}\label{Gammapgex}
	Recall that in Section \ref{Examples}, we illustrated Theorem \ref{BigInvTheorem2} for the 2-dimensional nonsymmorphic wallpaper group $\Gamma=pg$. In this section, we expand on this example and consider some subspaces of $L^2(\mathbb{R}^2)$ which are not symmorphic under $\Gamma=pg$.
	\begin{ex}
In Example \ref{Intropgex}, we saw that we could take the fundamental domain for $\Gamma^*$ to be the open rectangle 
\[
R = \left( -\tfrac12, \tfrac12\right) \times \left(-\tfrac12, 0 \right).
\]
We also saw that the the cross-section $\gamma$ can be taken to be $\gamma(\mathrm{id}) = [(0,0),\mathrm{id}]$ and $\gamma(\sigma) = [(\tfrac12,0),\sigma]$. Our theorem informs us that a closed subspace $V\subset L^2(\mathbb{R}^2)$ is invariant under shifts by $pg$ if and only if there is some measurable range function $r \mapsto K_r$ from $R$ into the closed subspaces of $l^2(\mathbb{Z}^2)$, uniquely determined almost everywhere by $V$, such that $V$ consists of all $f \in L^2(\mathbb{R}^2)$ such that, for some $F\in \int_R^\oplus l^2(\Pi)\otimes K_r \ud r$, 
	\begin{equation}\label{pgexformula}
	\widehat{f}(\chi_{[z,M]\cdot_{\idref{1}}r})	= 
	\left\{ \begin{array}{ll}
		\langle F(r),\delta_\sigma \otimes \delta_{z}\rangle &\text{if } M=\sigma\ \text{and}\ z_1\ \text{even}\\
		-\langle F(r),\delta_\sigma \otimes \delta_{z}\rangle &\text{if } M=\sigma\ \text{and}\ z_1\ \text{odd}\\
		\langle F(r),\delta_{\mathrm{id}}  \otimes \delta_{z}\rangle &\text{if } M=\mathrm{id}
	\end{array} \right. 
\end{equation}
for a.e. $r\in R$ and all $z\in \mathbb{Z}^2$. Given the measurable field $F\in \int_{R}^\oplus l^2(\Pi)\otimes K_r\ud r$, the above equation defines $\widehat{f}$ a.e. as a function of 3 variables, $r, z$, and $M$. Let $E$ be a measurable subset of $R$, as pictured below. 
	\begin{figure}[H]
	\pgfplotsset{
		axis lines=center, 
		grid,
		grid style={very thin, gray},
		xmin=-4,    xmax=4,   xtick={-4,...,4},
		x tick style={draw=none}, 
		ymin=-4,    ymax=4,   
		y tick style={draw=none},   
		x tick style={draw=none},  
		ytick={-4,...,4},
		restrict y to domain=-4:4, 
		ticklabel style={font=\scriptsize, fill=white, inner sep=2pt},
		domain=-1:11, samples=100,
		no marks, 
		every axis plot post/.append style={ultra thick,},
	}
	\centering
	\begin{tikzpicture}[scale=1]
	\begin{axis}[thick,
		xticklabels={-2,,-1,,0,,1,,2},
		yticklabels={-2,,-1,,0,,1,,2},
		xlabel={$x$},
		ylabel={$y$},
	every axis x label/.style={
		at={(ticklabel* cs:1.02)},
		anchor=west,
	},
	every axis y label/.style={
		at={(ticklabel* cs:1.02)},
		anchor=south,
	},
		]
		\fill[line width=1pt, color=RoyalBlue, opacity=.5] (-1, -1) -- (1,-1) -- (1,0) -- (-1,0) -- (-1,-1); 
		\draw[very thin, gray] (1, -1) -- (1,1);
		\draw[very thin, gray] (-1, -1) -- (-1,0);
		\draw[very thin, gray] (-1, 1) -- (1,1);
		\draw[very thin, gray] (-1, -1) -- (1,-1);
		\path[color=RoyalBlue, line width=.5pt, draw,use Hobby shortcut,closed=true]
		(-.7,-.7) .. (-.5,-.9) .. (-.3,-.6) .. (.8,-.4) .. (.2,-.2) .. (-.2,-.3) .. (-.6,-.5);
		\fill[color=RoyalBlue, opacity=.7, use Hobby shortcut,closed=true]
		(-.7,-.7) .. (-.5,-.9) .. (-.3,-.6) .. (.8,-.4) .. (.2,-.2) .. (-.2,-.3) .. (-.6,-.5);
		\node[color=black] at (-.5,.3) {\scriptsize $R$};	
		\node[color=black] at (.5,-.5) {\scriptsize $E$};	
	\end{axis}
\end{tikzpicture}
	\caption{The fundamental domain $R = \left( -\tfrac12, \tfrac12\right) \times \left(-\tfrac12, 0 \right)$ for $\Gamma^*$, with a measurable subset $E$.}
\end{figure}
Let $r \mapsto K_r$ be the measurable range function given as follows: let $\beta =  \delta_{(0,0)} + \delta_{(1,0)} \in l^2(\mathbb{Z}^2)$, and let 
\[
K_r = \left\{ \begin{array}{ll}
	\text{span}\, \beta &\text{if } r\in E\\
	\{0\} &\text{otherwise}
\end{array} \right.
\] 
In particular, for $F\in \int_R^\oplus l^2(\Pi)\otimes K_r \ud r$, we have $F=0$ on $E^c$. The picture below illustrates $\widehat{f}$ satisfying \eqref{pgexformula}:
	\begin{figure}[H]
	\pgfplotsset{
		axis lines=center, 
		grid,
		grid style={very thin, gray},
		xmin=-4,    xmax=4,   xtick={-4,...,4},
		x tick style={draw=none}, 
		ymin=-4,    ymax=4,   
		y tick style={draw=none},   
		x tick style={draw=none},  
		ytick={-4,...,4},
		restrict y to domain=-4:4, 
		ticklabel style={font=\scriptsize, fill=white, inner sep=2pt},
		domain=-1:11, samples=100,
		no marks, 
		every axis plot post/.append style={ultra thick,},
	}
	\centering
	\begin{tikzpicture}[scale=1, every node/.style={inner sep=0}]
		\begin{axis}[thick,
			xticklabels={-2,,-1,,0,,1,,2},
			yticklabels={-2,,-1,,0,,1,,2},
			xlabel={$x$},
			ylabel={$y$},
		every axis x label/.style={
			at={(ticklabel* cs:1.02)},
			anchor=west,
		},
		every axis y label/.style={
			at={(ticklabel* cs:1.02)},
			anchor=south,
		},
			]
			\draw[very thin, gray] (1, -1) -- (1,1);
			\draw[very thin, gray] (-1, -1) -- (-1,0);
			\draw[very thin, gray] (-1, 1) -- (1,1);
			\draw[very thin, gray] (-1, -1) -- (1,-1);
			\path[color=RoyalBlue, line width=.5pt, draw,use Hobby shortcut,closed=true]
			(-.7,-.7) .. (-.5,-.9) .. (-.3,-.6) .. (.8,-.4) .. (.2,-.2) .. (-.2,-.3) .. (-.6,-.5);
			\fill[color=RoyalBlue, opacity=.8, use Hobby shortcut,closed=true]
			(-.7,-.7) .. (-.5,-.9) .. (-.3,-.6) .. (.8,-.4) .. (.2,-.2) .. (-.2,-.3) .. (-.6,-.5);
			\path[color=RoyalBlue, line width=.5pt, draw,use Hobby shortcut,closed=true]
		(-.7,-.7) .. (-.5,-.9) .. (-.3,-.6) .. (.8,-.4) .. (.2,-.2) .. (-.2,-.3) .. (-.6,-.5);
		\node[color=black] at (2.3,-1.2) {\tiny $[(1,0),\mathrm{id}]\cdot E$};
		\node[color=black] at (2.3,1.2) {\tiny $[(1,0),\sigma]\cdot E$};	
			\node[color=black] at (-0.2,-1.2) {\tiny $E$};
		\node[color=black] at (-0.2,1.2) {\tiny $[(0,0),\sigma]\cdot E$};		
		\path[color=RoyalBlue, line width=.5pt, draw,use Hobby shortcut,closed=true]
		(1.3,-.7) .. (1.5,-.9) .. (1.7,-.6) .. (2.8,-.4) .. (2.2,-.2) .. (1.8,-.3) .. (1.4,-.5);
		\fill[color=RoyalBlue, opacity=.8, use Hobby shortcut,closed=true]
		(1.3,-.7) .. (1.5,-.9) .. (1.7,-.6) .. (2.8,-.4) .. (2.2,-.2) .. (1.8,-.3) .. (1.4,-.5);
			\path[color=RoyalBlue, line width=.5pt, draw,use Hobby shortcut,closed=true]
		(1.3,-.7) .. (1.5,-.9) .. (1.7,-.6) .. (2.8,-.4) .. (2.2,-.2) .. (1.8,-.3) .. (1.4,-.5);
			\begin{scope}[yscale=-1,xscale=1, shift={($(0,-8)-(0,0)$)}]
			\path[color=YellowOrange, line width=.5pt, draw,use Hobby shortcut,closed=true]
			(-.7,-.7) .. (-.5,-.9) .. (-.3,-.6) .. (.8,-.4) .. (.2,-.2) .. (-.2,-.3) .. (-.6,-.5);
			\fill[color=YellowOrange, opacity=.8, use Hobby shortcut,closed=true]
			(-.7,-.7) .. (-.5,-.9) .. (-.3,-.6) .. (.8,-.4) .. (.2,-.2) .. (-.2,-.3) .. (-.6,-.5);
			\path[color=YellowOrange, line width=.5pt, draw,use Hobby shortcut,closed=true]
			(1.3,-.7) .. (1.5,-.9) .. (1.7,-.6) .. (2.8,-.4) .. (2.2,-.2) .. (1.8,-.3) .. (1.4,-.5);
			\fill[color=white, opacity=0.8, use Hobby shortcut,closed=true]
			(1.3,-.7) .. (1.5,-.9) .. (1.7,-.6) .. (2.8,-.4) .. (2.2,-.2) .. (1.8,-.3) .. (1.4,-.5);
			\fill[color=white, opacity=0.8, pattern=north east lines,pattern color=YellowOrange, use Hobby shortcut,closed=true]
			(1.3,-.7) .. (1.5,-.9) .. (1.7,-.6) .. (2.8,-.4) .. (2.2,-.2) .. (1.8,-.3) .. (1.4,-.5);
		\end{scope}
					\end{axis}
	\end{tikzpicture}
	\caption{A graphical representation of possible functions $\widehat{f}$ satisfying \eqref{pgexformula}. We see that (i) $\widehat{f}$ is supported on the union of $E$ and the three translates of $E$ shown; (ii) $\widehat{f}\big([(1,0),\mathrm{id}]\cdot r \big) = \widehat{f}(r)$ for a.e. $r\in E$; and (iii) $\widehat{f}\big([(1,0),\sigma]\cdot r \big) = -\widehat{f}\big([(0,0),\sigma]\cdot r \big)$ for a.e. $r\in E$. Given the measurable range function $r\mapsto K_r$, Theorem \ref{BigInvTheorem2} says that the set $V$ of all $f\in L^2(\mathbb{R}^2)$ whose Plancherel transform satisfies properties (i)--(iii), is $\pi$-invariant.
	}
\end{figure}
Suppose that instead we take $\beta_1 =  \delta_{(0,0)} - \delta_{(1,0)} \in l^2(\mathbb{Z}^2)$, and let 
\[
H_r = \left\{ \begin{array}{ll}
	\text{span}\, \beta_1 &\text{if } r \in E\\
	\{0\} &\text{otherwise}
\end{array} \right.
\]
Then we would get the picture below:
	\begin{figure}[H]
	\pgfplotsset{
		axis lines=center, 
		grid,
		grid style={very thin, gray},
		xmin=-4,    xmax=4,   xtick={-4,...,4},
		x tick style={draw=none}, 
		ymin=-4,    ymax=4,   
		y tick style={draw=none},   
		x tick style={draw=none},  
		ytick={-4,...,4},
		restrict y to domain=-4:4, 
		ticklabel style={font=\scriptsize, fill=white, inner sep=2pt},
		domain=-1:11, samples=100,
		no marks, 
		every axis plot post/.append style={ultra thick,},
	}
	\centering
	\begin{tikzpicture}[scale=1, every node/.style={inner sep=0}]
		\begin{axis}[thick,
			xticklabels={-2,,-1,,0,,1,,2},
			yticklabels={-2,,-1,,0,,1,,2},
			xlabel={$x$},
			ylabel={$y$},
			every axis x label/.style={
				at={(ticklabel* cs:1.02)},
				anchor=west,
			},
			every axis y label/.style={
				at={(ticklabel* cs:1.02)},
				anchor=south,
			},
			]
			\draw[very thin, gray] (1, -1) -- (1,1);
			\draw[very thin, gray] (-1, -1) -- (-1,0);
			\draw[very thin, gray] (-1, 1) -- (1,1);
			\draw[very thin, gray] (-1, -1) -- (1,-1);
			\path[color=RoyalBlue, line width=.5pt, draw,use Hobby shortcut,closed=true]
			(-.7,-.7) .. (-.5,-.9) .. (-.3,-.6) .. (.8,-.4) .. (.2,-.2) .. (-.2,-.3) .. (-.6,-.5);
			\path[color=RoyalBlue, line width=.5pt, draw,use Hobby shortcut,closed=true]
			(-.7,-.7) .. (-.5,-.9) .. (-.3,-.6) .. (.8,-.4) .. (.2,-.2) .. (-.2,-.3) .. (-.6,-.5);
			\fill[color=RoyalBlue, opacity=.8, use Hobby shortcut,closed=true]
			(-.7,-.7) .. (-.5,-.9) .. (-.3,-.6) .. (.8,-.4) .. (.2,-.2) .. (-.2,-.3) .. (-.6,-.5);
			\node[color=black] at (2.3,-1.2) {\tiny $[(1,0),\mathrm{id}]\cdot E$};
			\node[color=black] at (2.3,1.2) {\tiny $[(1,0),\sigma]\cdot E$};	
			\node[color=black] at (-0.2,-1.2) {\tiny $E$};
			\node[color=black] at (-0.2,1.2) {\tiny $[(0,0),\sigma]\cdot E$};		
			\path[color=RoyalBlue, line width=.5pt, draw,use Hobby shortcut,closed=true]
			(-.7,-.7) .. (-.5,-.9) .. (-.3,-.6) .. (.8,-.4) .. (.2,-.2) .. (-.2,-.3) .. (-.6,-.5);
				\path[color=RoyalBlue, line width=.5pt, draw,use Hobby shortcut,closed=true]
			(1.3,-.7) .. (1.5,-.9) .. (1.7,-.6) .. (2.8,-.4) .. (2.2,-.2) .. (1.8,-.3) .. (1.4,-.5);
			\fill[color=white, opacity=0.8, use Hobby shortcut,closed=true]
			(1.3,-.7) .. (1.5,-.9) .. (1.7,-.6) .. (2.8,-.4) .. (2.2,-.2) .. (1.8,-.3) .. (1.4,-.5);
			\fill[color=white, opacity=0.8, pattern=north west lines,pattern color=RoyalBlue, use Hobby shortcut,closed=true]
			(1.3,-.7) .. (1.5,-.9) .. (1.7,-.6) .. (2.8,-.4) .. (2.2,-.2) .. (1.8,-.3) .. (1.4,-.5);
		\begin{scope}[yscale=-1,xscale=1, shift={($(0,-8)-(0,0)$)}]
		\path[color=YellowOrange, line width=.5pt, draw,use Hobby shortcut,closed=true]
					(-.7,-.7) .. (-.5,-.9) .. (-.3,-.6) .. (.8,-.4) .. (.2,-.2) .. (-.2,-.3) .. (-.6,-.5);
		\fill[color=YellowOrange, opacity=.8, use Hobby shortcut,closed=true]
					(-.7,-.7) .. (-.5,-.9) .. (-.3,-.6) .. (.8,-.4) .. (.2,-.2) .. (-.2,-.3) .. (-.6,-.5);
			\path[color=YellowOrange, line width=.5pt, draw,use Hobby shortcut,closed=true]
				(1.3,-.7) .. (1.5,-.9) .. (1.7,-.6) .. (2.8,-.4) .. (2.2,-.2) .. (1.8,-.3) .. (1.4,-.5);
	\fill[line width=.5pt, color=YellowOrange, opacity=.8,use Hobby shortcut,closed=true] (1.3,-.7) .. (1.5,-.9) .. (1.7,-.6) .. (2.8,-.4) .. (2.2,-.2) .. (1.8,-.3) .. (1.4,-.5);
			\end{scope}
		\end{axis}
	\end{tikzpicture}
	\caption{Possible functions $\widehat{f}$ with a different choice of range function $r\mapsto H_r$. With this new choice of range function, we see that (i) $\widehat{f}$ is supported on the union of $E$ and the three translates of $E$ shown; (ii) $\widehat{f}\big([(1,0),\sigma]\cdot r \big) = \widehat{f}\big([(0,0),\sigma]\cdot r \big)$ for a.e. $r\in E$; and (iii) $\widehat{f}\big([(1,0),\mathrm{id}]\cdot r \big) = -\widehat{f}(r)$ for a.e. $r\in E$. Given this new measurable range function $r\mapsto H_r$, Theorem \ref{BigInvTheorem2} says that the set $V$ of all $f\in L^2(\mathbb{R}^2)$ whose Plancherel transform satisfies properties (i)--(iii), is $\pi$-invariant.}
\end{figure}
 Comparing the pictures obtained from the two distinct range functions illustrates the one-to-one correspondence between measurable range functions and $\pi$-invariant subspaces.
	\end{ex}
		\begin{ex}
	We now consider a subspace of $L^2(\mathbb{R}^2)$ which is not invariant under shifts by symmetries of $pg$. Let $E$ be any bounded subset of $\mathbb{R}^2$, and let $V$ be the space of all functions supported on $E$. That is,
	\[
	V=\{f: f=0 \ \text{a.e.}\ \text{on}\ E^c\}.
	\]
Then $V$ is not $\pi$-invariant. For example, $f = \mathbbm{1}_E \in V$, but $\pi([y,\mathrm{id}])f = \mathbbm{1}_{[y,\mathrm{id}]\cdot E}\notin V$ for any $y\neq 0$. Thus, by our theorem, there exists no measurable range function $r \mapsto K_r$ for which $\widehat{f}$ satisfies equation \eqref{pgexformula}. 
		\end{ex}
		\begin{ex}
With $R$ being the fundamental domain for $\Gamma^*$, we let 
		\[
		\Phi(R) = \{\chi_y : y\in R\}.
		\]
 Let $V$ be the space of all functions whose Fourier transform is supported on $\Phi(R)$. That is, 
		\[
		V = \{f:\widehat{f}(\chi_y)=0\ \text{for a.e.}\ \chi_y\notin\Phi(R)\}.
		\]
Then $V$ is not $\pi$-invariant. This can be seen from the formula \ref{FT} for $\widehat{\pi}$.		
		\end{ex}
		\begin{ex}
			Let $\beta \in l^2(\mathbb{Z}^2)$ be the function 
			\[
			\beta(z_1,z_2) = \left\{ \begin{array}{ll}
				1 &\text{if } (z_1,z_2)= (0,0)\\
				0 &\text{otherwise}
			\end{array} \right.,
			\]
			and define $g \in  l^2(\Pi)\otimes l^2(\mathbb{Z}^2)$ to be $g = (\delta_{\mathrm{id}} + \delta_\sigma)\otimes \beta = \mathbbm{1}_\Pi \otimes \beta$. For each $r\in R$, let $W_r$ be the closed subspace of $l^2(\Pi)\otimes l^2(\mathbb{Z}^2)$ spanned by $g$. Thus $\{W_r\}_{r\in R}$ is a constant field of subspaces, and $r\mapsto W_r$ is a measurable range function from $R$ into the closed subspaces of $l^2(\Pi)\otimes l^2(\mathbb{Z}^2)$. Since $W_r\neq l^2(\Pi)\otimes K_r$ for any closed subspace $K_r\subset l^2(\mathbb{Z}^2)$, it follows that the range function $r\mapsto W_r$ does not give rise to a $pg$-invariant subspace of $L^2(\mathbb{R}^2)$ via the correspondence outlined in our theorem. That is, $U^{-1}\big(\int_R^\oplus W_r \ud r\big)$---where $U$ is the unitary defined in \ref{Udef}---is not $\pi$-invariant.
		\end{ex}
		\begin{ex}
Let 
\begin{equation}\label{R1R2}
R_1 = \left( -\tfrac12, 0 \right) \times \left(-\tfrac12, 0 \right)\ \ \text{and let}\ \ R_2 = \left( 0, \tfrac12 \right) \times \left(-\tfrac12, 0 \right).
\end{equation}
We illustrate these sets below.
	\begin{figure}[H]
	\pgfplotsset{
		axis lines=center, 
		grid,
		grid style={very thin, gray},
		xmin=-4,    xmax=4,   xtick={-4,...,4},
		x tick style={draw=none}, 
		ymin=-4,    ymax=4,   
		y tick style={draw=none},   
		x tick style={draw=none},  
		ytick={-4,...,4},
		restrict y to domain=-4:4, 
		ticklabel style={font=\scriptsize, fill=white, inner sep=2pt},
		domain=-1:11, samples=100,
		no marks, 
		every axis plot post/.append style={ultra thick,},
	}
	\centering
	\begin{tikzpicture}[scale=1]
		\begin{axis}[thick,
			xticklabels={-2,,-1,,0,,1,,2},
			yticklabels={-2,,-1,,0,,1,,2},
			xlabel={$x$},
			ylabel={$y$},
			every axis x label/.style={
				at={(ticklabel* cs:1.02)},
				anchor=west,
			},
			every axis y label/.style={
				at={(ticklabel* cs:1.02)},
				anchor=south,
			},
			]
	\fill[line width=1pt, color=RoyalBlue, opacity=.5] (-1, -1) -- (0,-1) -- (0,0) -- (-1,0) -- (-1,-1); 
	\fill[line width=1pt, color=YellowOrange, opacity=.5] (0, -1) -- (1,-1) -- (1,0) -- (0,0) -- (0,-1); 
			\draw[very thin, gray] (1, -1) -- (1,1);
			\draw[very thin, gray] (-1, -1) -- (-1,0);
			\draw[very thin, gray] (-1, 1) -- (1,1);
			\draw[very thin, gray] (-1, -1) -- (1,-1);
			\node[color=black] at (-.5,-.5) {\scriptsize $R_1$};	
			\node[color=black] at (.5,-.5) {\scriptsize $R_2$};	
		\end{axis}
	\end{tikzpicture}
	\caption{The two open sets $R_1$ and $R_2$ described above.}
\end{figure}
Let $F \in L^2\big(R;l^2(\Pi)\otimes l^2(\mathbb{Z}^2)\big)$ be the function
\[
F(r) = \left\{ \begin{array}{ll}
	\delta_{\mathrm{id}} \otimes \delta_{(0,0)}  &\text{if } r \in R_1\\
	-\delta_{\mathrm{id}} \otimes \delta_{(0,0)}  &\text{if } r \in R_2
\end{array} \right..
\]
Let $W = \langle F \rangle$ be the closed subspace spanned by $F$, and let
\[
G(r) = \delta_{\mathrm{id}} \otimes \delta_{(0,0)} \ \ \text{for all}\ r\in R.
\]
 Then 
\begin{align*}
	\langle F,G\rangle &= \int_R \langle F(r),G(r)\rangle\\
	&= \int_{R_1} \langle \delta_{\mathrm{id}} \otimes \delta_{(0,0)},\delta_{\mathrm{id}} \otimes \delta_{(0,0)} \rangle \ud r -\int_{R_2} \langle \delta_{\mathrm{id}} \otimes \delta_{(0,0)},\delta_{\mathrm{id}} \otimes \delta_{(0,0)} \rangle \ud r \\
	&= m(R_1) - m(R_2) \\
	&=0,
	\end{align*}
so that $G\perp W$. On the other hand, we have $\langle F(r),G(r)\rangle = 1$ for each $r \in R_1$. Thus $W$ violates the orthogonality condition of Proposition \ref{kappaOrth}. It follows that $V= U^{-1}(W)$ is not $\pi$-invariant. The orthogonality condition can be violated because the space $W$ is not of the form $\int_R^\oplus l^2(\Pi)\otimes K_r \ud r$ for any measurable range function $r \mapsto \{K_r\}$ taking values in the closed subspaces of $l^2(\mathbb{Z}^2)$.
		\end{ex}
	\chapter{The Central Decomposition}\label{CentralDecompSec}
	
	\section{The Structure of the Decomposition}
	
	Recall from Theorem \ref{centraldecomptheorem} that any unitary representation $(\pi,\mathscr{H})$ of a Type I, second countable, locally compact group $G$ on a separable Hilbert space $\mathscr{H}$ has a unique (up to a.e-equivalence) disintegration as a direct integral of factor representations. In this section we find a unitary equivalence between our decomposition $\tau = \int_{R}^\oplus U^r \otimes I \ud \mu(r)$ and a representation $\rho$ satisfying the properties listed in Theorem \ref{centraldecomptheorem}. We thus give the central decomposition of the natural crystal shift representation $\pi$. Recall from Theorem \ref{centraldecomptheorem} that if $\widehat{\Gamma}$ is given the Mackey Borel structure, then there is a standard measure on $\widehat{\Gamma}$, a measurable field of Hilbert spaces $\{\mathscr{H}_\zeta\}_{\zeta\in\widehat{\Gamma}}$, and a measurable field of (Type I) factor representations $\{\rho_\zeta\}_{\zeta\in \widehat{\Gamma}}$ on $\{\mathscr{H}_\zeta\}_{\zeta\in\widehat{\Gamma}}$ with the following properties:
	\begin{list}{(\arabic{ctr})}{\usecounter{ctr}}
		\setcounter{ctr}{0}
		\item\label{firstprop} we have $\rho_{\zeta}\sim \tilde{\zeta}\otimes I$ for some $\tilde{\zeta} \in \zeta$, for $\mu$-a.e. $\zeta \in \widehat{\Gamma}$.
		\item there is a unitary equivalence $U$ between $\pi$ and $\rho = \int_{\widehat{\Gamma}}^\oplus \rho_\zeta\ud\mu(\zeta)$. 
		\item $U$ transforms the center of $\pi(G)^-$ into the algebra $\mathscr{A}$ of diagonal operators on $\int_{\widehat{\Gamma}}^\oplus \mathscr{H}_\zeta\ud\mu(\zeta)$. Equivalently, by Lemmas \ref{closureofconjugate} and \ref{centerlemma}, $\mathscr{Z}(\rho(\Gamma)^-)= \mathscr{A}$.
		\item If $\mu^\prime$ and $\{\rho_\zeta^\prime\}_{\zeta\in \widehat{\Gamma}}$ satisfy the same properties, then $\mu$ is equivalent to $\mu^\prime$, and $\rho_\zeta$ is equivalent to $\rho_\zeta^\prime$ for $\mu$-a.e. $\zeta \in \widehat{\Gamma}$.
	\end{list}
By this theorem, it is possible to find a measurable field of Hilbert spaces $\{\mathscr{H}_\zeta\}_{\zeta\in\widehat{\Gamma}}$ and an explicit unitary $\tilde{\Psi}$ from $L^2\big(R;l^2(\Pi) \otimes l^2(\mathcal{L}^*)\big)$ onto $\int_{\widehat{\Gamma}}^\oplus \mathscr{H}_\zeta\ud\mu(\zeta)$ which intertwines $\tau$ with  $\rho = \int_{\widehat{\Gamma}}^\oplus \rho_\zeta\ud\mu(\zeta)$. We turn to finding this $\Psi$ and making the identification explicit.

	\section{The Subgroup Representation Topology}
 At this point it will be useful to discuss these \textit{subgroup representation pairs} and their topology. In what follows, $G$ will be a locally compact group.\label{CX1}
	\begin{defn}
Given a topological space $X$, we let $\mathcal{C}(X)$ denote the collection of all closed subsets of $X$. Let $\mathcal{K}(G)$ denote the set of all closed subgroups of $G$. 
	\end{defn}
\label{KG} Let $X$ be a locally compact Hausdorff space. There is a topology on $\mathcal{C}(X)$ called the \textit{compact-open (or Fell)} topology; this is described on p.~204 of \cite{TaylorKaniuth}. For a locally compact group $G$, we give $\mathcal{K}(G)$ the relative topology inherited from $\mathcal{C}(G)$. 
\begin{defn}\label{SG}
	Let $\mathcal{S}(G)$ be the set of all pairs $(K,\sigma)$ where $K\in \mathcal{K}(G)$, and $\sigma$ is a unitary representation of $K$. We call $\mathcal{S}(G)$ the set of \textit{subgroup representations} (sometimes called the set of \textit{subgroup-representation pairs}).
\end{defn}
There is a natural topology on $\mathcal{S}(G)$ called the \textit{subgroup representation topology}, or the \textit{inner hull kernel-topology}. This is defined on p.~429 of \cite{Fell} or pp.~221--222 \cite{TaylorKaniuth}.
\begin{defn}\label{Imap}
	We let $I$ denote the inducing map, which send an element $(K,\sigma)$ of $\mathcal{S}(G)$ to the equivalence class of the induced representation $\mathrm{ind}_K^G \sigma$.
\end{defn}
\begin{rem}\label{equivonrepparis}
	As discussed in Baggett (\cite{Baggett}, p.~175), there is an equivalence relation $\sim$ on $\mathcal{S}(G)$ such that $I(K,\sigma) = I(K^\prime,\sigma^\prime)$ if and only if $(K,\sigma) \sim (K^\prime,\sigma^\prime)$, so that the inducing map is well-defined. This equivalence is given by:  $(K,\sigma) \sim (K^\prime,\sigma^\prime)$ if and only if $K^\prime = gKg^{-1}$ for some $g\in G$ and $\sigma^\prime \sim g\cdot \sigma$, where $g\cdot\sigma (x) = \sigma(g^{-1}xg)$.
\end{rem}
\begin{defn}\label{AG}
	Let $\mathscr{A}(G)$ denote the set of all pairs $(K, \sigma)$ where $K \in \mathcal{K}(G)$ and $\sigma$ is an irreducible unitary representation of $K$.
\end{defn}
\begin{defn}\label{S0}
	Let $N$ be a nontrivial abelian closed normal subgroup of $G$. Define $\mathcal{S}_0(G)\subset \mathscr{A}(G)$ be the set of subgroup representation pairs $(K,\sigma)$ such that $K = G_\chi$ is the stabilizer for some element $\chi \in \widehat{N}$ and $\sigma$ is an irreducible unitary representation of $K$ whose restriction to $N$ is a multiple of $\chi$. We call $\mathcal{S}_0(G)$ the set of \textit{cataloguing pairs} for $\widehat{G}$.
\end{defn}
\begin{defn}\label{tildeSOG}
	Let $\widetilde{\mathcal{S}}_0(G) = \mathcal{S}_0(G)/\sim$, where $\sim$ is the equivalence relation described in Remark \ref{equivonrepparis}.
\end{defn}
Lemma \ref{InducedIrreducible} may now be restated as follows.
\begin{lemma}
	The induced representation $\mathrm{ind}_K^G \sigma$ is irreducible for each $(K,\sigma) \in \mathcal{S}_0(G)$.
	\end{lemma}
	It follows that $I(K,\sigma) \in \widehat{G}$ for each $(K,\sigma) \in \mathcal{S}_0(G)$.	The theorem which follows combines Theorem 6.39 of \cite{FollandAbstract}, the paragraph following this theorem, and Theorem 6.40 of the same.
	\begin{theorem}\label{BigInducedBij}
Let $G$ be a locally compact group with a nontrivial closed abelian normal subgroup $N$. Suppose further that $N$ satisfies the regularity condition described in \cite{FollandAbstract}, p.~196. Then the inducing map $I$ is a surjection from $\mathcal{S}_0(G)$ onto $\widehat{G}$. Moreover, $I$ becomes a bijection by identifying elements which are equivalent under the equivalence relation described in Remark \ref{equivonrepparis}.
	\end{theorem}
Lawrence Baggett's paper \cite{Baggett} builds on the work of Fell (\cite{Fell}) in describing the topology of $\widehat{G}$ in terms of cataloguing pairs $(K,\sigma) \in \mathcal{S}_0(G)$.\footnote{Note that the author is using different notation than Baggett; in particular, Baggett uses $Y$ for what we call $\mathcal{S}_0(G)$.} In particular, Conjecture 2 (p.~176) of Baggett's paper is valid for crystallographic groups, by a result of Raeburn (\cite{Raeburn}). We won't need this powerful result. We will, however, make use of the following proposition from p.~183 of \cite{Baggett}, which rests on a powerful result of Fell. The proof is outside the scope of this thesis, so we merely point the reader to the relevant theorems. The \textit{support} of a representation is described on p.~38 of \cite{TaylorKaniuth}.
\begin{prop}
	Let $(K_i, \sigma_i)_{i\in \mathscr{I}}$ be a net of elements of $\mathscr{A}(G)$ which converges to an element $(K, \sigma)\in \mathscr{A}(G)$. Suppose further that for each $i\in\mathscr{I}$, the representation $\mathrm{ind}_{K_i}^G \sigma_i$ is
	irreducible. Then the net $(I(K_i,\sigma_i))_{i\in\mathscr{I}}$ converges in $\widehat{G}$ to each $\rho \in \widehat{G}$ in the support of $\mathrm{ind}_K^G \sigma$.
\end{prop}
\begin{proof}
	This follows from Lemma 5.7, Proposition 5.11, and Corollary 5.40 in \cite{TaylorKaniuth}.
\end{proof}
\begin{cor}\label{Icont}
	Let $I$ denote the inducing map defined in \ref{Imap}. Then $I\vert_{\mathscr{A}(G)}$ is a continuous map from $\mathscr{A}(G)$ into $\widehat{G}$. 
\end{cor}
	\section[The Structure of \texorpdfstring{$\widehat{\Gamma}$}{Gamma-Hat}]{The Structure of \texorpdfstring{$\boldsymbol{\widehat{\Gamma}}$}{Gamma-Hat}}
	Recall from Theorem \ref{bigtheoremapp} that
\begin{equation}\label{GammaHatParam1}
	\widehat{\Gamma}=\bigcup_{\omega\in\Omega} \big\{[\mathrm{ind}_{\Gamma_{\omega}}^\Gamma \sigma]: [\sigma] \in \widehat{\Gamma_{\omega}}^\omega\big\}
\end{equation}
where $\widehat{\Gamma_{\omega}}^\omega$ denotes the equivalence classes of irreducible representations of $\Gamma_{\omega}$ which restrict to a multiple of $\chi_\omega^{\mathrm{T}}$ on $\mathrm{T}$. Here, $[\sigma] \in \widehat{\Gamma_{\omega}}^\omega$ is the equivalence class of an irreducible representation $\sigma$ of $\Gamma_\omega$ on a Hilbert space $\mathscr{H}(\sigma)$. Each $\mathrm{ind}_{\Gamma_{\omega}}^\Gamma \sigma$ is a unitary representation of $\Gamma$ on $l^2\big(\Gamma/\Gamma_\omega ; \mathscr{H}(\sigma)\big)$, as discussed in \cite{TaylorKaniuth}, \S 2.1. Moreover, for each $\omega \in \Omega$, we have that the quotient space $\Gamma/\Gamma_\omega$ is in bijection with the orbit $\Gamma(\chi_\omega^\mathrm{T})$, since $\Gamma_\omega$ is the stabilizer of $\chi_\omega^\mathrm{T}$ in $\Gamma$. Thus we may regard $\mathrm{ind}_{\Gamma_{\omega}}^\Gamma \sigma$ as a unitary representation of $\Gamma$ on $l^2\big(\Gamma(\chi_\omega^\mathrm{T}); \mathscr{H}(\sigma)\big)$. Equation \eqref{GammaHatParam1} says that each $\zeta \in \widehat{\Gamma}$ comes from an induced representation of an irreducible unitary representation $\sigma$ of a (closed) subgroup $\Gamma_\omega$ of $\Gamma$.
\begin{rem}
	Let $\Gamma$ be a crystallographic group. Then $\mathcal{S}_0(\Gamma)\subset \mathcal{S}(\Gamma)$ is the set of subgroup representation pairs of the form $(\Gamma_\nu, \sigma)$, where $\Gamma_\nu$ is the stabilizer of $\chi_\nu^\mathrm{T}$ in $\Gamma$ for $\nu \in \Omega_{\mathrm{T}^*}$, and $\sigma$ is an irreducible unitary representation of $\Gamma_\nu$ whose restriction to $\mathrm{T}$ is a multiple of $\chi_\nu^\mathrm{T}$. For each $\nu \in \Omega_{\mathrm{T}^*}$, let \label{Snu}
	\[
	S(\nu) = \Bigg\{\ \sigma:	\begin{aligned}
 \quad \sigma\ &\text{is an irreducible unitary representation of }\Gamma_\nu\\ 
&\text{whose restriction to }\mathrm{T}\ \text{is a multiple of } \chi_\nu^\mathrm{T}
	\end{aligned}\Bigg\}.
	\]
Then
	\[
	\mathcal{S}_0(\Gamma) = \{(\Gamma_\nu, \sigma):  \sigma \in S(\nu), \nu \in \Omega_{\mathrm{T}^*}\}.
	\]
\end{rem}
\begin{prop}\label{S0lemma}
		\[
	\mathcal{S}_0(\Gamma) = \{([0,L]\Gamma_\omega[0,L]^{-1}, [0,L]\cdot \rho) \in \mathcal{S}(\Gamma): \rho \in S(\omega), L\in \Pi, \omega \in \Omega\},
	\]
	where $[0,L]\cdot \sigma$ denotes the conjugation action defined in Remark \ref{equivonrepparis}.
	\end{prop}
	\begin{proof}
We first note that, for each $\omega \in \Omega$ and $L\in \Pi$, we have
\begin{equation}\label{actiononStab}
\Gamma_{L\cdot_{\idref{12}} \omega} = [0,L]\Gamma_\omega [0,L]^{-1}.
\end{equation}
Indeed, by Proposition \ref{actionsrel} we have
	\begin{align*}
		\Gamma_{L\cdot_{\idref{12}} \omega} &= \{[y,M] \in \Gamma: M\cdot_{\idref{6}} \chi_{L\cdot_{\idref{12}} \omega}^\mathrm{T} = \chi_{L\cdot_{\idref{12}} \omega}^\mathrm{T}\}\\
		&= \{[y,M] \in \Gamma:  M\cdot_{\idref{12}}(L\cdot_{\idref{12}} \omega) = L\cdot_{\idref{12}} \omega\}\\
		&= \{[y,M] \in \Gamma: (L^{-1}ML)\cdot_{\idref{12}} \omega= \omega\}\\
		&=\{[y,M] \in \Gamma: ([0,L]^{-1}[y,M][0,L])\cdot_{\idref{5}} \chi_\omega^\mathrm{T} = \chi_\omega^\mathrm{T}\}\\
		&=\{[0,L][y,M][0,L]^{-1} \in \Gamma: [y,M]\cdot_{\idref{5}} \chi_\omega^\mathrm{T} = \chi_\omega^\mathrm{T}\}\\
		&= [0,L]\Gamma_\omega [0,L]^{-1}.
	\end{align*}
We also show that
\begin{equation}\label{actiononStab1}
S(L\cdot_{\idref{12}} \omega) = \{[0,L] \cdot \rho : \rho \in S(\omega)\}.
	\end{equation}
To see this, let $\sigma \in S(L\cdot_{\idref{12}} \omega)$. Then since $\sigma$ is an irreducible representation of $\Gamma_{L\cdot_{\idref{12}} \omega}$, it is immediate from \eqref{actiononStab} that $\rho := [0,L]^{-1}\cdot \sigma$ is an irreducible representation of $\Gamma_\omega$. Moreover, since $\sigma\vert_{\mathrm{T}} = n\chi_{L\cdot_{\idref{12}} \omega}^\mathrm{T}$, we have, for $[y,\mathrm{id}]\in \mathrm{T}$, 
\begin{align*}
	\rho([y,\mathrm{id}]) &= ([0,L]^{-1}\cdot \sigma)([y,\mathrm{id}])\\
	&= \sigma([0.L]\cdot_{\idref{3}}[y,\mathrm{id}])\\
	&= n\chi_{L\cdot_{\idref{12}} \omega}^\mathrm{T}([0.L]\cdot_{\idref{3}}[y,\mathrm{id}])\\
	&= [n([0,L]^{-1}\cdot_{\idref{5}}\chi_{L\cdot_{\idref{12}}\omega}^\mathrm{T})]([y,\mathrm{id}])\\
	&= n \chi_\omega^\mathrm{T}([y,\mathrm{id}]),
\end{align*}
where we have invoked Proposition \ref{actionsrel} again. Thus, $\rho\vert_{\mathrm{T}} = n \chi_\omega^\mathrm{T}$, and hence $\sigma = [0,L]\cdot \rho$, with $\rho \in S(\omega)$. For the reverse inclusion, let $\sigma  = [0,L]\cdot \rho$, with $\rho \in S(\omega)$; then for $[y,\mathrm{id}]\in \mathrm{T}$, we have
\begin{align*}
	\sigma([y,\mathrm{id}])&= ([0,L]\cdot \rho)([y,\mathrm{id}])\\
	&= \rho([0,L]^{-1}\cdot_{\idref{3}}[y,\mathrm{id}])\\
	&= (n\chi_\omega^\mathrm{T})([0,L]^{-1}\cdot_{\idref{3}}[y,\mathrm{id}])\\
	&= [n([0,L]\cdot_{\idref{5}}\chi_\omega^\mathrm{T})]([y,\mathrm{id}])\\
	&=  n\chi_{L\cdot_{\idref{12}} \omega}^\mathrm{T}([y,\mathrm{id}]),
	\end{align*}
whence $\sigma\vert_\mathrm{T} =  n\chi_{L\cdot_{\idref{12}}\omega}^\mathrm{T}$. Therefore, $[0,L] \cdot \rho \in S(L\cdot_{\idref{12}}\omega)$ whenever $\rho \in S(\omega)$. Hence we have established \eqref{actiononStab1}.

Now, since $\Omega$ is a transversal for the action of $\Pi$ on $\Omega_{\mathrm{T}^*}$, it follows from  \eqref{actiononStab} and  \eqref{actiononStab1} that
\begin{align*}
\mathcal{S}_0(\Gamma) &= \{(\Gamma_\nu, \sigma):  \sigma \in S(\nu), \nu \in \Omega_{\mathrm{T}^*}\} \\
 &= 	\{(\Gamma_{L \cdot_{\idref{12}} \omega}, \sigma) \in \mathcal{S}(\Gamma):  \sigma \in S(L\cdot_{\idref{12}} \omega), L\in \Pi, \omega \in \Omega\} \\
&= 	\{([0,L]\Gamma_\omega[0,L]^{-1}, [0,L]\cdot \rho) \in \mathcal{S}(\Gamma): \rho \in S(\omega),  L\in \Pi,\omega \in \Omega\},
\end{align*}
where $[0,L]\cdot \rho$ denotes the conjugation action defined in Remark \ref{equivonrepparis}.
	\end{proof}
 It follows from Theorem \ref{BigInducedBij} that the elements of $\widehat{\Gamma}$ may be catalogued by the subgroup representation pairs $(\Gamma_\nu, \sigma) \in \mathcal{S}_0(\Gamma)$:
\begin{equation}\label{Gammahatnewparam}
	\widehat{\Gamma}= \{I(\Gamma_\nu, \sigma): [(\Gamma_\nu, \sigma)] \in \widetilde{\mathcal{S}}_0(\Gamma)\}.
\end{equation}
Moreover, it follows from the preceding proposition that
\begin{equation}\label{GammaHatParam2}
	\widehat{\Gamma}= \{I(\Gamma_\omega, \sigma): \sigma \in S(\omega), \omega \in \Omega\},
\end{equation}
which recaptures \eqref{GammaHatParam1}. Recall that
\[
\Omega_0 := \{\omega\in \Omega : \Pi_\omega = \{\mathrm{id}\}\}.
\]
We saw in \S \ref{inducedRep} that $\Gamma_{\omega}=\mathrm{T}$ if and only if $\omega \in \Omega_0$. From this we arrived at the decomposition 
\begin{equation}\label{paramGammahat1}
	\widehat{\Gamma}=  \big\{[U^\omega]: \omega\in \Omega_0 \big\} \cup   \bigcup_{\omega\in\Omega - \Omega_0} \big\{[\mathrm{ind}_{\Gamma_{\omega}}^\Gamma \sigma]: [\sigma] \in \widehat{\Gamma_{\omega}}^\omega\big\},
\end{equation}
Recall as well that $R\subset \Omega_0$, with $m(\Omega_0-R) =0$. We define an embedding $\iota: R \hookrightarrow \mathcal{S}_0(\Gamma)$ by 
\begin{equation}\label{iota}
\boxed{	\iota(r) = (\mathrm{T},\chi_r^\mathrm{T})}\,.
\end{equation}
\begin{lemma}\label{imageofOmega0}
	The image of $R$ under $\iota$ meets each equivalence class in $\mathcal{S}_0$ (for the equivalence relation $\sim$ described in Remark \ref{equivonrepparis}) at most once.
\end{lemma}
\begin{proof}
	Let $(\Gamma_{r_1},\chi_{r_1}^\mathrm{T})$ and $(\Gamma_{r_2},\chi_{r_2}^\mathrm{T})$ be elements of $\iota(R)$, with $(\Gamma_{r_1},\chi_{r_1}^\mathrm{T})\sim (\Gamma_{r_2},\chi_{r_2}^\mathrm{T})$. Then $\chi_{r_2}^\mathrm{T} \sim [0,L]\cdot \chi_{r_1}^\mathrm{T}$, for some $L\in\Pi$. For any $[y,\mathrm{id}]\in\mathrm{T}$, we have
	\[
	([0,L]\cdot \chi_{r_1}^\mathrm{T})([y,\mathrm{id}]) = \chi_{r_1}^\mathrm{T}([0,L]^{-1}\cdot_{\idref{3}}[y,\mathrm{id}]) = ([0,L]\cdot_{\idref{5}}\chi_r^\mathrm{T})([y,\mathrm{id}]) =\chi_{L\cdot_{\idref{12}} r_1}^\mathrm{T}([y,\mathrm{id}]);
	\]
	therefore,
	\[
	\chi_{r_2}^\mathrm{T} \sim \chi_{L\cdot_{\idref{12}} r_1}^\mathrm{T},
	\]
	and $\chi_{r_2}^\mathrm{T}=\chi_{L\cdot_{\idref{12}} r_1}^\mathrm{T}$ as elements of $\widehat{\mathrm{T}}$. Then by \eqref{ParamOfThat}, we must have $L\cdot_{\idref{12}} r_1=r_2 \in \Omega$; and by the definition of $\Omega$, it follows that $L\cdot_{\idref{12}} r_1=r_1$. Thus $r_1 = r_2$. 
\end{proof}
\begin{defn}\label{Psi}
We define $\Psi: R \to \widehat{\Gamma}$ to be the composition $I\circ \iota$, so that
\[
\Psi(r) = [U^r].
\]
\end{defn}
We set \label{Gamma0hat}
\[
\boxed{\widehat{\Gamma}_0 := \Psi(R)}\,.
\] 
We show that $\Psi$ is a Borel isomorphism onto its image $\widehat{\Gamma}_0$. For this we require some preliminary lemmas.
\begin{defn}\label{EX}
Let $X$ be a locally compact Hausdorff space. We define $\mathcal{E}(X)$ to be the set of all complex-valued functions $f$ defined on a closed subset $D(f)$ of $X$.
\end{defn}
\begin{rem}
Identifying functions $f\in \mathcal{E}(X)$ with their graphs, $\mathcal{E}(X)$ inherits a topology called the \textit{semicompact-open topology}---see \cite{TaylorKaniuth}, p.~224 for details. 
\end{rem}
\begin{defn}\label{EsG}
For a locally compact group $G$, we define $\mathcal{E}_s(G)$ to be the topological subspace of $\mathcal{E}(G)$ consisting of all $f\in \mathcal{E}(G)$ such that $D(f)$ is a closed subgroup of $G$.
	\end{defn}
The following is Lemma 5.30 in \cite{TaylorKaniuth}.
\begin{lemma}
Let $X$ be a locally compact Hausdorff space, and let $Y$ be a closed subset of $X$. The semicompact-open topology relativized to $\{f\in \mathcal{E}(X):D(f)=Y\}$ coincides with the topology of uniform convergence on compact subsets of $Y$.
\end{lemma}
\begin{cor}\label{topEsH}
	Let $G$ be a discrete group and let $H$ be a subgroup of $G$. Since any compact subset of a discrete group is finite, it follows from the preceding lemma that the semicompact-open topology relativized to  $\{f\in \mathcal{E}_s(G):D(f)=H\}$ coincides with the topology of pointwise convergence on $H$.
	\end{cor}
\begin{defn}
An element $\phi \in \mathcal{E}_s(G)$ is said to be of \textit{positive type associated with} $(K,\sigma) \in \mathcal{S}(G)$ if $D(\phi)=K$ and $\phi(x) = \langle \sigma(x)\xi,\xi\rangle$ for all $x\in K$, for some $\xi \in \mathscr{H}_\sigma$ (\textit{cf.} \S \ref{postype}). 
\end{defn}
The following (Theorem 5.36 in \cite{TaylorKaniuth}) gives us a useful characterization of convergence in $\mathcal{S}(G)$. 
\begin{theorem}\label{convergenceequiv}
	Let $(K_i,\sigma_i)_{i\in \mathscr{I}}$ be net in $\mathcal{S}(G)$ and let $(K,\sigma)\in \mathcal{S}(G)$. The following are equivalent.
	\begin{list}{(\roman{ctr1})}{\usecounter{ctr1}}
	\item $(K_i,\sigma_i) \to (K,\sigma)$ in $\mathcal{S}(G)$.
	\item For each finite sequence $\phi_1,\ldots,\phi_n$ of functions of positive type associated with $(K,\sigma)$ and each subnet $(K_{i_j}, \sigma_{i_j})_{j\in \mathscr{J}}$ of $(K_i,\sigma_i)_i$, there exists a further subnet $(K_{i_{j_k}}, \sigma_{i_{j_k}})_{k\in\mathscr{K}}$ such that, for each $k\in \mathscr{K}$ and each $1\leq \ell \leq n$, there exists a finite sum $S(\ell,k)$ of functions of positive type associated with $(K_{i_{j_k}}, \sigma_{i_{j_k}})$ such that the net $\big(S(\ell,k)\big)_{k\in \mathscr{K}}$ converges to $\phi_\ell$ in $\mathcal{E}_s(G)$.
	\end{list}
\end{theorem}
The following is Corollary 15.2 in \cite{Kechris}.
\begin{lemma}\label{BorelIsoLemma}
	Let $X$ and $Y$ be standard Borel spaces and let $f:X \to Y$ be Borel. 
	If $A \subset X$ is Borel and $f\vert_A$ is injective, then $f(A)$ is Borel and $f$ is a Borel 
	isomorphism of $A$ with $f(A)$.
\end{lemma}
\begin{lemma}\label{iotaBorelIso}
	Let $\Gamma$ be a crystal group and let $\mathcal{S}_0(\Gamma)$ have the relative topology inherited from the subgroup representation topology on $\mathcal{S}(\Gamma)$. Then $\iota$ defined in \eqref{iota} is a Borel isomorphism from $R$ onto its image in $\mathcal{S}_0(\Gamma)$.
\end{lemma}
\begin{proof} We first show that $\iota$ is injective, so that it is a bijection onto its image. If $\iota(r_1) = \iota(r_2)$, then $\chi_{r_1}^\mathrm{T}= \chi_{r_2}^\mathrm{T}$, which by our parametrization result, \eqref{ParamOfThat}, means that $r_1 = r_2$. Next we show that $\iota$ is continuous. Suppose that the net $\{r_i\}_{i\in \mathscr{I}}$ converges to $r$ in $R$, then $\iota(r_i) = (\mathrm{T},\chi_{r_i}^\mathrm{T})$ for all $i\in \mathscr{I}$, and $\iota(r) = (\mathrm{T},\chi_{r}^\mathrm{T})$. Note that any function $\phi$ of positive type associated with $(\mathrm{T},\chi_{r}^\mathrm{T})$ is given by
	\begin{equation}\label{postypecharT}
	\phi_{\xi,\xi}([y,\mathrm{id}]) = \langle \chi_{r}^\mathrm{T}([y,\mathrm{id}]) \xi,\xi\rangle = \chi_{r}^\mathrm{T}([y,\mathrm{id}])|\xi|^2\ \ \text{for}\ [y,\mathrm{id}]\in\mathrm{T}.
	\end{equation}
	Thus any function of positive type associated with $(\mathrm{T},\chi_{r}^\mathrm{T})$ is just $|\xi|^2\chi_{r}^\mathrm{T}$ for some $\xi \in \mathbb{C}$. Let $|\xi_1|^2 \chi_{r}^\mathrm{T}, \dots, |\xi_n|^2 \chi_{r}^\mathrm{T}$ be functions of positive type associated with $(\mathrm{T},\chi_{r}^\mathrm{T})$, where $\xi_\ell \in \mathbb{C}$ for $\ell=1,\dots,n$, and let $(\mathrm{T},\chi_{r_{i_j}}^\mathrm{T})_{j\in\mathscr{J}}$ be a subnet of $(\mathrm{T},\chi_{r_i}^\mathrm{T})_{i\in\mathscr{I}}$. Then for each $j\in \mathscr{J}$ and each $\ell=1,\dots,n$, the function $\psi(\ell,j)$, given by $\psi(\ell,j)([y,\mathrm{id}]) = \langle \chi_{r_j}^\mathrm{T}([y,\mathrm{id}])\xi_\ell, \xi_\ell\rangle$, is a function of positive type associated with $(\mathrm{T},\chi_{r_j}^\mathrm{T})$. Moreover, $\psi(\ell,j)$ converges to $|\xi_\ell|^2 \chi_{r}^\mathrm{T}$ in $\mathcal{E}_s(\Gamma)$, for each $\ell=1,\dots,n$. Indeed, since $r_j \to r$, it follows by the continuity of $\upsilon$\footnote{Recall that $\upsilon = \mathcal{C}\circ q\vert_{\Omega_{\mathrm{T}^*}}$, and each of these maps is continuous.} that $\chi_{r_j}^\mathrm{T} \to \chi_{r}^\mathrm{T}$ in $\widehat{\mathrm{T}}$; in particular, 
	\[
	\chi_{r_j}^\mathrm{T}([y,\mathrm{id}]) \to 	\chi_{r}^\mathrm{T}([y,\mathrm{id}])\ \ \text{for each}\ [y,\mathrm{id}] \in \mathrm{T},
	\]
	whence $\psi(\ell,j)$ converges to $|\xi_\ell|^2 \chi_{r}^\mathrm{T}$ pointwise, and hence in $\mathcal{E}_s(\Gamma)$, by Corollary \ref{topEsH}. It follows by Theorem \ref{convergenceequiv} that $(\mathrm{T},\chi_{r_i}^\mathrm{T}) \to (\mathrm{T},\chi_{r}^\mathrm{T})$ in $\mathcal{S}_0(\Gamma)$, i.e., $\iota$ is continuous.
	
It follows now from Lemma \ref{BorelIsoLemma} that $\iota(R)$ is Borel in $\mathcal{S}_0(\Gamma)$, and that $\iota$ is a Borel isomorphism onto its image.
\end{proof}
\begin{lemma}
	$\Psi$ is a bijection from $R$ onto $\widehat{\Gamma}_0$.
\end{lemma} 
\begin{proof}
Since $\widehat{\Gamma}_0$ was defined to be the image of $R$ under $\Psi$, it suffices to show that $\Psi$ is injective. Recall that $\Psi = I \circ \iota$, and $\iota$ was shown to be injective. Thus it suffices to show that $I\vert_{\iota(R)}$ is injective. But this follows immediately from Theorem \ref{BigInducedBij}, since $\iota(R)$ contains at most one element from each equivalence class in $\mathcal{S}_0$, by Lemma \ref{imageofOmega0}.
\end{proof}
\begin{prop}
 Let $\Gamma$ be a crystallographic group. Then $\widehat{\Gamma}_0$ is Borel in $\widehat{\Gamma}$, and $\Psi: R \to \widehat{\Gamma}_0$ is a Borel isomorphism.
\end{prop}
\begin{proof}
Since $R$ is Borel, it is, in particular, a standard measure space. Moreover, $\widehat{\Gamma}$ is standard by Theorem 7.6 in \cite{FollandAbstract}, since $\Gamma$ is type I. Thus by Lemma \ref{BorelIsoLemma} (with $R = A = X)$, it suffices to show that $\Psi$ is a Borel map from $R$ into $\widehat{\Gamma}$. This follows immediately from Lemmas \ref{iotaBorelIso} and \ref{Icont}.
\end{proof}
Since $\iota$ and $I\vert_{\mathcal{S}_0(\Gamma)}$ are continuous, $\Psi$ is in fact continuous. Since $\Psi$ is Borel measurable, we can give $\widehat{\Gamma}$ the image measure $\mu := m\circ \Psi^{-1}$.
\begin{cor}\label{dualspacedescription}
	We have
	\begin{equation}\label{dualspacedescriptioneqn}
		\widehat{\Gamma} = \{[U^{\Psi^{-1}(\zeta)}]: \zeta \in \widehat{\Gamma}_0\} \cup (\widehat{\Gamma}-\widehat{\Gamma}_0),
	\end{equation}
	with $\mu(\widehat{\Gamma}-\widehat{\Gamma}_0) = 0$.
\end{cor}
\begin{proof}
	We have
	\begin{align*}
		\widehat{\Gamma} &= \widehat{\Gamma}_0 \cup (\widehat{\Gamma}-\widehat{\Gamma}_0)\\
		&= \Psi(R) \cup (\widehat{\Gamma}-\widehat{\Gamma}_0)\\
		&= 	\big\{[U^r]: r\in R \big\}\cup (\widehat{\Gamma}-\widehat{\Gamma}_0)\\
		&= 	\{[U^{\Psi^{-1}(\zeta)}]: \zeta \in \widehat{\Gamma}_0\}\cup (\widehat{\Gamma}-\widehat{\Gamma}_0);
	\end{align*}
	and
	\[
	\mu(\widehat{\Gamma}- \widehat{\Gamma}_0) = m\big(\Psi^{-1}(\widehat{\Gamma}- \widehat{\Gamma}_0)\big) = m(\emptyset) = 0.
	\]	
\end{proof}
\begin{cor}
	The measure $\mu  =m\circ \Psi^{-1}$ is standard.
	\end{cor}

\section[Realizing the decomposition]{Realizing the decomposition}
Let $\Psi$ be as defined in \ref{Psi}. We define a mapping \label{tildePsi}
\[
\widetilde{\Psi}: L^2\big(R;l^2(\Pi)\otimes l^2(\mathcal{L}^*)\big) \to L^2\big(\widehat{\Gamma};l^2(\Pi)\otimes l^2(\mathcal{L}^*)\big)
\]
by
\[
(\widetilde{\Psi}F)(\zeta) = \left\{ \begin{array}{ll}
	(F\circ \Psi^{-1})(\zeta) &\text{if } \zeta \in \widehat{\Gamma}_0 \\
	0  &\text{if } \zeta \in \widehat{\Gamma}-\widehat{\Gamma}_0.
\end{array} \right.
\]
To see that this is well-defined, we must show that the vector-valued function $\widetilde{\Psi}F$ is weakly measurable for each $F \in L^2\big(R;l^2(\Pi)\otimes l^2(\mathcal{L}^*)\big)$. That is, we must show that 
\begin{equation}\label{psitildemeas}
	\zeta \mapsto \langle (\widetilde{\Psi}F)(\zeta) ,g\rangle
\end{equation}
is measurable for each $g \in l^2(\Pi)\otimes l^2(\mathcal{L}^*)$. First, define $F_g: R \to \mathbb{C}$ by $F_g(r) = \langle F(r),g\rangle$; then $F_g$ is measurable since $F$ is weakly measurable. Now we define $\widetilde{\Psi}_g: \widehat{\Gamma} \to \mathbb{C}$ by 
\[
\widetilde{\Psi}_g(\zeta) = \left\{ \begin{array}{ll}
	(F_g \circ \Psi^{-1})(\zeta) &\text{if } \zeta \in \widehat{\Gamma}_0\\
	0  &\text{if } \zeta \in \widehat{\Gamma}-\widehat{\Gamma}_0. 
\end{array}
\right.
\]
Since the mapping \eqref{psitildemeas} is just $\widetilde{\Psi}_g$, showing that $\widetilde{\Psi}F$ is weakly measurable amounts to showing that $\widetilde{\Psi}_g$ is measurable for each $g \in l^2(\Pi)\otimes l^2(\mathcal{L}^*)$. Being a composition of Borel functions, $F_g \circ \Psi^{-1}$ is Borel-measurable on $\widehat{\Gamma}_0$. Similarly, the constant $0$ function is Borel measurable on $\widehat{\Gamma}-\widehat{\Gamma}_0$. Therefore, by Lemma \ref{measonAB}, $\widetilde{\Psi}_g$ is measurable for each $g$, and hence $\widetilde{\Psi}F$ is weakly measurable.

We show presently that $\widetilde{\Psi}$ is an isometric isomorphism. We first note that $\widetilde{\Psi}$ is clearly linear; moreover, $\widetilde{\Psi}$ is an isometry, since 
\begin{align*}
	\|\widetilde{\Psi}F\|^2&= \int_{\widehat{\Gamma}} \|(\widetilde{\Psi}F)(\zeta)\|^2 \ud\mu(\zeta)\\
	&= \int_{\widehat{\Gamma}_0} \|(\widetilde{\Psi}F)(\zeta)\|^2 \ud\mu(\zeta)\\
	&= \int_{\widehat{\Gamma}_0} \|(\widetilde{\Psi}F)(\zeta)\|^2 \ud(m\circ \Psi^{-1})(\zeta)\\
	&= \int_{R} \|(\widetilde{\Psi}F)(\Psi(r))\|^2 \ud r\\
	&= \int_{R} \|F(r)\|^2 \ud r\\
	&= \|F\|^2.
\end{align*}
Moreover, $\widetilde{\Psi}$ is surjective: given $\varrho\in L^2\big(\widehat{\Gamma};l^2(\Pi)\otimes l^2(\mathcal{L}^*)\big)$, if we let $F = \varrho \circ \Psi$, then $\widetilde{\Psi}F =  F\circ \Psi^{-1} = \varrho$ and $F \in  L^2\big(R;l^2(\Pi)\otimes l^2(\mathcal{L}^*)\big)$, since $\|F\| = \|\widetilde{\Psi}F\| = \|\varrho\| <\infty$. Thus $\widetilde{\Psi}$ is an isometric isomorphism.

We of course have that
\[
L^2\big(\widehat{\Gamma};l^2(\Pi) \otimes l^2(\mathcal{L}^*)\big) = \int_{\widehat{\Gamma}}^\oplus \mathscr{H}_\zeta\ud \mu(\zeta),
\]
where $\mathscr{H}_\zeta \equiv l^2(\Pi) \otimes l^2(\mathcal{L}^*)$ is the constant field of Hilbert spaces. We now define a measurable field of operators $\{\rho_\zeta [x,L]\}_{\zeta \in \widehat{\Gamma}}$ by
\begin{equation}\label{rhozetaops}
	\rho_\zeta [x,L] = \left\{\begin{array}{ll} 
		U^{\Psi^{-1}(\zeta)}[x,L]\otimes I & \text{if } \zeta \in \widehat{\Gamma}_0 \\
		0 & \text{if } \zeta \in \widehat{\Gamma}-\widehat{\Gamma}_0,
	\end{array}	 \right.
\end{equation}
where $I$ is the identity operator on $l^2(\mathcal{L}^*)$. We check that this field is measurable. By Remark \ref{measopconstantfield}, this amounts to showing that
\begin{equation}\label{rhoopsmeas}
	\zeta \mapsto \langle \rho_\zeta [x,L] g,h\rangle
\end{equation}
is measurable for all $g, h \in l^2(\Pi)\otimes l^2(\mathcal{L}^*)$. We know that
\[
r \mapsto \langle (U^r[x,L] \otimes I) g,h\rangle
\]
is measurable from $R$ to $\mathbb{C}$, by Lemma \ref{Uchimeas1}. Define $U_{g,h}: R \to \mathbb{C}$, by 
\[
U_{g,h}(r) = \langle (U^r[x,L] \otimes I) g,h\rangle, 
\]
and define $\Psi_{g,h}: \widehat{\Gamma} \to \mathbb{C}$ by 
\[
\Psi_{g,h}(\zeta) = \left\{ \begin{array}{ll}
	(U_{g,h} \circ \Psi^{-1})(\zeta) &\text{if }\zeta \in \widehat{\Gamma}_0  \\
	0 &\text{if } \zeta \in \widehat{\Gamma}-\widehat{\Gamma}_0
\end{array}
\right.
\]
The map \eqref{rhoopsmeas} is just $\Psi_{g,h}$, so we must show the measurability of $\Psi_{g,h}(\zeta)$ for all $g, h \in  l^2(\Pi)\otimes l^2(\mathcal{L}^*)$. This measurability follows from Lemma \ref{measonAB}, since $U_{g,h} \circ \Psi^{-1}$ is measurable on $\widehat{\Gamma}_0$, and the constant zero function is measurable on $\widehat{\Gamma}-\widehat{\Gamma}_0$.

We define the unitary representation $\rho$ of $\widehat{\Gamma}$ on $L^2\big(\widehat{\Gamma};l^2(\Pi)\otimes l^2(\mathcal{L}^*)\big)$ by \label{rho1}
\[
\boxed{\rho = \int_{\widehat{\Gamma}}^\oplus \rho_\zeta \ud \mu(\zeta)}\,,
\]
where $\rho_\zeta$ is the representation of $\Gamma$ on $l^2(\Pi)\otimes l^2(\mathcal{L}^*)$ given by \eqref{rhozetaops}. We show that $\widetilde{\Psi} \tau \widetilde{\Psi}^{-1} = \rho$. Given $F \in  L^2\big(R;l^2(\Pi)\otimes l^2(\mathcal{L}^*)\big)$, $[x,L]\in \Gamma$, and $\bar{\zeta} \in \widehat{\Gamma}$, we have
\begin{align*}
	(\widetilde{\Psi}\tau[x,L]F)(\bar{\zeta}) &=  \left\{ \begin{array}{ll} 
		\tau[x,L]F(\Psi^{-1}\bar{\zeta}) &\text{if } \bar{\zeta} \in \widehat{\Gamma}_0 \\
		0 &\text{if } \bar{\zeta} \in \widehat{\Gamma}-\widehat{\Gamma}_0. 
	\end{array}
	\right.\\
	&=  \left\{ \begin{array}{ll} 
		(U^{\Psi^{-1}(\bar{\zeta})}[x,L]\otimes I )F(\Psi^{-1}\bar{\zeta}) &\text{if } \bar{\zeta} \in \widehat{\Gamma}_0 \\
		0 &\text{if } \bar{\zeta} \in \widehat{\Gamma}-\widehat{\Gamma}_0. 
	\end{array}
	\right.\\
	&=  \left\{ \begin{array}{ll} 
		(U^{\Psi^{-1}(\bar{\zeta})}[x,L]\otimes I )(\widetilde{\Psi}F)(\bar{\zeta}) &\text{if } \bar{\zeta} \in \widehat{\Gamma}_0 \\
		0 &\text{if } \bar{\zeta} \in \widehat{\Gamma}-\widehat{\Gamma}_0. 
	\end{array}
	\right.\\
	&= \rho_{\bar{\zeta}}[x,L](\widetilde{\Psi}F)(\bar{\zeta})\\
	&= \Big(\int_{\widehat{\Gamma}}^\oplus \rho_\zeta[x,L] \ud \mu(\zeta) \Big)(\widetilde{\Psi}F)(\bar{\zeta}).
\end{align*}
Thus, 
\[
\widetilde{\Psi}\tau[x,L]F = \Big(\int_{\widehat{\Gamma}}^\oplus \rho_\zeta[x,L] \ud \mu(\zeta) \Big)(\widetilde{\Psi}F)
\]
for all $F \in L^2\big(R;l^2(\Pi)\otimes l^2(\mathcal{L}^*)\big)$ and all $[x,L] \in \Gamma$, whence
\[
(\widetilde{\Psi}\circ\tau)[x,L] = \Big(\int_{\widehat{\Gamma}}^\oplus \rho_\zeta\ud \mu(\zeta) \Big)[x,L] \circ \widetilde{\Psi}
\]
Thus $\widetilde{\Psi}$ is an intertwining map for $\tau$ and $\rho =\int_{\widehat{\Gamma}}^\oplus \rho_\zeta\ud \mu(\zeta)$.

We now show that the direct integral decomposition $\int_{\widehat{\Gamma}}^\oplus \rho_\zeta\ud \mu(\zeta)$ is the central decomposition of $\pi$. First we note that $\rho_\zeta$ is a factor representation for all $\zeta \in \widehat{\Gamma}$. Recall from Lemma \ref{uchiirred} that $U^\omega$ is irreducible for all $\omega \in \Omega_0$. In particular, $U^r$ is irreducible for $r\in R$, and since $\Psi^{-1}(\widehat{\Gamma}_0)=R$, we have that $U^{\Psi^{-1}(\zeta)}$ is irreducible for all $\zeta \in \widehat{\Gamma}_0$. Therefore $\rho_\zeta$ is a factor representation for all $\zeta \in \widehat{\Gamma}_0$. And since $\rho_\zeta = 0$ is (trivially) a factor representation whenever $\zeta \in \widehat{\Gamma}- \widehat{\Gamma}_0$, we have that all the $\rho_\zeta$ are factor representations. Moreover, since 
\[
\zeta = \Psi(\Psi^{-1}(\zeta)) = [U^{\Psi^{-1}(\zeta)}]
\]
for each $\zeta \in \widehat{\Gamma}_0$, we have that $U^{\Psi^{-1}(\zeta)} \in \zeta$ for all $\zeta \in \widehat{\Gamma}_0$. It follows from this, together with Corollary \ref{dualspacedescription} and the fact that $\rho_\zeta =	U^{\Psi^{-1}(\zeta)} \otimes I$ whenever $\zeta \in \widehat{\Gamma}_0$, that our decomposition satisfies the first condition (\ref{firstprop}) of the central decomposition theorem. We will need the following lemmas:
\begin{lemma}
	Let $W$ be a closed, $\rho$-invariant subspace of $L^2\big(\widehat{\Gamma};l^2(\Pi)\otimes l^2(\mathcal{L}^*)\big)$ and let $\varrho_1 \in W$ and  let $\varrho_2 \in L^2\big(\widehat{\Gamma};l^2(\Pi)\otimes l^2(\mathcal{L}^*)\big)$. Then $\varrho_2 \in W^\perp$ if and only if $\varrho_1(\zeta) \perp \varrho_2(\zeta)$ for almost every $\zeta \in \widehat{\Gamma}$. 
\end{lemma}
\begin{proof}
	Let $\varrho_1 \in W$ and $\varrho_2 \in W^\perp$. Then $\widetilde{\Psi}^{-1}(W)$ and $[\widetilde{\Psi}^{-1}(W)]^\perp$ are $\tau$-invariant, and
	\[
	\langle \widetilde{\Psi}^{-1}\varrho_1,  \widetilde{\Psi}^{-1}\varrho_2\rangle = \langle \varrho_1, \varrho_2\rangle = 0.
	\]
	It follows by Proposition \ref{kappaOrth} that
	\[
	\langle \widetilde{\Psi}^{-1}\varrho_1(r),  \widetilde{\Psi}^{-1}\varrho_2(r)\rangle =0 \ \ \ \text{for a.e.}\ r\in R.
	\]
Then by an application of \ref{imagemeasurethm}, we have
\begin{align*}
	\int_{\widehat{\Gamma}} |\langle \varrho_1(\zeta),\varrho_2(\zeta) \rangle| \ud \mu(\zeta) &= \int_{\widehat{\Gamma}_0} |\langle \varrho_1(\zeta),\varrho_2(\zeta) \rangle| \ud(m\circ \Psi^{-1})(\zeta) \\
	&= \int_{R} |\langle \varrho_1(\Psi(r),\varrho_2(\Psi(r)) \rangle| \ud r\\
	&=  \int_{R} |\langle \varrho_1(\Psi(r),\varrho_2(\Psi(r)) \rangle| \ud r\\
	&= \int_{R} |\langle \widetilde{\Psi}^{-1}\varrho_1(r),  \widetilde{\Psi}^{-1}\varrho_2(r)\rangle| \ud r\\
	&= 0.
\end{align*}
Therefore, $\varrho_1(\zeta) \perp \varrho_2(\zeta)$ for almost every $\zeta \in \widehat{\Gamma}$. The converse is trivial. \end{proof}
\begin{lemma}\label{orthconditionforCentral}
	Let $W$ be a closed, $\tau$-invariant subspace of $L^2\big(\widehat{\Gamma};l^2(\Pi)\otimes l^2(\mathcal{L}^*)\big)$. Select a countable dense subset $\mathcal{A}$ of $L^2\big(\widehat{\Gamma};l^2(\Pi)\otimes l^2(\mathcal{L}^*)\big)$, and define a range function $J_W^\mathcal{A}:\widehat{\Gamma} \to \{\text{closed subspaces of}\ l^2(\Pi)\otimes l^2(\mathcal{L}^*)\}$ as in Theorem \ref{mainpointwisedecompthm}. Let $\varrho \in L^2\big(\widehat{\Gamma};l^2(\Pi)\otimes l^2(\mathcal{L}^*)\big)$. Then $\varrho \in W^\perp$ if and only if $\varrho(\zeta) \perp J_W^\mathcal{A}(\zeta)$ for almost every $\zeta \in \widehat{\Gamma}$. In particular, $W$ satisfies the orthogonality condition \eqref{orthCond}.
\end{lemma}
\begin{proof}
	The argument is identical to that of Theorem \ref{RangeFunctionOrthTheorem}.
\end{proof}
It follows from the above lemma and Definition \ref{mblerangefun} that each $\rho$-invariant subspace $W\subset L^2\big(\widehat{\Gamma};l^2(\Pi)\otimes l^2(\mathcal{L}^*)\big)$ has associated to it a unique (up to a.e.-equivalence) measurable range function $J_W$. The following is a direct analogue of the first part of Theorem \ref{MainDecompTheorem}, and it is proved in an identical fashion.
	\begin{theorem}\label{DecompTheoremforCentral}
		Let $W$ be a closed, $\rho$-invariant subspace of $L^2\big(\widehat{\Gamma};l^2(\Pi)\otimes l^2(\mathcal{L}^*)\big)$, and let $J_W$ denote the measurable range function associated with $W$. Let $P_W$ be the orthogonal projection of $L^2\big(\widehat{\Gamma};l^2(\Pi)\otimes l^2(\mathcal{L}^*)\big)$ onto $W$, and $P_J(\zeta)$ denote the orthogonal projection from $l^2(\Pi) \otimes l^2(\mathcal{L}^*)$ onto $J_W(\zeta)$. Then
		\begin{list}{(\alph{ctr1})}{\usecounter{ctr1}}
			\setcounter{ctr1}{0}
			\item Regarding $\int_{\widehat{\Gamma}}^\oplus J_W(\zeta)\ud \mu(\zeta)$ as a subspace of $L^2\big(\widehat{\Gamma};l^2(\Pi)\otimes l^2(\mathcal{L}^*)\big)$, we have
			\[
			W = \int_{\widehat{\Gamma}}^\oplus J_W(\zeta)\ud \mu(\zeta)\ \ \ \text{and}\ \ \ 	P_W =\int_{\widehat{\Gamma}}^\oplus P_{J}(\zeta)\ud \mu(\zeta).
			\]
			\item $J_W(\zeta)$ is invariant under $\rho_\zeta$ for almost every $\zeta \in \widehat{\Gamma}$.
		\end{list}
\end{theorem}
 Let $\mathscr{A}$ denote the diagonal algebra of operators on 
 \[
 \int_{\widehat{\Gamma}}^\oplus \mathscr{H}_\zeta\ud \mu(\zeta)=L^2\big(\widehat{\Gamma};l^2(\Pi)\otimes l^2(\mathcal{L}^*)\big);
 \]
 i.e.: $\mathscr{A}$ is the algebra of decomposable operators on $\int_{\widehat{\Gamma}}^\oplus \mathscr{H}_\zeta\ud \zeta$ of the form $T(\zeta) = g(\zeta)I_\zeta$, where $g \in L^\infty(\widehat{\Gamma},\mu)$ and $I_\zeta$ is the identity operator on $\mathscr{H}_\zeta$. If $\rho[x,L] \in \rho(\Gamma)$, then since $\rho[x,L]$ is decomposable as
\[
\rho[x,L] = \int_{\widehat{\Gamma}}^\oplus \rho_\zeta[x,L] \ud \zeta,
\]
it is immediate from Lemma \ref{directintegralprops} that $\rho[x,L]$ commutes with each element of $\int_{\widehat{\Gamma}}^\oplus \rho_\zeta(\Gamma)^\prime \ud\zeta$. Thus
\[
\rho(\Gamma) \subset \Big(\int_{\widehat{\Gamma}}^\oplus \rho_\zeta(\Gamma)^\prime \ud\zeta\Big)^\prime.
\]
It follows that
\[
\int_{\widehat{\Gamma}}^\oplus \rho_\zeta(\Gamma)^\prime \ud\zeta \subset \rho(\Gamma)^\prime.
\]
Let $P$ be a projection in $\rho(\Gamma)^\prime$. Then $P= P_W$ for some $\rho$-invariant subspace $W$, by Proposition \ref{invorth}. Thus $W$ satisfies the orthogonality condition, by Lemma \ref{orthconditionforCentral}. Hence by Theorem \ref{DecompTheoremforCentral}, we have
\[
W = \int_{\widehat{\Gamma}}^\oplus J_W(\zeta)\ud \zeta,
\]
where $J_W$ is the measurable range function associated with $W$. And by the same theorem, $P$ is decomposable as $\int_{\widehat{\Gamma}}^\oplus P_J(\zeta)\ud \zeta$, where $P_J(\zeta)$ denotes the orthogonal projection from $l^2(\Pi) \otimes l^2(\mathcal{L}^*)$ onto $J_W(\zeta)$; and furthermore $J_W(\zeta)$ is invariant under $\rho_\zeta$ for a.e. $\zeta$. Thus
\[
P \in \int_{\widehat{\Gamma}}^\oplus \rho_\zeta(\Gamma)^\prime \ud\zeta.
\]
Now since the projections are total subset of the von Neumann algebra $\rho(\Gamma)^\prime$ (Proposition \ref{ProjDense}), it follows that 
\[
\rho(\Gamma)^\prime \subset \int_{\widehat{\Gamma}}^\oplus \rho_\zeta(\Gamma)^\prime \ud\zeta;
\]
hence we have equality: $\rho(\Gamma)^\prime = \int_{\widehat{\Gamma}}^\oplus \rho_\zeta(\Gamma)^\prime \ud\zeta$. It follows by Theorem \ref{818Tak} that
\[
\rho(\Gamma)^{\prime\prime} = \int_{\widehat{\Gamma}}^\oplus \rho_\zeta(\Gamma)^{\prime\prime}  \ud\zeta.
\]
Finally, by Proposition \ref{820Tak}, we have
\[
\mathscr{Z}\big(\rho(\Gamma)^{\prime\prime}\big) = \int_{\widehat{\Gamma}}^\oplus  \mathscr{Z}\big( \rho_\zeta(\Gamma)^{\prime\prime} \big) \ud\zeta.
\]
But the center of $\rho_\zeta(\Gamma)^{\prime\prime}$ is $\mathbb{C}_{\mathscr{H}_\zeta}$ for every $\zeta \in \widehat{\Gamma}$, since each $\rho_\zeta$ is a factor representation. Hence $\mathscr{Z}\big(\rho(\Gamma)^{-}\big) = \mathscr{A}$, which is equivalent to showing that our sequence of transformations carries $\mathscr{Z}(\pi(\Gamma)^-)$ onto $\mathscr{A}$.

\chapter{Conclusion}

The goal we set out to accomplish in this thesis was to find a complete characterization of the closed subspaces of $L^2(\mathbb{R}^n)$ invariant under shifts by crystallographic groups. This was accomplished in Theorem \ref{BigInvTheorem2}. For this, we defined a natural crystallographic shift representation and showed that it was unitarily equivalent to a direct integral $\tau$ of factor representations on $l^2(\Pi)\otimes l^2(\mathcal{L}^*)$. We showed that any subrepresentation of $\tau$ disintegrates into a direct integral of subrepresentations of these factor representations, thereby ensuring that our classification is complete. We showed that $\tau$ is unitarily equivalent to a direct integral over the unitary dual $\widehat{\Gamma}$, and that this final decomposition is the central decomposition of $\pi$.

The construction of the fundamental set $\Omega_{\mathrm{T}^*}$ in Corollary \ref{paramdomaindef} can be done in such a way that $\Omega_{\mathrm{T}^*}$ is Borel; this follows from results in descriptive set theory. Likewise, the transversal $\Omega$ can be chosen to be Borel. These arguments rely on the axiom of choice. While not needed in this thesis, the author believes that it is possible to construct these sets explicitly, and in such a way that $R\subset \Omega$ and $\Pi R\subset \Omega_{\mathrm{T}^*}$. It may also be of interest to check that the set $\Omega_0$ of free points in $\Omega$ is Borel. We saw in Theorem \ref{bigtheoremcor} that the fundamental set $\Omega$ is divided into two pieces; a ``good'' piece, $\Omega_0$, on which we have a closed formula \eqref{Uchi1} for the induced representations; and a ``bad'' piece $\Omega - \Omega_0$ on which we have no such formula.

There are several lines of inquiry one may wish to pursue. We note that the underlying idea of our approach is that of decomposing $\mathbb{R}^n$ into $\dot\bigcup_{[z,M]\in \Gamma^*} [z,M]\Omega$ and then doing some further identifications. Thus in our approach, we work with a space isomorphic to $L^2(\Omega)\otimes l^2(\Gamma^*)$. Suppose that we instead select a fundamental set $F$ for the action of $\Gamma$ on $\mathbb{R}^n$, and write 
\[
\mathbb{R}^n = \dot\bigcup_{[x,M]\in \Gamma} [x,M]F,
\]
whence 
\[
L^2(\mathbb{R}^n) \cong \bigoplus_{[x,M]\in \Gamma} L^2([x,M]F) \cong l^2(\Gamma)\otimes L^2(F).
\]
Then
\[
\pi([y,L])f([x,M]z]) = f([y,L]^{-1}[x,M]z),
\]
from which it follows that  $\pi \sim \lambda \otimes I_{L^2(F)}$, where $\lambda$ is the left-regular representation of $\Gamma$ on $l^2(\Gamma)$. One could also work toward the direct integral decomposition from this direction. In \cite{Cstarcrystal}, Keith Taylor defines an analogue of the Plancherel transform for $C^*(\Gamma)$, the group $C^*$-algebra of $\Gamma$. He also gives a formula for it, and characterizes its range. These results might be useful for pursuing an alternate route to the decomposition in this thesis.

Another question is how our work relates to the ideas in Manning's 2012 Ph.D. thesis \cite{Manning}. In particular, can we realize some kind of crystallographic Daubechies wavelet, as he does, using the characterization we have found.

Finally, it was shown by MacArthur and Taylor in \cite{KeithJosh} that one can construct a multiwavelet with crystallographic shifts on the assumption that one has an appropriate generalized multiresolution analysis. Knowing the form of the $\Gamma$-invariant subspaces of $L^2(\mathbb{R}^n)$ is a major step towards constructing this GMRA, as the $V_0$ space where the finite scaling ensemble lives must be $\Gamma$-invariant. The author believes that the results of this thesis will aid significantly in the construction of these finite scaling ensembles, thus advancing the constructive approach to multiwavelets with crystal shifts. In particular, a logical next step is to attempt to reproduce for crystallographic shifts the results of Bownik \cite{Bownik} dealing with principal shift invariant and finitely generated shift-invariant subspaces.

\begin{appendices}
	
	\chapter{Relevant Group Actions}\label{GroupActions}
	We here include verification that the rules in Table \ref{TableOfActions} do indeed give group actions. We include this because some of these calculations are referenced elsewhere in this thesis.
	
	\section[Actions on \texorpdfstring{$\mathbb{R}^n$}{Euclidean space}, \texorpdfstring{$\mathrm{T}$}{T}, and \texorpdfstring{$\widehat{\mathrm{T}}$}{T-Hat}]{Actions on \texorpdfstring{$\boldsymbol{\mathbb{R}^n}$}{Euclidean space}, \texorpdfstring{$\boldsymbol{\mathrm{T}}$}{T}, and \texorpdfstring{$\boldsymbol{\widehat{\mathrm{T}}}$}{T-Hat}}
	
	\begin{list}{\textbf{Action \arabic{ctr}:}}{\usecounter{ctr}}
		\setcounter{ctr}{0}
		\item $\Gamma$ acts on $\mathbb{R}^n$: $[x,L]\cdot_1 u = L(u+x)$ is a group action.
		
		\textbf{Check:} $[0,\text{id}]\cdot_1 u = \text{id}(u+0) = u$. Moreover,
		\begin{align*} 
			[x,L]\cdot_1([y,M]\cdot_1 u) &= [x,L]\cdot_1\big(M(u+y)\big) \\
			&= [x,L]\cdot_1 (Mu + My) \\
			&= L(Mu + My + x)\\
			&= LMu + LMy + Lx \\
			&= LM(u+M^{-1}x + y) \\
			&= ([M^{-1}x+y,LM])\cdot_1 u \\
			&= ([x,L][y,M])\cdot_1 u,
		\end{align*} 
		so (i) does indeed give a group action. As a corollary, the translation subgroup $\mathrm{T}$ acts by translation on $\mathbb{R}^n$: $[x,\text{id}]\cdot_1 u = u+x$.
		
		\item \label{pionR} $\Pi$ acts on $\mathbb{R}^n$ by matrix multiplication: $L\cdot_2 u := Lu$.
		
		\textbf{Check:} Clearly $\text{id}\cdot_{\arabic{ctr}} u = u$; and 
		\[
		(LM)\cdot_{\arabic{ctr}} u = LMu = L\cdot_2(Mu) = L\cdot_{\arabic{ctr}}(M\cdot_{\arabic{ctr}} u).
		\]
		This shows that (ii) does indeed give a group action. 
		
		\item \label{gammaont} $\Gamma$ acts on $\mathrm{T}$ by conjugation:
		\begin{equation}\label{gammaonn}
			[x,L] \cdot_{\arabic{ctr}} [y,\text{id}] = [x,L][y,\text{id}][x,L]^{-1}. 
		\end{equation}
		That this gives a group action follows from the fact that $\mathrm{T}$ is normal in $\Gamma$, and is verified in many algebra texts.
		
		\item\label{Action4} \label{piont} Since $\mathrm{T}$ is abelian, $\Pi$ acts on $\mathrm{T}$ via the cross-section $\gamma$: 
		\[
		L\cdot_{\arabic{ctr}} [x,\text{id}] = \gamma(L)[x,\text{id}]\gamma(L)^{-1}. 
		\]
		Moreover, it can be shown that in this case we have
		\[
		L^{-1}\cdot_{\arabic{ctr}} [x,\text{id}] = \gamma(L)^{-1}[x,\text{id}]\gamma(L).
		\]
		We may express this group action using the notation of \eqref{gammaonn}:
		\[
		L\cdot_{\arabic{ctr}} [x,\text{id}] = \gamma(L)\cdot_{\idref{3}} [x,\text{id}].
		\]
		\textbf{Check:} Since $\mathrm{T}$ is abelian and $\gamma(\text{id}_{\Pi}) \in \mathrm{T}$, we have that $\text{id}\cdot_{\arabic{ctr}} [x,\text{id}] = [x,\text{id}]$ for all $[x,\text{id}]\in \mathrm{T}$. It can be shown that 
		\begin{equation}\label{cocycledef}
			\gamma(L)\gamma(M) = \gamma(LM)\alpha(L,M)
		\end{equation}
		for some $\alpha(L,M)\in \mathrm{T}$.  
		\begin{tcolorbox}
		\begin{defn}\label{cocycle}
	The map $\alpha: \Pi\times \Pi \to \mathrm{T}$ satisfying \eqref{cocycledef} is called a \textit{cocycle}.
		\end{defn}
	\end{tcolorbox}
	See \cite{TaylorKaniuth} and \cite{Cstarcrystal} for more about the cocycle. Thus we have 
		\begin{align*}
			(LM)\cdot_{\arabic{ctr}} [x,\text{id}] &= \gamma(LM)[x,\text{id}]\gamma(LM)^{-1}\\
			&= \gamma(L)\gamma(M)\alpha(L,M)^{-1}[x,\text{id}]\big(\gamma(L)\gamma(M)\alpha(L,M)^{-1}\big)^{-1}\\
			&= \gamma(L)\gamma(M)\alpha(L,M)^{-1}[x,\text{id}]\alpha(L,M)\gamma(M)^{-1}\gamma(L)^{-1}\\
			&= \gamma(L)\gamma(M)[x,\text{id}]\gamma(M)^{-1}\gamma(L)^{-1}\\
			&= \gamma(L)\cdot_{\idref{3}}\big(\gamma(M)\cdot_{\idref{3}} [x,\text{id}]\big)\\
			&= L\cdot_{\arabic{ctr}} (M\cdot_{\arabic{ctr}} [x,\text{id}]).
		\end{align*}
		
		\item \label{gammaonthat} $\Gamma$ acts on $\widehat{\mathrm{T}}$. The action of $\Gamma$ on $\widehat{\mathrm{T}}$ is given by 
		\[
		([x,L]\cdot_{\arabic{ctr}} \chi^{\mathrm{T}})([y,\text{id}]) = \chi^{\mathrm{T}}([x,L]^{-1}\cdot_{\idref{3}} [y,\text{id}])
		\]
		for all $[x,L] \in \Gamma$, $\chi^{\mathrm{T}}\in \widehat{\mathrm{T}}$, and $[y,\text{id}] \in \mathrm{T}$. Moreover, for any $[y,\text{id}]\in\mathrm{T}$, and any $[x,L]\in\Gamma$, we have
		\begin{equation}\label{actiongammaonN}
			([x,L][y,\text{id}])\cdot_{\arabic{ctr}}\chi^{\mathrm{T}} = [x,L]\cdot_{\arabic{ctr}}\chi^{\mathrm{T}}
		\end{equation}
		In particular,
		\begin{equation}\label{specialcase}
			[y,\text{id}]\cdot_{\arabic{ctr}}\chi^{\mathrm{T}} = \chi^{\mathrm{T}}
		\end{equation}
		for all $[y,\text{id}]\in \mathrm{T}$ and all $\chi^{\mathrm{T}}\in\widehat{\mathrm{T}}$. 
		
		\textbf{Check:} Since $\Gamma$ acts on $\mathrm{T}$, we have $[0,\text{id}]\cdot_{\idref{3}} [y,\text{id}] = [y,\text{id}]$ for each $[y,\text{id}]\in\mathrm{T}$. Thus
		\[
		\big([0,\text{id}]\cdot_{\arabic{ctr}}\chi^{\mathrm{T}}\big)([y,\text{id}]) = \chi^{\mathrm{T}}([0,\text{id}]^{-1}\cdot_{\idref{3}} [y,\text{id}]) = \chi^{\mathrm{T}}([y,\text{id}]),
		\]
		for all $[y,\text{id}] \in\mathrm{T}$, whence
		\[
		[0,\text{id}]\cdot_{\arabic{ctr}} \chi^{\mathrm{T}} = \chi^{\mathrm{T}}.
		\]
		Now it must be shown that 
		\begin{equation}\label{grpaction}
			([x,L][y,M])\cdot_{\arabic{ctr}} \chi^{\mathrm{T}} = [x,L]\cdot_{\arabic{ctr}}([y,M]\cdot_{\arabic{ctr}}\chi^{\mathrm{T}})
		\end{equation}
		for all $[x,L]$ and $[y,M]$ in $\Gamma$. This will follow again from the fact that $\Gamma$ acts on $\mathrm{T}$. For any $[z,\text{id}] \in\mathrm{T}$ we have
		\begin{align*}
			\big(([x,L][y,M])\cdot_{\arabic{ctr}} \chi^{\mathrm{T}}\big)([z,\text{id}]) &=
			\chi^{\mathrm{T}}\big(([x,L][y,M])^{-1}\cdot_{\idref{3}} [z,\text{id}]\big) \\
			&=\chi^{\mathrm{T}}\big(([y,M]^{-1}[x,L]^{-1})\cdot_{\idref{3}} [z,\text{id}]\big) \\
			&= \chi^{\mathrm{T}}\big([y,M]^{-1}\cdot_{\idref{3}}([x,L]^{-1}\cdot_{\idref{3}} [z,\text{id}])\big) \\
			&= ([y,M]\cdot_{\arabic{ctr}}\chi^{\mathrm{T}})\big( [x,L]^{-1}\cdot_{\idref{3}}[z,\text{id}]\big)\\
			&= \big([x,L]\cdot_{\arabic{ctr}} ([y,M]\cdot_{\arabic{ctr}}\chi^{\mathrm{T}})\big)([z,\text{id}]),
		\end{align*}
		whence equation \eqref{grpaction} holds. To prove \eqref{actiongammaonN} above, let $[y,\text{id}]$ and $[z,\text{id}]$ be in $\mathrm{T}$, and $[x,L]\in\Gamma$. We have
		\begin{align*}
			\big(([x,L][y,\text{id}])\cdot_{\arabic{ctr}}\chi^{\mathrm{T}}\big)([z,\text{id}])&= \chi^{\mathrm{T}}\big(([x,L][y,\text{id}])^{-1}\cdot_{\idref{3}} [z,\text{id}] \big)\\
			&= \chi^{\mathrm{T}}\big(([x,L][y,\text{id}])^{-1} [z,\text{id}] [x,L][y,\text{id}]\big)\\
			&= \chi^{\mathrm{T}}\big([y,\text{id}]^{-1} [x,L]^{-1}[z,\text{id}] [x,L][y,\text{id}]\big)\\
			&= \chi^{\mathrm{T}}([y,\text{id}])^{-1} \chi^{\mathrm{T}}([x,L]^{-1}[z,\text{id}] [x,L])\chi^{\mathrm{T}}([y,\text{id}])\\
			&= \chi^{\mathrm{T}}([x,L]^{-1}[z,\text{id}] [x,L])\\
			&=([x,L]\cdot_{\arabic{ctr}} \chi^{\mathrm{T}})([z,\text{id}]).
		\end{align*}
		\begin{rem}\label{moreaboutactiongammaonThat}
			Since
			\[
			[x,L]=[y,L][x-y,\mathrm{id}],
			\]
			it follows directly from \eqref{actiongammaonN} that
			\begin{equation}\label{moreaboutactiongammaonThat1}
				[x,L]\cdot_{\idref{5}} \chi^{\mathrm{T}} = 	[y,L]\cdot_{\idref{5}} \chi^{\mathrm{T}}. 
			\end{equation}
			In particular, 
			\begin{equation}\label{moreaboutactiongammaonThat2}
				[x,L]\cdot_{\idref{5}} \chi^{\mathrm{T}} = 	\gamma(L)\cdot_{\idref{5}} \chi^{\mathrm{T}}
			\end{equation}
			for all $[x,L] \in \Gamma$ and any choice of cross-section $\gamma$.
		\end{rem}
		
		\item \label{pionthat} $\Pi$ acts on $\widehat{\mathrm{T}}$ via the cross-section $\gamma$. This is given by 
		\[
		L\cdot_{\arabic{ctr}} \chi^{\mathrm{T}} = \gamma(L)\cdot_{\idref{5}} \chi^{\mathrm{T}},
		\]
		where $\gamma$ is the cross-section. This action is independent of the choice of cross-section. 
		
		\textbf{Check:} We have
		\[
		\text{id}\cdot_{\arabic{ctr}} \chi^{\mathrm{T}} = \gamma(\text{id})\cdot_{\idref{5}} \chi^{\mathrm{T}} = \chi^{\mathrm{T}}
		\]
		by \eqref{specialcase}, since $\gamma(\text{id})\in \mathrm{T}$. Moreover, for each $[x,\text{id}]\in\mathrm{T}$ we have 
		\begin{align*}
			&\big((LM)\cdot_{\arabic{ctr}}\chi^{\mathrm{T}}\big)([x,\text{id}])\\
			&= \big(\gamma(LM)\cdot_{\idref{5}}\chi^{\mathrm{T}}\big)([x,\text{id}]) \\
			&= \chi^{\mathrm{T}}\big(\gamma(LM)^{-1}\cdot_{\idref{3}} [x,\text{id}]\big)\\
			&= \chi^{\mathrm{T}}\big(\gamma(LM)^{-1}[x,\text{id}]\gamma(LM)\big)\\
			&= \chi^{\mathrm{T}}\big([\gamma(L)\gamma(M)\alpha(L,M)^{-1}]^{-1}[x,\text{id}]\gamma(L)\gamma(M)\alpha(L,M)^{-1}\big)\\
			&=\chi^{\mathrm{T}}\big(\alpha(L,M)\gamma(M)^{-1}\gamma(L)^{-1}[x,\text{id}]\gamma(L)\gamma(M)\alpha(L,M)^{-1}\big)\\
			&=\chi^{\mathrm{T}}\big(\alpha(L,M)\big)\chi^{\mathrm{T}}\big(\gamma(M)^{-1}\gamma(L)^{-1}[x,\text{id}]\gamma(L)\gamma(M)\big)\chi^{\mathrm{T}}\big(\alpha(L,M)\big)^{-1}\\
			&=\chi^{\mathrm{T}}\big(\gamma(M)^{-1}\gamma(L)^{-1}[x,\text{id}]\gamma(L)\gamma(M)\big)\\
			&=\chi^{\mathrm{T}}\big(\gamma(M)^{-1}\cdot_{\idref{3}}(\gamma(L)^{-1}\cdot_{\idref{3}} [x,\text{id}])\big)\\
			&= (\gamma(M)\cdot_{\idref{5}}\chi^{\mathrm{T}})\big(\gamma(L)^{-1}\cdot_{\idref{5}} [x,\text{id}]\big)\\
			&=\big(\gamma(L)\cdot_{\idref{5}}(\gamma(M)\cdot_{\idref{5}} \chi^{\mathrm{T}})\big)([x,\text{id}])\\
			&= \big(L\cdot_{\arabic{ctr}}(M\cdot_{\arabic{ctr}} \chi^{\mathrm{T}})\big)([x,\text{id}]),
		\end{align*}
		so that $(LM)\cdot_{\arabic{ctr}} \chi^{\mathrm{T}} = L\cdot_{\arabic{ctr}}(M\cdot_{\arabic{ctr}}\chi^{\mathrm{T}})$, and hence $\Pi$ acts on $\widehat{\mathrm{T}}$. To see that this action is independent of choice of cross-section, let $\gamma_1$ and $\gamma_2$ be two cross-sections. Then there exist $x_1,x_2 \in\mathbb{R}^n$ such that $\gamma_i(L) = [x_i,L]\in\Gamma$, for $i=1,2$. It follows immediately from \eqref{moreaboutactiongammaonThat1} that
		\[
		\gamma_1(L)\cdot_{\idref{5}} \chi^{\mathrm{T}} = \gamma_2(L)\cdot_{\idref{5}} \chi^{\mathrm{T}}  \ \ \ \text{for all}\ \chi^{\mathrm{T}}\in\widehat{\mathrm{T}},
		\]
		so that $\gamma_1$ and $\gamma_2$ give the same action.
	\end{list}
	
	\vspace{.3cm}
	
	\section[Actions on \texorpdfstring{$\mathcal{L}$}{L} and \texorpdfstring{$\widehat{\mathcal{L}}$}{L Hat}]{Actions on \texorpdfstring{$\boldsymbol{\mathcal{L}}$}{L} and \texorpdfstring{$\boldsymbol{\widehat{\mathcal{L}}}$}{L Hat} }\label{ActionsLLhat}
	
	Recall that $\mathrm{T} \cong \mathcal{L}$ via $\pi_1|_{\mathrm{T}}:\mathrm{T}\to\mathcal{L}$, where $\pi_1\vert_{\mathrm{T}}([x,\mathrm{id}]) = x$\label{pi1}. The above actions on $\mathrm{T}$ and $\widehat{\mathrm{T}}$ can be modified to yield group actions on $\mathcal{L}$ and $\widehat{\mathcal{L}}$. The following fact will be useful:
	\begin{rem}\label{actionfromiso}
		Let $f: \mathcal{S}_1 \to \mathcal{S}_2$ be a bijection of sets, and $G$ a group acting on $\mathcal{S}_1$, with the group action denoted by $\cdot_{\mathcal{S}_1}: G\times \mathcal{S}_1\to \mathcal{S}_1$. Then the operation $\cdot_{\mathcal{S}_2}: G\times {\mathcal{S}_2} \to {\mathcal{S}_2}$ defined by $g\cdot_{\mathcal{S}_2} y := f(g\cdot_{\mathcal{S}_1} (f^{-1}(y)))$ is a group action on ${\mathcal{S}_2}$. 
	\end{rem}
	We introduce five actions analogous to Actions \idref{2} through \idref{6} above. In what follows it will simplify notation to let $f = \pi_1|_{\mathrm{T}}$.
	\begin{list}{\textbf{Action \arabic{ctr}:}}{\usecounter{ctr}}
		\setcounter{ctr}{6}
		\item \label{pionl} There is an action of $\Pi$ on $\mathcal{L}$, which is just the restriction of Action \idref{2} to $\mathcal{L}\subset \mathbb{R}^n$:
		\[
		L\cdot_{\arabic{ctr}} y := L\cdot_{\idref{2}} y = Ly,\ \text{for}\ L\in \Pi\ \text{and}\ y \in\mathcal{L}.
		\]	
		\textbf{Check:}  We verify that the lattice $\mathcal{L}$ is invariant under $\cdot_{\arabic{ctr}}$ of $\Pi$. Let $y\in\mathcal{L}$ and $L\in\Pi$. Then $[x,\mathrm{id}]\in\mathrm{T}$ and $L\cdot_{\idref{4}} [y,\mathrm{id}] = \gamma(L) \cdot_{\idref{3}} [y,\mathrm{id}]  \in \mathrm{T}$. Letting $\gamma(L)=[x_L,L]$, we have
		\begin{align}
			L\cdot_{\idref{4}} [y,\mathrm{id}]&=[x_L,L]\cdot_{\idref{3}} [y,\text{id}] \nonumber \\
			&= [x_L,L][y,\text{id}][x_L,L]^{-1} \nonumber\\
			&= [x_L,L][y,\text{id}][-Lx_L,L^{-1}]\nonumber\\
			&= [x_L,L][Ly-Lx_L,L^{-1}]\nonumber\\
			&= [Lx_L + Ly - Lx_L,\text{id}]\nonumber\\
			&= [Ly,\text{id}];\label{conjbygamma}
		\end{align}
		so $[Ly,\text{id}] \in \mathrm{T}$ whenever $y\in\mathcal{L}$. Hence $L\cdot_{\arabic{ctr}} y \in \mathcal{L}$ whenever $y\in \mathcal{L}$.

		\item \label{gammaonl} Applying the above remark with $\mathcal{S}_1 = \mathrm{T}$ and $\mathcal{S}_2=\mathcal{L}$, and with $f= \pi_1|_{\mathrm{T}}: \mathrm{T}\to\mathcal{L}$ and $\cdot_{\idref{3}}$ the action of $\Gamma$ on $\mathrm{T}$ described above, we get an action of $\Gamma$ on $\mathcal{L}$:
		\[
		[x,L]\cdot_{\arabic{ctr}} y = f\big([x,L]\cdot_{\idref{3}}(f^{-1}(y)\big)\ \ \text{for all}\ [x,L]\in\Gamma\ \text{and}\ y\in\mathcal{L}.
		\]
		
		\textbf{Check:} This follows immediately from the remark above, the proof of which is trivial.
		
		\item \label{pionl2} There is an action of $\Pi$ on $\mathcal{L}$ via $\gamma$, given by
		\[
		L\cdot_{\arabic{ctr}}y = \gamma(L)\cdot_{\idref{8}}  y.
		\]
		\textbf{Check:} We have
		\begin{align*}
			\gamma(L)\cdot_{\arabic{ctr}} y &= \gamma(L)\cdot_{\idref{8}}  y\\
			&= f\big(\gamma(L)\cdot_{\idref{3}}(f^{-1}(y))\big)\\
			&= f\big(L\cdot_{\idref{4}}(f^{-1}(y))\big)
		\end{align*}
		It follows from Remark \ref{actionfromiso} that $\cdot_{\arabic{ctr}}$ is a group action.
	\end{list}
	
	For the next two actions, it will be helpful to introduce a bijective map $\vartheta: \widehat{\mathrm{T}} \overset{\cong}{\to}\widehat{\mathcal{L}}$, defined by $\vartheta(\chi^{\mathrm{T}}) = \chi^{\mathrm{T}} \circ f^{-1}$, where $f$ is defined above. Of course, $\vartheta^{-1}: \widehat{\mathcal{L}} \to \widehat{\mathrm{T}}$ is given by $\vartheta^{-1}(\chi^{\mathcal{L}})= \chi^{\mathcal{L}} \circ f$.
	
	\begin{list}{\textbf{Action \arabic{ctr}:}}{\usecounter{ctr}}
		\setcounter{ctr}{9}
		
		\item \label{gammaonlhat} There is an action of $\Gamma$ on $\widehat{\mathcal{L}}$, given by 
		\[
		([x,L]\cdot_{\arabic{ctr}} \chi^{\mathcal{L}})(y) = \chi^{\mathcal{L}}([x,L]^{-1}\cdot_{\idref{8}} y) 
		\]
		for all $[x,L]\in\Gamma, \chi^{\mathcal{L}}\in\widehat{\mathcal{L}}$, and $y\in\mathcal{L}$.
		
		\textbf{Check:} First observe that
		\begin{align}
			([x,L]\cdot_{\arabic{ctr}} \chi^{\mathcal{L}})(y) &=  \chi^{\mathcal{L}}([x,L]^{-1}\cdot_{\idref{8}} y) \nonumber\\
			&= \chi^{\mathcal{L}}\big(f\big([x,L]^{-1}\cdot_{\idref{3}} f^{-1}(y)\big)\big)\nonumber \\
			&=\chi^{\mathcal{L}}\big(f([x,L]^{-1}\cdot_{\idref{3}}[y,\mathrm{id}])\big)\nonumber \\
			&=\vartheta^{-1}(\chi^{\mathcal{L}})([x,L]^{-1}\cdot_{\idref{3}}[y,\mathrm{id}])\label{converttentothree} \\
			&=\big([x,L]\cdot_{\idref{5}}\vartheta^{-1}(\chi^{\mathcal{L}})\big)([y,\mathrm{id}])\label{convertten}
		\end{align}
		Thus, for any $[y,\mathrm{id}]\in\mathrm{T}$, 
		\begin{align*}
			([x,L] \cdot_{\idref{5}}\vartheta^{-1}(\chi^{\mathcal{L}}))([y,\mathrm{id}]) &= ([x,L]\cdot_{\arabic{ctr}} \chi^{\mathcal{L}})(y)\\
			&= ([x,L]\cdot_{\arabic{ctr}} \chi^{\mathcal{L}})\big(f([y,\mathrm{id}])\big)\\
			&=\vartheta^{-1}([x,L]\cdot_{\arabic{ctr}}\chi^{\mathcal{L}})([y,\mathrm{id}]),
		\end{align*}
		so that 
		\begin{equation}\label{derivedten}
			[x,L] \cdot_{\idref{5}}\vartheta^{-1}(\chi^{\mathcal{L}}) = \vartheta^{-1}([x,L]\cdot_{\arabic{ctr}}\chi^{\mathcal{L}}).
		\end{equation}
		We may write this as
		\begin{equation}\label{derivedten2}
			[x,L]\cdot_{\arabic{ctr}}\chi^{\mathcal{L}}=\vartheta\big([x,L] \cdot_{\idref{5}}\vartheta^{-1}(\chi^{\mathcal{L}})\big).
		\end{equation}
		It follows from Remark \ref{actionfromiso} that $\cdot_{\arabic{ctr}}$ is an action.
		
		\item \label{pionlhat}  There is an action of $\Pi$ on $\widehat{\mathcal{L}}$, given by
		\[
		L\cdot_{\arabic{ctr}} \chi^{\mathcal{L}} = \gamma(L)\cdot_{\idref{10}} \chi^{\mathcal{L}}.
		\]
		
		\textbf{Check:} We have
		\begin{align}\label{forpionlhat}
			(L\cdot_{\arabic{ctr}} \chi^{\mathcal{L}})(y) &= \big(\gamma(L)\cdot_{\idref{10}} \chi^{\mathcal{L}}\big)(y)\nonumber \\
			&=\big(\gamma(L)\cdot_{\idref{5}} \vartheta^{-1}(\chi^{\mathcal{L}})\big)([y,\mathrm{id}])\ \text{by equation}\ \eqref{convertten}\nonumber \\
			&= \big(L\cdot_{\idref{6}} \vartheta^{-1}(\chi^{\mathcal{L}})\big)([y,\mathrm{id}]).
		\end{align}
		Thus, for any $[y,\mathrm{id}]\in\mathrm{T}$, we have
		\begin{align*}
			\big(L\cdot_{\idref{6}} \vartheta^{-1}(\chi^{\mathcal{L}})\big)([y,\mathrm{id}]) &= 	(L\cdot_{\arabic{ctr}} \chi^{\mathcal{L}})(y)\\
			&= 	(L\cdot_{\arabic{ctr}} \chi^{\mathcal{L}})\big(f([y,\mathrm{id}])\big)\\
			&= 	\vartheta^{-1}(L\cdot_{\arabic{ctr}} \chi^{\mathcal{L}})([y,\mathrm{id}]), 
		\end{align*}
		so that
		\[
		L\cdot_{\idref{6}} \vartheta^{-1}(\chi^{\mathcal{L}}) = \vartheta^{-1}(L\cdot_{\arabic{ctr}}\chi^{\mathcal{L}}).
		\]
		Rewriting this as
		\[
		L\cdot_{\arabic{ctr}}\chi^{\mathcal{L}} = \vartheta\big(L\cdot_{\idref{6}}\vartheta^{-1}(\chi^{\mathcal{L}}) \big),
		\]
		it follows from Remark \ref{actionfromiso} that $\cdot_{\arabic{ctr}}$ is an action.
	\end{list}
	
	\section[The Action of \texorpdfstring{$\Pi$}{Pi} on \texorpdfstring{$\Omega_{\mathrm{T}^*}$}{Parameter Domain}]{The Action of \texorpdfstring{$\boldsymbol{\Pi}$}{Pi} on \texorpdfstring{$\boldsymbol{\Omega_{\mathrm{T}^*}}$}{Parameter Domain}}\label{ActionPiParamSection}
	
	\begin{list}{\textbf{Action \arabic{ctr}:}}{\usecounter{ctr}}
		\setcounter{ctr}{11}
		\item \label{pionparam} $\Pi$ acts on $\Omega_{\mathrm{T}^*}$ by
		\begin{equation}\label{pionparameq}
			L\cdot_{12} \nu := (q\vert_{\Omega_{\mathrm{T}^*}})^{-1}\big(q(L\nu)\big).
		\end{equation}
		\textbf{Check:} First, for any $\nu \in \Omega_{\mathrm{T}^*}$, we have
		\begin{align*}
			\mathrm{id} \cdot_{12} \nu &= q\vert_{\Omega_{\mathrm{T}^*}}^{-1}\big(q(\mathrm{id}\,\nu)\big)\\
			&= q\vert_{\Omega_{\mathrm{T}^*}}^{-1}\big(q(\nu)\big)\\
			&= q\vert_{\Omega_{\mathrm{T}^*}}^{-1}(\nu+\mathcal{L}^*)\\
			&= \nu.
		\end{align*}
		Second, for any $\nu \in \Omega_{\mathrm{T}^*}$, and $L,M\in\Pi$, we have
		\begin{align*}
			M\cdot_{12}(L\cdot_{12} \nu) &= \big(q\vert_{\Omega_{\mathrm{T}^*}}^{-1}\circ q\big)\Big(M\big[ (q\vert_{\Omega_{\mathrm{T}^*}}^{-1}\circ q)(L\nu)\big] \Big)\\
			&=\big(q\vert_{\Omega_{\mathrm{T}^*}}^{-1}\circ q\big)\Big(M\big[ (q\vert_{\Omega_{\mathrm{T}^*}}^{-1}(L\nu + \mathcal{L}^*)\big] \Big).
		\end{align*}
		Now  $q\vert_{\Omega_{\mathrm{T}^*}}^{-1}(L\nu + \mathcal{L}^*)$ belongs to the coset $L\nu + \mathcal{L}^*$, and therefore 
		\[
		M\big[ (q\vert_{\Omega_{\mathrm{T}^*}}^{-1}(L\nu + \mathcal{L}^*)\big]
		\]
		belongs to the coset $ML\nu + \mathcal{L}^*$. It follows that 
		\[
		q\Big(M\big[ (q\vert_{\Omega_{\mathrm{T}^*}}^{-1}(L\nu + \mathcal{L}^*)\big] \Big) = ML\nu + \mathcal{L}^*,
		\]
		and hence 
		\begin{align*}
			\big(q\vert_{\Omega_{\mathrm{T}^*}}^{-1}\circ q\big)\Big(M\big[ (q\vert_{\Omega_{\mathrm{T}^*}}^{-1}(L\nu + \mathcal{L}^*)\big] \Big) &= q\vert_{\Omega_{\mathrm{T}^*}}^{-1}( ML\nu + \mathcal{L}^*) \\
			&= q\vert_{\Omega_{\mathrm{T}^*}}^{-1}\big(q(ML\nu)\big)\\
			&= (ML)\cdot_{12}(\nu).
		\end{align*}
		Thus we have $M\cdot_{12}(L\cdot_{12} \nu) = (ML)\cdot_{12}(\nu)$ for all $\nu \in \Omega_{\mathrm{T}^*}$, which completes our verification that \eqref{pionparameq} gives a group action.
	\end{list}

	\chapter{Results on Weakly Measurable Functions}
	Let $(X,\mathcal{M}, \mu)$ be a measure space. 
	\begin{lemma}\label{WeakMeasurabilityImpliesNormMeasurability} If $\mathscr{H}$ is separable and $f: X \rightarrow \mathscr{H}$ is weakly measurable then $x \mapsto \|f(x)\|$ is a measurable real-valued function.
	\end{lemma}
	\begin{proof} Let $\{e_n\}_{n=1}^\infty$ be a basis for $\mathscr{H}$; by Parseval's identity,
		$$
		\|f(x)\|^2 = \sum_{n=1}^\infty |\langle f(x), e_n \rangle|^2.
		$$ 
		Since weak measurability means that $x \mapsto \langle f(x), e_n \rangle$ is measurable for each $n$, it follows that $x \mapsto \sum_{n=1}^N |\langle f(x), e_n \rangle|^2$ is measurable for each $N$. Since the pointwise limit of a sequence of measurable functions is measurable, we have that $x\mapsto \|f(x)\|^2$ is measurable, and therefore that $x\mapsto \|f(x)\|$ is measurable.
	\end{proof}		
	\begin{lemma}\label{LinearCombinationsWeaklyMeasurable} Constant functions are weakly measurable. Moreover, if $f:X\rightarrow \mathscr{H}$ and $g:X\rightarrow \mathscr{H}$ are weakly measurable, and $c\in\mathbb{C}$, then $cf$ and $f+g$ are weakly measurable. \end{lemma}
	\begin{proof} Given $\psi\in \mathscr{H}$, $x\mapsto \langle \psi, \phi\rangle$ is measurable for each $\phi\in \mathscr{H}$, since constant scalar-valued functions are measurable. Thus constant functions are weakly measurable. If $f$ is weakly measurable and $c\in\mathbb{C}$, then the weak measurability of $cf$ follows from the analogous property for scalar-valued measurable functions, since $\langle cf(x), \phi\rangle = c  \langle f(x),\phi\rangle$. If $f$ and $g$ are weakly measurable, then the weak  measurability of $f+g$ again follows from the analogous property for scalar-valued functions, since $\langle f(x) + g(x),\phi\rangle = \langle f(x),\phi\rangle + \langle g(x),\phi\rangle$.\end{proof}
	For any topological space $X$, we let $\mathcal{B}_{X}$ denote the Borel $\sigma$-algebra on $X$.
	
	\begin{prop}\label{KfWeaklyMeasurable} Let $\mathscr{H}$ and $\mathscr{K}$ be separable Hilbert spaces. If $f: X \rightarrow \mathscr{H}$ and $K: X \rightarrow \mathscr{B}(\mathscr{H},\mathscr{K})$ are weakly measurable, then $K(\cdot)f(\cdot): X \rightarrow \mathscr{K}$ is weakly measurable. \end{prop}
	\begin{proof}We first observe that if $f: X \rightarrow \mathscr{H}$ and $g: X\rightarrow \mathscr{H}$ are weakly measurable, then $x\mapsto \langle f(x),g(x)\rangle$ is measurable. To see this, observe that by the polarization identity we have
		\begin{multline*}
			\langle f(x),g(x)\rangle = \textstyle{\frac14} \big(\|f(x)+g(x)\|^2 - \|f(x)-g(x)\|^2 \\ + i\|f(x)+i g(x)\|^2-i\|f(x)-i g(x)\|^2\big),
		\end{multline*}
		and each of the four summands on the right-hand side is measurable.
		
		We now note that if $K: X\rightarrow \mathscr{B}(\mathscr{H},\mathscr{K})$ is weakly measurable, then so is $K^*: X\rightarrow \mathscr{B}(\mathscr{K},\mathscr{H})$, where $K^*(x)$ is the adjoint of $K(x)$. This is because for any $\psi \in \mathscr{H}$, $\phi \in \mathscr{K}$, and $x\in X$, we have
		$$
		\langle K^*(x)\phi, \psi\rangle = \overline{\langle \psi, K^*(x)\phi\rangle} = \overline{\langle K(x)\psi, \phi\rangle},
		$$
		which is measurable because $K$ is weakly measurable, and complex conjugation is measurable. Now $K^*$ being weakly measurable means, by definition, that $K^*(x) \phi$ is weakly measurable for each $\phi \in \mathscr{K}$, whence, by the first paragraph,
		$$
		x\mapsto \langle f(x), K^*(x)\phi\rangle
		$$
		is measurable for each $\phi\in \mathscr{K}$. But $\langle f(x), K^*(x)\phi\rangle = \langle K(x)f(x),\phi\rangle$, so we have proven that $K(x)f(x)$ is weakly measurable.
	\end{proof}
	
	\begin{defn} Let $f: X\to\mathscr{H}$. \hspace{2cm}
		\begin{list}{(\arabic{ctr})}{\usecounter{ctr}}
			\setcounter{ctr}{0}
			\item $f$ is called \textit{finitely-valued} if it is constant on each of a finite number is disjoint measurable sets $E_j$ and equal to $0$ on $X\setminus \bigcup E_j$.
			\item $f$ is said to be a \textit{simple function} if it is finite-valued and if $\mu(E_j)<\infty$ for each of the sets $E_j$ on which $f$ takes a nonzero value.
			\item $f$ is said to be \textit{countably-valued} if it assumes at most a countable set of values in $\mathscr{H}$, and if each such value is assumed on a measurable set.
			\item $f$ is called \textit{strongly measurable} if there exists a sequence of countably-valued functions converging almost everywhere in $X$ to $f$.
		\end{list}
	\end{defn}
	
	\begin{lemma}\label{StrongMeasConv}
		If $X$ is $\sigma$-finite, then $f: X\to \mathscr{H}$ is strongly measurable if and only if there exists a sequence of simple functions converging almost everywhere in $X$ to $f$.
	\end{lemma}
	\begin{proof}
		First note that one direction of this if and only if statement is trivial: if there is a sequence of simple functions converging to $f$ a.e., then since these functions are, \textit{a fortiori}, countably-valued, $f$ is strongly measurable. To prove the other direction, suppose that $f: X\to \mathscr{H}$ be strongly measurable. Since $X$ is $\sigma$-finite we can partition $X$ as a disjoint union of measurable sets: $X= \bigcup_{m=1}^\infty X_m$. Let $Y_n = \bigcup_{m=1}^n X_m$. Let $\{\psi_n\}_{n=1}^\infty$ be a sequence of countably-valued functions converging pointwise to $f$ on $X-N$, where $N$ is a null set. Since $\phi_n$ is countably-valued, we may write it as
		\[
		\psi_n = \sum_{m=1}^\infty a_m(n) \mathbbm{1}_{E_m(n)},
		\]
		where $E_m(n) = \{x\in X: \psi_n(x) = a_m(n)\}$. Note that each $E_m(n)$ is a measurable set, by definition; that the $E_m(n)$'s may be taken to be disjoint; and that $X= \bigcup_{m=1}^\infty E_m(n)$. Now define  
		\[
		\phi_n =\sum_{m=1}^n a_m(n)\mathbbm{1}_{E_m(n)\cap Y_n};
		\]
		i.e.: $\phi_n$ is the $n$-th partial sum of $\psi_n\mathbbm{1}_{Y_n}$. This gives a sequence $\{\phi_n\}_{n=1}^\infty$ of simple functions. We claim that $\phi_n$ converges to $f$ on $X-N$. To see this, let $x\in X-N$, and let $\epsilon >0$. There exists $N_1(x)\in \mathbb{N}$ such that $\|\psi_n(x) - f(x)\|<\epsilon$ whenever $n\geq N_1(x)$. Moreover, there exists $N_2(x)$ such that $x\in Y_n$ whenever $n\geq N_2(x)$, and there exists $N_3(x)$ such that $x\in \bigcup_{m=1}^n E_m(n)$ whenever $n\geq N_3(x)$. Thus, for $n\geq \max\{N_1(x),N_2(x),N_3(x)\}$, we have
		\[
		\|\phi_n(x)-f(x)\| = \|\psi_n(x)-f(x)\| <\epsilon.
		\]
		Thus $\phi_n$ converges pointwise to $f$ almost everywhere, which proves our claim. 
		
	\end{proof}
	\begin{theorem}\label{WeakStrong}
		If $\mathscr{H}$ is separable, then the notions of weak and strong measurability coincide.
	\end{theorem}
	\begin{proof}
		This is a corollary of theorem of B.J. Pettis. See \cite{HillePhillips}, p.~73, Corollary 2.
	\end{proof}

	\chapter{Technical Measure-theoretic Points and Lemmas}\label{Bappendix}
	\section{Restricting to subsets}
	\begin{defn}\label{relsigalg}
		Let $(X,\mathcal{M})$ be a measurable space, and let $E\subset X$ be any subset. Then $\mathcal{M}\vert_E = \{F\cap E: F\in\mathcal{M}\}$ is a $\sigma$-algebra on $E$, called the \textit{relative $\sigma$-algebra}.
	\end{defn}
	
	\begin{prop}\label{BorelSigmaAlg}
		Let $X$ be a topological space. Given any set $E\subset X$, we can give $E$ the Borel $\sigma$-algebra $\mathcal{B}_E$ generated by the open sets in $E$ (where $E$ is given the relative topology inherited from $X$). 	The Borel $\sigma$-algebra $\mathcal{B}_E$ coincides with the relative $\sigma$-algebra $\mathcal{B}\vert_E=\{F \cap E: F\in \mathcal{B}_X\}$.
	\end{prop}
	\begin{proof}
		By definition,
		\[
		\mathcal{B}_E = \sigma\{U: U\ \text{open in}\ E\} = \sigma\{V\cap E: V\ \text{open in}\ X\}.
		\]
		Moreover, since $V\cap E \in \mathcal{B}\vert_E$ whenever $V\subset X$ is open, we have
		\[
		\mathcal{B}_E=\sigma\{V\cap E: V\ \text{open in}\ X\} \subset  \mathcal{B}\vert_E.
		\]
		On the other hand, since the inclusion map $\iota: E\to X$ is continuous, $\iota^{-1}(F) \in \mathcal{B}_E$ for any $F\in\mathcal{B}_{X}$. But $\iota^{-1}(F) = F\cap E$, whence $\mathcal{B}\vert_E =\{F \cap E: F\in \mathcal{B}_{X}\} \subset \mathcal{B}_E$\footnote{Credit to \cite{inclusionBorel} for the trick with the inclusion map.}. 
	\end{proof}
	
	\begin{defn}\label{MbleonE}
		Let $(X,\mathcal{M})$ and $(Y,\mathcal{N})$ be measurable spaces, and let $f:X\to Y$ be $(\mathcal{M},\mathcal{N})$-measurable. Given $E\in\mathcal{M}$, we say that $f$ is \textit{measurable on $E$} if $f^{-1}(A)\cap E\in\mathcal{M}$ for all $A\in\mathcal{N}$.
	\end{defn}
	\begin{rem}
		With $(X,\mathcal{M})$, $(Y,\mathcal{N})$, $E$, and $f$ as in Definition \ref{MbleonE}, $f$ is measurable on $E$ if and only if $f\vert_E$ is $\mathcal{M}\vert_E$-measurable, where $\mathcal{M}\vert_E$ is the relative $\sigma$-algebra $\mathcal{M}\vert_E = \{F\cap E: F\in\mathcal{M}\}$. Indeed, given a measurable set $A\subset\mathcal{N}$, we have 
		\[
		(f\vert_E)^{-1}(A) = f^{-1}(A)\cap E.
		\] 
		Thus, if $f$ is measurable on $E$, $(f\vert_E)^{-1}(A) \in \mathcal{M}$, and hence 
		\[
		(f\vert_E)^{-1}(A)  = (f\vert_E)^{-1}(A) \cap E \in \mathcal{M}\vert_E.
		\]
		Conversely, if $f\vert_E$ is $\mathcal{M}\vert_E$-measurable, then since	\begin{align*}
			f^{-1}(A) \cap E &= \{x\in E: f(x)\in A\}\\
			&= \{x\in E: f\vert_E(x)\in A\}\\
			&= (f\vert_E)^{-1}(A),
		\end{align*}
		we have $f^{-1}(A)\cap E = (f\vert_E)^{-1}(A) \in\mathcal{M}\vert_E \subset\mathcal{M}$, whence $f$ is measurable on $E$.
	\end{rem}
	\begin{lemma}\label{measonAB}
		With $(X,\mathcal{M})$, $(Y,\mathcal{N})$, and $f$ as in  Definition \ref{MbleonE}, if $X=E\cup F$, where $E, F \in \mathcal{M}$, then $f$ is $(\mathcal{M},\mathcal{N})$-measurable if and only if it is measurable on $E$ and on $F$.
	\end{lemma}
	\begin{proof}
		If $f$ is $(\mathcal{M},\mathcal{N})$-measurable, then $f^{-1}(A) \in\mathcal{M}$ for all $A\in\mathcal{N}$. Thus $f^{-1}(A)\cap E \in\mathcal{M}$ and $f^{-1}(A)\cap F \in\mathcal{M}$, for all $A\in\mathcal{N}$, so that $f$ is measurable on $E$ and on $F$. Conversely, if $f^{-1}(A)\cap E \in\mathcal{M}$ and $f^{-1}(A)\cap F \in\mathcal{M}$, for all $A\in\mathcal{N}$, then
		\begin{align*}
			f^{-1}(A) &=f^{-1}(A) \cap X \\
			&= f^{-1}(A)\cap (E\cup F)\\
			&= [f^{-1}(A)\cap E] \cup [f^{-1}(A)\cap F] \in \mathcal{M}.
		\end{align*}
		Note: In the case where $(Y,\mathcal{N}) = (\mathbb{C},\mathcal{B}_\mathbb{C})$, this is Exercise 5 on p.~48 of \cite{FollandReal}.
	\end{proof}
	
	\begin{rem}\label{L2aevector}
		Let $\mathscr{H}$ be a separable Hilbert space.	Let $(X,\mathcal{M},\mu)$ be a measure space, and let $Y\in\mathcal{M}$ with $\mu(X-Y)=0$. Give $Y$ the relative $\sigma$-algebra $\mathcal{M}\vert_Y$. We show that $L^2(X;\mathscr{H})$ can be identified with $L^2(Y;\mathscr{H})$. Given $f\in L^2(Y;\mathscr{H})$, we note that $x\mapsto \langle f(x), e_n\rangle$ is $\mathcal{M}\vert_Y$-measurable by definition. We define $\tilde{f}$ on $X$ by
		\begin{equation}\label{ftilde}
			\begin{aligned}
				\tilde{f}(x) = \left\{ \begin{array}{ll}
					f(x) &\text{if } x\in Y \\
					0 &\text{if } x\in X-Y
				\end{array} \right. .
			\end{aligned}
		\end{equation}
		and define $f_n(x):=\langle \tilde{f}(x) , e_n\rangle$ on $X$, so that
		\[
		f_n(x) 
		= \left\{ \begin{array}{ll}
			\langle f(x), e_n\rangle  &\text{if } x\in Y \\
			0 &\text{if } x\in X-Y
		\end{array} \right. .
		\]
		By Lemma \ref{measonAB}, $f_n$ is $\mathcal{M}$-measurable precisely when it is measurable on both $Y$ and $X-Y$. Since $f_n\vert_Y$ is just $x\mapsto \langle f(x), e_n\rangle$, which is $\mathcal{M}\vert_Y$-measurable, $f_n$ is measurable on $Y$. And since $f_n\vert_{X-Y} \equiv 0$, it follows that $f_n$ is measurable on $X-Y$. Therefore $f_n$ is $\mathcal{M}$-measurable. It follows from Remark \ref{weakmeasITObases} that $\tilde{f}$ is a weakly measurable function from $X$ to $\mathscr{H}$. Therefore defining $\iota: L^2(Y;\mathscr{H}) \to L^2(X;\mathscr{H})$ by
		\[
		\iota([f]_{L^2(Y;\mathscr{H})}) = [\tilde{f}]_{L^2(X;\mathscr{H})}
		\] 
		gives a surjective isometry from $L^2(Y;\mathscr{H})$ onto $L^2(X;\mathscr{H})$. Indeed, given $g \in L^2(X;\mathscr{H})$, we have $g= \iota(g\vert_Y)$, since $\mu(X-Y)=0$; and $g\vert_Y$ is a weakly measurable $\mathscr{H}$-valued function on $(Y,\mathcal{M}\vert_Y)$, as we show presently. Since $g$ is weakly measurable on $X$, $x\mapsto \langle g(x),e_n\rangle$ is $\mathcal{M}$-measurable---and therefore measurable on $Y$---for all $n$. Hence $y\mapsto \langle g\vert_Y(y),e_n\rangle$ is $\mathcal{M}\vert_Y$-measurable for all $n$, which means that $g\vert_Y$ is weakly measurable.  Furthermore, $\|g\vert_Y\|^2 = \int_Y \|g(x)\|_\mathscr{H}^2 \ud\mu(x) = \int_X \|g(x)\|_\mathscr{H}^2 \ud\mu(x) <\infty$. This shows that $\iota$ is surjective. Lastly, we check that $\iota$ is an isometry:
		\begin{align*}
			\|\iota([f]_{L^2(Y;\mathscr{H})})\|^2 &= \|[\tilde{f}]_{L^2(X;\mathscr{H})}\|^2\\
			&= \int_X \|\tilde{f}(x)\|^2\ud \mu(x)\\
			&= \int_Y \|f(x)\|^2\ud \mu(x)+  \int_{X-Y} 0^2\ud \mu(x) \\
			&=\int_Y \|f(x)\|^2\ud \mu(x)\\
			&=\|[f]_{L^2(Y;\mathscr{H})}\|^2.
		\end{align*}
	\end{rem}
	We note that the proof of the above remark might be streamlined by proving it in the scalar case and noting that $L^2(X;\mathscr{H}) \cong L^2(X)\otimes \mathscr{H}$ (Theorem II.10(b) in \cite{ReedSimon}).
	
	\section{Using Complete Measure Spaces}
	\begin{prop}\label{BorelClosure}
		Euclidean space $\mathbb{R}^n$ is usually given the Lebesgue $\sigma$-algebra $\mathscr{L}$ and Lebesgue measure $m$, where $\mathscr{L}$ is the completion of $\mathcal{B}_{\mathbb{R}^n}$ with respect to $m$. Given any set $E\subset \mathbb{R}^n$, we can give $E$ the relative Borel $\sigma$-algebra $\mathcal{B}\vert_E$ (and this coincides with the Borel $\sigma$-algebra $\mathcal{B}_E$ generated by the open sets in $E$, by \ref{BorelSigmaAlg}). If $E$ is a Borel set, then the relative $\sigma$-algebra $\mathscr{L}\vert_E=\{F\cap E: F\in \mathscr{L}\}$ is the completion of $\mathcal{B}\vert_E$ with respect to $m$.
	\end{prop}
	\begin{proof} 
		Suppose $F\in\overline{\mathcal{B}\vert_E}$, so that $F=F_1\cup N_0$, where $F_1 \in \mathcal{B}\vert_E$ and $N_0\subset N$ for some $N\subset\mathcal{B}\vert_E$ with $m(N)=0$. Then $F_1=F_2\cap E$, for some $F_2\in \mathcal{B}_{\mathbb{R}^n}$; and $N =N_1\cap E$, for some $N_1\in \mathcal{B}_{\mathbb{R}^n}$. Since $N_0\subset N\subset E$, we have
		\begin{align*}
			F &= F_1\cup N_0 \\
			&= (F_2\cap E) \cup N_0\\
			&= (F_2 \cup N_0) \cap (E\cup N_0)\\
			&=  (F_2 \cup N_0) \cap E.
		\end{align*}
		Clearly, $F_2 \cup N_0 \in \mathscr{L}$, since $F_2 \in \mathcal{B}_{\mathbb{R}^n}$ and $N_0\subset N \in \mathcal{B}_{\mathbb{R}^n}$ (since $E\in\mathcal{B}_{\mathbb{R}^n}$) with $m(N) =0$; hence  $F\in \mathscr{L}\vert_E$. Therefore, $\overline{\mathcal{B}\vert_E} \subset\mathscr{L}\vert_E$. Conversely, if $F \in \mathscr{L}\vert_E$, then $F=F_1 \cap E$ for some $F_1 \in\mathscr{L}$. Since $\mathscr{L}$ is the completion of $\mathcal{B}_{\mathbb{R}^n}$, we have $F_1 = F_2 \cup N_0$, where $F_2 \in \mathcal{B}_{\mathbb{R}^n}$ and $N_0 \subset N$ for some $N\in \mathcal{B}_{\mathbb{R}^n}$ such that $m(N)=0$. Thus
		\begin{align*}
			F&=F_1\cap E\\
			&= (F_2\cup N_0)\cap E\\
			&= (F_2\cap E) \cup (N_0\cap E).
		\end{align*}
		Here $F_2\cap E \in\mathcal{B}\vert_E$; and $N_0\cap E \subset N\cap E \in \mathcal{B}\vert_E$, with $m(N\cap E) =0$, since $N\cap E \subset N$ and $m(N)=0$. Thus $F \in \overline{\mathcal{B}\vert_E}$, whence $\mathscr{L}\vert_E \subset\overline{\mathcal{B}\vert_E}$.
	\end{proof}
	
	\begin{prop}\label{SepBorelAlg}
		Let $X$ and $Y$ be separable metric spaces. Then the Borel $\sigma$-algebra $\mathcal{B}_{X\times Y}$ on $X\times Y$ coincides with the product $\sigma$-algebra $\mathcal{B}_X\times\mathcal{B}_Y$.
	\end{prop}
	\begin{proof}
		This is Proposition 1.5 in \cite{FollandReal}.
	\end{proof}
	\begin{lemma}\label{ProductClosure}
		Let $(X,\mathcal{M},\mu)$ be a complete measure space, and let $Y$ be a countable set with $\sigma$-algebra $\mathcal{P}(Y)$\label{PY} and counting measure $\nu$. Then the product measure $\mu\times \nu$ is complete on $\big(X\times Y, \mathcal{M}\times\mathcal{P}(Y)\big)$.
	\end{lemma}
	\begin{proof}
		We note that the $\sigma$-algebra $\mathcal{M}\times\mathcal{P}(Y)$ on $X\times Y$ is generated by rectangles $E\times F$, where $E \in \mathcal{M}$ and $F\subset Y$. Let $A\in \mathcal{M}\times\mathcal{P}(Y)$, and for each $y\in Y$ let $A^y= \{x\in X: (x,y)\in A\}$ denote the $y$-section of $A$. Each $A^y \in\mathcal{M}$ by Proposition 2.34(a) in \cite{FollandReal}. Moreover, $A$ can be written as
		\[
		A= \bigcup_{y\in Y} A^{y}\times \{y\}.
		\]
		Indeed, it is immediate that if $(x,y) \in A^y\times \{y\}$, then $x\in A^y$, and hence by definition, $(x,y) \in A$. Conversely, if $(x,y)\in A$, then $x\in A^y$, and hence $(x,y) \in  A^y\times \{y\}$. Now let $A \in  \mathcal{M}\times\mathcal{P}(Y)$ with $\mu\times\nu(A) =0$, and let $B\subset A$. Then
		\begin{equation}\label{Bdecomp}
			B= \bigcup_{y\in Y} B^y\times \{y\}
		\end{equation}
		where $B^y$ denotes the $y$-section of $B$. Now $B^y = \{x\in X: (x,y) \in B\} \subset A^y$, and $A^y \in \mathcal{M}$ with $\mu(A^y) = 0$ for a.e. $y\in Y$ (by Exercise 49(a) on p.~69 of \cite{FollandReal}). Since $Y$ has the discrete $\sigma$-algebra with counting measure, this means that $\mu(A^y) = 0$ for all $y\in Y$. Since $\mu$ is a complete measure on $(X,\mathcal{M})$, it follows that $B^y \in \mathcal{M}$ for all $y$. Thus \eqref{Bdecomp} shows that $B\in \mathcal{M}\times\mathcal{P}(Y)$. This concludes the proof that $\big(\mathcal{M}\times\mathcal{P}(Y),\mu\times\nu\big)$ is complete.
	\end{proof}
	\begin{cor}\label{completionofproduct}
		Let $(X,\mathcal{M},\mu)$ be a measure space, and let $Y$ be a countable set with $\sigma$-algebra $\mathcal{P}(Y)$ and counting measure $\nu$. Let $\overline{\mathcal{M}}$ be the completion of $\mathcal{M}$ with respect to $\mu$, and let $\overline{\mu}$ be the completion of $\mu$. Then $\overline{\mathcal{M}}\times\mathcal{P}(Y)$ is the completion of $\mathcal{M}\times\mathcal{P}(Y)$ with respect to $\mu\times\nu$, and $\overline{\mu\times \nu} = \overline{\mu}\times \nu$.
	\end{cor}
	\begin{proof} 
		Recall that 
		\begin{multline}\label{completeprod}
			\overline{\mathcal{M}\times\mathcal{P}(Y)} =\{E\cup F_0: E\in \mathcal{M}\times\mathcal{P}(Y),\ \text{and}\ F_0 \subset  F,\\ \text{where}\ F\in \mathcal{M}\times\mathcal{P}(Y)\ \text{with}\ (\mu\times \nu)(F)= 0\}.
		\end{multline}
		Given $E \cup F_0 \in \overline{\mathcal{M}\times\mathcal{P}(Y)}$ as in \eqref{completeprod}, we note that $\mathcal{M}\times\mathcal{P}(Y) \subset \overline{\mathcal{M}}\times\mathcal{P}(Y)$, and hence	$E \in \overline{\mathcal{M}}\times\mathcal{P}(Y)$. Now $F_0 \subset F$, with 
		\[
		(\overline{\mu}\times\nu)(F)=	(\mu\times\nu)(F) =0.
		\]
		Since $\overline{\mu}\times\nu$ is complete by Lemma \ref{ProductClosure}, it follows that $F_0 \in \overline{\mathcal{M}}\times\mathcal{P}(Y)$. Thus $E\cup F_0 \in \overline{\mathcal{M}}\times\mathcal{P}(Y)$, and hence we have shown 
		\[
		\overline{\mathcal{M}\times\mathcal{P}(Y)} \subset \overline{\mathcal{M}}\times\mathcal{P}(Y).
		\]
		On the other hand, let $E\times F$ be a measurable rectangle in $\overline{\mathcal{M}}\times\mathcal{P}(Y)$. Then $E = A \cup B_0$, where $A\in \mathcal{M}$ and $B_0\subset B$ with $\mu(B) =0$. Then 
		\[
		E\times F = (A\cup B_0) \times F = (A\times F) \cup (B_0\times F).
		\]
		It is clear that $A\times F \in \mathcal{M}\times \mathcal{P}(Y)$, and $B_0 \times F \subset B\times F \in \mathcal{M}\times\mathcal{P}(Y)$, with $(\mu\times\nu)(B\times F) = \mu(B)\nu(F) =0$. Thus $B_0 \times F \in \overline{\mathcal{M}\times\mathcal{P}(Y)}$. It follows that 
		\[
		E\times F = (A\times F) \cup (B_0\times F) \in \overline{\mathcal{M}\times\mathcal{P}(Y)}.
		\]
		Since $E\times F$ was an arbitrary measurable rectangle in $\overline{\mathcal{M}}\times\mathcal{P}(Y)$, and this $\sigma$-algebra is generated by such measurable rectangles, we have
		\[
		\overline{\mathcal{M}}\times\mathcal{P}(Y)\subset \overline{\mathcal{M}\times\mathcal{P}(Y)}.
		\]
		Lastly, we note that by Theorem 1.9 in \cite{FollandReal}, there is a \textit{unique} extension $\overline{\mu\times\nu}$ of $\mu\times\nu$ to a complete measure on $\overline{\mathcal{M}\times\mathcal{P}(Y)}$. Since $\overline{\mu\times\nu}$ and $\overline{\mu} \times \nu$ are both complete measures on $\overline{\mathcal{M}\times\mathcal{P}(Y)} = 	\overline{\mathcal{M}}\times\mathcal{P}(Y)$, it follows that they are equal.
	\end{proof}
	
	\section{Images of measures}
	\begin{defn}\label{imegameas}
		Let $(X,\mathcal{M})$ and $(Y,\mathcal{N})$ be measurable spaces, and let $f: X\to Y$ be an $(\mathcal{M},\mathcal{N})$-measurable function. For any (nonnegative) measure $\mu$ on $(X,\mathcal{M})$, we obtain a measure on $(Y,\mathcal{N})$ by
		\[
		E \mapsto \mu\big(f^{-1}(E)\big),\ \ \ E\in \mathcal{N}.
		\]
		This measure is called the \textit{image of the measure} $\mu$ under the mapping $f$, and is denoted by $\mu\circ f^{-1}$.
	\end{defn}
	\begin{theorem}[Theorem 3.6.1 in \cite{Bogachev}; or Proposition 2.6.8 in \cite{Cohn}]\label{imagemeasurethm}
		Let $(X,\mathcal{M})$ and $(Y,\mathcal{N})$ be measurable spaces, and let $f: X\to Y$ be an $(\mathcal{M},\mathcal{N})$-measurable function. Let $\mu$ be a measure on $X$. An $\mathcal{N}$-measurable function $g$ on $Y$ is integrable with respect to the measure $\mu\circ f^{-1}$ precisely when the function $g\circ f$ is integrable with respect to $\mu$. In addition, one has
		\begin{equation}\label{imageofmeasure}
			\int_Y g(y) \ud(\mu\circ f^{-1})(y) = \int_X g(f(x))\ud\mu(x).
		\end{equation}
	\end{theorem}
	\begin{proof}
		We first show that $g\circ f$ is $\mathcal{M}$-measurable whenever $g$ is $\mathcal{N}$-measurable. Let $F\subset \mathbb{C}$ be a Borel set. Then $g^{-1}(F)\in\mathcal{N}$. Moreover, since
		\begin{align*}
			(g\circ f)^{-1}(F) &= \{x\in X: (g\circ f)(x)\in F\}\\
			&=  \{x\in X: g(f(x)) \in F\}\\
			&=\{ x\in X: f(x) \in g^{-1}(F)\}\\
			&= f^{-1}\big(g^{-1}(F)\big).
		\end{align*}
		and $f$ is $(\mathcal{M},\mathcal{N})$-measurable, $(g\circ f)^{-1}(F)\in\mathcal{M}$. Hence $g\circ f$ is $\mathcal{M}$-measurable  whenever $g$ is $\mathcal{N}$-measurable
		
		Let $E\in\mathcal{N}$. Then 
		\[
		\int_Y \mathbbm{1}_{E}\ud(\mu\circ f^{-1})(y) =(\mu\circ f^{-1})(E)= \mu\big(f^{-1}(E)\big) = \int_X \mathbbm{1}_{f^{-1}(E)}\ud\mu(x).
		\]
		Moreover, $\mathbbm{1}_{f^{-1}(E)} = \mathbbm{1}_E\circ f$, and hence
		\[
		\int_Y \mathbbm{1}_{E}\ud(\mu\circ f^{-1})(y) =\int_X \mathbbm{1}_E\circ f\ud\mu(x),
		\]
		so that \eqref{imageofmeasure} holds for indicators of measurable sets. By linearity, it follows that \eqref{imageofmeasure} holds for simple functions. Let $g$ be a nonnegative $\mathcal{N}$-measurable function. It follows from Proposition 2.10(a) in \cite{FollandReal} that there exists a sequence $\{\phi_n\}$ of simple functions such that $0\leq \phi_1\leq\phi_2 \leq \cdots \leq g$, and $\phi_n \to g$ pointwise. It then follows from the Monotone Convergence Theorem (2.14 in \cite{FollandReal}) that 
		\[
		\int_Y g \ud(\mu\circ f^{-1})(y) = \lim_{n\to\infty} \int_Y  \phi_n \ud(\mu\circ f^{-1})(y) =  \lim_{n\to\infty} \int_X \phi_n \circ f\ud\mu(x) = \int_X g \circ f\ud\mu(x),
		\]
		where both sides may be infinite. Thus \eqref{imageofmeasure} holds for all nonnegative measurable functions. In particular, if $g$ is an $\mathcal{N}$-measurable function, applying the above to $|g|$ gives
		\[
		\int_Y |g| \ud(\mu\circ f^{-1})(y)  = \int_X |g| \circ f\ud\mu(x) =\int_X |g \circ f| \ud\mu(x),
		\]
		so that $g\circ f$ is integrable with respect to $\mu$ precisely when $g$ is integrable with respect to $\mu\circ f^{-1}$. Finally, it follows from the linearity of the integral that \eqref{imageofmeasure} holds for all $\mathbb{C}$-valued, $\mathcal{N}$-measurable functions $g$, since any such function can be expressed as a linear combination of nonnegative measurable functions.
	\end{proof}
	\begin{cor}\label{L2isom}
		Let $\mu$ be a measure on $(X,\mathcal{M})$, and let $f: X\to Y$ be an $(\mathcal{M},\mathcal{N})$-isomorphism, i.e., an $(\mathcal{M},\mathcal{N})$-measurable bijection whose inverse is $(\mathcal{N},\mathcal{M})$-measurable. Then 
		\[
		\Phi: L^2(Y,\mathcal{N},\mu\circ f^{-1})\to L^2(X,\mathcal{M},\mu)
		\]
		given by
		\[
		\Phi(g) = g\circ f
		\]
		is an isometric isomorphism, with inverse $h\mapsto h\circ f^{-1}$.
	\end{cor}
	\begin{proof}
		Let $g$ be an $\mathcal{N}$-measurable function. Applying theorem \ref{imagemeasurethm} to $|g|^2$, we see that $g$ is square-integrable with respect to $\mu\circ f^{-1}$ precisely when $g\circ f$ is square-integrable with respect to $\mu$, and moreover,
		\[
		\|g\|^2 =\int_Y |g|^2 \ud(\mu\circ f^{-1})(y)  =\int_X |g \circ f|^2 \ud\mu(x) = \|g\circ f\|^2= \|\Phi(g)\|^2.
		\]
		Thus $\Phi$ is an isometry. Let $h\in L^2(X,\mathcal{M},\mu)$. We show that $h\circ f^{-1}\in L^2(Y,\mathcal{N},\mu\circ f^{-1})$. First, $h\circ f^{-1}$ is an $\mathcal{N}$-measurable function, because
		\begin{align*}
			(h\circ f^{-1})^{-1}(F) &= \{x\in X: (h\circ f^{-1})(x)\in F\}\\
			&=  \{x\in X: h(f^{-1}(x)) \in F\}\\
			&=\{ x\in X: f^{-1}(x) \in h^{-1}(F)\}\\
			&= (f^{-1})^{-1}\big(h^{-1}(F)\big)\in\mathcal{N},
		\end{align*}
		test
		because $f^{-1}$ is $(\mathcal{N},\mathcal{M})$-measurable, and $h^{-1}(F) \in \mathcal{M}$ for any Borel set $F\subset\mathbb{C}$. Thus we may apply theorem \ref{imagemeasurethm} to $|h\circ f^{-1}|^2$ to get 
		\[
		\|h\circ f^{-1}\|^2 = \int_Y |h\circ f^{-1}|^2 \ud(\mu\circ f^{-1})(y)  =\int_X |h\circ f^{-1} \circ f|^2 \ud\mu(x) = \|h\|^2 <\infty,
		\]
		so that $h\circ f^{-1}\in L^2(Y,\mathcal{N},\mu\circ f^{-1})$. Since $h = \Phi(h\circ f^{-1})$, we have shown that $\Phi$ is a surjective isometry, and $\Phi^{-1}(h)= h\circ f^{-1}$.
	\end{proof}
	
	\begin{theorem}[Image of Measure---Extended]\label{imagemeasurethmE}\label{imagemeasureextended}
		Let $\mu$ be a complete measure on $(X,\mathcal{M})$, and let $f: X\to Y$ be $(\mathcal{M},\mathcal{N})$-measurable. Let $\overline{\mathcal{N}}$ denote the completion of $\mathcal{N}$ with respect to $\mu\circ f^{-1}$. An $\overline{\mathcal{N}}$-measurable function $g$ on $Y$ is integrable with respect to the completed measure $\overline{\mu\circ f^{-1}}$ precisely when the function $g\circ f$ is integrable with respect to $\mu$. In addition, one has
		\begin{equation}\label{imageofmeasureEX}
			\int_Y g(y) \ud(\overline{\mu\circ f^{-1}})(y) = \int_X g(f(x))\ud\mu(x).
		\end{equation}
	\end{theorem}
	\begin{proof}
		The proof is very similar to the proof of Theorem \ref{imagemeasurethm}; we include it nonetheless, for completeness. We first show that $g\circ f$ is $\mathcal{M}$-measurable whenever $g$ is $\overline{\mathcal{N}}$-measurable. By Proposition 2.12 in \cite{FollandReal}, there exists an $\mathcal{N}$-measurable function $g_0$ such that $g=g_0$ on the complement of some set $N \in\mathcal{N}$ with $(\mu\circ f^{-1})(N) =0$. It follows now from Theorem \ref{imagemeasurethm} that $g_0\circ f$ is $\mathcal{M}$-measurable. Moreover, since $\mu\big(f^{-1}(N)\big) =0$, and $g_0\circ f = g\circ f$ on $X-f^{-1}(N)$, it follows from Proposition \ref{MeasAEresult} that $g\circ f$ is $\mathcal{M}$-measurable.
		
		Let $E\in\overline{\mathcal{N}}$; thus $E=F\cup N_0$, where $F \in\mathcal{N}$ and $N_0\subset N_1$, where $N_1 \in\mathcal{N}$ with $(\mu\circ f^{-1})(N_1) =0$. Since $f^{-1}(N_0)\subset f^{-1}(N_1)$, and $\mu\big(f^{-1}(N_1)\big) = 0$, it follows from the completeness of $\mu$ that  $f^{-1}(N_0) \in \mathcal{M}$ with $\mu\big(f^{-1}(N_0)\big) =0$. Since $f^{-1}(E) = f^{-1}(F)\cup f^{-1}(N_0)$, we see that 
		\begin{align*}
			(\overline{\mu\circ f^{-1}})(E) &= 	(\overline{\mu\circ f^{-1}})(F\cup N_0) \\
			&= (\mu\circ f^{-1})(F) \\
			&=  \mu\big(f^{-1}(F) \big)\\
			&= \mu\big(f^{-1}(F) \cup f^{-1}(N_0) \big)\\
			&= \mu\big( f^{-1}(F \cup N_0) \big)\\
			&= \mu\big(f^{-1}(E)\big).
		\end{align*}
		Thus $\overline{\mu\circ f^{-1}} = \mu\circ f^{-1}$. Thus
		\[
		\int_Y \mathbbm{1}_{E}\ud(\overline{\mu\circ f^{-1}})(y) =(\overline{\mu\circ f^{-1}})(E)= \mu\big(f^{-1}(E)\big) = \int_X \mathbbm{1}_{f^{-1}(E)}\ud\mu(x).
		\]
		Moreover, $\mathbbm{1}_{f^{-1}(E)} = \mathbbm{1}_E\circ f$, and hence
		\[
		\int_Y \mathbbm{1}_{E}\ud(\overline{\mu\circ f^{-1}})(y) =\int_X \mathbbm{1}_E\circ f\ud\mu(x),
		\]
		so that \eqref{imageofmeasureEX} holds for indicators of measurable sets. By linearity, it follows that \eqref{imageofmeasureEX} holds for simple functions. Let $g$ be a nonnegative $\overline{\mathcal{N}}$-measurable function. It follows from Proposition 2.10(a) in \cite{FollandReal} that there exists a sequence $\{\phi_n\}$ of simple functions such that $0\leq \phi_1\leq\phi_2 \leq \cdots \leq g$, and $\phi_n \to g$ pointwise. It then follows from the Monotone Convergence Theorem (2.14 in \cite{FollandReal}) that 
		\[
		\int_Y g \ud(\overline{\mu\circ f^{-1}})(y) = \lim_{n\to\infty} \int_Y  \phi_n \ud(\overline{\mu\circ f^{-1}})(y) =  \lim_{n\to\infty} \int_X \phi_n \circ f\ud\mu(x) = \int_X g \circ f\ud\mu(x),
		\]
		where both sides may be infinite. Thus \eqref{imageofmeasureEX} holds for all nonnegative measurable functions. In particular, if $g$ is an $\overline{\mathcal{N}}$-measurable function, applying the above to $|g|$ gives
		\[
		\int_Y |g| \ud(\overline{\mu\circ f^{-1}})(y)  = \int_X |g| \circ f\ud\mu(x) =\int_X |g \circ f| \ud\mu(x),
		\]
		so that $g\circ f$ is integrable with respect to $\mu$ precisely when $g$ is integrable with respect to $\overline{\mu\circ f^{-1}}$. Finally, it follows from the linearity of the integral that \eqref{imageofmeasureEX} holds for all $\mathbb{C}$-valued, $\overline{\mathcal{N}}$-measurable functions $g$, since any such function can be expressed as a linear combination of nonnegative measurable functions.
	\end{proof}
	\begin{cor}\label{imagemeasExcor}
		Let $\mu$ be a complete measure on $(X,\overline{\mathcal{B}_X})$, and let $f: X \to Y$ be a Borel isomorphism (i.e. a Borel bijection whose inverse is Borel also). Then 	
		\[
		\Phi: L^2(Y,\overline{\mathcal{B}_Y},\overline{\mu\circ f^{-1}})\to L^2(X,\overline{\mathcal{B}_X},\mu)
		\]
		given by
		\[
		\Phi(g) = g\circ f
		\]
		is an isometric isomorphism, with inverse $h\mapsto h\circ f^{-1}$.
	\end{cor}
	\begin{proof}
		The proof that $\Phi$ is an isometry follows from Theorem \ref{imagemeasurethmE}, just as in Corollary \ref{L2isom}. Indeed, 	given $g$ an $\overline{\mathcal{N}}$-measurable function, applying theorem \ref{imagemeasurethmE} to $|g|^2$, we see that $g$ is square-integrable with respect to $\overline{\mu\circ f^{-1}}$ precisely when $g\circ f$ is square-integrable with respect to $\mu$, and moreover,
		\[
		\|g\|^2 =\int_Y |g|^2 \ud(\overline{\mu\circ f^{-1}})(y)  =\int_X |g \circ f|^2 \ud\mu(x) = \|g\circ f\|^2= \|\Phi(g)\|^2.
		\]
		Given $h\in  L^2(X,\overline{\mathcal{B}_X},\mu)$, we show that $h\circ f^{-1} \in  L^2(Y,\overline{\mathcal{B}_Y},\overline{\mu\circ f^{-1}})$. By Proposition 2.12 in \cite{FollandReal}, there exists a Borel set $N\subset X$ of measure zero, and a $\mathcal{B}_X$-measurable function $h_0$, such that $h=h_0$ on $X-N$. Thus $h=h_0$ in $L^2(X,\overline{\mathcal{B}_X},\mu)$, and by Basic version of Theorem, there exists $g \in L^2(Y,\mathcal{B}_Y,\mu\circ f^{-1})$ such that $\Phi(g) = h$.

		Thus $h\circ f^{-1} = h_0\circ f^{-1}$ on $Y-f(N)$, and $f(N)$ is a Borel set (being the preimage of $N$ under $f^{-1}$, since $f$ is a bijection), with $(\mu\circ f^{-1})\big(f(N)\big)=0$. It follows from Proposition \ref{MeasAEresult} that $h\circ f^{-1}$ is $\overline{\mathcal{B}_Y}$-measurable. As in Corollary \ref{L2isom}, applying Theorem \ref{imagemeasurethmE} to $|h\circ f^{-1}|^2$, we have
		\[
		\|h\circ f^{-1}\|^2 = \int_Y |h\circ f^{-1}|^2 \ud(\overline{\mu\circ f^{-1}})(y)  =\int_X |h\circ f^{-1} \circ f|^2 \ud\mu(x) = \|h\|^2 <\infty,
		\]
		so that $h\circ f^{-1} \in L^2(Y,\overline{\mathcal{B}_Y},\overline{\mu\circ f^{-1}})$. Since $h = \Phi(h\circ f^{-1})$, we have shown that $\Phi$ is a surjective isometry with $\Phi^{-1}(h)= h\circ f^{-1}$.
	\end{proof}
	\begin{lemma}\label{imagemeascomp}
		Let $(X,\mathcal{M})$, $(Y,\mathcal{N})$, and  $(Z,\mathcal{O})$ be measurable spaces, and let $f: X\to Y$ be an $(\mathcal{M},\mathcal{N})$-measurable function, and let $g: Y\to Z$ be an $(\mathcal{N},\mathcal{O})$-measurable function. Suppose that $\mu$ is a measure on $(X,\mathcal{M})$. We can give $(Y,\mathcal{N})$ the image measure $\mu\circ f^{-1}$, and give $(Z,\mathcal{O})$ the image measure $(\mu\circ f^{-1}) \circ g^{-1}$. Then $f\circ g$ is $(\mathcal{M},\mathcal{O})$-measurable, and 
		\[
		(\mu\circ f^{-1}) \circ g^{-1} = \mu\circ (f\circ g)^{-1}.
		\]
	\end{lemma}
	\begin{proof}
		The first assertion is trivial (see Proposition 2.6.1 in \cite{Cohn}). Given any $A\in \mathcal{O}$, 
		\begin{align*}
			[(\mu\circ f^{-1}) \circ g^{-1}](A) &= (\mu \circ f^{-1})\big(g^{-1}(A)\big)\\
			&= \mu\big( f^{-1}\big(g^{-1}(A)\big)\big)\\
			&=\mu\big((g\circ  f)^{-1}(A)\big)\\
			&= [\mu\circ (f\circ g)^{-1}](A),
		\end{align*}
		so that $(\mu\circ f^{-1}) \circ g^{-1} = \mu\circ (f\circ g)^{-1}$.
	\end{proof}
	
	\section[Proof that the image measure \texorpdfstring{$m\circ \upsilon^{-1}\,$}{}is Haar]{Proof that the image measure \texorpdfstring{$\boldsymbol{m\circ \upsilon^{-1}}\,$}{}is Haar}\label{MOUHaar}
	Recall that $\upsilon$ is given by the composition of mappings
	\[
	\begin{tikzcd} 
		\Omega_{\mathrm{T}^*}\arrow[to=1-3, bend right=30, "\upsilon"{name=D}] \arrow[r, "X_1"{name=X1}]
		& \mathbb{R}^n/\mathcal{L}^* \arrow[r, "\mathcal{C}"{name=C}]
		& \widehat{\mathrm{T}} \, ,
	\end{tikzcd}
	\]
	where
	\[
	X_1(\omega) = \omega+\mathcal{L}^*,
	\]
	and 
	\[
	\mathcal{C}(\omega+\mathcal{L}^*)= \chi_\omega^\mathrm{T},
	\]
	where $\chi_\omega^\mathrm{T}([y,\mathrm{id}]) = e^{2\pi i\omega\cdot y}$. Consider natural quotient map
	\[
	q: \mathbb{R}^n \to \mathbb{R}^n/\mathcal{L}^*
	\]
	given by $q(x) = x+\mathcal{L}^*$. Then $X_1$ is just the restriction of $q$ to $\Omega_{\mathrm{T}^*}$. By definition of the quotient topology on $\mathbb{R}^n/\mathcal{L}^*$, $q$ is \textit{strongly continuous}, meaning that $U \subset \mathbb{R}^n/\mathcal{L}^*$ is open if and only if the inverse image $q^{-1}(U)$ is open in $\mathbb{R}^n$.
	\begin{defn}[Cf. \cite{DixmierC}, p.~397]
		The quotient space $\mathbb{R}^n/\mathcal{L}^*$ has the \textit{quotient Borel structure}: i.e.: the $\sigma$-algebra $\mathcal{B}_{\mathbb{R}^n/\mathcal{L}^*}$ given by 
		\[
		\mathcal{B}_{\mathbb{R}^n/\mathcal{L}^*}  := \{E \in \mathbb{R}^n/\mathcal{L}^*: q^{-1}(E) \in \mathcal{B}_{\mathbb{R}^n}\}
		\]
	\end{defn}
	Note that $\mathcal{B}_{\mathbb{R}^n/\mathcal{L}^*}$ is indeed a $\sigma$-algebra, by Corollary 1.2.9, p.~6, \cite{Bogachev}.
	\begin{prop}
		The Borel $\sigma$-algebra induced by the quotient topology on $\mathbb{R}^n/\mathcal{L}^*$ is contained in the quotient Borel structure.
	\end{prop}
	\begin{proof}
		Let $\mathbb{R}^n/\mathcal{L}^*$ be given the quotient topology and let $\tilde{\mathcal{B}}$ denote the associated Borel $\sigma$-algebra on $\mathbb{R}^n/\mathcal{L}^*$. If $U \in \mathbb{R}^n/\mathcal{L}^*$ is open, then $q^{-1}(U)$ is open in $\mathbb{R}^n$, so that $U \in \mathcal{B}_{\mathbb{R}^n/\mathcal{L}^*}$. It follows that
		\[
		\tilde{\mathcal{B}}= \sigma\{U\subset \mathbb{R}^n/\mathcal{L}^*: U\ \text{is open}\} \subset \mathcal{B}_{\mathbb{R}^n/\mathcal{L}^*}.
		\]
	\end{proof}

	\begin{lemma}[Exercise 10 on p.~27 of \cite{FollandReal}]
		 Given a measure space $(X,\mathcal{M},\mu)$ and $E\in\mathcal{M}$, $\mu\vert_E(A) = \mu(A\cap E)$ defines a measure on $\mathcal{M}$.
		\end{lemma}
	\begin{lemma}[Exercise 7 on p.~220 of \cite{FollandReal}]
	Given a $\sigma$-finite Radon measure $\mu$, on $(X,\mathcal{B}_X)$, and $E\in \mathcal{B}_X$, the restricted Borel measure $\mu\vert_E$ given by $\mu\vert_E(A) = \mu(F\cap A)$ defines a Radon measure on $\mathcal{B}_X$. 
	\end{lemma}		
	\begin{lemma}\label{radonmeasaeequal} 
	If $N\in\mathcal{M}$ is a null set, then $\mu\vert_{X-N} = \mu$.
	\end{lemma}
		\begin{proof}
		For any $F\in\mathcal{M}$, we have
		\begin{align*}
			\mu(F) &= \mu\big(F\cap [N\,\dot{\cup} \,(X-N)]\big)\\
			&= \mu\big( [F\cap N]\, \dot{\cup}\, [F\cap (X-N)]\big)\\
			&= \mu(F\cap N) + \mu\big(F\cap (X-N)\big)\\
			&=\mu\big(F\cap (X-N)\big)\\
			&=\mu\vert_{X-N}(F).
		\end{align*}
	\end{proof}
	\begin{lemma}\label{imageofaeequalfns}
		If $(Y,\mathcal{N})$ is a measurable space, and $f,g: X\to Y$ are $(\mathcal{M},\mathcal{N})$-measurable functions such that $f =g$ on $\mu$-a.e., then $\mu\circ f^{-1} = \mu \circ g^{-1}$.
		\end{lemma}
\begin{proof}
	If $f =g$ on $X-N$, where $N\in\mathcal{M}$ is $\mu$-null, then for each $A\in\mathcal{N}$ we have
		\begin{align*}
			(\mu\circ f^{-1})(A) &= (\mu\vert_{X-N} \circ f^{-1})(A)  \\
			&=\mu\vert_{X-N} \big(f^{-1}(A)\big)\\
			&=\mu \big(f^{-1}(A)\cap(X-N)\big)\\
			&=\mu \big(g^{-1}(A)\cap(X-N)\big)\\
			&= (\mu\circ g^{-1})(A);
		\end{align*}
		whence $\mu\circ f^{-1}=\mu\circ g^{-1}$.
	\end{proof}
	\begin{lemma}\label{imageofhomeoradon}
		Let $(X,\mathcal{B}_X)$ and $(Y,\mathcal{B}_Y)$ be Borel spaces, and let $\mu$ be a Radon measure on $(X,\mathcal{B}_X)$. If $f: X\to Y$ is a homeomorphism, then $\mu\circ f^{-1}$ is a Radon measure on $(Y,\mathcal{B}_Y)$.
	\end{lemma}
	\begin{proof}
		First note that $\mu\circ f^{-1}$ is finite on compact sets, since $K\subset Y$ is compact if and only if $f^{-1}(K)\subset X$ is compact, and $\mu$ is Radon on $(X,\mathcal{B}_X)$. We show that $\mu\circ f^{-1}$ is outer regular on Borel sets, i.e.:
		\[
		(\mu\circ f^{-1})(E) = \inf\{(\mu\circ f^{-1})(U):U\supset E,\ U\ \text{open in }Y\}\\
		\]
		for all $E\in \mathcal{B}_Y$. Let $E\in \mathcal{B}_Y$. Since $f$ is a homeomorphism,
		\begin{multline*}
			\{V:  f^{-1}(E)\subset V\ \text{and } V\ \text{open in }X\} \\
			= 	\{f^{-1}(U):f^{-1}(E)\subset f^{-1}(U)\ \text{and } U\ \text{open in }Y\}. 
		\end{multline*}
		Now since $E\subset U$ if and only if $f^{-1}(E)\subset f^{-1}(U)$, it follows that
		\begin{align*}
			&\inf\{(\mu\circ f^{-1})(U):U\supset E,\ U\ \text{open in }Y\}\\
			&= \inf\{\mu\big(f^{-1}(U)\big):f^{-1}(U)\supset f^{-1}(E),\ U\  \text{open in }Y\}\\
			&= \inf\{\mu(V):V\supset f^{-1}(E),\ V\ \text{open in }X\}\\
			&= \mu\big(f^{-1}(E)\big)\\
			&= (\mu\circ f^{-1})(E),
		\end{align*}
		where we have used the fact that $\mu$ is Radon in the second-to-last equality. Finally, we show that $\mu\circ f^{-1}$ is inner regular on open sets, i.e.:
		\[
	 (\mu\circ f^{-1})(U) = \sup\{(\mu\circ f^{-1})(K):K\subset U,\ K\ \text{compact in } Y\}
		\]
		for all open $U\subset Y$. Let $U\subset Y$ be open. Then since $f$ is a homeomorphism,
		\begin{multline*}
			\{f^{-1}(K):f^{-1}(K)\subset f^{-1}(U)\ \text{and } K\ \text{compact in }Y\} \\ = \{F:  F\subset f^{-1}(U)\ \text{and } F\ \text{compact in }X\}.
		\end{multline*}
		It follows that
		\begin{align*}
			&\sup\{(\mu\circ f^{-1})(K):K\subset U,\ K\ \text{compact in } Y\}\\
			&= \sup\{\mu\big(f^{-1}(K)\big):f^{-1}(K)\subset f^{-1}(U),\ K\  \text{compact in }Y\}\\
			&= \sup\{\mu(F):F\subset f^{-1}(U),\ F\ \text{compact in }X\}\\
			&= \mu\big(f^{-1}(U)\big)\\
			&= (\mu\circ f^{-1})(U),
		\end{align*}
		where we have used the fact that $\mu$ is Radon in the second-to-last equality.
	\end{proof}
	\begin{lemma}[From Theorem B.14.4 in \cite{OSearcoid}]
		Suppose that $X$ and $Y$ are sets and $f:X\to Y$. Suppose $C\subset X$, $D\subset Y$. Then
		\begin{list}{(\arabic{ctr1})}{\usecounter{ctr1}}\label{inclusionsets}
			\item $C\subset f^{-1}(f(C))$ with equality if $f$ is injective;
			\item $f(f^{-1}(D))\subset D$ with equality if $f$ is surjective.
		\end{list}
	\end{lemma}
	\begin{lemma}\label{innerRegOpen}
		Let $(X,\mathcal{B}_X)$ and $(Y,\mathcal{B}_Y)$ be Borel spaces, and let $\mu$ be a Radon measure on $(X,\mathcal{B}_X)$. Let $f: X\to Y$ be an $(\mathcal{B}_X,\mathcal{B}_Y)$-measurable function satisfying
		\begin{tcolorbox}	
			\begin{center}
				\begin{tabular}{ l c }
					(1) &	$f$ is surjective;		\\						
					(2) &	$f$ is a \textit{proper map}, i.e.:
				\end{tabular}
				\begin{equation}\label{radonhyp2}
					K\subset Y\ \text{is compact} \implies f^{-1}(K)\subset X\ \text{is compact;}
				\end{equation}
				\begin{tabular}{lc}
					(3) & $f$ is continuous.
				\end{tabular}
			\end{center}
		\end{tcolorbox}
		Then the image of $\mu$ under $f$ is inner regular on open sets, i.e.:
		\[
		(\mu\circ f^{-1})(U) = \sup\{(\mu\circ f^{-1})(K):K\subset U,\ K\ \text{compact}\}
		\]
		for all open $U\subset Y$.
	\end{lemma}
	\begin{proof}
		Let $U\subset Y$ be open. By assumption, $f^{-1}(K)$ is compact in $X$ whenever $K$ is compact in $Y$. It follows that
		\begin{multline*}
			\{f^{-1}(K):f^{-1}(K)\subset f^{-1}(U)\ \text{and } K\ \text{compact in }Y\} \\ \subset \{F:  F\subset f^{-1}(U)\ \text{and } F\ \text{compact in }X\}.
		\end{multline*}
		Therefore,
		\begin{multline}\label{supmuK}
			\sup\{\mu\big(f^{-1}(K)\big):f^{-1}(K)\subset f^{-1}(U)\ \text{and } K\ \text{compact in }Y\} \\ \leq \sup\{\mu(F):  F\subset f^{-1}(U)\ \text{and } F\ \text{compact in }X\}.
		\end{multline}
		For the opposite inequality, we claim that for any compact $F\subset f^{-1}(U)$, there exists a compact $K\subset Y$ such that $f^{-1}(K)\subset f^{-1}(U)$ and $\mu(F)\leq \mu\big(f^{-1}(K)\big)$. Proving this claim will establish that \eqref{supmuK} holds with equality. Let  $F\subset f^{-1}(U)$ be compact and let $K= f(F)$. Note that $f^{-1}(K) \in \mathcal{B}_X$, since $f$ is $(\mathcal{B}_X,\mathcal{B}_Y)$-measurable, and $K$ is compact (being the continuous image of a compact set) and hence Borel. Moreover, $F \subset f^{-1}\big(f(F)\big) = f^{-1}(K)$, by Lemma \ref{inclusionsets}, and hence $\mu(F) \leq \mu\big(f^{-1}(K)\big)$. Finally, since $F\subset f^{-1}(U)$, we have $K = f(F) \subset f\big(f^{-1}(U)\big)\subset U$, by Lemma \ref{inclusionsets}, and hence $f^{-1}(K)\subset f^{-1}(U)$, establishing our claim. Since $f$ is surjective, $K\subset U$ if and only if $f^{-1}(K) \subset f^{-1}(U)$ (by Lemma \ref{inclusionsets} again). It follows that
		\begin{align*}
			&\sup\{(\mu\circ f^{-1})(K):K\subset U,\ K\ \text{compact in } Y\}\\
			&= \sup\{\mu\big(f^{-1}(K)\big):f^{-1}(K)\subset f^{-1}(U),\ K\  \text{compact in }Y\}\\
			&= \sup\{\mu(F):F\subset f^{-1}(U),\ F\ \text{compact in }X\}\\
			&= \mu\big(f^{-1}(U)\big)\\
			&= (\mu\circ f^{-1})(U),
		\end{align*}
		where we have used the fact that $\mu$ is Radon in the second-to-last equality.
	\end{proof}
		\begin{lemma}\label{closureofparamdomcompact}
		$\overline{\Omega_{\mathrm{T}^*}}$ is compact.
	\end{lemma}
	\begin{proof}
		Recall that $\overline{\Omega_{\mathrm{T}^*}} = \Pi\overline{R}=\cup_{L\in \Pi} [0,L]\overline{R}$, and $\overline{R}$ is compact by Theorem 6.6.9 in\cite{Ratcliffe}, since $R$ is convex and locally finite. It suffices to show that $[0,L]\overline{R}$ is compact for each $L\in \Pi$, for a finite union of compact spaces is compact. Let $\mathcal{U} = \{U_i: i \in I\}$ be an open cover for $[0,L]\overline{R}$. Then $\{ [0,L^{-1}]U_i: i\in I\}$ is an open cover for $\overline{R}$, by Lemma \ref{translatesopen}. Thus, choosing a finite subcover $\{ [0,L^{-1}]U_i: i=1,\dots,n\}$ for $\overline{R}$, we see that $\{U_i: i=1,\dots,n\}$ is a finite subcover of $\mathcal{U}$ for $[0,L]\overline{R}$, and hence $[0,L]\overline{R}$ is compact.
	\end{proof}
	\begin{cor}\label{innerregimage}
		Let $m$ denote Lebesgue measure on $\mathbb{R}^n$. The image measure $m \circ q\vert_{\Omega_{\mathrm{T}^*}}^{-1}$ is inner regular on open sets, where $q\vert_{\Omega_{\mathrm{T}^*}}$ denotes the restriction of the natural quotient map to $\Omega_{\mathrm{T}^*}$.
	\end{cor}
	\begin{proof}
		We first prove that the image measure $m \circ q\vert_{\overline{\Omega_{\mathrm{T}^*}}}^{-1}$ is inner regular on open sets, where $q\vert_{\overline{\Omega_{\mathrm{T}^*}}}$ denotes the restriction of the natural quotient map to $\overline{\Omega_{\mathrm{T}^*}}$.	By definition, $q$ is continuous, from which it follows that $q\vert_{\overline{\Omega_{\mathrm{T}^*}}}$ is continuous (giving $\overline{\Omega_{\mathrm{T}^*}}$ the relative topology induced by $\mathbb{R}^n$). It follows that for any compact $K\subset Y$, preimage of $K$ under $q\vert_{\overline{\Omega_{\mathrm{T}^*}}}$ is closed; and since $\overline{\Omega_{\mathrm{T}^*}}$ is compact, the preimage of $K$ under $q\vert_{\overline{\Omega_{\mathrm{T}^*}}}$ is compact (being a closed subset of a compact set). Thus $q\vert_{\overline{\Omega_{\mathrm{T}^*}}}$ is proper. Thus, giving $\mathbb{R}^n$ and $\mathbb{R}^n/\mathcal{L}^*$ the usual $\sigma$-algebras, it follows from Lemma \ref{innerRegOpen} that $m \circ q\vert_{\overline{\Omega_{\mathrm{T}^*}}}^{-1}$ is inner regular on open sets. The result now follows from Lemma \ref{radonmeasaeequal}, since $\overline{\Omega_{\mathrm{T}^*}} - \Omega_{\mathrm{T}^*} \subset \Pi\overline{R} - \Pi R$, which has measure zero, by Lemma \ref{paramdomainminusPiRNull}.
	\end{proof}
	\begin{lemma}\label{outerregimage}
		Let $m$ denote Lebesgue measure on $\mathbb{R}^n$. The image measure $m \circ q\vert_{\Omega_{\mathrm{T}^*}}^{-1}$ is outer regular on Borel sets, where $q\vert_{\Omega_{\mathrm{T}^*}}$ denotes the restriction of the natural quotient map to $\Omega_{\mathrm{T}^*}$.
	\end{lemma}
	\begin{proof}
		Let $E \subset \mathbb{R}^n/\mathcal{L}^*$ be Borel. Since $q\vert_{\Omega_{\mathrm{T}^*}}$ is continuous, 
		\begin{multline*}
			\{q\vert_{\Omega_{\mathrm{T}^*}}^{-1}(U):q\vert_{\Omega_{\mathrm{T}^*}}^{-1}(E)\subset q\vert_{\Omega_{\mathrm{T}^*}}^{-1}(U)\ \text{and } U\ \text{open in }\mathbb{R}^n/\mathcal{L}^*\} \\
			\subset	\{V:  q\vert_{\Omega_{\mathrm{T}^*}}^{-1}(E)\subset V\ \text{and } V\ \text{open in }\Omega_{\mathrm{T}^*}\}. 
		\end{multline*}
		Indeed, $q\vert_{\Omega_{\mathrm{T}^*}}^{-1}(U)$ is open in $\Omega_{\mathrm{T}^*}$ whenever $U$ open in $\mathbb{R}^n/\mathcal{L}^*$. Therefore,
		\begin{multline}\label{supmuE}
			\inf\{m\big(q\vert_{\Omega_{\mathrm{T}^*}}^{-1}(U)\big) :q\vert_{\Omega_{\mathrm{T}^*}}^{-1}(E)\subset q\vert_{\Omega_{\mathrm{T}^*}}^{-1}(U)\ \text{and } U\ \text{open in }\mathbb{R}^n/\mathcal{L}^*\} \\
			\geq \inf\{m(V):  q\vert_{\Omega_{\mathrm{T}^*}}^{-1}(E)\subset V\ \text{and } V\ \text{open in }\Omega_{\mathrm{T}^*}\}. 
		\end{multline}
		For the reverse inequality, we claim that for any open $V \subset \Omega_{\mathrm{T}^*}$ with $q\vert_{\Omega_{\mathrm{T}^*}}^{-1}(E)\subset V$, there exists an open $U\subset \mathbb{R}^n/\mathcal{L}^*$ such that $q\vert_{\Omega_{\mathrm{T}^*}}^{-1}(E)\subset q\vert_{\Omega_{\mathrm{T}^*}}^{-1}(U)$ and $m\big(q\vert_{\Omega_{\mathrm{T}^*}}^{-1}(U)\big) \leq m(V)$. Proving this claim will establish that \eqref{supmuE} holds with equality. Let $V \subset \Omega_{\mathrm{T}^*}$ be open (in $\Omega_{\mathrm{T}^*}$), with $q\vert_{\Omega_{\mathrm{T}^*}}^{-1}(E)\subset V$. Let $V_0 = V\cap \Pi R$, where $R\subset \Omega_{\mathrm{T}^*}$ is the fundamental domain for $\Gamma^*$. Let $U = q(V_0)$. We show that $U$ is open in $\mathbb{R}^n/\mathcal{L}^*$. We claim that 
		\begin{equation}\label{qinverse}
			q^{-1}(U) = \bigcup_{z\in\mathcal{L}^*} V_0+z.
		\end{equation}
		Indeed, if $x\in q^{-1}(U)$, then since $\Omega_{\mathrm{T}^*}$ is a transversal for $\mathbb{R}^n/\mathcal{L}^*$, there exists a unique $y\in \Omega_{\mathrm{T}^*}$ such that $x= y +z_0$ for some $z_0\in \mathcal{L}^*$. Then 
		\[
		q\vert_{\Omega_{\mathrm{T}^*}}(y)= q(x) \in U = q\vert_{\Omega_{\mathrm{T}^*}}(V_0),
		\]
		whence 
		\[
		y \in q\vert_{\Omega_{\mathrm{T}^*}}^{-1}\big(q\vert_{\Omega_{\mathrm{T}^*}}(V_0)\big)=  V_0,
		\]
		since $q\vert_{\Omega_{\mathrm{T}^*}}$ is injective. Thus $x\in V_0 + z_0 \in  \bigcup_{z\in\mathcal{L}^*} V_0+z$. Conversely, if $x = y +z_0$ for some $y \in V_0$ and $z_0 \in \mathcal{L}^*$, then $q(x) = q(y) \in q(V_0)$, whence $x\in q^{-1}\big(q(V_0)\big) = q^{-1}(U)$, establishing our claim \eqref{qinverse}. Since $V$ is open in $\Omega_{\mathrm{T}^*}$, there exists an open $\mathcal{U} \in\mathbb{R}^n$ such that $V= \mathcal{U} \cap \Omega_{\mathrm{T}^*}$. Then $V_0 = V \cap \Pi R = \mathcal{U} \cap \Pi R$ is open in $\mathbb{R}^n$, whence $q^{-1}(U)$ is open in $\mathbb{R}^n$, by \eqref{qinverse}. It follows by the strong continuity of $q$ that $U$ is open in $\mathbb{R}^n/\mathcal{L}^*$.
		
		Since $q\vert_{\Omega_{\mathrm{T}^*}}$ is injective, $q\vert_{\Omega_{\mathrm{T}^*}}^{-1}(U) = q\vert_{\Omega_{\mathrm{T}^*}}^{-1}\big(q\vert_{\Omega_{\mathrm{T}^*}}(V_0)\big) = V_0$; in particular, $m\big(q\vert_{\Omega_{\mathrm{T}^*}}^{-1}(U)\big) = m(V_0) = m(V)$, because $\Omega_{\mathrm{T}^*} -\Pi R$ is of measure zero. Finally, since $q\vert_{\Omega_{\mathrm{T}^*}}$ is surjective, it follows that if $q\vert_{\Omega_{\mathrm{T}^*}}^{-1}(E)\subset V$, then $E= q\vert_{\Omega_{\mathrm{T}^*}}\big(q\vert_{\Omega_{\mathrm{T}^*}}^{-1}(E)\big)\subset q\vert_{\Omega_{\mathrm{T}^*}}(V) =U$, and hence $q\vert_{\Omega_{\mathrm{T}^*}}^{-1}(E)\subset q\vert_{\Omega_{\mathrm{T}^*}}^{-1}(U)$, establishing our claim. Finally, since $q\vert_{\Omega_{\mathrm{T}^*}}$ is surjective, $E\subset U$ if and only if $q\vert_{\Omega_{\mathrm{T}^*}}^{-1}(E) \subset q\vert_{\Omega_{\mathrm{T}^*}}^{-1}(U)$. It follows that
		\begin{align*}
			&\inf\{(m\circ q\vert_{\Omega_{\mathrm{T}^*}}^{-1})(U):U\supset E,\ U\ \text{open in }\mathbb{R}^n/\mathcal{L}^*\}\\
			&= \inf\{m\big(q\vert_{\Omega_{\mathrm{T}^*}}^{-1}(U)\big):q\vert_{\Omega_{\mathrm{T}^*}}^{-1}(U)\supset q\vert_{\Omega_{\mathrm{T}^*}}^{-1}(E),\ U\  \text{open in }\mathbb{R}^n/\mathcal{L}^*\}\\
			&= \inf\{m(V):V\supset q\vert_{\Omega_{\mathrm{T}^*}}^{-1}(E),\ V\ \text{open in }\Omega_{\mathrm{T}^*}\}\\
			&= m\big(q\vert_{\Omega_{\mathrm{T}^*}}^{-1}(E)\big)\\
			&= (m\circ q\vert_{\Omega_{\mathrm{T}^*}}^{-1})(E),
		\end{align*}
		where we have used the fact that $m$ is Radon in the second-to-last equality.
	\end{proof}
	\begin{cor}\label{lambdaofq}
		Let $m$ denote Lebesgue measure on $\mathbb{R}^n$. The image measure $m \circ q\vert_{\Omega_{\mathrm{T}^*}}^{-1}$ is a Radon measure, where $q\vert_{\Omega_{\mathrm{T}^*}}$ denotes the restriction of the natural quotient map to $\Omega_{\mathrm{T}^*}$.
	\end{cor}
	\begin{proof}
		Since $m \circ q\vert_{\Omega_{\mathrm{T}^*}}^{-1} = m \circ q\vert_{\overline{\Omega_{\mathrm{T}^*}}}^{-1}$, it follows that $m \circ q\vert_{\Omega_{\mathrm{T}^*}}^{-1}$ is finite on compact sets. Indeed, given any compact $K\subset  \mathbb{R}^n/\mathcal{L}^*$, we have $q^{-1}(K) \cap \overline{\Omega_{\mathrm{T}^*}}$ compact, whence
		\[
		(m \circ q\vert_{\overline{\Omega_{\mathrm{T}^*}}}^{-1})(K) = m\big(q\vert_{\overline{\Omega_{\mathrm{T}^*}}}^{-1}(K) \big)=m\big(q^{-1}(K)\cap \overline{\Omega_{\mathrm{T}^*}}\big) <\infty.
		\]
		The necessary regularity properties are given by Corollary \ref{innerregimage} and Lemma \ref{outerregimage}.
	\end{proof}
	\begin{cor}
		The image measure $m \circ \upsilon^{-1}$ is a Radon measure on $(\widehat{\mathrm{T}},\mathcal{B}_{\widehat{\mathrm{T}}})$, where $\mathcal{B}_{\widehat{\mathrm{T}}}$ is the Borel $\sigma$-algebra on $\widehat{\mathrm{T}}$ and $\upsilon$ is given by the following composition:
		\[
		\begin{tikzcd} 
			\Omega_{\mathrm{T}^*}\arrow[to=1-3, bend right=30, "\upsilon"{name=D}] \arrow[r, "q\vert_{\Omega_{\mathrm{T}^*}}"{name=qrest}]
			& \mathbb{R}^n/\mathcal{L}^* \arrow[r, "\mathcal{C}"{name=C}]
			& \widehat{\mathrm{T}} \, .
		\end{tikzcd}
		\]	
	\end{cor}
	\begin{proof}
		Corollary \ref{lambdaofq} establishes that $m \circ q\vert_{\Omega_{\mathrm{T}^*}}^{-1}$ is a Radon measure on $\mathbb{R}^n/\mathcal{L}^*$. Since $\mathcal{C}$ is a homeomorphism, it follows from Lemma \ref{imageofhomeoradon} that $(m \circ q\vert_{\Omega_{\mathrm{T}^*}}^{-1})\circ \mathcal{C}^{-1}$ is a Radon measure on $\widehat{\mathrm{T}}$. Moreover, it follows from Lemma \ref{imagemeascomp} that $m\circ \upsilon^{-1}=(m \circ q\vert_{\Omega_{\mathrm{T}^*}}^{-1})\circ \mathcal{C}^{-1}$, and hence $m\circ \upsilon^{-1}$ is a Radon measure on $\widehat{\mathrm{T}}$.
	\end{proof}
	We now show that $\Omega_{\mathrm{T}^*}$ can be given a group structure, and that $\upsilon$ is a group homomorphism. The group structure on $\Omega_{\mathrm{T}^*}$ is essentially taken from $\mathbb{R}^n/\mathcal{L}^*$: since $\Omega_{\mathrm{T}^*}$ is a transversal for $\mathbb{R}^n/\mathcal{L}^*$, each element $x+ \mathcal{L}^*\in \mathbb{R}^n/\mathcal{L}^*$ can be written uniquely as $\bar{x} + \mathcal{L}^*$ with $\bar{x}\in \Omega_{\mathrm{T}^*}$. Letting $p$ denote the inverse of $q\vert_{\Omega_{\mathrm{T}^*}}$, we have by definition that $\bar{x}= p(x+\mathcal{L}^*)$. Thus, given $\omega, \nu \in \Omega_{\mathrm{T}^*}$, we can define $\omega + \nu =  p(\omega + \nu + \mathcal{L}^*)$. It is easily verified that this makes $\Omega_{\mathrm{T}^*}$ into a group. It follows that $q\vert_{\Omega_{\mathrm{T}^*}}$ is a group homomorphism: given $\omega, \nu \in \Omega_{\mathrm{T}^*}$, 
	\begin{align*}
		q\vert_{\Omega_{\mathrm{T}^*}}(\omega+\nu) &= 	q\vert_{\Omega_{\mathrm{T}^*}}\big(p(\omega+\nu+\mathcal{L}^*)\big) \\
		&= \omega+\nu+\mathcal{L}^*\\
		&= (\omega+\mathcal{L}^*) + (\nu+\mathcal{L}^*)\\
		&= 	q\vert_{\Omega_{\mathrm{T}^*}}(\omega) + q\vert_{\Omega_{\mathrm{T}^*}}(\nu).
	\end{align*}
	Moreover, $\mathcal{C}$ is a homomorphism, since
	\begin{align*}
		\mathcal{C}\big( (\omega+ \mathcal{L}^*) +(\nu+ \mathcal{L}^*)\big) &= 	\mathcal{C}(\omega+\nu + \mathcal{L}^*)\\
		&= \chi_{\omega+\nu}^\mathrm{T} \\
		&= \chi_{\omega}^\mathrm{T} \chi_{\nu}^\mathrm{T}\\
		&= \mathcal{C}(\omega+\mathcal{L}^*) +\mathcal{C}(\nu+\mathcal{L}^*).
	\end{align*}
	It follows that $\upsilon = \mathcal{C}\circ 	q\vert_{\Omega_{\mathrm{T}^*}}$ is a group homomorphism.
	\begin{lemma}\label{scalingHaar}
		The image measure $m\circ \upsilon^{-1}$ is left-invariant, and hence, a Haar measure on $\widehat{\mathrm{T}}$. Moreover,
		\begin{equation}\label{HaarEquation}
			m\circ \upsilon^{-1} = |\Pi|m(R) \mu_{\widehat{\mathrm{T}}},
		\end{equation}
		where the measure $\mu_{\widehat{\mathrm{T}}}$ on $\widehat{\mathrm{T}}$ is the dual of the counting measure $\mu_{\mathrm{T}}$ on $\mathrm{T}$.
	\end{lemma}
	\begin{proof}
		Given $\chi_\omega^\mathrm{T}\in \widehat{\mathrm{T}}$ and $E\subset \widehat{\mathrm{T}}$, we have 
		\begin{align*}
			\chi_\omega^\mathrm{T}E &= \{\chi_\omega^\mathrm{T}\chi_\nu^\mathrm{T} : \chi_\nu^\mathrm{T} \in E \}\\
			&= \{\upsilon(\omega)\upsilon(\nu) : \nu\in \upsilon^{-1}(E) \}\\
			&= \{\upsilon(\omega+\nu) : \nu\in \upsilon^{-1}(E) \}\\
			&= \upsilon\big(\omega+ \upsilon^{-1}(E)\big).
		\end{align*}
		Since $\upsilon$ is a bijection, it follows that $\upsilon^{-1}(\chi_\omega^\mathrm{T}E) = \omega + \upsilon^{-1}(E)$. Now, since $m$ is a Haar measure, we have
		\begin{align*}
			(m\circ \upsilon^{-1})(\chi_\omega^\mathrm{T}E ) & = m\big(\upsilon^{-1}(\chi_\omega^\mathrm{T}E )\big)\\
			&= m\big(\omega + \upsilon^{-1}(E)\big)\\
			&= m\big(\upsilon^{-1}(E)\big)\\
			&= (m\circ \upsilon^{-1})(E).
		\end{align*}
		This proves that $m\circ \upsilon^{-1}$ is left-invariant. For the second claim, note that 
		\[
		(m\circ \upsilon^{-1})(\widehat{\mathrm{T}}) = m\big(\upsilon^{-1}(\widehat{\mathrm{T}})\big) = m(\Omega_{\mathrm{T}^*}) = |\Pi|m(R).
		\]
		On the other hand, $\mu_{\widehat{\mathrm{T}}}(\widehat{\mathrm{T}}) =1 = 	\tfrac{1}{|\Pi|m(R)}(m\circ \upsilon^{-1})(\widehat{\mathrm{T}})$. Since the Haar measure is unique up to scaling, the result follows.
	\end{proof}
	\begin{tcolorbox}
		\begin{rem}\label{HaarCompletionRem}
			Having established that \eqref{HaarEquation} holds on all Borel sets, it follows that
			\begin{equation}\label{mofupsiloninverse}
			\overline{m\circ \upsilon^{-1}} = |\Pi|m(R) \overline{\mu_{\widehat{\mathrm{T}}}}.
			\end{equation}
			As with the Lebesgue measure $m$, we drop the overline bar and simply write $\mu_{\widehat{\mathrm{T}}}$ for the complete Haar measure on $\widehat{\mathrm{T}}$.
		\end{rem}
	\end{tcolorbox}
		\begin{lemma}\label{restrictedqBoreliso}
		Let $R$ be the fundamental domain for $\Gamma^*$ defined in \ref{Rdomain}. Then the restricted quotient map $X_1=q\vert_{\Pi R}$ is a Borel isomorphism from $(\Pi R,\mathcal{B}_{\Pi R})$ onto its image $q(\Pi R)$ in $(\mathbb{R}^n/\mathcal{L}^*,\mathcal{B}_{\mathbb{R}^n/\mathcal{L}^*})$. 
	\end{lemma}	 
	\begin{proof}
		First note that $q\vert_{\Pi R} = X_1\vert_{\Pi R}$, and we have already shown in Remark \ref{qrestrictedbij} that $X_1$ is bijective.	Since $q$ is continuous, it is Borel measurable, and hence measurable on $\Pi R$. Now let $f:q(\Pi R) \to \Pi R$ denote the inverse of $q\vert_{\Pi R}$. Let $E\subset \Pi R$ be Borel, so that in particular, $E\in \mathcal{B}_{\mathbb{R}^n}$. Since $f$ is a bijection, the preimage $f^{-1}(E)$ coincides with $f^{-1}(E)$---the image of $E$ under $f^{-1}$---for any $E \subset \Pi R$. Hence
		\[
		f^{-1}(E) = q\vert_{\Pi R}(E)= q(E).
		\]
		We claim that
		\begin{equation}\label{qinverseclaim}
			q^{-1}\big(q(E)\big) = \bigcup_{x\in\mathcal{L}^*} E+x.
		\end{equation}
		To see this, note that if $y\in	q^{-1}\big(q(E)\big)$, then $q(y)\in q(E)$, and hence $q(y) = q(y_0)$ for a (unique) $y_0\in E$. In other words, $y+\mathcal{L}^* = y_0+\mathcal{L}^*$, so that $y = y_0 + z$ for a unique $z \in\mathcal{L}^*$. Thus $y\in  E + z$, and therefore the left-hand side of \eqref{qinverseclaim} is contained in the right-hand side. For the reverse inclusion, suppose that $y \in E + x$ for some $x\in\mathcal{L}^*$. Then $y=y_0+x$, for some $y_0\in E$, whence $q(y) = q(y_0) \in q(E)$; in other words, $y\in q^{-1}\big(q(E)\big)$. Therefore \eqref{qinverseclaim} holds. Observe that since $E$ is Borel in $\mathbb{R}^n$, and $\mathcal{B}_{\mathbb{R}^n}$ is invariant under translations, it follows from \eqref{qinverseclaim} that $q^{-1}\big(q(E)\big)$ is Borel. Thus, by the definition of the Borel structure on $\mathbb{R}^n/\mathcal{L}^*$, we have $f^{-1}(E) \in \mathcal{B}_{\mathbb{R}^n/\mathcal{L}^*}$, and therefore $f = q\vert_{\Pi R}^{-1}$ is Borel measurable. This proves that $q\vert_{\Pi R}$ is a Borel isomorphism onto its image.
	\end{proof}	
	\begin{rem}\label{qpir}
	We note in particular that it follows from \eqref{qinverseclaim} in the above proof that $q(\Pi R)$ is open in $\mathbb{R}^n/\mathcal{L}^*$.
	\end{rem}
	Recall that we defined $\upsilon_0 = \mathcal{C}\circ q\vert_{\Pi R}$. Then since $\upsilon=\upsilon_0$ a.e., it follows from Lemma \ref{imageofaeequalfns} and \eqref{mofupsiloninverse} that
	\begin{cor}\label{mofupsiloninverse1}
		\[
\mu_{\widehat{\mathrm{T}}} = \frac{1}{|\Pi|m(R)} \overline{m\circ \upsilon_0^{-1}}
\]
	\end{cor}
	
\end{appendices}

\bibliographystyle{plain}
\bibliography{Draft1Bib}{}

\end{document}